\documentclass[11pt]{article}
\usepackage{amsmath, exscale, epsfig,amssymb, 
amsthm,
 makeidx}

\usepackage{titlesec}
\titlespacing*{\section}{0pt}{0.4\baselineskip}{0.3\baselineskip}
\titlespacing*{\subsection}{0pt}{0.4\baselineskip}{0.3\baselineskip}

\usepackage[pdftex]{hyperref} %%not compatible with showkeys
\hypersetup{
    unicode      = false,     % non-Latin characters in AcrobatÃÂ¢ÃÂÃÂs bookmarks
    pdftoolbar   = true,      % show AcrobatÃÂ¢ÃÂÃÂs toolbar?
    pdfmenubar   = true,      % show AcrobatÃÂ¢ÃÂÃÂs menu?
    pdffitwindow = true,      % page fit to window when opened
    pdfnewwindow = true,      % links in new window
    colorlinks   = true,      % false: boxed links; true: colored links
    linkcolor    = blue,      % color of internal links
    citecolor    = red,      % color of links to bibliography
    filecolor    = blue,      % color of file links
    urlcolor     = green       % color of external links
}

\small
\normalsize
\usepackage[utf8]{inputenc}

\usepackage{setspace}

\usepackage[T1]{fontenc}
\usepackage[english]{babel}
\usepackage{enumerate,vmargin}

\usepackage{times}
\usepackage{amsfonts}
\usepackage{amsbsy}
\usepackage{amscd}
\usepackage{multicol}
\allowdisplaybreaks

\usepackage[all]{xy}

\usepackage{stmaryrd}
\usepackage{graphicx}
\usepackage{paralist}
\usepackage{amsfonts}
\usepackage{amssymb,mathrsfs}
\setmarginsrb{2.9cm}{2.6cm}{2.9cm}{1.7cm}{0cm}{0mm}{0cm}{9mm}

\def\build#1_#2^#3{\mathrel{\mathop{\kern 0pt#1}\limits_{#2}^{#3}}}
\def\noi{{\noindent}}

\def\cq{$\hfill \square$}
\def\un{{\bf 1}}

\newcommand{\fdelta}{{\boldsymbol{\delta}}}

\newcommand{\bbT}{\mathbb{T}}

\newcommand{\bD}{\mathbb{D}}

\newcommand{\bE}{{\bf E}}

\newcommand{\bs}{{\bf s}}

\newcommand{\bM}{\mathbb{M}}
\newcommand{\bN}{\mathbf{N}}
\newcommand{\bbN}{\mathbb{N}}

\newcommand{\bP}{{\bf P}}
\newcommand{\bR}{\mathbf{R}}

\newcommand{\bbU}{\mathbb{U}}

\newcommand{\cI}{\mathcal{I}}

\newcommand{\cM}{\mathcal{M}}

\def\be{\begin{equation}}

\def\ee{\end{equation}}
\def\ba{\begin{eqnarray*}}
\def\ea{\end{eqnarray*}}

\def\noi{\noindent}

\newcommand{\lgeo}{[\![}
\newcommand{\rgeo}{]\!]}
\def\cqfd{ \hfill $\blacksquare$ }

\newcommand{\eqo}{\! = \! }
\newcommand{\geqo}{\! \geq \! }
\newcommand{\leqo}{\! \leq \! }
\newcommand{\ino}{\! \in \! }

\newcommand{\bbR}{\mathbb{R}}
\newcommand{\leko}{\! < \! }
\newcommand{\geko}{\! > \! }

\newcommand{\deHaus}{\delta_{\mathtt{Haus}}}
\newcommand{\dePro}{\delta_{\mathtt{Pro}}}

\newcommand{\dGHP}{\fdelta_{\mathtt{GHP}}}

\newcommand{\epp}{\varepsilon}

\def\btt{\mathbf{t}}
\def\fftree{t}
\def\cW{\mathcal{W}}
\def\baa{\mathbf{a}}
\def\bbb{\mathbf{b}}
\def\bC{\mathbf{C}}
\def\bD{\mathbf{D}}
\def\bCqq{\mathbf C^{_0}_{^{d}}}
\def\bCpl{\mathbf C^{_0}_{^{1}}}

\newcommand{\bbZ}{\mathbb Z}

\newcommand{\bmu}{\boldsymbol \mu}
\newcommand{\bzeta}{\boldsymbol \zeta}

\def\ccL{\mathscr{L}}

\def\bdd{\mathbf{d}}
\def\bgg{\mathbf{g}}

\def\bgam{{\boldsymbol{\gamma}}}
\def\bsigma{{\boldsymbol{\sigma}}}
\def\bdelta{{\boldsymbol{\delta}}}
\def\bxxi{{\boldsymbol{\xi}}}
\def\bxxi{\overline{\xi}}
\newcommand{\ccF}{\mathscr F}

\def\gamun{\gamma_{_{1}}}
\newcommand{\fSigma}{\boldsymbol{\Sigma}}
\newcommand{\bSigma}{\boldsymbol{\Sigma}}
\newcommand{\fbeta}{\boldsymbol{\beta}}
\def\Psimu{\Psi_{\! \mu}}
\newcommand{\diffe}{\partial}
\def\bS{\mathbf{S}}

\def\MMT{\mathbf{MT}}
\def\bbd{\mathbf{d}}

\newtheoremstyle{thmstyl}
{3.5pt} % space above
{2.5pt} % space below
{\em} % Body font
{} % Indent amount
{\bfseries} % Theorem head font
{.} % Punctuation after theorem head
{.5em} % space after theorem head
{} % theorem head spec 
\theoremstyle{thmstyl}

\newtheorem{theorem}{Theorem}[section]

\newtheorem{lemma}[theorem]{Lemma}
\newtheorem{proposition}[theorem]{Proposition}

\newtheoremstyle{dfstyl}
{3.5pt} % space above
{2.5pt} % space below
{} % Body font
{} % Indent amount
{\bfseries} % Theorem head font
{.} % Punctuation after theorem head
{.5em} % space after theorem head
{} % theorem head spec 
\theoremstyle{dfstyl}

\newtheorem{definition}[theorem]{Definition}

\newtheorem{remark}[theorem]{Remark}

\newcommand{\eqnsection}{
\renewcommand{\theequation}{\arabic{section}.\arabic{equation}}
    \makeatletter
    \csname  @addtoreset\endcsname{equation}{section}
    \makeatother}
\eqnsection

\begin{document}

\setlength{\abovedisplayskip}{5pt}
\setlength{\belowdisplayskip}{5pt}

\title{\textbf{Scaling limits of critical branching \\ random walks and their snakes} }
\date{}

%%%%%
%%%%% Limites d'échelle des marches aléatoires branchantes critiques via leurs serpents. 
%%%%%

\author{Thomas \textsc{Duquesne}
\thanks{{\scriptsize LPSM, Sorbonne Universit\'e, 4 place Jussieu, F-75252 Paris Cedex 05, France
Email: thomas.duquesne@sorbonne-universite.fr}}
\and  Fael \textsc{Rebei}
\thanks{{\scriptsize LPSM, Sorbonne Universit\'e, 4 place Jussieu, F-75252 Paris Cedex 05, France
Email: rebei@lpsm.paris}}
%\and Other \textsc{Person}
%\thanks{other institution
%Email: blibli@bloblo.fr }
}

\maketitle

\vspace{-7mm}

\begin{abstract} 

We study scaling limits of $\mathbb R^d$-valued branching random walks (BRWs for short) in the following setting:
 indexing-trees are Galton-Watson trees with critical offspring distribution which are assumed, when appropriately rescaled, to converge to a critical $\psi$-L{\'e}vy tree; conditional on the indexing tree, the jumps of the BRW may exhibit dependence within the same sibling group, but distinct sibling groups are independent, and their distributions are centered but may vary (and for instance depend on the whole indexing tree), with a typical value of the walk nevertheless remaining within the domain of attraction of a mixture of Gaussian distributions. Under these assumptions and a necessary additional assumption, which we call the \emph{discrete Sheu assumption}, we show that rescaled discrete snakes converge functionally to the $\psi$-Brownian snake introduced by Le Gall \& Le Jan (1998) and D.~\& Le Gall (2002). The method uses a coupling result for BRWs of independent interest. An application is given concerning the scaling limits of the range of BRWs that take their values in the $\mathtt b$-ary tree to the reflected $\psi$-Brownian cactus which is a variant of the Brownian cactus introduced by Curien, Le Gall \& Miermont (2013). 
{\small 

\medskip

\noi
{\textbf{Keywords} $\, $ Branching random walk; Branching process; Brownian snake; Continuum random tree; Coupling; Discrete snake; Galton–Watson tree; L{\' e}vy tree; Limit theorem

\smallskip

\noi
}
%60B05 = Probability measures on topological spaces
%53C23= Global geometric and topological methods ({\`a} la Gromov); differential geometric analysis on metric spaces
%60J80 Branching processes (Galton-Watson, birth-and-death, etc.)
%60G50= Sums of independent random variables; random walks
%60G51 =Processes with independent increments; L ́evy processes
%60F17 =Functional limit theorems; invariance principles

}

\end{abstract}
\section{Introduction}
\label{Introdevsec}
Branching random walks (BRW for short) combine a branching phenomenon with spatial motion: they can be viewed as a process $\smash{(S_u)_{u\in \btt}}$ where $\btt$ is the random genealogical tree of a population and $S_u$ is the position of the individual $u\ino \btt$. The dynamic is specified by the joint distribution of the jumps $\xi_u$ made by individual $u$ relative to its direct parent $\smash{\overleftarrow{u}}$, that is, $\smash{\xi_u = S_u \! -\! S_{\overleftarrow{u}}}$. BRWs have been the subject of intensive research since the 1970s. They are closely related to numerous models in physics and statistical mechanics (traveling waves and FKPP equations, the GREM, Mandelbrot's multiplicative cascades, Gaussian free fields): we refer to 
Shi \cite{Shi15} for an overview of this topic when the ambient space is $\bbR^d$; we also refer to the works of Gouëzel, Huerter, Lalley \& Sellke
\cite{LalSel97, LalHue00, Lal06, LalGou13} for the study of BRWs in hyperbolic spaces and, more specifically, to 
Liggett \cite{Lig96} and to Benjamini \& Müller \cite{BenMue12} when the ambient space is a tree.

In this article, we consider scaling limits of sequences of BRWs $\, \smash{\mathbf S_n\eqo (S_{n,u})_{u\in \btt_n}}$ when the space is $\smash{\bbR^d}$ (except for an application to BRWs with values in the $\mathtt b$-ary tree) and 
when the $\btt_n$ are Galton-Watson trees (subject to various conditionings, or forests) with offspring  distributions $\mu_n$ (GW($\mu_n$)-trees, for short), which are assumed to be critical. The jumps are supposed to be $(i)$ centered, that is, 
$\smash{\bE \big[ S_{n,u} \! -\! S_{n,\overleftarrow{u} }\,  \big|\,  \btt_n \big]\eqo 0}$, $(ii)$ independent by group of siblings (i.e., within the same sibling group, the jumps may be dependent, but jumps from different sibling groups are independent: see Definition \ref{sibinddef}), $(iii)$ their laws, conditional on $\btt_n$, may vary and depend, for instance, on $u\ino \btt_n$, but these inhomogeneities remain moderate (i.e., controlled by a moment assumption) so that the jumps stay in the domain of attraction 
of $\smash{\bbR^d}$-valued mixtures of Gaussian laws: see, more precisely, the set of assumptions \textbf{Sib-Ind}$_{_{\, }}$($\smash{\baa, \bbb, \bS_\cdot, G, C_\cdot, \fbeta}$) stated between Remarks \ref{siblindep} and \ref{iidass}. 

Our asymptotic regime yields two renormalization sequences $\smash{(a_n)_{n\in \bbN}}$ and $\smash{(b_n)_{n\in \bbN}}$ satisfying:
\begin{equation}
\label{renormab} 
\textrm{{\small \texttt{Norm} ($\baa, \bbb$):}}  \qquad  \qquad  a_n \xrightarrow[n\to \infty]{\; } \infty \quad  \textrm{and} \quad  \lambda_n\! :=\! \frac{b_n}{a_n}    \xrightarrow[n\to \infty]{\; } \infty \; .
\end{equation}  
(The sequence $\smash{(\lambda_n)_{n\in \bbN}}$ is introduced to simplify notation). Here, $b_n$ is a \emph{size parameter} for the tree 
$\btt_n$, i.e., $\smash{\frac{1}{b_n} \# \btt_n}$ converges in distribution on $\smash{\bbR_+^*}$,  
and $\lambda_n$ is a \emph{height parameter} for $\btt_n$: if $\mathbf u_n$ is an individual chosen uniformly at random 
in $\btt_n$ and if $|\mathbf u_n|$ denotes its height, i.e., its generation, then $\smash{\frac{1}{\lambda_n}  |\mathbf u_n|}$ converges in distribution on $\smash{\bbR_+^*}$. Accordingly, $\smash{\frac{1}{\sqrt{\lambda_n}}S_{n, \mathbf u_n}}$ converges to a Gaussian vector, or a mixture of Gaussian vectors.  

In this context, let us briefly review the limit theorems concerning the trees $\btt_n$, starting with their Galton-Watson processes $\smash{(Z^{n}_k)_{k\in \bbN}}$ (here, $\smash{Z^n_k}$ is the number of individuals situated at height $k$ in $\btt_n$). As shown in Grimvall \cite{Gr74}, the only admissible limits in law for $\smash{(\frac{1}{a_n} Z^n_{^{\lfloor \lambda_n s\rfloor}} )_{s\in \bbR_+}}$ are the Continuous State Space Branching Processes (CSBP for short) introduced by Jirina \cite{Ji} and Lamperti \cite{La1,La2,La3}. A CSBP $\smash{(Z_s)_{s\in \bbR_+}}$ is a Feller Markov process that can be viewed as a spectrally positive Lévy process $\smash{(X_s)_{s\in \bbR_+}}$, which is time-changed through the Lamperti transform (see Helland \cite{He} for continuity properties of this transform on the Skorokhod space): the distribution of the CSBP $Z$ is characterized by its branching mechanism $\psi$, which is the Laplace exponent of the Lévy process $X$. It is necessarily of the Lévy–Khintchine form recalled in (\ref{LK}).
Assuming that $\smash{\frac{1}{a_n} Z^n_{\lfloor \lambda_n \rfloor}\! \to \! Z_1}$ in distribution, the appropriately renormalized trees $\btt_n$ themselves converge (in the sense of their finite dimensional marginal laws) 
to a random metric space called the $\psi$-Lévy tree; this continum tree is defined via a process $\smash{(H_{ s})_{s\in \bbR_+}}$ called the $\psi$-height process (introduced in Le Gall \& Le Jan \cite{LGLJ98} and further studied in D.~\& Le Gall \cite{DuLG02}), which is a local-time functional of $X$: see (\ref{defH}). 
More precisely, we denote by $(C_s (\btt_n))_{s\in \bbR_+}$ the contour process of $\btt_n$ (see Definition \ref{contsnadef}):  informally, it is the distance from the root at time $s$ of a particle originating from the root of $\btt_n$ that travels clockwise through this tree at unit speed, backtracking as little as possible (the trees $\btt_n$ 
are chosen to be ordered and can therefore be drawn in the clockwise-oriented half-plane). 
Assuming that $\smash{\frac{1}{a_n} Z^n_{\lfloor \lambda_n \rfloor}\! \to \! Z_1}$, we have convergence 
\begin{equation}
\label{cvContintro}  
\big(C^{(n)}_s\big)_{s\in \bbR_+} \xrightarrow[n\to \infty]{\; } (H_{ s/2})_{s\in \bbR_+}  
\end{equation}  
in the sense of the finite dimensional marginal laws, 
where $\smash{C^{(n)}_s\! :=\! \frac{1}{\lambda_n} C_{b_n s} (\btt_n)}$, $\smash{s\in \bbR_+}$, and where 
$\smash{(H_s)_{s\in \bbR_+}}$ is the $\psi$-height process. As shown in D.~\& Le Gall \cite{DuLG02}, it admits 
a continuous version if and only if $\psi$ satisfies the almost-sure extinction condition of the corresponding CSBP, 
the so-called Grey condition: $\smash{\int^\infty_1 dz/ \psi (z) \leko \infty}$ (see Grey \cite{Grey}). In this case, and 
under the necessary and sufficient assumption that the height distributions of the $\btt_n$ are tight in the scale 
$\lambda_n$, which amounts to the assumption \texttt{Hght}$_{^{\,}}$($\baa, \bbb, \bmu$) in (\ref{Hghtabmu}), 
convergence (\ref{cvContintro}) holds in distribution in $\smash{\bC(\bbR_+, \bbR)}$.  As in the discrete case, we 
think of $H_s$ as the height in the $\psi$-Lévy tree of a particle at time $s$ that travels clockwise through the tree. 
In particular, the distance $\smash{d_H(s,s')}$ between the vertices visited at times $s$ and $s'$ is given by 
$\smash{H_s + H_{s'} - 2m_H(s,s')}$, where $\smash{m_H(s,s') := \min_{s\wedge s'\leq r\leq s\vee s'} H_r}$ 
(see Section \ref{traceapplsec}). 

In the case of a quadratic branching mechanism $\smash{ \psi (z)\eqo z^2}$, the $\psi$-Lévy tree 
is (a variant of) Aldous' CRT \cite{Al91}, which is also called the Brownian tree, and when 
$\mu_n\eqo \mu$, $n\ino \bbN$, has a moment of order $2$ the functional convergence (\ref{cvContintro}) 
was proven in Aldous \cite{Al93}.  
In these cases, the empirical measures of the BRWs also converge: let $\smash{ \mathcal Z^n_{^k}(dy)
\eqo \sum \delta_{^{\! S_{n,u}/\sqrt{\lambda_n}}} (dy)}$, where the sum ranges over the $u \ino \btt_n$ of height 
$k$; then $\smash{(\frac{1}{{a_n}} \mathcal Z^n_{^{\lfloor b_n s \rfloor}}(dy))_{s\in \bbR_+} \! \to \! 
(\mathcal Z_s (dy))_{s\in \bbR_+}}$ in distribution for processes taking values in the space of finite measures on 
$\smash{ \bbR^d}$. Here, $\smash{(\mathcal Z_s (dy))_{s\in \bbR_+}}$ is a Dawson-Watanabe superprocess, 
a process whose Markov dynamics is closely related to the semilinear PDEs 
$\smash{\partial_t v\! -\! \frac{_1}{2}\Delta v\eqo v^2}$, for which it provides a probabilistic representation: 
see Perkins \cite{Per02}, Etheridge \cite{Et}, and Le Gall \cite{LG99} for a detailed exposition of 
Dawson-Watanabe superprocesses. They were generalized to general branching mechanisms $\psi$ by 
Dynkin in \cite{Dy0, Dy92, Dy1, DynCRM, dynkin2002diffusions} and are called 
$\psi$-superprocesses 
(here we restrict ourselves to Brownian spatial motions), and their Markov dynamics are related to the PDEs 
$\smash{\partial_t v\! -\! \frac{_1}{^2}\Delta v\eqo \psi(v)}$, for which they also provide a probabilistic representation. 

The $\psi$-Brownian snakes $\smash{(W_{\! s}(\cdot))_{s\in \bbR_+}}$ facilitate the study of $\psi$-superprocesses 
by closely linking them to their underlying Lévy tree, as in the case of BRWs. In the case of quadratic branching mechanisms, Brownian snakes were introduced by Le Gall in \cite{LG1,LG93bis, LG94, LG94bis, LG95} and were later 
extended to general $\psi$-mechanisms in Le Gall \& Le Jan \cite{LGLJ98bis} and in D.~\& Le Gall \cite{DuLG02}. Informally, $\smash{W_{\! s} (\cdot) \ino \bC(\bbR_+, \bbR^d)}$ represents the historical path, i.e., the path of the ancestors of the individual in the $\psi$-Lévy tree that is visited at time $s$ (let us call it the $s$-individual for simplicity): $\smash{W_{\! s} (\cdot)}$ therefore has the same distribution as $\smash{B_{\cdot \wedge H_s}}$, where $B$ is a $\smash{\bbR^d}$-valued  Brownian motion, the position of the $s$-individual being 
$\smash{\widehat{W}_{\! s} \! :=\!  W_{\! s} (H_s)}$; for distinct $\smash{s, s' \in \mathbb{R}_+}$, $\smash{W_{\! s}(\cdot)}$ and $\smash{W_{\! s'} (\cdot)}$ coincide up to time $\smash{m_H(s, s')}$, which is the height of the most recent common ancestor of the $s$- and $s'$-individuals, and $\smash{W_{\! s} (\cdot + m_H(s, s')) \! - \! W_{\! s}  (m_H (s,s'))}$ and $\smash{W_{\! s'} (\cdot + m_H (s,s')) \! -\! W_{\! s}  (m_H (s,s'))}$ are two independent Brownian motions stopped at the respective times 
$\smash{H_s\! -\! m_H(s,s')}$ and $\smash{H_{s'} \! -\! m_h(s,s')}$: see (\ref{Brosnadefbis}). Even when $H$ is continuous, $W$ does not necessarily have a continuous version: as shown in D.~\& Le Gall \cite{DuLG02}, a necessary and sufficient condition for the continuity of the $\psi$-Brownian snake is the so-called Sheu condition of compactness of the support of the $\psi$-superprocess: $\smash{\int_1^\infty ds/ (\int_0^s \psi (r) dr)^{1/2}} $ $\smash{ \leko \infty}$ (see Sheu \cite{Sh1} or Hesse \& Kyprianou \cite{HesKyp14}). 

As in the continuous setting, we can define the discrete snakes $\smash{(W_{\! s}(\mathbf S_n, r))_{s,r\in \bbR_+}}$ 
associated with the BRWs $\smash{\mathbf S_n \eqo (S_{n, u})_{u\in \btt_n}}$ (see Definition \ref{contsnadef}). 
Discrete snakes appear in the metric bijections at the heart of the constructions of planar maps and their limits 
(see Addario-Berry \& Albenque \cite{AddarioBerryAlbenque2017}, Bouttier, Di Francesco \& Guitter 
\cite{Bouttier2004}, Chassaing \& Schaeffer \cite{ChassSchaef_2004}, Janson \& Stef{\'a}nsson 
\cite{janson2015scaling}, Le Gall \cite{le_gall_2013}, Le Gall \& Miermont \cite{LeGallMiermont2011}, 
Marckert \& Miermont \cite{marckert2007invariance}, 
Marckert \& Mokkadem \cite{MaMo03, MarckertMokkadem2006}, 
Miermont \cite{miermont2013brownian}, 
Schaeffer \cite{Schaeffer1998}) 
but also in the study of other objects in random geometry such as cactuses 
(see Curien, Le Gall \& Miermont \cite{CuLGMi13}, Le Gall \cite{LG15}) or the trace of tree-valued BRWs 
(see D.~, Khanfir, Lin \& Torri \cite{DuKhLiTo22}, D.~\& Khanfir \cite{duquesne2026scaling}) or other 
combinatorial constructions (see  Addario-Berry, Donderwinkel, Goldschmidt \& Mitchell \cite{ABDoGoMi25+}, 
Gittenberger \cite{Gittenberger2003} Janson \& Marckert \cite{JanMar05}, Marckert \cite{marckert2008lineage}, 
Marckert \& Mokkadem \cite{marckert2003depth}, Marzouk \cite{Mar20}). 
The purpose of this article is to establish a limit theorem analogous to (\ref{cvContintro}) for renormalized discrete 
snakes $\smash{W^{_{(n)}}_{\! s}(r)\eqo \frac{1}{\sqrt{\lambda_n}} W_{b_n s}(\mathbf S_n, \frac{1}{\lambda_n} r)} $, 
$\smash{r, s \ino \mathbb{R}_+}$, for locally centered, independent by group of siblings, 
“moderately inhomogeneous” jump distributions, under the convergence assumption (\ref{cvContintro}) 
of the genealogy, and under the additional following assumption that we call the \emph{discrete Sheu condition}: 
\begin{equation}
\label{Sheudis} 
\textrm{{\small \texttt{Sheu}$_{^{\,}}$($\baa, \bbb, \bmu$)}:} \qquad \quad   \lim_{y\to \infty}
 \limsup_{n\to \infty}\int_{y}^{a_n}\!\!\!\! 
 \frac{\mathrm ds}{\sqrt{\! \int_0^s\! \psi_n (r) \mathrm dr\, }} =0 \; , 
\end{equation} 
where we have defined $\smash{\psi_n (z)\eqo b_n \big( g_{\mu_n} \! \big(1\! -\! \tfrac{z}{a_n} \big) \! -\! 1+ 
\tfrac{z}{a_n} \big)}$, $\smash{z \ino [0, a_n]}$, with $\smash{g_{\mu_n} (r) \eqo \sum_{k \in \mathbb{N}} r^k \mu_n (k)}$, $\smash{r \ino [0, 1]}$. 
Under these assumptions, Theorem \ref{maincvsnake} states that the following convergence holds in 
distribution on $\smash{\bC(\bbR_+, \bbR)\times \bC(\bbR_+, \bC(\bbR_+, \bbR^d))}$
\begin{equation}
\label{cvsnakeintro} 
\Big( \big(C^{_{(n)}}_s\big)_{s\in \bbR_+}, \big(W^{_{(n)}}_{\! s} (\cdot) \big)_{s\in \bbR_+}\Big)  \xrightarrow[n\to \infty]{\; } 
\big( (H_{s/2})_{s\in \bbR_+}  , (\sqrt{\! \fbeta}.W_{\! s/2} (\cdot) ) \big)_{s\in \bbR_+} , 
\end{equation}
where $\smash{(W_s (\cdot))_{s\in \bbR_+}}$ is a $\psi$-Brownian snake associated with $H$ as above, 
and where $\smash{\sqrt{\!  \fbeta \, }}$ is the square root of a random covariance matrix $\fbeta$ that depends only on $H$. 

The difficulty in (\ref{cvsnakeintro}) does not lie in the convergence of finite-dimensional marginal distributions, but rather in the tightness: in fact, our assumptions allow the limiting trees to have infinite branch points, and an extreme value phenomenon can conflict with tightness if the assumption (\ref{Sheudis}) is not satisfied, even when the jumps are simply i.i.d. with a symmetric exponential  distribution, or even a uniform distribution on $\{ -1, 1\} $ (the assumption (\ref{Sheudis}) is, moreover, necessary in these cases). The proof of (\ref{cvsnakeintro}) proceeds in two distinct steps: first, we prove (\ref{cvsnakeintro}) for specific, relatively explicit cases (Theorem \ref{Sheuexplain}), and then we extend this convergence using a coupling argument of general relevance and independent interest (see Theorem \ref{brwcouplingth}, which applied to our cases in Theorem \ref{extenscv}). 
Up to logarithmic factors due to the coupling, the assumptions of Theorem \ref{maincvsnake} proving  (\ref{cvsnakeintro}), are in some sense optimal: see Section \ref{assumpsec} for a detailed discussion. 
The proof approach differs from previous works of Janson \& Marckert \cite{JanMar05}, Marckert \& Miermont \cite{marckert2007invariance}, Marzouk \cite{Mar20} and Addario-Berry, Donderwinkel, Goldschmidt \& Mitchell \cite{ABDoGoMi25+}.
With regard to centered snakes only, Theorem \ref{maincvsnake} extends these previous works concerning 
the branching mechanism (which is quadratic, or more generally stable in the case of Marzouk \cite{Mar20}), the assumptions regarding the independence of jumps and their possible inhomogeneities in distribution; however, this theorem only recovers, up to a logarithmic factor, the optimal moment assumptions obtained in those articles for stable mechanisms and for i.i.d. or homogeneous jumps. We refer the reader to Section \ref{assumpsec} for a more detailed discussion on our assumptions and on the connections with previous works.

In application of Theorem \ref{maincvsnake}, we prove, at Theorem \ref{applicactus}, that the range of critical BRWs with values in the $\mathtt b$-ary tree converges to the reflected $\psi$-Brownian cactus, which extends the result obtained in 
D., Khanfir, Lin \& Torri \cite{DuKhLiTo22} concerning the cases of stable branching mechanisms: see Section \ref{traceapplsec}.

\smallskip

\noi
\textbf{Organization of the article.} Sections \ref{presentsec}, \ref{4casessec}, and \ref{statementssec} provide a detailed overview of the probabilistic objects, the assumptions, and the cases under consideration, as well as a precise statement of the main results; Section \ref{assumpsec} examines the assumptions (interdependence, optimality), provides examples, and discusses 
previous works related to our article; Section \ref{traceapplsec} presents an application to the convergence of the range of BRWs with values in the $\mathtt b$-ary tree to the $\psi$-reflected Brownian cactus.
Section \ref{discrobjsec} precisely defines the various encodings of discrete objects; Section \ref{Brosnasec} introduces Brownian snakes and proves a technical result on their oscillations (Lemma \ref{osc1lemma}). Section \ref{prelimBRWsec} is devoted to preliminary estimates on BRWs (in particular, Section \ref{couplingsec} states and proves the coupling results). 
The proofs of the main theorems are given in Section \ref{Thm1pfsec}. The proof of the application to BRWs with values in the $\mathtt b$-ary tree is given in Section \ref{pfsecapplicactus}. In Appendix \ref{unifboundedsec}, we state and prove a convergence result for the case of i.i.d. jumps with uniformly bounded support. This result is optimal and cannot be fully derived from the previous theorems due to logarithmic factors arising from coupling. Nevertheless, its proof, which is largely analytical, is lengthy and disconnected from the proofs of the other results in the article; therefore, it appears only in the full version of the article pre-published on arXiv, and not in the version submitted to a journal.

\smallskip

\noi
\textbf{Acknowledgments.} One of the two authors would like to thank, on the one hand, Nicolas Fournier for suggesting, a good while ago now, to consider coupling to extend the convergence of our snakes, and, on the other hand, Robin Khanfir 
for carefully  reading an earlier version of the manuscript 
(of course, any errors or typos in the current article are solely our own).

\subsection{A brief presentation of the random objects considered in this paper}
\label{presentsec}

\noi
\textbf{$\;$ Branching random walks and their encoding processes.} 
We will use the following common abbreviations: e.g.~for \emph{example given}, i.e.~for \emph{id est}, 
i.i.d.~for \emph{independent and identically distributed}, iff for \emph{if and only if}, r.v. for \emph{random variable}, RW for \emph{Random Walk}, w.l.o.g~for \emph{without loss of generality}. A few others 
are introduced as we go along. Although, in its strict sense, the term \emph{Polish} refers only to a topology, we find convenient to use it as a synonym for a \emph{complete and separable metric space}. 
We use notations $\bbN\! :=\! \{ 0,1,2, \ldots \}$, for the set of natural numbers, including $0$, $\smash{\bbN^*\!:=\! \bbN\backslash \{ 0\}}$, $\smash{\bbR_+ \! :=\! [0, \infty)}$ and $\smash{\bbR_+^*\! :=\! \bbR_+\backslash \{ 0\}}$. Finally, unless explicitly stated otherwise, 
all r.v.s are assumed to be defined on the same probability space $\smash{(\Omega, \mathscr F, \bP)}$, 
which is sufficiently rich to carry as many independent r.v.s as needed.

In this article we study scaling limits of $\smash{\bbR^d}$-valued critical \emph{branching random walks} (BRWs for short) which we generically denote by $\smash{\bS\eqo (S_u)_{u\in \btt}}$. The indexing tree $\btt$ is assumed to have a \emph{root} which is denoted by $\varnothing$ and to be \emph{ordered}: namely if we interpret $\btt$ as the \emph{family tree} of a population whose $\varnothing$ is the \emph{ancestor}, then it means that siblings have a relative birth-order (see Definition \ref{Ulamtree} for a more precise framework on rooted ordered trees). For all $u\ino \btt$, we denote by $\smash{k_u (\btt)}$ the \emph{number of children of $u$}, i.e., its outdegree, and if $\smash{k_u(\btt)\eqo 0}$, we call it a \emph{leaf} of $\btt$. 
For all $u ,v\ino \btt$, we denote by $\lgeo u,  v\rgeo$ the $\smash{d_{\mathtt{gr}}}$-shortest path joining $u$ to $v$, where $\smash{d_{\mathtt{gr}}}$ stands for the graph distance on $\btt$. Therefore $\lgeo \varnothing , u\rgeo$ is the ancestral line of $u$. It is convenient to use the notation 
$\smash{ |u|:=d_{\mathtt{gr}} (\varnothing, u)}$ for the \emph{height of $u$}. 
We also denote by $u\wedge v$ the \emph{most recent common ancestor of $u$ and $v$}. Namely, $\lgeo \varnothing , u\rgeo\cap \lgeo \varnothing , v\rgeo\eqo \lgeo \varnothing , u \wedge v\rgeo$. 

For any individual $u\ino \btt$, $\smash{S_u\ino \bbR^d}$ denotes its \emph{spatial position}. We restrict ourselves to BRWs starting from the origin, i.e.~$\smash{S_\varnothing \eqo 0}$, so that they are fully characterized by their \emph{jumps}, which are given for all 
$u\ino \btt \backslash \{\varnothing \}$ by $\smash{\xi_u \eqo S_u\! -\! S_{\overleftarrow{u}}}$, where here $\smash{\overleftarrow{u}}$ stands for the \emph{direct parent of $u$}.

 Scaling limits of BRWs are obtained by means of two encoding processes: \emph{the contour process} $\smash{(C_{\! s} (\btt))_{s\in \bbR_+}}\! $ of $\btt$ and the \emph{discrete snake} $\smash{(W_{\! s}(\bS, r))_{r,s\in \bbR_+}\!} $ of $\bS$. They are defined as follows.   
\begin{definition}
\label{contsnadef} We assume that $\btt$ is finite and that $(S_u)_{u\in \btt}$ is a $\bbR^d$-valued BRW. 
\begin{compactenum}

\smallskip

\item[$(a)$] The \emph{contour exploration} of $\btt$ is the sequence of vertices $\smash{(v_k)_{0\leq k\leq 2( \# \btt -1)}}$ 
that is recursively defined as follows: \emph{$v_0\eqo \varnothing$; $v_{k+1}$ is the first child (relative to birth-order) of $v_k$ that does not belong to $\{ v_0, v_1, \ldots, v_k\}$ if there is one; if not and if furthermore $v_k \! \neq \! \varnothing$, then 
$\smash{v_{k+1}\eqo \overleftarrow{v}_{\!\! k}}$; otherwise $\smash{v_{k+1}\eqo \varnothing}$ and $\smash{k+1 \eqo 2(\# \btt \! -\! 1)}$, necessarily}.

\smallskip

\item[$(b)$]  For all $\smash{k\ino \{ 0, \ldots ,2(\# \btt \! -\! 1)\}}$, we set $\smash{C_k (\btt)\eqo |v_k|}$, the height of $v_k$, and we conveniently set 
$\smash{C_{k} (\btt)\eqo 0}$ for all integers $\smash{k\geqo 2(\# \btt\! -\! 1)}$.  
We next extend $\smash{C_\cdot (\btt)}$ continuously on $\smash{\bbR_+}$ by setting for all $\smash{s\ino \bbR_+}$: 
$$ C_{s}(\btt) \eqo C_{\lfloor s\rfloor} (\btt)+ \{ s\} (C_{\lceil s\rceil} (\btt) \! -\! C_{\lfloor s \rfloor} (\btt)), $$ 
where $\smash{\lceil s \rceil \eqo \lfloor s \rfloor +1}$ and $\smash{\{ s \}\eqo s\! -\! \lfloor s \rfloor}$. 
The resulting continuous function $\smash{(C_s(\btt))_{s\in \bbR_+}}$ is called the \emph{contour process} of $\btt$.

\smallskip

\item[$(c)$] For all $\smash{k\ino \{ 0, \ldots ,2(\# \btt \! -\! 1)\}}$ and all $\smash{\ell \ino \bbN}$, let $\smash{v_k(\ell)\ino \lgeo \varnothing , v_k \rgeo}$ be such that $\smash{|v_k(\ell)|\eqo \ell \! \wedge \! |v_k|}$. 
We set $\smash{W_{\! k} (\bS, \ell)\eqo S_{\! v_k(\ell)}}$. We first fix $k$ and we extend $\smash{W_{\! k} (\bS, \cdot )}$ continuously by setting 
$\smash{W_{\! k} (\bS,r)\!  = \!    W_{\! k} (\bS,\lfloor r \rfloor ) + \{ r\} \big( W_{\! k} (\bS,\lceil r \rceil ) \! -\!  W_{\! k} (\bS,\lfloor r \rfloor ) \big)}$ for all $\smash{r\ino \bbR_+}$. We extend $W_{\! \cdot} (\bS, \cdot)$ continuously by setting 
$$ \forall r,s\ino \bbR_+, \quad W_{\! s} (\bS, r) = \left\{ \begin{array}{ll}
  W_{\! \lceil s \rceil }  \big(\bS ,  r\! \wedge \! C_{\! s} (\fftree) \big) & \textrm{if $\; |v_{\lceil s \rceil } | \eqo 1+ |v_{\lfloor s \rfloor }|$,}\\
 W_{\! \lfloor s \rfloor}  \big( \bS , r\!  \wedge \! C_{\! s} (\fftree) \big)  & \textrm{if $\; |v_{\lfloor s \rfloor} | \eqo 1+ |v_{ \lceil s \rceil }|$.}
\end{array} \right. $$
The process $\smash{W_{\! \cdot} (\bS, \cdot )}$ is the \emph{snake associated with} $\bS$. We also introduce 
$\smash{\widehat{W}_{\! s} (\bS) \eqo  W_{\! s} (\bS, C_{\! s} (\fftree))}$, $s\ino \bbR_+$, and we call the $\bbR^d$-valued process 
$\smash{(\widehat{ W}_{\! s} (\bS))_{s\in \bbR_+^*}}\! $ the \emph{endpoint process}. \cq 
\end{compactenum}
\end{definition}
The snake $\smash{(W_{\! s}(\bS, \cdot))_{s\in \bbR_+}}\! $ associated with the BRW $\smash{(S_u)_{u\in \btt}}$ 
encodes the historical paths of the individuals of $\btt$ which are visited in contour order and which are interpolated in a convenient way to obtain a continuous process. 
Note that for all $s\ino \bbR_+$, $W_{\! s}(\bS, \cdot)$ belongs to $\mathbf C(\bbR_+, \bbR^d)$, the space of continuous functions from $\bbR_+$ to $\bbR^d$. Then 
 $W_\cdot (\bS, \cdot)$  is a \emph{snake with lifetime process} $C_\cdot (\btt)$ according to the following general definition 
(see Definition \ref{snadef} for more details): \emph{let $h\! :\! \bbR_+\! \to \!  \bbR_+$ be lower semicontinuous and such that $h(0)\eqo 0$; for all $s\ino \bbR_+$ let $w_s (\cdot) \ino \mathbf C(\bbR_+, \bbR^d)$; 
$w$ is a \emph{snake with lifetime process $h$} if}
\begin{equation}    
\label{snakedef}
\forall s,s',r\ino \bbR_+, \quad w_s(r)\! = \!  \widehat{w}_s , \; r\ino \big[ h(s), \infty\big) \quad \textrm{\emph{and}} \quad w_{s} (r)\eqo w_{s'} (r), \quad r\ino \big[ 0, m_h(s,s')\big], 
\end{equation}
\emph{where we have set $\smash{\widehat{w}_s \eqo w_s(h(s))}$ and $\smash{m_h(s,s')\eqo \min_{r\in [s\wedge s', s\vee s']} h(r)}$, which is a finite quantity since $h$ is lower semicontinuous. 
The $\smash{\bbR^d}$-valued function $\smash{(\widehat{w}_s)_{s\in \bbR_+}}$ is called the \emph{endpoint process of the snake $w$}.}

We shall also consider \emph{contour processes of trees with possibly variable lifespans}. Namely, let $T\eqo (\btt, (\ell_u)_{u\in \btt})$ be such a tree: here $\ell_u \ino \bbR_+^*$ is the lifespan of $u$ and we 
introduce $\zeta_u \eqo \sum_{v\in \lgeo \varnothing , u\rgeo} \ell_v$ and $\zeta^*_u \eqo \zeta_{\overleftarrow{u}}$, 
which are respectively the \emph{death-time} and the \emph{birth-time} of $u$, with the convention that 
$\zeta^*_{\varnothing}\eqo 0$. We denote by $w_1, \ldots, w_q$ the leaves of $\btt$ listed in increasing contour order.  
Then, the \emph{contour process of $T$} is the continuous piecewise affine function $s \ino\bbR_+ \! \mapsto \mathscr C_s(T)$ whose slopes are $\pm 1$, with $q$ local maxima and which takes the successive values $ 0$, $\zeta_{w_1}$, $\zeta_{w_1 \wedge w_2}$, $\zeta_{w_2}$, $\ldots$, $\zeta_{w_{q-1}\wedge w_q}$, $\zeta_{w_q}$, $0$. Note that $\mathscr C_\cdot (T)$ completely encodes $T$. We refer to Definition \ref{Contlifespan} $(a)$ for more details. 

\smallskip

The tree $\btt$ can be coded by two additional processes: namely, its \emph{height process} $\smash{(H_{\! s} (\btt))_{s\in \bbR_+}}$ and its \emph{\L{}ukasiewicz path} $\smash{(V_{\! s} (\btt))_{s\in \bbR_+}}\! $ that are defined as follows.

\begin{definition}
\label{heightlukadef} We assume that $\btt$ is a finite ordered rooted tree. 

\begin{compactenum}

\smallskip

\item[$(a)$] The \emph{depth-first exploration of $\btt$} $(u_j)_{0\leq j<\# \btt}$ is the contour exploration without repetition: namely, $u_j\! :=\! v_{m_j}$ where $m_j\eqo \inf \big\{ k\ino \bbN: \# \{ v_0, v_1, \ldots, v_k \} \eqo j+1\big\} $. 
It is convenient to set $u_{\# \btt} \eqo \varnothing$. 

\smallskip

\item[$(b)$] For all $0\leqo l\leko \# \btt$, we set $V_{0}(\btt)\eqo 0$ and $V_{l+1} (\btt) \! = \! V_{l}(\btt)  + k_{u_{l}} (\btt)\! -\! 1 $. We extend $V_\cdot (\btt)$ to $\bbR_+$ by setting $V_s (\btt)\eqo V_{\lfloor s \rfloor} (\btt)$ if $s\ino [0, \# \btt +1)$ and $V_{s} (\btt) \eqo -1$ for all $s\ino [\# \btt +1 , \infty)$. The process $V_\cdot (\btt)$ is the \emph{\L{}ukasiewicz path} of $\btt$.  

\smallskip

\item[$(c)$] For all $0\leqo l\leqo \# \btt$, we set $H_l(\btt) \eqo |u_l|$, the height of $u_l$. We extend $H(\btt)$ continuously on $[0, \#\btt]$ as follows: for all integers $0\leqo l\leko \#\btt $ and all $s\ino [0,1]$, we set 
\end{compactenum}
\begin{equation}
\label{ctrvshght}
\!\! H_{l+s}(\btt)  \eqo \left\{  \begin{array}{ll}
\!\! H_l (\btt)\! - \! \big(H_l (\btt) \! -\! |u_l \! \wedge \! u_{l+1}| \big) ((2s)\! \wedge \! 1) \!  + (2s\! -\! 1 )_+ \un_{\{ l\neq \# \btt-1\}}  \!\! & \!\!  \textrm{if $|u_{l+1}| \! \leq \! |u_l |$ } \\
\!\! H_l (\btt) + s   \!\! &\!\!  \textrm{if  $|u_{l+1}| \geko | u_l | $.}
\end{array} \right.
\end{equation}
\begin{compactenum}
\item[]$H_\cdot (\btt)$ is the \emph{height process} of $\btt$. \cq 
\end{compactenum}
\end{definition}
We observe that 
$\smash{s\ino \bbR_+\! \mapsto \! H_s (\btt)}$ is a continuous function, which fully encodes $\btt$. Moreover, there is an increasing bijection $\smash{\phi_{\btt} \! :\! \bbR_+ \! \to \! \bbR_+}$ such that $\smash{H_{s} (\btt)\eqo C_{\! \phi_{\btt} (s)} (\btt)}$ for all $\smash{s\ino\bbR_+}$. We refer to Remark \ref{contord} $(b)$ for more details. 

Similarly, we define the \emph{height snake} by $\smash{\cW_s (\bS, \cdot)\eqo W_{\! \phi_{\btt} (s)} (\bS, \cdot)}$, $\smash{s\ino \bbR_+}$, which is a snake with lifetime process $\smash{H_\cdot (\btt)}$. Here, $\smash{\cW_{j} (\bS, \cdot)}$ is simply the path of the ancestral line of $u_j$.  
 
\smallskip

\noi  
\textbf{$\;$ BRWs indexed by forests.} BRWs indexed by forests are finite or infinite sequences of BRWs indexed by trees as previously introduced. We denote them 
$\smash{\bS \eqo \big( (S_u(p))_{u\in \btt (p)} ; 1\leqo p \leqo N\big)}$, with $\smash{N\ino \bbN^*\! \cup \! \{ \infty\}}$, and with an obvious convention if $N\eqo \infty$. The jumps of $\bS$ are meant to be the sequences of the jumps of the $\bS (p)$: namely, $\smash{\big( (\xi_u(p))_{u\in \btt (p)\backslash \{ \varnothing\} } ; 1\leqo p \leqo N \big)}$. It is convenient to use the notation $\btt$ instead of $\smash{(\btt (p) ; 1\leqo p \leqo N)}$ and to use the notation $\smash{(S_u)_{u\in \btt}}$ instead of $\smash{\big( (S_u(p))_{u\in \btt (p)} ; 1\leqo p \leqo N\big)}$ 
(see Definition \ref{forestdef} and Remark \ref{contord} for more details). 

\emph{The processes encoding $\bS$ are obtained by concatenating the processes that encode the $\bS(p)$}. Namely, 
for all $1\leqo p \leqo N$, we set 
$\smash{\mathtt{r}_p\eqo \sum_{1\leq q\leq p} \# \btt (q)} $ and $\smash{\mathtt{r}_{0} \eqo 0}$. Then for all $\smash{s\ino [\mathtt{r}_{p-1}, \mathtt{r}_p]}$, we set 
$$ C_{\!2s}(\btt)\eqo C_{\!2(s-\mathtt{r}_{p-1})} (\btt (p)) \quad \textrm{and} \quad W_{\! 2s}(\bS, \cdot)\eqo W_{\! 2(s-\mathtt{r}_{p-1})} (\bS (p), \cdot).$$ 
We also set $\smash{H_{\! s}(\btt)\eqo H_{\!s-\mathtt{r}_{p-1}} \! (\btt (p))}$, 
$\smash{\cW_{\! s}(\bS, \cdot)\eqo \cW_{\! s-\mathtt{r}_{p-1}}\!  (\bS (p), \cdot)}$ and $\smash{V_{\! s}(\btt)\eqo V_{\!s-\mathtt{r}_{p-1}} \!(\btt (p))\! -\! p+1}$.  
  
\smallskip

Similarly, let $T$ $\eqo$ $\smash{\big( T(p)\eqo \big( \btt(p), (\ell_{u} (p))_{ u\in \btt (p) } \big)\, ; 1\leqo p \leqo N \big)}$ be a forest of trees with variable lifespans. The contour process of the foret $T$, which is denoted by $\smash{(\mathscr C_s(T))_{s\in \bbR_+}}$, is also defined as the concatenation of the contour processes 
$\smash{(\mathscr C_s(T(p)))_{s\in \bbR_+}}$. Namely,  for all $1\leqo p \leqo N$ we set $\smash{\mathtt{r}'_p\eqo \sum_{1\leq q\leq p} \sum_{u\in \btt(q)} \ell_u(q)}$, $\smash{\mathtt{r}'_{0} \eqo 0}$ and 
$\smash{\mathscr C_{\!s}(T)\eqo \mathscr C_{\!s-2\mathtt{r}'_{p-1}} (T(p))}$, for all $\smash{s\ino [2\mathtt{r}'_{p-1}, 2\mathtt{r}'_p]}$. See Definition \ref{lifespanforest} for more details.  
  
  \smallskip

\noi
\textbf{$\;$ Rescalings.} As previously mentionned we study  limits of \emph{sequences of BRWs} $\smash{\bS_n \eqo (S_{n,u})_{u\in \btt_n}}$, $n\ino \bbN$, whose 
genealogies, when suitably rescaled, converge to continuum random trees and whose spatial motions are, roughly speaking, in the domain of attraction of Gaussian laws. More precisely, we fix two renormalization sequences $\smash{\baa\eqo (a_n)_{n\in \bbN}}$ and $\smash{\bbb\eqo (b_n)_{n\in \bbN}}$ which satisfy \texttt{Norm} ($\baa, \bbb$) as in (\ref{renormab}).
We then rescale the encoding processes as follows: for all $n\ino \bbN$ and all $\smash{s, r\ino \bbR_+}$, 
\begin{equation}
\label{renormCW} 
C^{_{(n)}}_{s} \!\! \eqo \tfrac{1}{\lambda_n} C_{b_ns} (\btt_n) , \;  W^{_{(n)}}_{s} \! (r) \! \eqo  \tfrac{1}{\sqrt{\lambda_n}} 
W_{\! b_ns} (\bS_n, \lambda_n r), \;   H^{_{(n)}}_{s}\!\!   \eqo  \tfrac{1}{\lambda_n} H_{b_ns} (\btt_n), \; V^{_{(n)}}_{s} \!\! 
\eqo \tfrac{1}{a_n} V_{b_ns} (\btt_n).
\end{equation}  
 Under specific assumptions on the indexing trees and on the jumps of the BRWs that are discussed further, we prove 
the convergence of $\smash{ \big(C^{_{(n)}}_{\cdot }, W^{_{(n)}}_{\cdot })}$ in law in 
$\smash{\bC^{_0}_{^1}\times \mathbf C (\bbR_+, \bC^{_0}_{^d})}$ equipped with the product topology. 
Here, we use the shorthand 
$\smash{\bC^{_0}_{^d}}$ for $\smash{\mathbf C (\bbR_+, \bbR^d)}$ and for any Polish metric space $\smash{(E, d_E)}$, $\smash{\mathbf C (\bbR_+ , E)}$ stands for the space of continuous functions from $\smash{\bbR_+}$ to $E$: it is also Polish when equipped with the distance $\smash{\delta (f,g)\eqo \sum_{p\geq 1} 2^{-p} \big( 1\wedge \max_{s\in [0, p]} d_E (f(s),g(s))\big)}$, which induces the topology of uniform convergence on every compact intervals. 
 
 \smallskip
 
\noi
\textbf{$\;$ Galton-Watson forests.} We mostly focus on cases where the indexing trees $\btt_n$ are Galton-Watson trees (or variants of such trees) that converge to L{\'e}vy trees. Before stating our results we need to recall the main limit-theorems on rescaled Galton-Watson trees on which we rely. 
To that end we fix a sequence $\smash{\bmu\eqo (\mu_n)_{n\in \bbN}}$ of laws on $\bbN$ (the \emph{offspring distributions}) that are \emph{non trivial} and \emph{critical}: 
\begin{equation}
\label{nontricri}
\forall n\ino \bbN, \quad \mu_n(0)+\mu_n(1) <1 \quad \textrm{and} \quad \sum_{k\in \bbN} k\mu_n (k)= 1\; .
\end{equation}
For all $n\ino \bbN$, we denote by $\smash{\tau_{n}^{_{\infty}}}$ a sequence $\smash{(\tau_n(p))_{p\in \bbN}}$ of independent 
Galton-Watson trees with offspring distribution $\mu_n$ (GW($\mu_n$)-trees for short). Its \L{}ukasiewicz path 
$\smash{(V_k(\tau_{n}^{_{\infty}}))_{k\in \bbN}}$ is then a RW whose jump distribution $\nu_n$ is 
left-skipfree: $\smash{\nu_n(k)\eqo \mu_n (k+1)}$, $\smash{k\ino \bbN\cup \{ -1\}}$ (see e.g.~Le Gall \& Le Jan \cite{LGLJ98}). 
By (\ref{nontricri}), it is also a critical RW in the sense that $\smash{\sum_{k\in \bbZ} k\nu_n(k)\eqo 0}$. 
 
 \smallskip

\noi
\textbf{$\;$ Scaling limits of \L{}ukasiewicz paths.} The time- and space-continuous analogues of critical left-skipfree 
RWs are the 
\emph{critical, spectrally positive (i.e., without negative jump) L{\' e}vy processes} $\smash{X\eqo}$ $\smash{ (X_s)_{s\in \bbR_+}}$. Here critical means that for all 
$\smash{s\ino \bbR_+}$, $\smash{X_s}$ is integrable and $\smash{\bE [X_s] \eqo 0}$. The law of such a process $X$ is characterised by 
its \emph{Laplace exponent} $\smash{\psi (\lambda)\eqo \log \bE \big[ \exp (-\lambda X_1)\big]}$, $\smash{\lambda \ino \bbR_+}$, 
which necessarily belongs to $\mathscr L$, a class of functions which is defined as follows. 
\begin{definition}
\label{defLaplclass} We denote by $\ccL$ the space of critical \emph{L{\'e}vy-Khintchine} Laplace exponents $\psi$, which are of the form:
\begin{equation}
\label{LK}
\forall z\ino \bbR_+, \quad \psi (z) \eqo \beta_\psi z^2 + \int_{\bbR_+^*} \!\!\!  (e^{-zx} \! -\! 1+ zx)\,  \pi_\psi (dx).
\end{equation}
Here $\smash{\beta_\psi \ino \bbR_+}$ is the \emph{Brownian component} and the \emph{L{\'e}vy measure} $\smash{\pi_\psi} $ satisfies 
the following integrability condition: 
$\smash{\int_{\bbR^*_+}  ( x \wedge x^2) \pi_\psi (dx) \leko\infty}$. \cq 
\end{definition}

\noi
Most of the times $X$ is supposed to have infinite variation sample paths, which is equivalent to assume (see e.g.~Bertoin \cite{Be}) that 
\begin{equation}
\label{varinf}
\tag{{\small \texttt{Var}$_{\infty}$($\psi$)}} \textrm{either} \quad \beta_\psi >0 \quad \textrm{or} \quad  \int_{(0, 1)} \!\!\!\!  x\pi_\psi (dx) =\infty . 
\end{equation}

Our main assumption is the convergence of the rescaled \L{}ukasiewicz paths: namely, we assume that there are two renormalization sequences $\smash{\baa\eqo (a_n)_{n\in \bbN}}$ and $\smash{\bbb\eqo (b_n)_{n\in \bbN}}$ that tend to $\infty$, such that the following convergence holds in law in $\bbR$: 
\begin{equation*}
\label{Luka} 
\tag{{\small \texttt{\L{}uka}$_{\,}$($\baa, \bbb, \bmu,\psi$)}} \tfrac{1}{a_n} V_{\lfloor b_n \rfloor } \big( \tau_{n}^{_\infty}\big)    \xrightarrow[n\to \infty]{\; } X_1  .
\end{equation*}  
By standard results, it is equivalent to the convergence $\smash{V^{_{(n)}}_\cdot \!\!  \to \! X}$ in law in $\smash{\bD(\bbR_+, \bbR)}$,  
the space of c\`adl\`ag functions equipped with the Skorokhod topology (see Jacod \& Shiryaev \cite{JaSh02}, Corollary 3.6, Chapter VII, Section 3.a, p.~415). 
\texttt{\L{}uka}$_{^{\,}}$($\smash{\baa, \bbb, \bmu,\psi}$) translates into analytic conditions on $a_n, b_n$ and $\mu_n$: see the proof of Lemma \ref{Lukagrowth} for more details. 
Let us mention that Lemma \ref{Lukagrowth} $(iii)$ asserts that \texttt{Var}$_{\infty}$($\psi$) and \texttt{\L{}uka}$_{^{\,}}$($\baa, \bbb, \bmu,\psi$) imply $\smash{\lim_{n\to \infty} \lambda_n \eqo \infty}$, i.e., \texttt{Norm}$_{^{\,}}$($\baa, \bbb$).

 \smallskip
 
 \noi
 \textbf{$\;$ Scaling limits of GW-Markov chains.} \texttt{Var}$_{\infty}$($\psi$) and \texttt{\L{}uka}$_{^{\,}}$($\smash{\baa, \bbb, \bmu,\psi}$) (implying \texttt{Norm}$_{^{\,}}$($\baa, \bbb$)) entail the convergence of the corresponding Galton-Watson Markov chains as proved by Grimvall in \cite{Gr74}, Theorem 3.1 p.~1030 and Theorem 3.4 p.~1040. 
More precisely, we fix $\smash{x\ino \bbR_+^*}$, for all $n\ino \bbN$ we set $\smash{\tau_n^{{x}}\eqo (\tau_n(p))_{1\leq p \leq \lfloor a_n x\rfloor} }$ and we denote by $\smash{Z^{_{(n)}}_{^k} (x)}$ the number of individuals at height $k$ in $\smash{\tau_n^{{x}}}$: 
$\smash{(Z^{_{(n)}}_{^k} (x))_{k\in \bbN}}$ is thus a GW-Markov chain with offspring distribution $\mu_n$ and initial value 
$\smash{\lfloor a_n x \rfloor}$. Then under \texttt{Var}$_{\infty}$($\psi$) and \texttt{\L{}uka}$_{^{\,}}$($\smash{\baa, \bbb, \bmu,\psi }$), the following convergence holds in law in $\smash{\bD(\bbR_+, \bbR)}$: 
 $$ \Big(\tfrac{1}{a_n} Z^{_{(n)}}_{{\lfloor \lambda_n s \rfloor}} (x)\Big)_{s\in \bbR_+} \xrightarrow[n\to \infty]{\; }
 \big(Z^{{x}}_{s} \big)_{s\in \bbR_+}, $$
where $\smash{Z^{{x}}_{\cdot}}$ is a \emph{continuous state branching process} with \emph{branching mechanism} $\psi$ (a CSBP($\psi$) for short) with initial value $\smash{Z^{{x}}_{0}\eqo x}$. Namely, $\smash{Z^x_\cdot }$ is a $\smash{\bbR_+}$-valued Markov process whose transitions probabilities are characterised by the following: for all $\smash{s,s', \lambda \ino \bbR_+}$, $\smash{\bE [e^{-\lambda Z^x_{{s+s'}} } | Z^x_{{s'}}] } $ $\eqo$  $\smash{e^{-u_{s} (\lambda)Z^x_{s'} }}$ a.s.~where $\smash{u_{s} (\lambda)}$ is the unique solution to the equation $\smash{s\eqo \int_{u_s(\lambda)}^{\lambda} \frac{\mathrm d z }{ \psi (z)}}$. Since $\smash{\psi'(0^+)\eqo 0}$, $\smash{(Z^x_s)_{s\in \bbR_+}}$ is necessarily a martingale and $\smash{\lim_{s\to \infty} Z^x_s\eqo 0}$. Note that $0$ is an absorbing state for CSBPs. We recall here \emph{Grey's criterion} (see Grey \cite{Grey}) on extinction in finite time, which plays a role in what follows. 
\begin{equation*}
\label{Grey} 
\tag{{\small \texttt{Grey}$_{^{\,}}$($\psi$)}}  \bP \big( \exists s\ino \bbR_+ : Z^x_s\eqo  0 \big) \eqo 1 \; \Longleftrightarrow \int_{1}^\infty
 \frac{\mathrm d z}{\psi (z)} <\infty .
\end{equation*} 
Note that  \texttt{Grey}($\psi$) implies  \texttt{Var}$_{\infty}$($\psi$). We refer e.g.~to Bingham \cite{Bi76} for more details on CSBPs. 

\smallskip

\noi
\textbf{$\;$ Scaling limits of height and contour processes.} Under 
\texttt{Var}$_{\infty}$($\psi$)
and \texttt{\L{}uka}$_{^{\,}}$($\baa, \bbb, \bmu,\psi$), Theorem 2.2.1 in D.\& Le Gall \cite{DuLG02} asserts that finite dimensional marginal laws of the height processes $\smash{H^{_{(n)}}_{\cdot}}$ (rescaled as in (\ref{renormCW})) converge in law to those of a limiting process $\smash{(H_s)_{s\in \bbR_+}}$, which is called the \emph{$\psi$-height process} and which is derived from $X$ as follows: for all $\smash{s\ino \bbR_+}$, the following limit holds in probability 
\begin{equation}
\label{defH}
H_s= \lim_{\epp \to 0^+} \frac{1}{\epp} \int_0^s \!\! \un_{\{ X_r -\inf_{u\in [r,s]} X_u < \epp\}} \, \mathrm d r\; , 
\end{equation}
(see Section 1.2 in D.~\& Le Gall \cite{DuLG02}). The $\psi$-height process always admits a lower semicontinuous modification which enjoys the intermediate values property. Furthermore, it admits a continuous modification iff $\psi$ satisfies \texttt{Grey}$_{^{\,}}$($\psi$) (see Theorem 1.4.3 in D.~\& Le Gall \cite{DuLG02}). In this case, we get the following functional convergence of $\smash{H^{_{(n)}}_{\cdot}}$ to $H$: 
we suppose that $\psi$ satisfies \texttt{Grey}$_{^{\,}}$($\psi$) and we introduce the following assumption.
\begin{equation}
\label{Hghtabmu}
\textrm{{\small \texttt{Hght}$_{^{\,}}$($\baa, \bbb, \bmu$)}} \qquad \quad 
\exists x\ino \bbR_+^*\; \, \textrm{such that for all $\delta\ino \bbR^*_+$,} \quad  \liminf_{n\to \infty}\bP\big( Z^{_{(n)}}_{\lfloor \lambda_n \delta \rfloor} (x)\eqo 0 \big) \geko 0 \; .
\end{equation} 
Then \texttt{\L{}uka}$_{^{\,}}$($\baa, \bbb, \bmu,\psi$), \texttt{Hght}$_{^{\,}}$($\baa, \bbb, \bmu$) and \texttt{Grey}$_{^{\,}}$($\psi$) are equivalent to 
the convergence 
\[ ( C^{_{(n)}}_{\cdot}\! ,H^{_{(n)}}_{\cdot}\! ,V^{_{(n)}}_{\cdot})\!  \longrightarrow \! \big( (H_{s/2})_{s\in \bbR_+} ,  (H_{s})_{s\in \bbR_+} ,  (X_{s})_{s\in \bbR_+} \big)\]
in law in $\smash{(\bC^{_0}_{^1})^2 \! \times \! \bD (\bbR_+, \bbR)}$ (see Theorem 2.3.1 and Corollary 2.5.1 in D.~\& Le Gall \cite{DuLG02}). The quantity $x$ does not play a significant role in  \texttt{Hght}$_{^{\,}}$($\baa, \bbb, \bmu,$), which is actually equivalent to $\smash{\limsup_{n\to \infty} a_n \big(1\! -\! g_{{\mu_n}}^{_{\circ \lfloor \delta \lambda_n \rfloor}} (0) \big) \leko \infty}$ for all $\delta\ino \bbR^*_+$, where $\smash{g^{\circ m}_{\mu_n}}$ is the 
$m$-th iteration of the generating function of $\mu_n$ given by 
\begin{equation}
\label{genemuenn}g_{\mu_n} (r) \eqo \sum_{k\in \bbN} r^k \mu_n (k), \quad r\ino [0, 1]\; .
\end{equation} This condition is not easy to check directly in general. 
We mention here from Lemma 7.4 in Broutin, D.~\& Wang \cite{BrDuWa21} a more pratical condition involving the discrete branching mechanism $\psi_n \! :\! [0, a_n] \! \to \! \bbR_+$, which is defined by 
\begin{equation}
\label{psidisdef}
\forall \lambda \ino [0, a_n], \quad \psi_n (\lambda)\eqo b_n \big( g_{\mu_n} \! \big(1\! -\! \tfrac{\lambda}{a_n} \big) \! -\! 1+ 
\tfrac{\lambda}{a_n} \big) .
\end{equation}
Namely, \texttt{Var}$_{\infty}$($\psi$) and  \texttt{\L{}uka}$_{^{\,}}$($\baa, \bbb, \bmu, \psi$) imply that $\psi_n\!  \to \! \psi$ in $\bC^{_0}_{^1}$ (see Lemma \ref{Lukagrowth} ($iv$)). 
Furthermore, under these assumptions, the following condition
\begin{equation*}
\label{Greydis} 
\tag{{\small \texttt{Grey}$_{^{\,}}$($\baa, \bbb, \bmu$)}}. \!\!\!\!\!\!\!\!\!\!  \!\!\!\!\!\!\!\!\!\!  \!\!\!\!   \lim_{y\to \infty} \limsup_{n\to \infty}\int_{y}^{a_n}\!\!\!\! 
 \frac{ds}{\psi_n (s)} =0 , 
\end{equation*}  
is equivalent to  \texttt{Grey}$_{^{\,}}$($\psi$), \texttt{Hght}$_{^{\,}}$($\baa, \bbb, \bmu$) and $\lim_{n\to\infty} b_n / (a_n \log a_n )\eqo \infty$. We refer to 
Theorem \ref{Greyexplain} for a more detailed statement. 

\smallskip

\noi
\textbf{$\; $ Brownian snakes.} 
We next introduce the limiting processes 
for discrete snakes. Namely, \emph{Brownian snakes}, which are defined as follows: \emph{a snake $\smash{(W_{\! s} (\cdot))_{s\in \bbR_+}}$ with (possibly random) lifetime process $h$, as defined in (\ref{snakedef}), is a $\smash{\bbR^d}$-valued 
Brownian snake if conditionally given $h$, for all real numbers $s'\geqo s\geqo 0$},  
\begin{align}
(a): \; & \textrm{\emph{$W_{\! \cdot \wedge s}$ is independent of 
$W_{\! s'}  \big( \! \cdot + m_h(s,s')\big)\! -\! W_{\! s'}   \big(  m_h(s,s')\big)$;}} \nonumber \\
(b):\;  & \textrm{\emph{$W_{\! s'}  \big( \! \cdot + m_h(s,s')\big)\! -\! W_{\! s'}   \big(  m_h(s,s')\big)$ is distributed as $B_{\cdot\,  \wedge ( h(s')-m_h(s,s'))}$,}}
\label{Brosnadefbis} 
\end{align}
\emph{where here $\smash{(B_r)_{r\in \bbR_+}}$ 
stands for a $d$-dimensional standard Brownian motion starting at the origin.} 
Note that (\ref{Brosnadefbis}) characterizes the finite dimensional marginal laws of $W$ (see Definition \ref{Brosnadef}).

 Even if $h$ is continuous, there is in general no continuous version of either 
$W$ or $\smash{\widehat{W}}$. We recall here Theorem 4.5.2 from D.~\& Le Gall \cite{DuLG02} which asserts that 
there is a continuous version of a $\bbR^d$-valued Brownian snake whose lifetime process is a $\psi$-height process $H$ iff $\psi$ satisfies \emph{Sheu's condition} on compactness of the range of $\psi$-super-Brownian motion (see Sheu \cite{Sh1} or Hesse \& Kyprianou \cite{HesKyp14}): 
\begin{equation*}
\label{Sheu} 
\tag{{\small \texttt{Sheu}$_{^{\,}}$($\psi$)}}   \int_{1}^\infty\!\!\!\! 
 \frac{\mathrm d s}{\sqrt{\! \int_0^s\! \psi (r) \mathrm d r\, }} <\infty .
\end{equation*}  
Note that  \texttt{Sheu}$_{^{\,}}$($\psi$) implies \texttt{Grey}$_{^{\,}}$($\psi$) but e.g.~there are Laplace exponents 
$\psi$ of the form (\ref{LK}) such that $\smash{\psi(\lambda)\!  \asymp_\infty \! \lambda (\log \lambda)^2}$, which therefore satisfy \texttt{Grey}$_{^{\,}}$ ($\psi$) \emph{but not} \texttt{Sheu}$_{^{\,}}$($\psi$). Here $\smash{f(\lambda)\!  \asymp_\infty \! g(\lambda)}$ means that there are $\smash{\lambda_0, c\ino [1,\infty)}$ such that $\smash{0\leqo c^{-1} f(\lambda) \leqo g(\lambda) \leqo cf(\lambda) }$ for all $\lambda\ino [\lambda_0, \infty)$.

\subsection{The four cases considered in this article.} 
\label{4casessec}

We fix two renormalization sequences $\smash{\baa\eqo (a_n)_{n\in \bbN}}$ and $\smash{\bbb\eqo (b_n)_{n\in \bbN}}$, a sequence of non-trivial critical offspring distributions $\smash{\bmu\eqo (\mu_n)_{n\in \bbN}}$ and $\psi$, a function of the form (\ref{LK}) which satisfies at least  \texttt{Var}$_{\infty}$($\psi$). We want to obtain scaling limits of sequences of BRWs $\smash{\bS_n \eqo (S_{n,u})_{u\in \btt_n}}$, $n\ino \bbN$, via their encoding processes $\smash{C^{_{(n)}}_\cdot}$,  $\smash{H^{_{(n)}}_\cdot}$,  $\smash{V^{_{(n)}}_\cdot}$ and $\smash{W^{_{(n)}}_\cdot}$, which are normalized as in (\ref{renormCW}). More precisely in this article we investigate joint  convergences of $\smash{(C^{_{(n)}}_\cdot, V^{_{(n)}}_\cdot, W^{_{(n)}}_\cdot)}$ to limiting processes $\smash{(C,Y,W)}$ in the four following cases. 

\smallskip

\noi
$\bullet$ $\textbf{Case (0)}$.
In this case, we take $\smash{\btt_n\eqo \tau^{_\infty}_n\eqo (\tau_n (p))_{p\in \bbN^*}}$ as a sequence of independent GW($\mu_n$)-trees. In the case where $\smash{\btt_n}$ is equipped with variable lifespans, we denote 
the resulting $\smash{\bbR_+}$-marked tree by  
$\smash{\mathcal T_n\eqo \big( \mathcal T_{n} (p) \eqo (\tau_n (p), (\ell_{n,u} (p) )_{u\in \tau_n(p)}) \big)_{\! p\in \bbN^*}}$ and we assume that \emph{conditionally given the $\smash{\tau_n (p)}$, the lifespans $\smash{\ell_{n,u} (p)}$ are independent exponential r.v.s with parameter $1$. }
  
We denote by $X$ a spectrally positive L{\'e}vy process with Laplace exponent $\psi$, we denote by $H$ its associated height process (which is always chosen to be continuous if \texttt{Grey}$_{^{\,}}$($\psi$) holds) and we also denote by $\mathbf{W}$ a $\smash{\bbR^d}$-valued Brownian snake with lifetime process $H$ (which is always chosen to be continuous if \texttt{Sheu}$_{^{\,}}$($\psi$) is satisfied). Then, $(C,Y,W)$ is defined by $\smash{C_{\! 2s} \eqo H_{\! s}}$, $\smash{Y_{\! s} \eqo X_{\! s}}$ and $\smash{W_{\! 2s} (\cdot) \eqo \mathbf{W}_{\! s} (\cdot) }$ for all $\smash{s\ino \bbR_+}$. \cq

\smallskip

\noi
$\bullet$ $\textbf{Case (1)}$.
In this case we take $\smash{\btt_n\eqo (\tau_n(p))_{1\leq p\leq \lfloor a_n x\rfloor}}$, where $\smash{x\ino \bbR_+^*}$ (with the same notations as above). 
As in $\textbf{Case (0)}$ if $\smash{\btt_n}$ is equipped with variable lifespans then 
the resulting marked tree is denoted by $\smash{\mathcal T^{x}_n\eqo \big( \mathcal T_{n} (p) \eqo (\tau_n (p), (\ell_{n,u} (p) )_{u\in \tau_n (p)} ) ;  1\leqo p \leqo \lfloor a_n x \rfloor \big) }$ and we assume that \emph{conditionally given the $\smash{\tau_n (p)}$, the lifespans $\smash{\ell_{n,u} (p)}$ are independent exponential r.v.s with parameter $1$.}
Then for all $\smash{s\ino \bbR_+}$, we set $\smash{C_{\! 2s} \eqo H_{\! s\wedge \varsigma_{-x}}}$, $\smash{Y_{\! s} \eqo X_{\! s\wedge \varsigma_{- x}}}$ and $\smash{W_{\! 2s} (\cdot) \eqo \mathbf{W}_{ \! s\wedge \varsigma_{- x}} (\cdot)}$ where 
\begin{equation}
\label{infireach}
\varsigma_{-x}\eqo \inf \big\{ s\ino \bbR_+: X_s\leko  -x\big\} \; .
\end{equation}   
Here we recall from e.g.~Bertoin \cite{Be} that $\smash{(\varsigma_{-x})_{x\in \bbR_+}}$ is a subordinator whose Laplace exponent is $\smash{\psi^{-1}}$, the inverse of $\psi$.    \cq 
   
\smallskip

To discuss the next cases, we introduce the excursion measure of the height process. 
Here we assume that $X$ has infinite variation sample paths and for all $\smash{s\in \bbR_+}$, we set $\smash{I_s}$ 
$\eqo$  $\smash{\inf_{r\in [0, s]} X_r}$, which is the \emph{infimum process of $X$}. Basic results on fluctuation theory 
(see e.g.~Bertoin \cite{Be}, Chapters VI.1 and VII.1) entail that $X\! -\! I$ is a strong Markov process in $\smash{\bbR_+}$ 
and that $0$ is regular for 
$\bbR_+^*$ and recurrent with repect to this Markov process. Moreover $-I$ is a local time at $0$ for $X\! -\! I$ 
(see Theorem VII.1 \cite{Be}). We denote by $\bN$ the corresponding 
\emph{excursion measure} of $X\! -\! I$ above $0$. We denote by $\smash{(l_j, r_j)}$, $j\ino  \cI$, the excursion intervals of 
$X\! -\! I$ above $0$ and by $\smash{X^j = X_{(l_j + \cdot )\wedge r_j}-I_{l_j}}$, $j\ino \cI$, the corresponding excursions. 
Then the point measure $\smash{\sum_{j\in \cI} \delta_{(-I_{l_j}, X^j)}}$ is a \emph{Poisson point measure on $\smash{\bbR_+\! \times \! \bD(
\bbR_+, \bbR)}$ with intensity 
$dx \otimes \bN$}.

We next recall from D.~\& Le Gall \cite{DuLG02}, Section 1.2, that $H$ and $X\! -\! I$ share the same excursion intervals.  Namely, 
$\smash{\{ s\ino \bbR_+\! : \! H_{\! s} \geko 0\}\eqo \{ s\ino \bbR_+\! : \! X_{\! s} \geko I_{s}\}\eqo \bigcup_{j\in \cI} (l_j, r_j)}$. 
It allows to define 
$H$ under $\bN$ as an adapted function $H(X)$ of the excursion $X$. 
We then assume \texttt{Sheu} ($\psi$) and that $\mathbf{W}$ is a continuous version of $\bbR^d$-valued Brownian snake with lifetime process $H$. Elementary arguments imply that under $\bN$, there is also a continuous $d$-dimensional snake starting at $0$. 
More precisely, the processes $\smash{(H^j, \mathbf{W}^j)\! :=\! (H_{(l_j+s )\wedge r_j}, \mathbf{W}_{(l_j+s )\wedge r_j} (\cdot) )_{s\in \bbR_+}}$, 
$j\ino  \cI$, are 
the excursions of $(H, \mathbf{W})$ out of the null function and the point measure $\smash{\sum_{j\in \cI} \delta_{(-I_{l_j}, X^j, H^j \! ,\,  \mathbf{W}^j)}}$ is distributed as a Poisson point measure on $\bbR_+ \!\times \bD (\bbR_+, \bbR)\times \bC^{_0}_{^1}\! \times \! \bC(\bbR_+, \bC^{_0}_{^d})$ with intensity $dx  \otimes \bN (dXdHd\mathbf{W})$ (with a slight abuse of notation). We refer to  D.~\& Le Gall \cite{DuLG02}, Chapters 1 and 4, for more details. 

\smallskip

\noi
$\bullet$ $\textbf{Case (2)}.$ Here $\smash{c\ino \bbR_+^*}$ and $\smash{\btt_n}$ is distributed as a single 
GW($\mu_n$)-tree conditionned to have total height $\smash{\, \geqo \lambda_n c}$. If $\btt_n$ is equipped with 
variable lifespans, then we shall denote the resulting marked tree by 
$\smash{\mathcal T_n\eqo (\btt_n ,(\ell_{n,u} )_{u\in \btt_n})}$ and we assume that 
\emph{conditionally given $\smash{\btt_n}$, the lifespans $\smash{\ell_{n,u}}$ are independent exponential r.v.s with parameter $1$. }

We then define $(C,Y, W)$ as follows. We denote by $\smash{H^{\geq c}}$ the first excursion of $H$ above $0$ to hit $c$. More precisely, we introduce $\smash{\rho_c\eqo \inf \big\{ s\ino \bbR_+: H_{\! s}\eqo c\big\}}$ and 
\begin{equation}
\label{sigellerrc}
 \ell_c \eqo \sup \big\{ s \ino [0, \rho_c]: H_{\! s}\eqo 0\big\}\; \textrm{and} \; r_c \eqo \inf
\big\{ s \ino [ \rho_c, \infty): H_{\! s}\eqo 0\big\}, 
\end{equation}
which are a.s.~finite. Then, we set $\smash{H^{\geq c} \eqo H_{\! (\ell_c + \cdot )\wedge r_c}}$. We define $Y$ as  $\smash{X_{\! (\ell_c +\cdot ) \wedge r_c} \! -\! I_{ \ell_c}}$  which is the excursion of $X\! -\! I$ above $0$ straddling $\smash{\rho_c}$. We next define the processes $C$ and $W$ 
by setting $\smash{\smash{C_{\! 2s}\eqo H^{\geq c}_{\! s}}}$ and $\smash{W_{\! 2s} (\cdot) \eqo \mathbf{W}_{\! (\ell_c + s)\wedge r_c} (\cdot) }$ for all $\smash{s\ino \bbR_+}$. 
Standard arguments on Poisson point processes entail 
\begin{equation}
\label{selectexcu}
(Y_s,C_{2s}, W_{2s})_{s\in \bbR_+} \; \textrm{under $\bP$} \; \overset{\textrm{(law)}}{=} \; (X, H , \mathbf{W})  \; \textrm{under $\bN (\, \cdot \, | \, \max_{s\in \bbR_+} H_s \geko c)$.} 
\end{equation}
We recall that \texttt{Sheu} ($\psi$) implies \texttt{Grey} ($\psi$) which implies that $v(c) \! :=\! 
\bN (\max_{s\in \bbR_+} H_s \geko c) \ino \bbR_+^*$. Moreover we recall that $v$ satisfies the equation $\smash{\int_{v(c)}^\infty\frac{ \mathrm dz }{ \psi (z)} \eqo c}$.  
\cq

\smallskip

\noi
$\bullet$ $\textbf{Case (3)}$. Here $\smash{\psi (\lambda)\eqo \lambda^\alpha}$, $\smash{\lambda\ino \bbR_+}$, with $\alpha \ino (1, 2]$ (note that \texttt{Grey}$_{^{\,}}$($\psi$) and \texttt{Sheu}$_{^{\,}}$($\psi$) hold in this case). We refer to this case as to the \emph{$\alpha$-stable case}. We also assume that for all $n\ino \bbN$, $\smash{\mu_n \eqo \mu}$ and that $\mu$ is \emph{aperiodic}. It implies that $\smash{\bP (\# \tau\eqo n) \geko 0}$ for all sufficiently large $n$, where $\tau$ stands here for a GW($\mu$)-tree. We next assume that $\smash{\btt_n}$ is distributed as $\tau$ under $\smash{\bP (\, \cdot \, | \, \# \tau\eqo n)}$. We also take $\smash{b_n\eqo n}$, $\smash{n\ino \bbN^*}$. If $\smash{\btt_n}$ is equipped with variable lifespans, then we denote the resulting marked tree by $\smash{\mathcal T_n\eqo (\btt_n (\ell_{n,u} )_{u\in \btt_n})}$ and we assume that \emph{conditionally given $\smash{\btt_n}$, the lifespans $\smash{\ell_{n,u}}$ are independent exponential r.v.s with parameter $1$. }

The processes $(C,Y, W)$ are then derived from the \emph{normalized excursions} of $H$ and of $X\! -\! I$ above $0$, i.e., excursions conditioned to last one unit of time. As already mentioned, $\bN$-a.e.~$X$ and $H$ have the same lifetime $\zeta$: namely, $\bN$-a.e.~$ \zeta \leko \infty$, $\smash{X_s\eqo X_0\eqo H_0\! =\! H_{s}\! =\! 0}$ for all $s \ino [\zeta, \infty)$ and $\smash{H_s}$ and $\smash{X_s \! >\! 0}$ for all $s \! \in \! (0, \zeta)$. The scaling property of $X$ and $H$ in stable cases 
allows to define the \emph{normalized excursion measure} $\smash{\bN (\, \cdot \, | \, \zeta\eqo 1)}$, 
which admits e.g.~the following representation: we set $\smash{\bgg\eqo \sup \{ s\ino [0, 1]\! : \! H_{\! s} \eqo 0 \}}$, 
$\smash{\bdd \eqo \inf \{ s\ino [1, \infty)\! : \! H_{\! s} \eqo 0 \}}$, $\smash{\bzeta\eqo \bdd \! -\! \bgg}$ and 
$\smash{H^{\mathtt{nr}}_s \eqo \bzeta^{-\frac{\alpha-1}{\alpha}} H_{\! \bgg + \bzeta s}}$ for all $s\ino [0, 1]$. Then $\smash{H^{\mathtt{nr}}}$ under $\bP$ has the same law as $H$ under $\smash{\bN (\, \cdot \, | \, \zeta\eqo 1)}$. More generally we define 
$\smash{(Y,C, W)}$ by $\smash{C_{\! 2s}\eqo H^{\mathtt{nr}}_{\! s}}$, $\smash{W_{\! 2s} (r)\eqo \bzeta^{\frac{\alpha-1}{2\alpha}} \mathbf{W}_{\! \bgg + \bzeta s } \big( \bzeta^{-\frac{\alpha-1}{\alpha}} r\big)}$ and $\smash{Y_{\! s} \eqo \bzeta^{-\frac{1}{\alpha}} (X_{\! \bgg +\bzeta s }\! -\! I_{\! \bgg})}$. Then 
\begin{equation}
(Y_s,C_{2s}, W_{2s})_{s\in \bbR_+} \; \textrm{under $\bP$} \; \overset{\textrm{(law)}}{=} \; (X, H , \mathbf{W})  \; \textrm{under $\bN (\, \cdot \, | \, \zeta\eqo 1)$.} 
\end{equation}
We refer to D.~\& Le Gall \cite{DuLG02}, and D.~\cite{Du2} for more details. \cq

\smallskip

In all of these cases, $W$ is always a Brownian snake with lifetime process $C$.

\smallskip

\emph{The limit theorems on the trees $\btt_n$ on which we rely in this article can be summarised as follows}: we assume \texttt{Sheu}$_{^{\,}}$($\psi$) (and thus  \texttt{Grey}$_{^{\,}}$($\psi$)), \texttt{Norm}$_{^{\,}}$($\baa, \bbb$), \texttt{\L{}uka}$_{^{\,}}$($\baa, \bbb, \bmu, \psi$) and \texttt{Hght}$_{^{\,}}$($\baa, \bbb, \bmu$). 

\smallskip
\noi
$-$  Then in $\textbf{Case (i)}$, $\textbf i\ino \{ 0, 1,2\}$, 
Theorems 2.3.1, Corollary 2.5.1 and Proposition 2.5.2 in D.~\& Le Gall \cite{DuLG02} show that the following convergence holds in law in $\smash{(\bC^{_0}_{^1})^2\! \times\!  \bD(\bbR_+, \bbR)}$: 
\begin{equation}
\label{cvcodgen}
\big( C^{_{(n)}}_\cdot\! \! , H^{_{(n)}}_\cdot \! ,  V^{_{(n)}}_\cdot \big) \xrightarrow[n\to \infty]{\; } (C,H,  Y) .
\end{equation}
This convergence holds jointly in distribution with the convergence $\smash{\# \btt_n /b_n\! \to \! \varsigma_{-x}}$ in $\textbf{Case (1)}$ and with the convergence $\smash{\# \btt_n /b_n\! \to \! r_c\! -\! \ell_c}$ in $\textbf{Case (2)}$. 

\noi
$-$ In $\textbf{Case (3)}$
where $b_n\eqo n$, there exists a sequence $\smash{(a_n)_{n\in \bbN}}$ satisfying \texttt{\L{}uka}$_{^{\,}}$($\smash{\baa, \bbb, \bmu, \psi}$) iff  $\mu$ is in the domain of attraction of an $\alpha$-stable law. 
We set $\smash{\Psi_{\! \mu} (r)\eqo g_\mu (1\! -\! r) \! -\! (1\! -\! r)}$, $r\ino [0, 1]$ (namely, $\smash{\psi_n (\lambda) \eqo n \Psi_{\! \mu} (\lambda/a_n)}$, $\smash{\lambda\ino [0, a_n]}$). Then by standard results on stable laws, 
\begin{equation}
\label{genedomain} 
\exists (a_n)_{n\in \bbN} \;\,  \textrm{s.t.} \; \textrm{\texttt{\L{}uka}$_{^{\,}}$($\baa, \bbb, \bmu, \psi$)}\; \textrm{holds}  \Longleftrightarrow \Psi_{\! \mu} (r) \sim_{0^+} \! C_\alpha r^\alpha L(1/r) 
\end{equation}
where $L$ is a slowly varying function at $\infty$ and where $\smash{C_\alpha\eqo \frac{\alpha -1}{\Gamma (2-\alpha)} }$ if 
$\alpha \ino (1, 2)$ and $\smash{C_2\eqo 1}$ if $\alpha \eqo 2$. Moreover, if one of the two equivalent conditions 
in (\ref{genedomain}) are satisfied, then $\smash{a_n^\alpha \sim nL(a_n)}$. For more details see 
e.g.~Bingham, Goldies \& Teugels \cite{BiGoTe}, Sections 8.1 and 8.4. 

Consequently in $\textbf{Case (3)}$, \texttt{Norm}$_{^{\,}}$($\baa, \bbb$) is always met. Moreover we prove in 
Lemma \ref{controlgrey} that 
\texttt{Grey}$_{^{\,}}$($\baa, \bbb, \bmu$), and thus \texttt{Hght}$_{^{\,}}$($\baa, \bbb, \bmu$), are always satisfied too. Thus, 
Theorem 3.1 in D.~\cite{Du2} asserts that (\ref{cvcodgen}) holds true. 

\subsection{Statement of the main results}
\label{statementssec}

We first prove a limit theorem for specific cases and we extend it thanks to a coupling result. 

\smallskip

\noi
\textbf{$\; $Scaling limit of Brownian snakes indexed by GW-trees.} Our first result concerns a semi-discrete model: 
the $1$-dimensional Brownian snakes whose lifetime process is the contour process of a GW-forest with 
i.i.d.~exponential lifespans. 
More precisely, we fix two renormalization sequences $\smash{(a_n)_{n\in \bbN}}$, $\smash{(b_n)_{n\in \bbN}}$ that satisfy 
\texttt{Norm} ($\smash{\baa, \bbb}$) and we consider a sequence of finite trees or forests of finite trees $\smash{\btt_n}$, $\smash{n\ino \bbN^*}$. 
Each $\smash{\btt_n}$ is 
equipped with variable lifespans. We denote the resulting $\smash{\bbR_+}$-marked tree by $\smash{\mathcal T_n}$ and we assume 
that \emph{conditionally given $\smash{\btt_n}$, the lifespans are independent exponential r.v.s with parameter $1$.} We denote by 
$\smash{(\mathscr C_s(\mathcal T_n))_{s\in \bbR_+}}\! $ the contour process of $\smash{\mathcal T_n}$ and we denote by 
$\smash{(\mathscr W_{s} (\mathcal T_n, r))_{r, s\in \bbR_+}}\! $ a process whose conditional law given $\smash{\mathcal T_n}$ is that of a 
$1$-dimensional Brownian snake starting at $0$ and whose lifetime process is 
$\smash{(\mathscr C_{s} (\mathcal T_n))_{s\in \bbR_+}}\! $. We rescale $\smash{\mathscr C_\cdot (\mathcal T_n)}$ and $\smash{\mathscr W_{\cdot } (\mathcal T_n, \cdot)}$ as in (\ref{renormCW}) and 
to simplify notation we set: 
$$ \forall s, r \ino \bbR_+, \qquad \mathscr C_{s}^{_{(n)}} \eqo \tfrac{1}{\lambda_n}
 \mathscr C_{b_ns} (\mathcal T_n)\quad \textrm{and} \quad \mathscr W^{_{(n)}}_{s} (r)
= \tfrac{1}{\sqrt{\lambda_n}} \mathscr W_{b_n s} (\mathcal T_n, \lambda_n r) \; .$$
We recall Assumption \texttt{Sheu}$_{^{\,}}$($\smash{\baa, \bbb, \bmu}$) from (\ref{Sheudis}). We first prove the following scaling limit. 
\begin{theorem}
\label{Sheuexplain} We fix $x,c ,\beta \ino \bbR_+^*$. 
Let $\psi\ino \mathscr L$ satisfy \texttt{Var}$_{\infty}$($\psi$). Let $\smash{(a_n)_{n\in \bbN}}$, 
$\smash{(b_n)_{n\in \bbN}}$, $\smash{(\mu_n)_{n\in \bbN}}$, $\smash{\bS_n \eqo (S_{n, u})_{u\in \btt_n}}$, 
$\smash{(\mathcal T_n)_{n\in \bbN}}$, $n\ino \bbN$, and $(C,Y,W)$ be as in $\emph{\textbf{Case (i)}}$, $ \mathbf i \ino \{ 0,1,2\}$. 
Let $\smash{\mathscr C^{_{(n)}}_\cdot}$ and $\smash{\mathscr W^{_{(n)}}_\cdot}$ 
be as above. We also recall that $\smash{C^{_{(n)}}_{\cdot}}$, $\smash{V^{_{(n)}}_{\cdot}}$ and 
$\smash{W^{_{(n)}}_{\cdot}}$ are rescaled as in (\ref{renormCW}) and we recall from (\ref{psidisdef}) the definition of $\smash{\psi_n}$. 
\begin{compactenum}

\smallskip

\item[$(a)$] We assume \emph{\texttt{\L{}uka}$_{^{\,}}$($\baa, \bbb, \bmu, \psi$)} and 
\emph{\texttt{Grey}$_{^{\,}}$($\baa, \bbb, \bmu $)}.

\smallskip

\item[$(b)$] We assume that conditionally given $\smash{\btt_n}$, the jumps of $\smash{\bS_n}$ are $\bbR$-valued i.i.d.~r.v.s with law $\smash{\bgam (dx)\eqo \frac{1}{\sqrt{2\beta}}e^{-(|x|\sqrt{2})/\sqrt{\beta}} dx}$. 

\smallskip

\end{compactenum}
Then the following holds true. 
\begin{compactenum}

\smallskip

\item[$(i)$] \emph{\texttt{Sheu}$_{^{\,}}$($\baa, \bbb, \bmu$)} implies \emph{\texttt{Sheu}$_{^{\,}}$($\psi$)}
and $\smash{b_n/ (a_n (\log  a_n)^2)\! \to \! \infty}$. 

\smallskip

\item[$(ii)$] In $\emph{\textbf{Case (0)}}$, 
\emph{\texttt{Sheu}$_{^{\,}}$($\smash{\baa, \bbb, \bmu }$)} 
is \emph{equivalent} to the tightness of the laws of the endpoint processes $\smash{(\widehat{\mathscr W}^{_{(n)}}_\cdot)_{n\in \bbN}}$, in $\smash{\bC^{_0}_{^1}}$. Moreover, it implies that 
$\smash{\lim_{n\to \infty} ( \mathscr C^{_{(n)}}_\cdot \! , C^{_{(n)}}_{\cdot}\! ,V^{_{(n)}}_{\cdot} , \mathscr W^{_{(n)}}_{\cdot}  )\eqo ( C_\cdot , C_\cdot, Y_\cdot, W_\cdot )}$ in law in $\smash{(\bC^{_0}_{^1})^2 \! \times \! \bD (\bbR_+, \bbR)\times \mathbf C (\bbR_+, \bC^{_0}_{^1})}$.

\smallskip

\item[$(iii)$]  In $\emph{\textbf{Case (i)}}$, $\mathbf i \ino \{ 0,1,2\}$, 
\emph{\texttt{Sheu}$_{^{\,}}$($\smash{\baa, \bbb, \bmu }$)} is \emph{equivalent} to the tightness of the laws of the endpoint processes 
$\smash{(\widehat{W}^{_{(n)}}_\cdot)_{n\in \bbN}}$ in $\smash{\bC^{_0}_{^1}}$. Moreover, it implies that 
$\smash{\lim_{n\to \infty} ( C^{_{(n)}}_{\cdot}\! ,V^{_{(n)}}_{\cdot} , W^{_{(n)}}_{\cdot}  )}$ $\eqo$ $\smash{ ( C_\cdot, Y_\cdot, \sqrt{\beta} W_\cdot )}$ in law in $\smash{\bC^{_0}_{^1} \! \times \! \bD (\bbR_+, \bbR) \! \times \! \mathbf C (\bbR_+, \bC^{_0}_{^1})}$. 
\end{compactenum}
\end{theorem}
\noi
\textbf{Proof.} See Section \ref{ThmSheuexplainPfsec}. \cqfd

\begin{remark}
\label{heightstuff} As already mentioned, \texttt{\L{}uka}$_{^{\,}}$($\baa, \bbb, \bmu, \psi$) and  \texttt{Var}$_{\infty}$($\psi$) imply \texttt{Norm}$_{^{\,}}$($\baa, \bbb$) and we recall that \texttt{Grey}$_{^{\,}}$($\baa, \bbb, \bmu$) implies 
\texttt{Hght}$_{^{\,}}$($\baa, \bbb, \bmu$) and \texttt{Grey}$_{^{\,}}$($\psi$). \cq 
\end{remark}

\noi
$\; $\textbf{Extension to mildly dependent and inhomogeneous jumps.}  We extend the previous limit theorem 
to sequences of BRWs whose jumps satisfy the following weaker independence assumption. 
 \begin{definition} 
\label{sibinddef} 
Let $\smash{\bS\eqo (S_u)_{u\in \btt}}$ be a $\smash{\bbR^d}$-valued 
BRW whose jumps are $\smash{(\xi_u)_{u\in \btt \backslash \{ \varnothing\}}}$. We denote by 
$\smash{\mathtt{Lf} (\btt)\eqo \{ u\ino \btt: k_u(\btt)\eqo 0\}}$ the set of \emph{leaves} of $\btt$ and for all 
$\smash{u\ino \btt \backslash \mathtt{Lf} (\btt)}$ we set $\smash{\bxxi_u \eqo (\xi_{u\ast [j]} )_{1\leq j\leq k_u(\btt)} }$, 
where $\smash{u\ast [j]}$ stands for the $j$-th child of $u$. We then say that the jumps of $\bS$ are 
\emph{independent by group of siblings} if the following holds true. 
\begin{equation}
\label{sibinddef}
\hspace{1mm} \textrm{Conditionally given $\btt$, the $\smash{\bxxi_u}$, $\smash{u\ino  \btt \backslash \mathtt{Lf} (\btt)}$, 
are independent. \hspace{10mm} \cq } 
\end{equation}
\end{definition}
\begin{remark}
\label{siblindep} $(a)$ Of course independence by group of siblings includes the simpler cases where conditionally 
given $\btt$ the jumps $\smash{(\xi_u)_{u\in \btt \backslash \{ \varnothing\}}}$ are independent. 

\smallskip

\noi
$(b)$ Note that the conditional law of $\smash{\bxxi_u}$ may depend on $\btt$ and on $u$, and that the jumps of two siblings 
are not necessarily independent conditionally given $\btt$. However, if $\smash{u\ino \mathtt{Lf} (\btt)}$, then the r.v.s 
$\smash{\xi_v}$, $\smash{v\! \in\,  \rgeo \varnothing , u \rgeo}$ are independent conditionally given $\btt$ (with possibly distinct 
conditional laws). \cq 
\end{remark}

To control any inhomogeneities of the laws of the jumps, 
we use a coupling argument based on results due to Sakhanenko in \cite{Sak91} (Corollary 12 p.~81) and to $\smash{\textrm{Koml{\'o}s, Major \& Tusn{\'a}dy}}$ in \cite{KomMajTus76} (Theorem 4). This coupling argument on BRWs, which is of independent interest, is stated in Theorem \ref{brwcouplingth}, Section \ref{couplingsec}. It holds under a moment assumption that is expressed in terms of a \emph{moment gauge function}, i.e., a \emph{continuous and even function $G\!:\!  \bbR \! \to \! \bbR_+^*$, nondecreasing on $\bbR_+$ and such that for some $\kappa , x_0\ino \bbR_+^*$ the following holds: }
\begin{eqnarray}
\label{momgaudef}
\textrm{$x\ino [x_0, \infty)$} \!\! \! \!\!   & \mapsto&  \!\! \! \!\!  \textrm{\emph{$x^{-2 -\kappa} G(x)$ is nondecreasing}}  \\
& \textrm{and} & \textrm{\emph{$x\ino [x_0, \infty)\mapsto x^{-1} \log G(x)$ is nonincreasing.}} \nonumber
\end{eqnarray}
Such functions range from power functions $\smash{x\! \mapsto \! |x|^c}$, with $\smash{c\geko 2}$ to exponential ones $\smash{x\mapsto e^{c|x|}}$ with $\smash{c\ino \bbR_+^*}$. We shall assume that the jumps of $\bS$ have $G$-moments conditionally given $\btt$: 
\begin{equation}
\label{momjumdef}
\textrm{$\bP$-a.s.} \quad M(G, \bS)\! :=\! \max_{u\in \btt \backslash \{ \varnothing \} }\!\!  \bE \big[ G( \lvert \xi_u \rvert )  \big| \, \btt \big] <\infty \; . 
\end{equation}
Here $\smash{\lvert\,  \cdot \, \rvert}$ stands for the canonical Euclidean norm on $\smash{\bbR^d}$,$\langle \cdot, \cdot \rangle$ for the canonical scalar product and $\smash{(\mathtt{e}_i)_{1\leq i\leq d}}$ for the canonical basis. We shall also restrict ourselves 
to BRWs whose jumps are conditionally centered: 
\begin{equation}
\label{concendef}
\textrm{$\bP$-a.s.} \quad \forall u \ino \btt \backslash \{ \varnothing \}  \; : \quad \bE \big[ \xi_u  \big| \, \btt \big] =0 \; .
\end{equation}
Control on jumps also involves their conditional covariance matrix $\smash{\beta (\bS, u)\! :=\! (\beta_{i,j} (\bS,u))_{1\leq i,j\leq d}}$ given by 
\begin{equation}
\label{concovdef}
\beta_{i,j} (\bS,u)= \bE \big[\langle \mathtt{e}_i, \xi_u \rangle \langle \mathtt{e}_j, \xi_u \rangle   \big| \, \btt \big] , \quad 1\leqo i,j\leqo d, \; u\ino  \btt \backslash \{ \varnothing \} ,
\end{equation}
which is well-defined since $\smash{\lvert \xi_u \rvert}$ has a conditional moment of order $2$ by (\ref{momgaudef}) and (\ref{momjumdef}).

For a $d\! \times \! d$ matrix $A$ with real entries, we denote by 
$\smash{\lVert A\rVert_2 \eqo \max \{ |A.x|\, ; \, x\ino \bbR^d \! : \! |x|\eqo 1 \}}$ 
the operator norm associated with the Euclidean norm. We denote by $\mathrm{Sym}^{_+}_{^d}$, 
the space of symmetric spectrally nonnegative $d \times d$ matrices, equipped with the norm topology. 
For all $\smash{\beta\ino \mathrm{Sym}^{_+}_{^d}}$, 
we denote by 
$\smash{\sqrt{\beta} \ino \mathrm{Sym}^{_+}_{^d}}$ the square root given by 
$\smash{\sqrt{\beta} \! :=\!\! \sqrt{2\lVert \beta \rVert_2}\sum_{n\in \bbN} c_n  
\big( \mathtt{Id}\! -\! \frac{1}{{2\lVert \beta \rVert_2}} \beta)^n}$, the sum being absolutely convergent: here, $\smash{c_0\eqo 1}$ and 
$\smash{c_n\eqo -|\binom{1/2}{n}|}$, $\smash{n\ino \bbN^*}$, are the coefficients of the power serie expansion of $\smash{z\ino [0, 1]\! \mapsto \! \sqrt{1\! -\! z}}$. In this article, we do not use any continuity property of such a square root, but only its measurability.

  We now fix $\smash{\bS_n\eqo (S_{n,u})_{u\in \btt_n}}$, $\smash{n\ino \bbN}$, a sequence of $\bbR^d$-valued BRWs, two normalization sequences $\smash{(a_n)_{n\in \bbN}}$ and $\smash{(b_n)_{n\in \bbN}}$, a moment gauge function $G$ as in (\ref{momgaudef}), a nonnegative continous process $\smash{(C_{\! s })_{s\in \bbR_+}}$, which plays the role of a lifetime process of a snake, and a r.v.~$\smash{\fbeta\! :\! \Omega \! \to \! \mathrm{Sym}_{^d}^{_+}}$. We introduce the following set of assumptions. 
 
\smallskip

\noindent 
$\,$ \textbf{Sib-Ind}$_{_{\, }}$($\baa, \bbb, \bS_\cdot, G, C_\cdot, \fbeta$) :
 
\smallskip
 
\begin{compactenum} 
\item[$\mathbf{(1)}$] We assume \texttt{Norm}$_{^{\,}}$($\baa, \bbb$) and we also assume for all $c\ino \bbR_+^*$ that 
\begin{equation*}
\label{GrowthabG} 
\tag{{\small \texttt{Growth}$_{^{\,}}$($\baa, \bbb, G$)}}. 
\!\!\!\!\!\!\!\!\!\!  \!\!\!\!\!\!\!\!\!\!  \!\!\!\!  \! \lim_{n\to \infty} b_n /G \Big( c \sqrt{\lambda_n}/\log b_n \big) \eqo 0\; .
\end{equation*} 

\smallskip

\item[$\mathbf{(2)}$] For all $n\ino \bbN$, $\btt_n$ is a.s.~finite and conditionally given $\btt_n$, the jumps of $\bS_n$ are independent by group of siblings as in (\ref{sibinddef}) and conditionally integrable and centered as in (\ref{concendef}).

\smallskip

\item[$\mathbf{(3)}$] For all $n\ino \bbN$, $M(\bS_n, G) \leko \infty$ a.s.~and $\lim_{y\to \infty} \limsup_{n\to \infty} \bP \big( M(\bS_n, G) \geko y \big) \eqo 0$.

\smallskip

\item[$\mathbf{(4)}$] For all $\smash{n\ino \bbN}$, there is a $\smash{\btt_n}$-measurable r.v.~$\smash{\fbeta_n\! :\! \Omega \! \to \! \mathrm{Sym}^{_+}_{^d}}$ such that $\smash{\bdelta_n \! \to \! 0}$ in probability where we have set $\smash{\bdelta_n\eqo \max_{u\in \btt_n \backslash \{ \varnothing \}} 
\lVert \fbeta_n \! -\! \beta( \bS_n, u)  \rVert_\infty}$, with $\smash{ \beta( \bS_n, u)}$ as in (\ref{concovdef}). 

\smallskip
 
\item[$\mathbf{(5)}$] 
$\smash{(C^{_{(n)}}_{\cdot} \!, \fbeta_n) \! \to \! (C_\cdot, \fbeta)}$ holds in law in $\smash{\bC^{_0}_{^1} \! \times \! \mathrm{Sym}^{_+}_{^d}}$, where $\smash{C^{_{(n)}}_{\cdot}}$ stands for the contour process of $\smash{\btt_n}$, rescaled as in (\ref{renormCW}).

\smallskip
 
\item[$\mathbf{(6)}$] The laws of the r.v.s $\smash{\# \btt _n / b_n}$, $n\ino \bbN$, are tight in $\smash{\bbR_+}$. 

\end{compactenum}

\begin{remark}
\label{iidass}  When the jumps of the $\smash{\bS_n}$ are i.i.d., the previous set of assumptions simplifies and becomes  
\textbf{i.i.d.}$_{_{\, }}$($\smash{\baa, \bbb, \bS_\cdot, G, C_\cdot, \beta}$): 
$\mathbf{(1')}\! =\! \mathbf{(1)}$,  $\mathbf{(6')}\! =\! \mathbf{(6)}$ and 
\begin{compactenum} 

\smallskip

\item[$\mathbf{(2')}$] For all $\smash{n\ino \bbN}$, $\smash{\btt_n}$ is a.s.~finite and conditionally given $\smash{\btt_n}$, the jumps of $\smash{\bS_n}$ are i.i.d.~r.v.s whose law $\smash{\bgam_n}$ is deterministic integrable and centered.

\smallskip

\item[$\mathbf{(3')}$] $\smash{\sup_{n\in \bbN } \int_{\bbR^d} G( \lvert y\rvert ) \bgam_n (dy)\leko \infty}$.

\smallskip

\item[$\mathbf{(4')}$] $\smash{\lVert \beta_n \! -\! \beta \rVert_\infty \! \to \! 0} $ where $\smash{\beta_n}$ is the deterministic matrix  $\smash{\big(\int_{\bbR^d} y_iy_j \bgam_n (dy) ; 1\leqo i,j\leqo d \big)}$.  

\smallskip

\item[$\mathbf{(5')}$] $\smash{C^{_{(n)}}_{\cdot} \! \to \! C_\cdot}$ 
in law in $\smash{\bC^{_0}_{^1}} $. \cq 
\end{compactenum}
\end{remark}

The following theorem (based on the coupling Theorem \ref{brwcouplingth}) basically asserts the following: under the previously introduced assumptions (independence by group of siblings and a moment condition) \emph{if rescaled snakes converge in one case, then convergence holds in all cases}.  
\begin{theorem}
\label{extenscv}
 Let $\smash{(a_n)_{n\in \bbN}}$ and $\smash{(b_n)_{n\in \bbN}}$ be two renormalization sequences satisfying \emph{\texttt{Norm}$_{^{\,}}$($\smash{\baa, \bbb}$)}. Let $\smash{(\btt_n)_{n\in \bbN}}$ be a sequence of a.s.~finite random trees whose rescaled contour processes $\smash{C^{_{(n)}}_{\cdot}}$ are as in (\ref{renormCW}). Let $\smash{(C_s)_{s\in \bbR_+}}\, $ be a continuous nonnegative process. 
We assume that there are two moment gauge functions $\smash{G_{\! \ast}}$ and $G$, two r.v.s~$\smash{\fbeta_{\! \ast}\!  : \! \Omega \! \to \! \bbR_+}$ and $\smash{\fbeta\! : \! \Omega \! \to \! \mathrm{Sym}^{_+}_{^d}}$, and two sequences of BRWs, 
$\smash{\bS^\ast_n\eqo (S^\ast_{n,u})_{u\in \btt_n}}$ and $\smash{\bS_n \eqo (S_{n,u})_{u\in \btt_n}}$, $n\ino \bbN$, that satisfy the following. 
\begin{compactenum}

\smallskip

\item[$(a)$] The BRWs $\smash{(\bS^{ \ast}_{n})_{n\in \bbN}}$ are $\bbR$-valued and satisfy \emph{\textbf{Sib-Ind}$_{_{\, }}$($\smash{\baa, \bbb, \bS^\ast_\cdot, G_{\! \ast}, C_\cdot, \fbeta_{\! \ast}}$)}. 

\smallskip

\item[$(b)$] $\bP$-.a.s.~$\smash{\fbeta_{\! \ast} \in \bbR_+^*}$.

\item[$(c)$] The laws of the endpoint processes $\smash{(\widehat{W}_{\cdot}^{_{\ast,(n)}})_{n\in \bbN}}$, are tight in $\smash{\bC^{_0}_{^1}}$, where $\smash{W_{\cdot}^{_{\ast,(n)}}\!\!}$ stands for the snake of $\smash{\bS^\ast_n}$ rescaled as in (\ref{renormCW}).

\smallskip

\item[$(d)$] The BRWs $\smash{(\bS_{n})_{n\in \bbN}}$ are $\smash{\bbR^d}$-valued and satisfy \emph{\textbf{Sib-Ind}$_{_{\, }}$($\smash{\baa, \bbb, \bS_\cdot, G, C_\cdot, \fbeta}$)}.
\end{compactenum}

\smallskip

\noi
We denote by $\smash{W^{_{(n)}}_{\cdot}\! }$ the rescaled snake of $\smash{\bS_n}$ as in (\ref{renormCW}). 
Then, there is a continuous version of a $\smash{\bbR^d}$-valued Brownian snake $W$ with lifetime process $C$.
and $\smash{(C^{_{(n)}}_\cdot\! ,  W^{_{(n)}}_\cdot) \! \to \! (C \, , \, \sqrt{ \fbeta}.W )}$  in law in $\smash{\bC^{_0}_{^1} \! \times \! \mathbf C (\bbR_+, \bC^{_0}_{^d})}$. 
\end{theorem}
\noi
\textbf{Proof.} See Section \ref{Thms23pfsec}. \cqfd 

\smallskip

Remark \ref{iidass} and Theorem \ref{extenscv} allow to extend Theorem \ref{Sheuexplain} $(iii)$ 
as follows. 
\begin{theorem}
\label{maincvsnake} 
We fix $\smash{x, c \ino \bbR_+^*}$ and $\smash{\alpha \ino (1, 2]}$. 
Let $\smash{\psi\ino \mathscr L}$ satisfy \emph{\texttt{Var}$_{\infty}$($\psi$)}. Let $\smash{(a_n)_{n\in \bbN}}$, 
$\smash{(b_n)_{n\in \bbN}}$, $\smash{(\mu_n)_{n\in \bbN}}$, $\smash{\bS_n \eqo (S_{n, u})_{u\in \btt_n}}$ and 
$(C,Y,W)$ be as in any $\emph{\textbf{Case (i)}}$, $ \mathbf i \ino \{ 1,2,3\}$.
Recall that $\smash{C^{_{(n)}}_{\cdot}}$, $\smash{V^{_{(n)}}_{\cdot}}$ and $\smash{W^{_{(n)}}_{\cdot}}$ 
are rescaled as in (\ref{renormCW}). 
We make the following assumptions. 

\begin{compactenum}

\smallskip

\item[$(a)$] We assume \emph{\texttt{\L{}uka}$_{^{\,}}$($\smash{\baa, \bbb, \bmu, \psi}$)} and \emph{\texttt{Sheu}$_{^{\,}}$($\smash{\baa, \bbb, \bmu} $)}.
\smallskip

\item[$(b)$] We suppose that there is $\smash{\fbeta\! : \! \Omega \! \to \! \mathrm{Sym}^{_+}_{^d}}$ and a moment gauge function $G$ as in (\ref{momgaudef}) such that the sequence of BRWs $\smash{(\bS_n)_{n\in \bbN}}$ satisfies 
\emph{\textbf{Sib-Ind}$_{_{\, }}$($\smash{\baa, \bbb, \bS_\cdot, G, C_\cdot, \fbeta}$)} $(\mathbf{1}$-$\mathbf{5})$.  

\smallskip

\end{compactenum}

\noi
Then, the following limit holds in law in 
$\smash{\bC^{_0}_{^1} \! \times \! \mathbf C (\bbR_+, \bC^{_0}_{^d})}$.
\begin{equation}
\label{jcvbis}
\big( C^{_{(n)}}_\cdot\! ,  W^{_{(n)}}_\cdot \big)  \xrightarrow[n\to \infty]{\; }\big( C \, , \,  \sqrt{\fbeta}.W \big) .
\end{equation}
Furthermore, if $\smash{\lim_{n\to \infty}\big( C^{_{(n)}}_\cdot\! ,  V^{_{(n)}}_\cdot , \fbeta_n \big) \eqo  \big( C,V,  \fbeta \big) }$ holds in law in $\smash{\bC^{_0}_{^1} \! \times \! \bD (\bbR_+, \bbR)\times \mathrm{Sym}^{_+}_{^d}}$, then 
(\ref{jcvbis}) holds jointly with $\smash{ V^{_{(n)}}_\cdot \! \! \to \! Y}$. 
\end{theorem}
\noi
\textbf{Proof.} See Section \ref{Thmmainsnapfsec}. \cqfd 
\begin{remark}
\label{afterthm3}
$(a)$ Let us suppose that the assumptions of Theorem \ref{maincvsnake} hold true in $\textbf{Case (1)}$ for all $\smash{x\ino \bbR_+^*}$ with $\smash{\fbeta}$ and $G$ independent of $x$. Then (\ref{jcvbis}) holds in $\textbf{Case (0)}$ by standard arguments, and if in $\textbf{Case (1)}$, $\smash{\lim_{n\to \infty}\big( C^{_{(n)}}_\cdot\! ,  V^{_{(n)}}_\cdot , \fbeta_n \big) \eqo  \big( C,V,  \fbeta \big) }$ holds in law in $\smash{\bC^{_0}_{^1} \! \times \! \bD (\bbR_+, \bbR)\times \mathrm{Sym}^{_+}_{^d}}$, for all $\smash{x\ino \bbR_+^*}$, then (\ref{jcvbis}) holds jointly with $ \smash{V^{_{(n)}}_\cdot \! \! \to \! Y}$ in $\textbf{Case (0)}$. 

\smallskip

\noi
$(b)$ As mentioned right after (\ref{cvcodgen}), in $\textbf{Cases (1)}$, Assumptions 
\texttt{\L{}uka}$_{^{\,}}$($\baa, \bbb, \bmu, \psi$) alone, and in $\textbf{Cases (2)}$, 
\texttt{\L{}uka}$_{^{\,}}$($\baa, \bbb, \bmu, \psi$) combined with (\texttt{Hght}$_{^{\,}}$($\baa, \bbb, \bmu $) $+$ \texttt{Grey}$_{^{\,}}$($\psi$)) (which are implied by 
\texttt{Sheu}$_{^{\,}}$($\baa, \bbb, \bmu $) by Lemma \ref{controlpsi})) actually imply Assumption $(\mathbf{6})$ in \textbf{Sib-Ind}$_{_{\, }}$($\baa, \bbb, \bS_\cdot, G, C_\cdot, \fbeta$). This obviously holds true in $\textbf{Case (3)}$ too.

\smallskip

\noi
$(c)$ Let us mention that Lemma \ref{sheustable} shows that \texttt{\L{}uka}$_{^{\,}}$($\baa, \bbb, \bmu, \psi$) in $\textbf{Case (3)}$ implies \texttt{Sheu}$_{^{\,}}$($\baa, \bbb, \bmu $).

\smallskip

\noi
$(d)$ Recall that there is a time-change $\smash{\phi_{n}}$, which allows to derive discrete height processes and height-snakes from  contour processes and snakes: namely, $\smash{H^{_{(n)}}_{\! s} \eqo C^{_{(n)}}_{^{\! \phi_n(s)}}}$ and $\smash{\cW^{_{(n)}}_{\! s} \eqo W^{_{(n)}}_{^{\! \phi_n(s)}}}$, for all $\smash{s\ino \bbR_+}$. Since $\smash{\phi_n\! \to \! 2\mathrm{Id}_{\bbR_+} }$ 
 in probability in $\smash{\bC^{_0}_{^1}}$ (see (\ref{controphit}) for more details), then (\ref{jcvbis}) holds in law jointly with 
 $\smash{(H^{_{(n)}}_\cdot\! , \cW^{_{(n)}}_\cdot) \! \! \to \! ((C_{\! 2s})_{s\in \bbR_+} , (W_{\! 2s})_{s\in \bbR_+})}$.

\smallskip
 
\noi 
$(e)$ In $\textbf{Case (3)}$
and when the jumps of the BRWs are $\bbR$-valued independent r.v.s with a fixed deterministic law $\bgam$, Theorem \ref{maincvsnake} has been obtained by Janson \& Marckert \cite{JanMar05} when $\alpha\eqo 2$ and by Marzouk \cite{Mar20} in the other stable cases, with a necessary and sufficient condition in terms of the tail of $\bgam$ (see Section \ref{assumpsec} for a more precise discussion on this point). 
We actually do not provide an independent proof of this result: we only extend it thanks to Theorem \ref{extenscv}.

\smallskip
 
\noi 
$(f)$ Let $\psi \ino \ccL$ satisfy \texttt{Sheu}$_{^{\,}}$($\psi$) (and thus \texttt{Var}$_{\infty}$($\psi$)). Let the 
$a_n$, $b_n$, $\mu_n$, $\smash{\bS_n \eqo (S_{n, u})_{u\in \btt_n}}$ and $(C,Y,W)$ be as in $\textbf{Case (0)}$.
We assume \texttt{\L{}uka}$_{^{\,}}$($\baa, \bbb, \bmu, \psi$) and we assume that 
conditionally given $\btt_n$, the jumps of $\smash{\bS_n}$ are $\bbR$-valued and i.i.d.~r.v.s whose deterministic law $\smash{\bgam_n}$ 
satisfies the following: there are $\smash{c, \beta\ino \bbR_+^*}$ such that for all $n\ino \bbN$, 
\begin{equation}
\label{cenunii}
\bgam_n ([-c,c]) \eqo 1, \quad 
\int_{\bbR}\!  y \bgam_n(dy)\eqo 0 \quad \textrm{and} \quad \beta_n \! :=\! \int_{\bbR}\! y^2 \bgam_n(dy) \xrightarrow[n\to \infty]{\; } \beta.
\end{equation}
Then, $\smash{\sup_{n\in \bbN } \int_{\bbR} G( \lvert y\rvert ) \bgam_n (dy)\leko \infty}$ with $\smash{G(x)\eqo \exp (|x|\, )}$ and the condition \texttt{Growth}$_{^{\,}}$($\smash{\baa, \bbb, G}$) holds iff $\smash{b_n /(a_n  (\log a_n)^4) \! \to \! \infty}$. 
In this case the BRWs $\smash{(\bS_{n})_{n\in \bbN}}$ satisfy \textbf{i.i.d.}$_{_{\, }}$($\smash{\baa, \bbb, \bS_\cdot, G, C_\cdot, \beta}$) and Theorem \ref{maincvsnake} implies the convergence (\ref{jcvbis}) as noticed in the above Remark \ref{afterthm3} $(a)$. 
However, let us mention that we can prove actually the following slightly better result. 
\begin{compactenum}

\smallskip

\item[$(i)$] Under the assumptions \texttt{\L{}uka}$_{^{\,}}$($\baa, \bbb, \bmu, \psi$), (\ref{cenunii}) and \texttt{Sheu}$_{^{\,}}$($\baa, \bbb, \bmu$) only, the convergence $\smash{( C^{_{(n)}}_{\cdot}\! ,V^{_{(n)}}_{\cdot} , W^{_{(n)}}_{\cdot}  )}$  $\to $ $\smash{( C_\cdot, Y_\cdot, \sqrt{\beta} W_\cdot )}$ holds in law in $\smash{\bC^{_0}_{^1} \! \times \! \bD (\bbR_+, \bbR)\times \mathbf C (\bbR_+, \bC^{_0}_{^1})}$.

\smallskip

\item[$(ii)$] If $\smash{\bgam_n \eqo \frac{1}{2} (\delta_{-1} \! + \delta_{1})}$, for all $n\ino \bbN$ , 
and $\smash{\liminf_{n} \mu_n (1) \geko 0}$, then \texttt{Sheu}$_{^{\,}}$($\smash{\baa, \bbb, \bmu}$) is equivalent to the tightness in $\smash{\bC^{_0}_{^1}}$ of the laws of the endpoint processes $\smash{(\widehat{W}^{_{(n)}}_\cdot)_{n\in \bbN}}$. 

\smallskip

\end{compactenum}
This result is stated more precisely in Theorem \ref{unifbounded}, in Appendix 
\ref{unifboundedsec}. Its proof does not rely on the coupling result stated in Theorem \ref{brwcouplingth} and it is rather long. 
As already mentioned, this statement is better than the result obtained by applying Theorem \ref{maincvsnake}. Indeed, Proposition \ref{excntrex} below shows that for each $\psi$ as above, there are 
$a_n, b_n, \mu_n$ such that \texttt{\L{}uka}$_{^{\,}}$($\baa, \bbb, \bmu, \psi$) and \texttt{Sheu}$_{^{\,}}$($\baa, \bbb, \bmu$) hold true but $\smash{\limsup_{n\to \infty} b_n /(a_n  (\log a_n)^4)\leko \infty}$: thus, Theorem \ref{maincvsnake} does not apply in these cases.  \cq  
\end{remark}

\subsection{On the assumptions of Theorems \ref{Sheuexplain}, \ref{extenscv} and \ref{maincvsnake}, examples and related works}
\label{assumpsec}

Let us first discuss in detail each of the assumptions in Theorems \ref{Sheuexplain}, \ref{extenscv} and \ref{maincvsnake}.
To that end, we fix $\psi \ino \ccL$. We also fix 
two renormalization sequences $\smash{(a_n)_{n\in \bbN}}$ and $\smash{(b_n)_{n\in \bbN}}$ and a sequence $\smash{(\mu_n)_{n\in \bbN}}$ of offspring distributions satisfying (\ref{nontricri}) (i.e., they are non-trivial and critical). 

\smallskip

\noi
$\; $\textbf{Sharpness of Assumptions} \texttt{Grey}$_{^{\,}}$($\baa, \bbb, \bmu$) \textbf{and}  
\texttt{Sheu}$_{^{\,}}$($\baa, \bbb, \bmu$). As already mentioned, Theorem \ref{Greyexplain}, 
Section \ref{PfsecThmGreyexplain}, asserts the following: under \texttt{Grey}$_{\infty}$ ($\psi$) and \texttt{\L{}uka}$_{^{\,}}$($\baa, \bbb, \bmu, \psi$), 
\texttt{Grey}$_{^{\,}}$($\baa, \bbb, \bmu$) is a necessary and sufficient condition for the convergence of the contour 
processes $\smash{\mathscr C^{_{(n)}}_\cdot}$ of GW($\mu_n$)-forests with independent exponentially distributed lifespans 
with parameter $1$. Moreover, Theorem \ref{Sheuexplain} shows, under \texttt{Var}$_{\infty}$ ($\psi$) and \texttt{\L{}uka}$_{^{\,}}$($\baa, \bbb, \bmu, \psi$),  
that 
\texttt{Sheu}$_{^{\,}}$($\baa, \bbb, \bmu$) is a necessary and sufficient condition for the convergence of 
$1$-dimensional Brownian snakes $\smash{\mathscr W^{_{(n)}}_\cdot}$ with lifetime processes 
$\smash{(\mathscr C^{_{(n)}}_{ s})_{s\in \bbR_+}}$ or for the convergence of the snake of $\bbR$-valued BRWs with 
i.i.d.~jumps with law $\smash{\frac{1}{2} e^{-|x|} dx}$. 
Furthermore, as explained in Remark \ref{afterthm3} $(d)$, under \texttt{Var}$_{\infty}$ ($\psi$) and \texttt{\L{}uka}$_{^{\,}}$($\baa, \bbb, \bmu, \psi$) and 
if $\smash{\liminf_{n\to \infty} \mu_n (1) \geko 0}$, then 
\texttt{Sheu}$_{^{\,}}$($\baa, \bbb, \bmu$) is a necessary and sufficient condition for the convergence of snakes of 
$\bbR$-valued BRWs with i.i.d.~jumps with law $\smash{\tfrac{1}{2} (\delta_{-1}\! + \delta_{1})}$.

\smallskip

\noi
$\; $\textbf{About the growth of $b_n$ in terms of $a_n$}. 
Here we remind that $\psi \ino \mathscr L$ is fixed. Let us take a closer look to \texttt{Growth}$_{^{\,}}$($\baa, \bbb, G$) which connects $G$ to the growth of $\lambda_n\eqo b_n /a_n$ in terms of $a_n$. We observe that the weaker $b_n$ grows with $a_n$, the stronger is the moment assumption that we need to make in Theorem \ref{maincvsnake} to obtain a convergence of the snakes. To discuss the growth of $b_n$, it is convenient to introduce a function $\phi \! :\! \bbR_+ \! \to \bbR_+$ such that $b_n \eqo  \phi (a_n)$ and such that it faithfully interpolates the values of $b_n$ is terms of $a_n$. 

On one hand, we notice that \texttt{\L{}uka}$_{^{\,}}$($\baa, \bbb, \bmu, \psi$), \texttt{Var}$_{\infty}$($\psi$), \texttt{Grey}$_{^{\,}}$($\baa, \bbb, \bmu$) and \texttt{Sheu}$_{^{\,}}$($\baa, \bbb, \bmu$) \emph{actually impose minimal growth conditions of $b_n$ in terms of $a_n$}. Namely, Lemma \ref{Lukagrowth} ($iii$), Theorem \ref{Greyexplain}  
and Theorem \ref{Sheuexplain} show 
under \texttt{\L{}uka}$_{^{\,}}$($\baa, \bbb, \bmu, \psi$) that the following holds true. 

\begin{compactenum}

\smallskip

\item[$-$] Assumption \texttt{Var}$_{\infty}$($\psi$) implies $\lim_{x\to \infty}\phi (x)/x\eqo \infty$,

\smallskip

\item[$-$] Assumption \texttt{Grey}$_{^{\,}}$($\baa, \bbb, \bmu$)  implies $\lim_{x\to \infty} \phi (x)/(x\log x) \eqo \infty$,

\smallskip

\item[$-$] Assumption \texttt{Sheu}$_{^{\,}}$($\baa, \bbb, \bmu$)  implies $\lim_{x\to \infty}\phi (x)/(x(\log x)^2)\eqo \infty$.
\end{compactenum}

\smallskip

\noindent
On the other hand, Lemma \ref{Lukagrowth} ($i$) and Remark \ref{Lukagrowthrem} show that these assumptions \emph{do not impose any maximal growth condition}. Namely, for any function $\phi$ such that $\phi(x)/x\! \to \!  \infty$ (resp.~$\phi(x)/(x\log x) $, $\phi(x)/(x(\log x)^2)\! \to \!  \infty$) we can find $\baa, \bbb, \bmu$ such that $b_n \eqo \phi (a_n) $ and such that  
\texttt{\L{}uka}$_{^{\,}}$($\baa, \bbb, \bmu, \psi$) holds true and \texttt{Var}$_{\infty}$($\psi$) too (resp.~\texttt{Grey}$_{^{\,}}$($\baa, \bbb, \bmu$) too, \texttt{Sheu}$_{^{\,}}$($\baa, \bbb, \bmu$) too). 

However if $b_n /a_n^{2}$ is unbounded, then Lemma \ref{Lukagrowth} $(ii)$ asserts that $\smash{\limsup_{n\to \infty} 
 \mu_n (1)\eqo 1}$ and in that case \emph{$\smash{(b_n)_{n\in \bbN}}$ can not be seen as the genuine time-renormalization sequence for the RWs $\smash{(V^{{n}}_{k})_{k\in \bbN}}$} (indeed the RWs $\smash{(V^{{n}}_{k})_{k\in \bbN}}$ have geometric holding times with parameters $\mu_n (1)$, which are therefore longer and longer). 
From now on in this section, we assume that $\smash{(b_n)_{n\in \bbN}}$ is the \emph{true time-renormalization}. Namely, we assume 
 \begin{equation}
 \label{righttime}
q:=\limsup_{n\to \infty} 
 \mu_n (1)<  1.
\end{equation} 
 But Lemma \ref{Lukagrowth} $(ii)$ implies $\limsup_{n\to \infty} b_n/a_n^2 \leqo 4\beta_\psi / (1 \! -\!  q)$. 
 This is why, under \texttt{\L{}uka}$_{^{\,}}$($\baa, \bbb, \bmu, \psi$) and \texttt{Var}$_{\infty}$($\psi$), it is natural to restrict our attention to functions $\phi$ such that $\phi(x)/x\! \to \!  \infty$ and $\phi(x)\eqo \mathcal O_\infty (x^2)$. 
More precisely, to provide examples it is convenient to take $\phi$ in the class of functions $\ccL$ (see Definition \ref{defLaplclass}) which satisfy \texttt{Var}$_{\infty}$($\phi$). Indeed, in this case, we get 
$\lim_{x\to \infty}\phi(x)/x\eqo \infty$ and $\lim_{x\to \infty}  \phi(x)/x^2\eqo \beta_\phi$. 

Let us mention that the set of functions $ f \ino \ccL$ satisfying \texttt{Var}$_{\infty}$($f$) is sufficiently rich. It contains for instance functions $f_{\alpha, c}$ such that 
$f_{\alpha, c} (x)\asymp_\infty x^{\alpha} (\log x)^c$ with $(\alpha, c)$ belonging to
$\big( \{ 1\} \times (1, \infty) \big)\cup \big((1,2) \! \times \! \bbR \big)\cup \big( \{ 2\} \times (-\infty, 0]\big)$ (see Lemma \ref{examplphi} for more details). 
If $b_n \! \asymp a_n (\log a_n)^c$, then \texttt{Growth}$_{^{\,}}$($\baa, \bbb, G$) holds for moment gauge functions of the form $G(x)\eqo \exp (c' x^{2/(c-2)})$ (here with $c\geko 4$, necessarily). If 
$b_n \! \asymp f_{\alpha, c} (a_n)$ with $(\alpha, c) \ino (1,2)\! \times \bbR$, then  \texttt{Growth}$_{^{\,}}$($\baa, \bbb, G$) holds for moment gauge functions of the form $G(x)\eqo x^{2\alpha/(\alpha-1)} (\log x)^{c'}$ with $c'\geko 2(\alpha  - c)/(\alpha-1)$ (the same remains true if $\alpha\eqo 2$ and $c\leqo 0$).

The following proposition provides examples showing that none of the statements in Theorems \ref{Sheuexplain}, \ref{extenscv} and \ref{maincvsnake} is empty. 
This proposition also shows (under \texttt{\L{}uka}$_{^{\,}}$($\baa, \bbb, \bmu, \psi$) and (\ref{righttime})) that \emph{the growth of $b_n$ as a function of $a_n$} $-$ and thus Assumption \texttt{Growth}$_{^{\,}}$($\baa, \bbb, G$) $-$ \emph{depends essentially neither on $\psi$ nor on} \texttt{Sheu}$_{^{\,}}$($\baa, \bbb, \bmu$).

\begin{proposition} 
\label{excntrex} Let $\psi, \phi \ino \mathscr L$ satisfy \emph{\texttt{Var}$_{\infty}$($\psi$)} and \emph{\texttt{Var}$_{\infty}$($\phi$)}. 
Suppose that $\beta_\phi \leqo \beta_\psi$ and let $q \ino (0, 1)$. Let $\mathbf a \eqo (a_n)_{n\in \bbN}$ be a fixed sequence of positive real numbers such that $\lim_{n\to \infty} a_n \eqo \infty$. 
Below $ \mathbf b\eqo (b_n)_{n\in \bbN}$ stands for a sequence of positive real numbers 
such that $ \lim_{n\to \infty}b_n\eqo \infty$ 
and $\bmu \eqo (\mu_n)_{n\in \bbN}$ stands for a sequence of offspring distributions 
that satisfies (\ref{nontricri}). Then, the following holds true. 
\begin{compactenum}

\smallskip

\item[$(i)$] If \emph{\texttt{Sheu}$_{^{\,}}$($\psi$)} and \emph{\texttt{Sheu}$_{^{\,}}$($\phi$)} hold true, then there are 
$\mathbf b$ and $\bmu$ satisfying  \emph{\texttt{\L{}uka}$_{^{\,}}$($\baa, \bbb, \bmu, \psi$)}, $\lim_{n\to \infty} \mu_n (1) \eqo q$, $b_n \! \asymp_\infty\!  \phi(a_n)$ and \emph{\texttt{Sheu}$_{^{\,}}$($\baa, \bbb, \bmu$)}.

\smallskip

\item[$(ii)$] If \emph{\texttt{Grey}$_{^{\,}}$($\psi$)} and \emph{\texttt{Grey}$_{^{\,}}$($\phi$)} hold true but not \emph{\texttt{Sheu}$_{^{\,}}$($\phi$)}, then there are $\mathbf b$ and $\bmu$ satisfying  \emph{\texttt{\L{}uka}$_{^{\,}}$($\baa, \bbb, \bmu, \psi$)}, $\lim_{n\to \infty} \mu_n (1) \eqo q$, $b_n \! \asymp_\infty\!  \phi(a_n)$, \emph{\texttt{Grey}$_{^{\,}}$($\baa, \bbb, \bmu$) } but not \emph{\texttt{Sheu}$_{^{\,}}$($\baa, \bbb, \bmu$)}.

\smallskip

\item[$(iii)$] If \texttt{Grey}$_{^{\,}}$($\phi$) does not hold, 
then there are $\mathbf b$ and $\bmu$ that satisfy  \emph{\texttt{\L{}uka}$_{^{\,}}$($\baa, \bbb, \bmu, \psi$)}, $\lim_{n\to \infty}$ $ \mu_n (1) \eqo q$, $b_n \! \asymp_\infty\!  \phi(a_n)$ but not  \emph{\texttt{Grey}$_{^{\,}}$($\baa, \bbb, \bmu$)}. 
\end{compactenum}
\end{proposition}
\noi
\textbf{Proof}: see Section \ref{Proofsecexcntrex}. \cqfd

\smallskip

\noi
$\; $\textbf{About Assumption} \texttt{Growth}$_{^{\,}}$($\baa, \bbb, G$). 
We now discuss how sharp \texttt{Growth}$_{^{\,}}$($\baa, \bbb, G$) is in Theorem \ref{maincvsnake}.
We fix $\psi\ino \ccL$ and $(a_n)_{n\in \bbN}$ such that $\lim_{n\to \infty} a_n\eqo \infty$. 
We also fix $\alpha \ino (1, 2)$ and $q\ino (0, 1)$. 
We first assume \texttt{Sheu}$_{^{\,}}$($\psi$) and to simplify the discussion, it is also convenient to assume that $(a_n)$ increases and does not go too fast to  $\infty$. Namely, we suppose that there is $\kappa \ino (1,\infty)$ such that $1\leko a_{n+1} / a_n \leko \kappa $, $n\ino \bbN$. 

By Proposition \ref{excntrex} $(i)$ we can find $(b_n)_{n\in \bbN}$ and $(\mu_n)_{n\in \bbN}$ such that \texttt{\L{}uka}$_{^{\,}}$($\baa, \bbb, \bmu, \psi$) holds true, $\lim_{n\to \infty} \mu_n (1) \eqo q\ino (0, 1)$, 
 $b_n \! \asymp \!  a_n^\alpha$ and \texttt{Sheu}$_{^{\,}}$($\baa, \bbb, \bmu$). To simplify again, we next consider 
 $\bS_n \eqo (S_{n, u})_{u\in \btt_n}$, $n\in \bbN$, a sequence of real-valued BRWs such that 
 conditionally given $\btt_n$, the jumps of $\bS_n$ are i.i.d.~r.v.s~with the same deterministic law $\bgam$ which is independent of $n$ and which satisfies $\beta \! :=\! \int_{\bbR} x^2 \bgam (dx) \leko \infty$ and 
 $\int_{\bbR} x \bgam (dx)\eqo 0$. Moreover, we assume that the $\btt_n$ and $(C,Y,W)$ are as in $\textbf{Case (0)}$. Let $\xi$ be a r.v.~distributed according $\bgam$. 

We first note that \texttt{Growth}$_{^{\,}}$($\baa, \bbb, G$) is satisfied by moment gauge functions $G$ that satisfy 
$\smash{G(x)}$ $\smash{ \asymp_\infty}$ $\smash{ x^{2\alpha/(\alpha-1)} (\log x)^{\epp +(2\alpha/(\alpha-1))}}$, where $\epp \ino \bbR_+^*$.   
If $\bE [G (|\xi|)] \leko \infty$, then Theorem \ref{maincvsnake} applies and the joint convergence (\ref{jcvbis}) of the contour process and of the discrete snakes holds true. 
On the other hand, we observe that the sole assumption that the laws of $\smash{\widehat{W}^{_{(n)}}_\cdot}$ are tight in $\bC^{_0}_{^1}$, implies that $\smash{b_n \bP \big(|\xi | \geko c\sqrt{\lambda_n} \big)\! \to \! 0}$ for all $c\ino \bbR_+^*$. 
Since $(a_n)$ does not increase too fast, it implies that $\smash{\bP ( |\xi | \geko x) \eqo o_\infty (x^{-2\alpha/(\alpha-1)})}$, whereas the assumption 
$\bE [G (|\xi|)] \leko \infty$ implies that $\smash{\bP ( |\xi | \geko  x)} $  $\eqo$ $\smash{ \mathcal O_\infty (x^{-2\alpha/(\alpha-1)} (\log x)^{-\epp -(2\alpha/(\alpha-1)) }}$. 
So this shows that the moment assumption in Theorem \ref{maincvsnake}, which is controlled by 
\texttt{Growth}$_{^{\,}}$($\baa, \bbb, G$), is sharp up to logarithmic factors (that are due to coupling).

\smallskip 
 
\noi
\textbf{Related works.}  Apart from Marzouk \cite{Mar20} who considers the cases of fixed offspring distributions 
in the domain of attraction of a stable law, the already known results on the scaling limits of discrete snakes concern fixed 
offspring distributions $\mu_n\eqo \mu$, $n\ino \mathbb{N}$, which have a moment of order $2$, and 
the Galton-Watson tree is as in \textbf{Case (3)}, that is, $\btt_n$ is distributed as a GW($\mu$)-tree conditioned on having $n$ vertices: Chassaing \& Schaeffer \cite{ChassSchaef_2004} consider the case of the critical geometric offspring distribution and i.i.d. uniform jumps on $\{ 1, -0, 1\}$; Marckert \& Mokkadem \cite{MaMo03} also consider the critical geometric offspring distribution but with independent, homogeneous, centered jumps with moments of order $>6$; Gittenberger \cite{Gittenberger2003} extends this result to general offspring distributions with a moment of order $2$ and jumps with moments of order $>8$; Janson \& Marckert \cite{janson2015scaling} consider the case of a general reproduction law with a moment of order $2$ and i.i.d.~centered jumps satisfying the hypothesis $\smash{\bP (|\xi| \geko x)\eqo  o_\infty(x^{-4})}$, which is the optimal result. 
This result was extended by Marzouk \cite{Mar20} to the case of a fixed offspring distribution $\mu$ 
in the domain of an $\alpha$-stable distribution $\alpha \ino (1, 2]$ (and still in \textbf{Case (3)}), 
for i.i.d.~centered jumps satisfying the hypothesis $\smash{\bP (|\xi | \geko  x)\eqo  o_\infty (x^{-2\alpha/(\alpha -1)})}$, which is optimal. 

In Marckert \cite{marckert2008lineage} and Addario-Berry, Donderwinkel, Goldschmidt \& Mitchell \cite{ABDoGoMi25+}, 
the offspring distribution is fixed and has a moment of order $2$; jumps are independent by group of siblings 
and exhibit inhomogeneities. More precisely, the following assumptions hold:  for every $k\ino \bbN^*$, 
let $\nu_{k}(dz)$ be a distribution on $\bbR^k$; for every $u\ino \btt_n\backslash \mathtt{Lf} (\btt_n)$, let 
$(\xi_{[u]\ast i})_{1\leq i\leq k_u(\btt_n)}$ denote the jumps of the children of $u$ (recall that $k_u(\btt_n)$ is the number of children of $u$); then, 
\begin{compactenum}
\item[$(a)$] conditionally on $\btt_n$, the distribution of $(\xi_{[u]\ast i})_{1\leq i\leq k_u(\btt_n)}$ is $\nu_{k_u(\btt_n)}(dz)$; 
\item[$(b)$] jumps have a global variance: $\smash{\beta\! :=\! \sum_{k\in \bbN^*} \mu(k) \int_{\bbR^k}\nu_k (dz) |z|^2  
\leko \infty}$, and they are globally centered: 
$\smash{\sum_{k\in \bbN^*} \mu(k) \sum_{1\leq i\leq k} \overline{\nu}_{k,i} \eqo 0}$ where $\smash{\overline{\nu}_{k,i}\! :=\! \int_{\bbR^k}\nu_k (dz) \langle z, \mathtt e^k_i \rangle }$ and $ \smash{\mathtt e^{_k}_i}$ stands for the $i$-th canonical vector of $\bbR^k$.  
\end{compactenum}
Marckert  \cite{marckert2008lineage} shows under $(a)$, $(b)$, and under the additional assumptions that the offspring distribution has a finite support in $\{ 1, \ldots, K\}$ and that there exists $p \geko 4$ such that 
$$\max_{1\leq k\leq K}\max_{1\leq i\leq k} \int_{\bbR^k}\nu_k (dz) \big|  \langle z, \mathtt e^k_i  \rangle -\overline{\nu}_{k,i}\big| ^p  < \infty, $$
the joint functional convergence of the discrete snake and the contour function of $\btt_n$ holds, that is, (\ref{jcvbis}). 
Addario-Berry, Donderwinkel, Goldschmidt, \& Mitchell \cite{ABDoGoMi25+} extend this functional limit theorem under $(a)$, $(b)$, and under the additional, weaker assumptions that the offspring distribution $\mu$ has a moment of order $3$ and that 
$\smash{ \sum_{k\in \bbN^*} \mu(k) \nu_k (\{ z\ino \bbR^k : |z| \geko y\}) \eqo o_\infty(y^{-4})}$.
This last assumption is, moreover, necessary. 

These results are obtained by a fine combinatorial analysis and 
by precisely calibrated probabilistic arguments of uniformization.
 Our approach is less precise:   
while Theorem \ref{maincvsnake} allows for wider jump inhomogeneities and 
requires no more than a moment of order $2$ for the offspring distribution, however, the coupling used in its proof requires, on the one hand, stronger moments (a.s.~$\smash{\max_{u\in \btt_n\backslash \{ \varnothing \}} \int \nu_{k_u(\btt_n)} ( dz) |z|^4 (\log_+ |z|)^{4+\epp} \leko \infty}$) and, on the other hand, a stricter local control of the expectation and variance: specifically, in \textbf{Sib-Ind}$_{_{\, }}$($\baa, \bbb, \bS_\cdot, G, C_\cdot, \fbeta$) 
a.s.~for all $u\in \btt_n \backslash \{ \varnothing \} $ and 
$\smash{\bE [ \xi_u | \btt_n] \eqo 0}$ and $\smash{\max_{u\in \btt_n\backslash \mathtt{Lf} (\btt_n)} \max_{1\leq i\leq k_u(\btt_n) }| \beta \! -\! 
 \int \nu_{k_u(\btt_n)} ( dz) \langle z,  \mathtt e_i \rangle^2| \leko \infty}$. 

\subsection{Application to scaling limits of the range of BRWs on regular trees} 
\label{traceapplsec}
As an application of Theorem \ref{maincvsnake}, we prove a limit theorem for the range of $\bbT_{\! \mathtt b}$-valued 
BRWs $(\Upsilon_{\! n,u})_{u\in \btt_n}$. More precisely: 
\begin{compactenum}

\smallskip

\item[$-$] ${ \bbT_{\! \mathtt b} \! :=\!  \bigcup_{n\in \bbN} \{ 1, \ldots, \mathtt b\}^n}$ is the \emph{$\mathtt b$-ary rooted tree}, i.e., the set of finite words written in the alphabet $\{ 1, \ldots, \mathtt{b}\}$, $\mathtt b\geqo 2$. Here, it is 
convenient to denote by $\mathtt o$ the root of $\bbT_{\! \mathtt b}$, which is the empty word. 
We denote by $d_{\mathtt{gr}}$ the \emph{graph-distance} on $ \bbT_{\! \mathtt b}$. 
For $x\eqo (j_1, \ldots, j_n)\ino   \bbT_{\! \mathtt b} \backslash \{\mathtt o\}$, we denote by 
$\smash{\overleftarrow{x}\! :=\! (j_1, \ldots, j_{n-1})}$ the neighbour of $x$ that is closest to $\mathtt o$ ($\smash{\overleftarrow{x}\eqo\mathtt o}$ if $n\eqo 1$).

\smallskip

\item[$-$] For all $n\ino \bbN$, $\btt_n$ is as in Theorem \ref{maincvsnake}, $\mathbf{Case}(2)$, i.e., a GW($\mu_n$)-tree conditioned to be higher than $c\lambda_n$, where $c\ino \bbR_+^*$ is fixed.  

\smallskip

\item[$-$] The BRWs a.s.~start at the root, i.e., $\Upsilon_{\! n,\varnothing } \eqo \mathtt o$. Conditionally given $\btt_n$, the jumps of the BRW $(\Upsilon_{\! n,u})_{u\in \btt_n}$ are distributed as those of a 
\emph{nearest neighbour null recurrent RW} on $\bbT_{\! \mathtt{b}}$, whose transition probabilities 
are defined as follows: for all $j\ino \{ 1, \ldots, \mathtt b\}$, 
$\smash{\mathtt p(\mathtt o, (j))\eqo \frac{1}{\mathtt b}}$, (the RW is reflected at $\mathtt o$); for all 
$\smash{x,y\ino \bbT_{\! \mathtt b}\backslash \{ \mathtt o\} }$ such that $\smash{\overleftarrow{y}\eqo x}$, $\smash{\mathtt p (x,y)\eqo \frac{1}{2\mathtt b}}$ and 
$\smash{\mathtt p (x, \overleftarrow{x}) \eqo \frac12}$ (thus, $\smash{\mathtt p(x,z)\eqo 0} $ if $\smash{d_{\mathtt{gr}} (x,z)\geqo 2}$).

\smallskip

\end{compactenum}

\noi
Namely, for all finite ordered rooted tree $\mathtt t$ and for all $x_u \ino \bbT_{\! \mathtt{b}}$, $u\ino \mathtt t$,
\begin{equation}
\label{lawBRW}
\bP \big( (\Upsilon_{\! n,u})_{u\in \mathtt t} \eqo (x_u)_{u\in \mathtt t} \, ;\,  \btt_n \eqo \mathtt t\big)= \bP ( \btt_n \eqo \mathtt t) 
\, \un_{\{ x_{\varnothing} = \mathtt o\}}\!\!\!  \prod_{u\in \mathtt t \backslash \{ \varnothing\} }\!\!  \mathtt p (x_{\overleftarrow{u}}, x_u) . 
\end{equation}
We want to prove a convergence in law as $n\! \to \! \infty$, of the following rescaled random metric spaces 
\begin{equation}
\label{Rbndef} \big( \mathcal R_{ n} , \tfrac{1}{\sqrt{\lambda_n}} d_{\mathtt{gr}} ,  \mathtt o, \tfrac{1}{b_n}\mathbf m_{\mathtt{occ}}^{(n)} \big) \quad \textrm{where} \quad \mathcal R_{ n}\! :=\! \big\{\Upsilon_{\! n,u} ; u\ino \btt_n \big\} \quad \textrm{and} \quad 
\mathbf m_{\mathtt{occ}}^{(n)}= \sum_{u\in \btt_n} \delta_{\Upsilon_{\! n,u}} .
\end{equation}
Here, $\mathcal R_{ n}$ is the \emph{range} of the BRW $(\Upsilon_{\! n, u})_{u\in \btt_n}$, $\mathbf m_{{\mathtt{occ}}}^{_{(n)}}$ is its \emph{occupation measure} and the convergence holds 
weakly on the space $\mathbb M$ of \emph{isometry classes of pointed measured compact metric spaces} equipped with the \emph{Gromov-Hausdorff-Prokhorov distance} 
$\bdelta_{\mathtt{GHP}}$ that makes it a Polish space, as proved in Abraham, Delmas \& Hoscheit in \cite{ADHoscheit}, Theorem 2.5. (See (\ref{GHP}) in Definition \ref{GHPpointdef} for a precise definition of $\bdelta_{\mathtt{GHP}}$). 

   To define the limiting space, let us first briefly introduce Lévy trees, which are the scaling limits of $(\btt_n, \frac{1}{\lambda_n} d_{\mathtt{gr}}, \varnothing, \frac{1}{b_n} \#_{\btt_n})$ where $d_{\mathtt{gr}}$ stands for the graph-distance on $\btt_n$ and $\#_{\btt_n}\! :=\! \sum_{u\in \btt_n} \delta_u$ for the counting measure on $\btt_n$. 
Here we work under the assumptions of Theorem \ref{maincvsnake}, $\mathbf{Case}(2)$. In particular the branching mechanism 
$\psi \ino \mathscr L$ satisfies $\texttt{Grey}$($\psi$), 
and the limiting height process $H$ has law $\bN(\, \cdot \,  | \max_{s\in \bbR_+} H_s \geko c)$, where $\bN$ stands for the excursion measure above $0$ of the $\psi$-height process. The lifetime of $H$ is denoted by $\zeta$. We define the contour process by $C_s\eqo H_{s/2}$, $s\ino \bbR_+$ (its lifetime is therefore $2\zeta$) and to define the \emph{$\psi$-Lévy tree higher than $c$}, we introduce 
$$ \forall s_1, s_2 \ino [0, 2\zeta] , \quad  d_C (s_1, s_2) \! = \! C_{s_1}+ C_{s_2} -\, 2\!\!\!\!\!\! \!\! \!\! \min_{\quad s_1\wedge s_2 \leq s \leq 
s_1 \vee s_2} \!\! \!\! \!\!\!\! C_s.
$$ 
We easily check that $d_C$ is a pseudo-metric on $[0, 2\zeta]$ and we introduce the relation $\sim_C $ on 
$[0, 2\zeta]$ by setting $s_1 \! \sim_C\!  s_2$ iff $d_C (s_1, s_2)\! = \! 0$; 
it is an equivalence relation and the $\psi$-Lévy tree higher than $c$ is defined as the quotient space 
$ T_C \! = \! [0, 2\zeta ] / \! \sim_C $, equipped with the distance induced by $d_C$, that we keep denoting $d_C$. 
We denote by $\mathtt{proj}  \! : \! [0, 2\zeta]\!  \rightarrow \! T_C$ the canonical projection. Note that $\mathtt{proj}$ is continuous. Therefore $T_C$ is compact and connected. 
Moreover, $T_C$ is a \emph{$\bbR$-tree}, namely, a metric space such that all pairs of points $\sigma_1, \sigma_2\ino T_C$ are joined by a unique arc denoted by $\lgeo \sigma_1, \sigma_2 \rgeo_{T_C}$, 
which turns out to be a geodesic (see Definition \ref{Rtreedef} for more details). 
We set {$r_C\! := \! \mathtt{proj} (0)$} that is viewed as the root of $T_C$ and we equip $T_C$ with the measure $\mu_C$ that is the image of the Lebesgue measure on $[0, 2\zeta]$ via $\mathtt{proj}$, namely, $\smash{\int_{T_C} \! f \, d\mu_C\! = \! \int_0^{2\zeta} \! f(\mathtt{proj}(s)) \, ds}$, for all 
measurable and bounded $f\! : \! T_C \! \rightarrow \! \bbR$. The convergence (\ref{cvcodgen}) 
then implies the following weak convergence on $(\mathbb M, \dGHP)$: 
\begin{equation}
\label{cvtolev}
\big( \btt_n , \tfrac{1}{\lambda_n} d_{\mathtt{gr}} , \varnothing, \tfrac{2}{b_n} \#_{\btt_n}  \big) 
\xrightarrow[n\to \infty]{}\; (T_C, d_C, r_C, \mu_C) \; .
\end{equation}

The limiting space of (\ref{Rbndef}) is then constructed as follows. Under the assumptions of Theorem \ref{maincvsnake}, $\psi $ satisfies $\texttt{Sheu}$($\psi$) and there is a continuous version $\smash{(W_{\! s})_{s\in [0, 2\zeta]}}$ of the $1$-dimensional Brownian snake with lifetime $C$ and whose initial value is the null function. By the snake property, it makes sense to define 
a $\smash{T_C}$-indexed process $\smash{(\mathtt W (\sigma))_{\sigma \in T_C}}$ by setting $\smash{\mathtt W (\mathtt{proj} (s)) \! :=\! \widehat{W}_{\! s}}$, for all $s\ino [0, 2\zeta]$. Conditionally given $C$, $\smash{\sigma \!  \mapsto\!  W (\sigma)}$ is thus a continuous centered Gaussian process whose conditional covariance is characterized by $\smash{\bE \big[ |\, \mathtt W (\sigma_1)\!  \! -\! \mathtt W (\sigma_2) |^2 \big| C \big]\eqo d_C(\sigma_1, \sigma_2)}$, for all $\smash{\sigma_1, \sigma_2 \ino T_C}$. Then, we set 
\begin{equation*}
\forall \sigma_1, \sigma_2\ino T_C, \quad d_{C, |W|} ( \sigma_1, \sigma_2) = |\, \mathtt W (\sigma_1)  |+ | \, \mathtt W (\sigma_2)| -2 \!\!\!\! \!\! \!   \min_{ \quad \sigma \in \lgeo \sigma_1, \sigma_2 \rgeo_{T_C}} \!\! \!\! \!\! |\,  \mathtt W (\sigma) | , 
\end{equation*}    
where $\smash{\lgeo \sigma_1, \sigma_2 \rgeo_{T_C} }$ is the unique arc joining $\sigma_1$ to 
$\sigma_2$ in $\smash{T_C}$. Lemma 4.2 in D., Khanfir, Lin \& Torri \cite{DuKhLiTo22} asserts that 
$\smash{d_{C, |W|}}$ is a pseudo-metric on $\smash{T_C}$. We define the equivalence relation 
$\smash{\sim_{C, |W|} }$ on $\smash{T_C}$ by setting $\smash{\sigma_1 \sim_{C, |W|}  \sigma_2}$ iff 
$\smash{d_{C, |W|} ( \sigma_1, \sigma_2) \! = \! 0}$. We define $\smash{T_{C, |W|}}$ as the quotient metric 
space $\smash{T_C / \! \sim_{C, |W|}}$ equipped with the induced metric 
$\smash{d_{C, |W|}}$; we denote by $\smash{\mathtt{proj}'\! : \! T_C \! \rightarrow \! T_{C, |W|}}$ the canonical projection. 
Since it is continuous, $\smash{T_{C, |W|}}$ is compact and connected. Moreover $\smash{(T_{C, |W|}, d_{C, |W|})}$ is a 
$\bbR$ tree (see Lemma 4.2 \cite{DuKhLiTo22} for a proof, this result is more precisely recalled in Proposition \ref{contsnadist}). We next set $\smash{r_{C, |W|}\! = \! \mathtt{proj}' (r_C)}$ that is viewed as the root of $\smash{T_{C, |W|}}$ 
and we equip $\smash{T_{C, |W|}}$ with the measure $\smash{\mu_{C, |W|}}$ that is the image 
of $\smash{\mu_C}$ via $\smash{\mathtt{proj}'}$, i.e., $\smash{\int_{T_{C, |W|}} \! f \, d\mu_{C, |W|}\! = \! \int_{T_C} \! f(\mathtt{proj}'(\sigma)) \, \mu_C (d\sigma) }$, for all bounded measurable $\smash{f\! : \! T_{C, |W|} \! \rightarrow \! \bbR}$. 

This kind of space has been introduced for a quadratic branching mechanism 
by Curien, Le Gall \& Miermont \cite{CuLGMi13} (see also Le Gall \cite{LG15} for a different purpose) who coined the name \emph{Brownian cactus}. So we call $(T_{C, |W|}, d_{C,|W|})$ the \emph{reflected Brownian cactus with branching mechanism $\psi$}.  
As an application of Theorem \ref{maincvsnake}, we prove the following limit theorem. 
\begin{theorem}
\label{applicactus} We keep the notations and assumptions of Theorem \ref{maincvsnake}, $\mathbf{Case}(2)$. Let $(\Upsilon_{\! n,u})_{u\in \btt_n}$ be a $\bbT_{\! \mathtt b}$-valued BRW whose law is given by (\ref{lawBRW}). Recall from (\ref{Rbndef}), its range $\mathcal R_{n}$ and its occupation measure $\mathbf m_{{\mathtt{occ}}}^{_{(n)}}$, and recall that $(T_{C,|W|}, d_{C, |W|}, r_{C, |W|}, \mu_{C,|W| })$ stands for the \emph{reflected Brownian cactus with branching mechanism $\psi$}. Then 
 the following limit holds weakly on $(\mathbb M, \delta_{\mathtt{GHP}})$: 
 \begin{equation}\label{cvcactus} \big( \mathcal R_{ n} , \tfrac{1}{\sqrt{\lambda_n}} d_{\mathtt{gr}} ,  \mathtt o, \tfrac{1}{b_n}\mathbf m^{(n)}_{\mathtt{occ}}\big) \xrightarrow[n\to \infty]{\;} \big( T_{\! C, |W| }, d_{C, |W|}, r_{ C, |W|}, \mu_{ C, |W|} \big) 
\end{equation} 
jointly with (\ref{cvtolev}).  
\end{theorem}
\noi
\textbf{Proof.} See Section \ref{pfsecapplicactus}.  \cqfd

\smallskip

A straightforward adaptation makes the convergence work also in $\mathbf{Case}(1)$. We emphasize that $\mathbf{Case}(3)$ has been proved in D., Khanfir, Lin \& Torri \cite{DuKhLiTo22}, Theorem 1.2. 
It is also important to mention here that the proof of Theorem \ref{applicactus} relies to a substantial extent on general arguments developed in \cite{DuKhLiTo22}, which establish a connection between, on the one hand, the geometry of the image $\mathcal R_{n}$ of the $\bbT_{\! \mathtt b}$-valued BRWs and, on the other hand, $\bbR$-valued discrete snakes indexed by the same GW-tree $\btt_n$. 
These results, which are summarized in detail but without proof in Section \ref{pfsecapplicactus}, make it possible to 
deduce relatively simply the convergence (\ref{cvcactus}) from that of the $\bbR$-valued BRWs stated in Theorem \ref{maincvsnake}.  

\section{BRWs and their encoding processes, GW-trees}
\label{discrobjsec}
\noi
\textbf{Rooted ordered trees.} In this section we recall from e.g.~Neveu \cite{Ne} Ulam's formalism of rooted ordered trees. We also discuss the connection between the various processes encoding these trees and BRWs. 
We first introduce the set of finite words written with positive integers: 
\begin{equation}
\label{UUlamdef}
\bbU\! = \! \bigcup_{n\in \bbN} (\bbN^*)^n\; .
\end{equation}
Here, $(\bbN^*)^0$ stands for $\{ \varnothing \}$; $\bbU$ is totally ordered by the \textit{lexicographical order} $\leq_{\mathtt{lex}}$ (the strict order is denoted by $<_\mathtt{lex}$). 
It is convenient to think of $\bbU$ as a family tree whose $\varnothing$ is the ancestor. 

Let $u\! = \! [j_1, \ldots , j_n]\ino \bbU\backslash \{ \varnothing\}$. 
We use the notation $|u|\! = \! n$ for the \textit{height} of $u$ (with the convention $|\varnothing |\! = \! 0$) and we 
set $\overleftarrow{u}\! = \! [j_1, \ldots, j_{n-1}]$ that is interpreted as the (direct) \textit{parent of $u$}. More generally, for all $p\ino \bbN^*$, we set $u_{|p} \eqo  [j_1, \ldots , j_{n\wedge p}]$ and $u_{|0}\eqo \varnothing$. Namely, $u_{|p}$ is the \emph{ancestor of $u$ at the generation $p$}. We observe that $u_{|p}\eqo u$ if $p\geqo |u|$ and that 
$u_{|p}\eqo \overleftarrow{u}$ if $p\eqo |u|\! -\! 1$. 
We let $v\! = \! [k_1, \ldots, k_m]\ino \bbU$ and we denote by $u\ast v$ the concatenated word $[j_1, \ldots, j_n, k_1, \ldots, k_m]$, with the convention $u\ast \varnothing \eqo \varnothing \ast u\eqo u$. 
For all $u, v\ino \bbU$, we then denote by $u \wedge v$ the \emph{most recent common ancestor} of $u$ and $v$ which is defined by $u\wedge v\eqo u_{|p} \eqo v_{|p}$ where $p\eqo \max \{ q\ino \bbN:   u_{|q} \eqo v_{|q}\}$.

We view $\bbU$ as a graph whose undirected edges are $\{ u, \overleftarrow{u} \}$, 
$u\ino \bbU\backslash \{ \varnothing\}$. For all $u,v\ino \bbU$, $\lgeo u, v\rgeo  \! \subset \! \bbU$ is 
the shortest path joining $u$ to $v$ in $\bbU$ and we 
use the following notations 
$\, \rgeo u , v\rgeo \!   = \! \lgeo u , v\rgeo  \backslash \{ u\}$, $\lgeo u,v\lgeo  \,  = \! \lgeo u, v\rgeo  \backslash \{ v\}$ and $\, \rgeo u,v\lgeo  \,  = \! \lgeo u,v \rgeo  \backslash \{u,v \}$. 
We equip $\bbU$ with the \emph{genealogical} order $\preceq$. Namely $u\preceq v$ if $u\ino \lgeo \varnothing , v \rgeo$.

\begin{definition}
\label{Ulamtree} $(a)$ A \textit{rooted ordered tree} is a subset $  \fftree\! \subset \! \bbU$ satisfying the following conditions.
\begin{compactenum}
\item[($a_1$)] $\varnothing \ino \fftree$.
\item[($a_2$)] If $u\ino  \fftree\backslash \{ \varnothing\}$ 
then $\overleftarrow{u}\ino  \fftree$. 
\item[($a_3$)] For all $u\ino  \fftree$, there exists $k_u ( \fftree)\ino \bbN\cup \! \{ \infty\}$ such that: $u \! \ast \!  [j]\ino  \fftree$  iff  $1\! \leq \! j \! \leq \! k_u ( \fftree)$. 
\end{compactenum}
Here, $k_u ( t)$ is interpreted as the \textit{number of children of $u$} and if 
$1\! \leq \! j \! \leq \! k_u ( t)$, then $u\ast\!  [j] $ is the \emph{$j$-th child of $u$.} If $k_u(t)\eqo 0$, $u$ is called a \emph{leaf} of $t$ and we denote by $\mathtt{Lf} (t)$ the set of leaves of $t$. We denote by $\bbT$ the set of rooted ordered trees. 

\smallskip

\noi
$(b)$ Let $ t\ino \bbT$ and $u\ino t$. Then 
$\theta_u  t\!  = \! \{ v\ino \bbU: u\ast v\ino  t \}\ino \bbT$ is the \emph{subtree of $ t$ stemming from $u$}.

\smallskip

\noindent
$(c)$ A \emph{forest} is a sequence of trees 
$t\eqo (t(p))_{1\leq p\leq N}$, $N\ino \bbN^{*}\! \cup \!\{ \infty\}$. It is sometimes convenient to view it as a single tree $t'\ino \bbT$ with $k_\varnothing(t')\eqo N$ such that  $t(p)\eqo \theta_{[p]}t'$ for all $1\leqo p\leqo N$. \cq 
\end{definition}

\begin{remark}
\label{contord} Let $t\ino \bbT$. Recall from Definition \ref{contsnadef} $(a)$ its contour exploration $(v_k)_{0\leq k < 2(\# t-1)}$. 

\smallskip

\noi
$(a) $ Let $(u_l)_{0\leq l<\# t}$ be the $t$-valued sequence that is recursively defined as follows:  
we set $u_0\eqo \varnothing$; if $l+1 \leko \#t$, then $u_{l+1}$ is the $\leq_{\mathtt{lex}}$-least element of the (necessarily nonempty) set $\{ u\ino t \! : \! u_l \! <_{\mathtt{lex}} \! u\}$. 
We call $(u_l)_{0\leq l<\# t}$ the \emph{depth-first exploration of $t$}, which corresponds to Definition \ref{heightlukadef} $(a)$. 

\smallskip

\noi $(b)$ In general the depth-first and contour explorations may not visit the whole tree. In this article we shall restrict to trees $t\ino \bbT$ such that 
$\# \theta_{[j]} t\leko \infty$ for all $1\leqo j\leqo k_\varnothing (t)$ (with  possibly $k_\varnothing (t)\eqo \infty$). Then for such trees we get $t\eqo \{ u_l; 0\leqo j\leko \# t\}\eqo  \{ v_k; 0\leqo k\leko 2(\#t \! -\! 1)\}$ 
(the tree associated with a forest of finite trees satisfies this condition for instance). 
In this case, the contour exploration crosses each edge of $t$
exactly twice (upwards first and then downwards). Namely, for all $u\ino t\backslash \{ \varnothing\}$, there is 
a unique pair of integers $k,k^\prime$ such that  
$k\leko k^\prime$ and $(\overleftarrow{u},u)\eqo (v_k, v_{k+1})\eqo (v_{k^\prime+1} , v_{k^\prime})$. \cq 
\end{remark}
 
 \smallskip

\noi
\textbf{Encoding processes of BRWs.} Let $t\ino \bbT$ be finite and let $S\eqo (S_u)_{u\in t}$ be a $\bbR^d$-valued BRW such that $S_\varnothing \eqo 0$. We recall from Definitions \ref{contsnadef} and \ref{heightlukadef} that 
$C_\cdot (t)$, $H_\cdot (t)$, $V_\cdot (t)$ and $W_\cdot (S, \cdot)$ stand for the contour process, the 
height process, the \L{}ukasiewicz path of $t$ and the snake of $S$. 
\begin{definition}
\label{forestdef} Let $t\eqo (t(p))_{1\leq p\leq N}$ be a forest of \emph{finite} trees. Let 
$S\eqo (S(p)\eqo (S_u(p))_{u\in t(p)})_{1\leq p\leq N}$ be $\bbR^d$-valued BRWs starting at the origin, $N\ino \bbN^* \cup \{ \infty\}$. 
We define the contour process $(C_{\! s} (t))_{s\in \bbR_+}$, 
the height process $(H_{\! s} (t))_{s\in \bbR_+}$, the \L{}ukasiewicz path $(V_{\! s} (t))_{s\in \bbR_+}$ of the forest $t$ 
and the snake $(W_{\! s} (S, r))_{r,s\in \bbR_+}$ of $S$ in the following way: we set $\mathtt{r}_{0} \eqo 0$ and $\mathtt{r}_p\eqo \sum_{1\leq q\leq p} \# t(q) $, $1\leqo p\leqo N$; 
for all $s\ino [\mathtt{r}_{p-1}, \mathtt{r}_p]$, we set 
$$C_{\!2s}(t)\eqo C_{\!2(s-\mathtt{r}_{p-1})} (t (p)), \quad  H_{\! s}(t)\eqo H_{\!s-\mathtt{r}_{p-1}} \! (t (p)), \quad V_{\! s}(t)\eqo V_{\!s-\mathtt{r}_{p-1}} \!(t (p))\! -\! p +1, 
$$
and for all $r\ino \bbR_+$, we set $W_{\! 2s}(S, r)\eqo W_{\! 2(s-\mathtt{r}_{p-1})} (S (p), r)$. If 
$N\leko \infty$, then for all $s\ino [\mathtt{r}_N, \infty)$ we also set $C_{\! s } (t)$ $\eqo$  $H_{\! s} (t)\eqo 0$, $V_{\! s}(t)\eqo -N$ and $W_{\! s} (S,r)\eqo 0$, $r\ino \bbR_+$. \cq 
\end{definition} 

\begin{remark}
\label{contord} Let $t\eqo (t(p))_{1\leq p\leq N}$ and $S\eqo \big(S(p)\eqo (S_u(p))_{u\in t(p)}\big)_{1\leq p\leq N}$ be as 
above. Let $t'\ino \bbT$ be the tree associated to the forest $t$ 
as in Definition \ref{Ulamtree} $(c)$. We define $S'\eqo (S'_u)_{u\in t'}$ by setting 
$S'_\varnothing\eqo 0$ and for all $1\leqo p\leqo N$, 
$S'_{[p]\ast u} \eqo S_u(p)$, $u\ino t(p)$. 
In what follows, we denote by $(v_k)_{0\leq k \leq 2(\# t'-1)}$ and $(u_l)_{0\leq l < \#t'}$ the contour exploration and the depth-first exploration of $t'$. 

\medskip

\noi
$(a)$ The respective encoding processes of $(t, S)$ and $(t', S')$ are closely related. Namely, 
observe for all $s,r\ino \bbR_+$ that 
$$C_s (t)\eqo (C_{s+1} (t') -1)_+  \quad \textrm{and} \quad W_{\! s} (S, r)= W_{\! s+1} (S', r+1), $$
and for all integers $0\leqo l \leko  \# t'$, $H_l (t)\eqo (H_{l+1} (t') -1)_+$ and 
$V_{l+1} (t)\eqo V_l (t)+ k_{u_{l+1}} (t) \! -\! 1$. Therefore results on encoding processes of forests can be derived from similar results on trees.

\smallskip

\noi
$(b)$ Contour and height processes are closely related. To simplify, let us assume first that $N\eqo \infty$. 
We observe that 
\begin{equation}
\label{Klcontour}
\forall l\ino \bbN,  \quad K_l:= 2l \! -\! H_{l} (\fftree) \eqo \inf \{ k\ino \bbN: v_{k+1}\eqo u_{l+1}\}
\end{equation}
which increases with $l$. Thus for all $l\ino \bbN$, 
$C_{K_l} (t) \eqo H_{l} (t)$ and for all $s\ino [K_l, K_{l+1} ]$, we get
\begin{equation}
\label{ctrvshght_prime}
C_s (t) = \left\{  \begin{array}{ll}
(H_{l}(t) \! -\! s +K_l)_+ &  \textrm{if $s\ino [K_l, K_{l+1} \! -\! 1)$ } \\
(H_{l+1}(t) \! -\! K_{l+1} +s )_+ & \textrm{if $s\ino [K_{l+1} \! -\! 1, K_{l+1}] $.}
 \end{array} \right.
\end{equation}
See D.~\& Le Gall \cite{DuLG02}, Section 2.4, for more details. 
We can derive $H_\cdot (t)$ from $C_\cdot (t)$ by a time-change $\phi_t \! : \!\bbR_+  \! \to \! \bbR_+ $ that is defined as follows: 
we set $\phi_t(0)\eqo 0$ and for all $l\ino \bbN$ and all $s\ino [0, 1]$, we set $\phi_t(l+s)\eqo K_l+s$ if $|u_{l+1}| \geko |u_l|$ and otherwise we set 
$$\phi_t(l+s) \eqo  K_l + (K_{l+1} \! -\! 1 \! -\!  K_l ) 
((2s)\! \wedge \! 1) + (2s\! -\! 1)_+ \! \eqo   2(l+s) \! -\! H_l(t)+ (H_l(t)\! -\! H_{l+1} (t)) ((2s)\! \wedge \! 1).$$
Then, $\phi_t$ is an increasing continuous function such that 
\begin{equation}
\label{controphit}
\forall s\ino \bbR_+ , \forall l\ino \bbN^*, \quad H_s (t)\eqo C_{\! \phi_t(s)} (t) \quad \textrm{and} \quad \max_{r\in [0, l]} \big| \tfrac{1}{2} \phi_t(r) \! -\! r\big| \leq \tfrac{1}{2} \!\!\!\!\! \sup_{\quad r\in [0, l] }  \!\!\!\!\!H_r (t) \; .
\end{equation}
When the number $N$ of trees in the forest $t$ is finite, the height process $H_\cdot (t)$ is also obtained from $C_\cdot (t)$ by a continuous increasing time-change
$\phi_t \! : \! \bbR_+\! \to \! \bbR_+$ and (\ref{controphit}) holds true too.  \cq 
\end{remark}

\noi
\textbf{Contour and height processes of $\bbR_+$-marked-trees.}
We next define contour processes (and also height processes) 
of trees equipped with variable lifespans. Namely, let $t\ino \bbT$ and $\ell_u \ino \bbR_+$, $u\ino t$.  
Here, $\ell_u$ is interpreted as the \emph{lifespan of $u$} and we 
introduce $\zeta_u \eqo \sum_{v\in \lgeo \varnothing , u\rgeo} \ell_v$ and $\zeta^*_u \eqo \zeta_{\overleftarrow{u}}$, 
which are respectively the \emph{death-time} and the \emph{birth-time} of $u$, with the convention that 
$\zeta^*_{\varnothing}\eqo 0$. 
\begin{definition}
\label{Contlifespan} Let $T\eqo (t, (\ell_u)_{u\in t})$ be as above. 
Let $(v_k)_{0\leq k\leq 2 (\# t-1)}$ and $(u_l)_{0\leq l< \# t }$ be resp.~the contour and the depth-first explorations of $t$. 

\smallskip

\noi
$(a)$ We define $(\Lambda_T(k))_{0\leq k\leq 2\#t}$ by 
$$  \Lambda_T(0)\eqo 0, \; \Lambda_T(1)\eqo \zeta_{v_0} \quad \textrm{and} \quad  \Lambda_T (k+1)\eqo \Lambda_T (k)+ |\zeta_{v_k} \!\! -\! \zeta_{v_{k-1}}|, \quad  1\leqo k \leqo 2(\#t \! -\! 1),  $$   
and $\Lambda_T (2\# t)\eqo \Lambda_T (2\# t\! -\! 1) + \ell_\varnothing$, if $\#t \leko \infty$.

\smallskip

\noi
$(b)$ We define $(L_T (l))_{0\leq l\leq \#t}$ by $L_T(0)\eqo 0$ and $L_T (l)\eqo \sum_{0\leq j< l} \ell_{u_j}$, for all $1\leqo l\leqo \# t$. 
We also set $\mathtt r'\eqo L_T(\# t)$. Note that $\Lambda_T(2\#t)\eqo 2 \mathtt r '$.

\smallskip

\noi

\noi
$(c)$ The \emph{contour of $T$} is the unique 
continuous function $(\mathscr C_s(T))_{s\in \bbR_+\! }$ that satisfies the following conditions: 
$\mathscr C_0 (T)\eqo 0$; $\mathscr C_{\Lambda_T(k)} (T)\eqo \zeta_{v_{k-1}}$, for all $1\leqo k\leko 2\#t $, and 
$\mathscr C_{\Lambda_T(2\# t)} (T)\eqo 0$; $\mathscr C_\cdot (T)$ is affine on $[\Lambda_T (k), \Lambda_T(k+1)]$, for all $0\leqo k\leko2 \# t$;  $\mathscr C_{s} (T)\eqo 0$ for all $s\ino [2\mathtt r'\! , \infty)$.

\smallskip

\noi
$(d)$ The \emph{height process of $T$} is the unique c{\`a}dl{\`a}g function $(\mathscr H_s(T))_{s\in \bbR_+\! }$ that satisfies the following conditions: $\mathscr H_0 (T)\eqo 0$; $\mathscr H_{L_T (l) +s} (T)\eqo \zeta^*_{u_l} +s$, for all $0\leqo l\leko \# t $ and $s\ino [0, \ell_{u_l})$; $\mathscr H_s (T)\eqo 0$, for all $s\ino [ \mathtt r' \! , \infty)$. \cq 
\end{definition}

In the cases where the contour exploration visits the whole tree $t$ and if $\ell_u \ino \bbR_+^*$, then $\mathscr C_\cdot (T)$ or $\mathscr H_\cdot (T)$ completely encode $T$. The height process of $T$ is only used in a technical argument. The fact that it is not continuous is not a problem for our purpose. We next define the contour and the height processes of $\bbR_+$-marked forests. 
\begin{definition}
\label{lifespanforest} Let $\smash{T\eqo \big( T(p)\eqo (t(p), (\ell_u(p))_{u\in t(p)} ) \big)_{1\leq p\leq N}}$ be a forest of finite $\bbR_+$-marked trees, $N$ being possibly infinite.
The contour and height processes of the foret $T$, denoted by 
$(\mathscr C_s(T))_{s\in \bbR_+}$ and $(\mathscr H_s(T))_{s\in \bbR_+}$, are obtained by concatenation of the processes 
$\mathscr C_\cdot (T(p))$ and $\mathscr H_\cdot (T(p))$ and similarly for $\Lambda_T $ and $L_T $.  
Namely, we define $\Lambda_T(0)\eqo L_T(0)\eqo \mathtt r_0\eqo \mathtt{r}'_{0} \eqo 0$ and 
for all $1\leqo p \leqo N$, $\smash{\mathtt r_p\eqo \sum_{1\leq q\leq p} \# t(q)}$, $\smash{\mathtt{r}'_p\eqo \sum_{1\leq q\leq p} \sum_{u\in \btt(q)} \ell_u(q)}$, 
$$ \forall s\ino [\mathtt{r}'_{p-1}, \mathtt{r}'_p], \quad \mathscr C_{2s}(T)\eqo \mathscr C_{2s-2\mathtt{r}'_{p-1}}\!  (T(p)) \quad \textrm{and} \quad \mathscr H_{\!s}(T)\eqo \mathscr H_{\!s-\mathtt{r}'_{p-1}}\!  (T(p)), $$ 
and $\mathscr C_{2s}(T)\eqo \mathscr H_{\!s}\eqo 0$, for all $s\ino [\mathtt r'_N, \infty)$ 
(with an obvious convention if $N\eqo \infty$ and 
$\mathtt r_N'\eqo \infty$). Similarly, we set $\Lambda_T (k+ 2\mathtt r_{p-1}) \eqo \Lambda_T(2\mathtt r_{p-1} )+ \Lambda_{T(p)} (k)$ for all $0\leqo k\leqo 2\# t(p)$, 
and $L_T (l+ \mathtt r_{p-1}) \eqo L_T(\mathtt r_{p-1} )+ L_{T(p)} (l)$ for all $0\leqo l\leqo \# t(p)$ and 
$\Lambda_T$ and $L_T$ are constant on resp.~$[2\mathtt r_N, \infty)$ and $[\mathtt r_N, \infty)$ (with an obvious convention if $N\eqo \infty$ and thus $\mathtt r_N\eqo \infty$). 
Note that $\Lambda_T(2\mathtt r_p)\eqo 2L_T(\mathtt r_p)\eqo 2\mathtt r_p'$\cq
\end{definition}

In the following remark, we discuss the links between $\mathscr C (T)$, $\mathscr H (T)$, $C (T)$, $H (t)$, $\Lambda_T$ and $L_T$. 
\begin{remark}
\label{contlifespanrem} Let $\smash{T\eqo \big( T(p)\eqo (t(p), (\ell_u(p))_{u\in t(p)} ) \big)_{1\leq p\leq N}}$ be as in Definition \ref{lifespanforest} and let $t'\ino \bbT$ be the tree associated with the forest $t\eqo (t(p))_{1\leq p\leq N}$ 
as in Definition \ref{Ulamtree}. We define $T'\eqo (t', (\ell'_u)_{u\in t'})$ similarly by setting 
$$\ell'_\varnothing \! =\!  0, \quad \ell'_{[p]\ast u}  \! =\! \ell_u(p), \; 1\leqo p\leqo N, \, u\ino t(p), \quad \textrm{and} \quad \zeta'_u = \zeta'^*_u+ \ell'_{u}= \!\!  \!   \sum_{v\in \lgeo \varnothing, u\rgeo }\!\!    \ell'_v, \;\,  u\ino t' \backslash \{ \varnothing\}. $$
In what follows, we denote by $(v_k)_{0\leq k \leq 2(\# t'-1)}$ and $(u_l)_{0\leq l < \#t'}$ the contour and depth-first explorations of $t'$. 
When $t'$ is finite (which happens iff $N$ is finite), we observe that $\# t'\! -\! 1\eqo \mathtt r_N$. In this case, it is convenient to set $u_{\# t'}\eqo \varnothing$. 

\smallskip

\noi
$(a)$ Let $\Lambda_T$ be as in Definition \ref{lifespanforest}. We first observe that $\Lambda_T (k+1)\! -\! \Lambda_T(k) \eqo |\zeta'_{v_{k+1}}\!\!\! -\! \zeta'_{v_k}|$, for all $0\leqo k \leko 2\mathtt r_N$. 
We extend $\Lambda_T$ continuously on $\bbR_+$ 
by setting $\Lambda_T (s)\eqo 
\Lambda_T(\lfloor s \rfloor ) + \{ s\} (\Lambda_T(\lceil s \rceil )\! -\! \Lambda_T(\lfloor s \rfloor ))$ for $s\ino \bbR_+$.  
We easily check that 
\begin{equation}
\label{timechangecontourlifespan}
\forall s \ino[0, 2\mathtt r_N) \quad \mathscr C_{\Lambda_T(s)} (T)\eqo \zeta'_{v_{\lfloor s \rfloor }} + 
\{ s\} (\zeta'_{v_{\lceil s \rceil }}\! -\! \zeta'_{v_{\lfloor s \rfloor }} ) .
\end{equation}
In the next sections, the lifespans are most of the time i.i.d.~$\texttt{expo} (1)$ r.v.s, so that $C (t)$ and 
$\mathscr C(T)$ are close. More precisely, for all $0\leqo k \leqo 2\mathtt r_N$, 
(\ref{timechangecontourlifespan}) entails
\begin{equation}
\label{compCTCt1}
\max_{r\in [0, k]} \big| \mathscr C_{\Lambda_T (r)} (T)\! -\! C_{r}(t) \big| \leq 2+
 \max \Big\{ \big| \zeta'_{v_i} \! \! -\! 1\! -\! (|v_i|\! -\! 1)_+ \big|\,  ; \,  1\leqo i\leqo k  \Big\} .
\end{equation}
\emph{Indeed}, we recall from Remark \ref{contord} $(a)$ that $C_k(t)\eqo (C_{k+1} (t')\! -\! 1)_+\eqo (|v_{k+1}| \! -\! 1)_+$ and 
$|C_{r} \! -\! C_{(r-1)_+} (t)| \leqo 1$, for all $r\ino \bbR_+$.

\smallskip

\noi
$(b)$ We next want to compare $\Lambda_T $ with $2L_T$. To that end, we first note 
that 
\begin{equation}
\label{LTident} L_T(l)\eqo \sum_{0\leq j\leq l} \ell'_{u_j} , \quad 0\leqo l\leqo \mathtt r_N\eqo \# t'\! -\! 1.
\end{equation}
We extend $L_T$ continuously on $\bbR_+$ 
by setting $L_T (s)\eqo 
L_T(\lfloor s \rfloor ) + \{ s\} (L_T(\lceil s \rceil )\! -\! L_T(\lfloor s \rfloor ))$, $s\ino \bbR_+$. 
We recall $K_l$ from (\ref{Klcontour}) and $\phi_t$ from Remark \ref{contord} $(b)$. For all $0\leqo l\leko \mathtt r_N$, by definition of $K_l$, $v_{1+K_l}\eqo u_{l+1}$ and $\smash{v_{K_l}\eqo \overleftarrow{u}_{\! l+1}}$. Then, we check that 
$\smash{\Lambda_T (1+\phi_t(l))\eqo\Lambda_T (1+K_l)\eqo 2L_T(l+1) \! -\! \zeta'_{u_{l+1}}}$.
Thus $\smash{\Lambda_T (1+K_l)\! -\!  \Lambda_T (K_l)\eqo 
|\zeta'_{u_{l+1}}\!\!\! -\! \zeta'_{\overleftarrow{u}_{\! l+1}}|\eqo \ell'_{u_{l+1}}\eqo L_T(l+1)\! -\! L_T(l)}$. 
Therefore we get 
\begin{equation}
\label{basicconstat}
\Lambda_T (K_l)\eqo \Lambda_T (\phi_t(l)) \eqo  2L_T(l) \! -\! (\zeta'_{u_{l+1}} \! \! \! -\ell'_{u_{l+1}})= 2L_T(l)\! -\! \zeta'^*_{u_{l+1}} . 
\end{equation} 
Since $L_T$ and $\Lambda_T \circ \phi_t $ are nondecreasing, we get 
\begin{eqnarray*} 
\forall s\ino [l, l+1], \quad \zeta'_{u_{l+2}}\!\!\! -\! \ell'_{u_{l+2}} \!\! \! \!\! - 2\ell'_{u_{l+1}} \! \!\!= 2L_T(l)  \!\!\! \!\! &- & \!\!\! \!\!\Lambda_T(\phi_t(l\! +\! 1))\leq  2L_T(s)\! -\! \Lambda_T(\phi_t(s))\\
\!\!\! & \leqo & \!\!\! 2L_T(l\! +\! 1)\! -\! \Lambda_T(\phi_t(l)) \eqo \ell'_{u_{l+1}}\!\! + \zeta'_{u_{l+1}} .
\end{eqnarray*}
Namely, $\smash{\max_{s\in [0, l]}|2L_T(s)\! -\! \Lambda_T(\phi_t(s))| \leqo 2\max_{1\leq j\leq l}\zeta'_{u_{j}}}$. 
By definition of $K_j$ and by (\ref{timechangecontourlifespan}), we get 
$\smash{ \zeta'_{u_{j+1}}\eqo \zeta'_{v_{1+K_j}}\eqo \mathscr C_{\Lambda_T (1+K_j)} (T)}$. Since $\smash{\Lambda_T (1+K_j) \leqo 2L_T (j+1)}$, we get for all $0\leqo l\leko \mathtt r_N$, 
\begin{equation}
\label{controlLambda}
\max_{s\in [0, l]} \big|2L_T(s)\! -\!  \Lambda_T (\phi_t(s))  \big| \leq 2\!\!\!\!\!\! \!\!\!\!\!\!\!\!  \max_{\quad \quad s\in [0, 2L_T(l)]} \!\!\!\!\!\! \!\!\!\!\!\! \mathscr C_s(T) .
\end{equation}
This extends to all $l\in \bbN$, since $ \Lambda_T $ and $2L_T$ are constant and both equal to $2\mathtt r'_N$  on resp.~$[2\mathtt r_N, \infty)$ and $[\mathtt r_N, \infty)$.

\smallskip

\noi
$(c)$ We view $\smash{\mathscr C(T)}$ as a time-change of $\smash{\mathscr H(T)}$. As already noticed, 
$\smash{v_{1+K_l}\eqo u_{l+1}}$ and $\smash{v_{K_l}\eqo \overleftarrow{u}_{\! l+1}}$. By (\ref{timechangecontourlifespan}), we then get 
$\smash{\mathscr C_{\Lambda_T(K_l)}(T)\eqo \zeta'_{v_{K_l}} \!\! \eqo \zeta'_{ \overleftarrow{u}_{\! l+1}}  \!\!\!\!\eqo \zeta'^*_{u_{l+1}}}$. We also see that 
\begin{equation}
\label{cHindetof}
\mathscr H_{L_T(l)+s}\eqo \zeta'^*_{u_{l+1}} \!\!\!+s, \quad s\ino [0, \ell'_{u_{l+1}})\eqo [0, L_T(l+1)\! -\! L_T(l)).
\end{equation}
By definition, $\smash{\mathscr C(T)}$ is affine on 
$\smash{[\Lambda_T(K_l), \Lambda_T(1+K_l)]}$ and, as already noticed, $\smash{\mathscr C_{\Lambda_T(1+K_l)}(T) }$ $ \eqo$ $ \smash{\zeta'_{u_{l+1}}\eqo \zeta'^*_{u_{l+1}}+ \ell'_{u_{l+1}}}$ and $\smash{\zeta'_{v_{1+K_l}} \!\!\!\! -\! \zeta'_{v_{K_l}} \eqo \ell'_{u_{l+1}}}$. 
Thus $\smash{\Lambda_T(1+K_l)\eqo \Lambda_T(K_l)+ \ell'_{u_{l+1}}}$, and we get $\smash{\mathscr C_{\Lambda_T(K_l)+s}(T)\eqo \zeta'^*_{u_{l+1}}\!\! +s\eqo\mathscr H_{L_T(l)+s}} $, for all $\smash{s\ino [0, \ell'_{u_{l+1}})}$. 
For all $0\leqo l\leko \mathtt r_N$, and for all 
$s\ino [L_T(l), L_T(l+1))$, we set 
\begin{equation} 
\label{defcapK}\mathcal K_T(s)\eqo s\! -\! L_T(l)+\Lambda_T(K_l)
\end{equation}
 and for all $s\ino [\mathtt r_N, \infty)$, we set $\mathcal K_T(s) \eqo 2s$. Then, we have proved that 
\begin{equation}
\label{ccCvsccH}
\forall s\ino \bbR_+, \quad \mathscr C_{\mathcal K_T(s)} (T)= \mathscr H_s(T).  
\end{equation} 
We now compare $\mathcal K_T(s)$ to $2s$. 
Let $l \leko \mathtt r_N$ and $s\ino [L_T(l), L_T(l+1))$. We note that $\smash{2s\! -\! \mathcal K_T(s)}$ $\eqo$ $\smash{s\! -\! L_T(l)+ 2L_T(l)\! -\! \Lambda_T(K_l)\eqo  s\! -\! L_T(l)+ \zeta'^*_{u_{l+1}}\!\! \ino [0,  \zeta'_{u_{l+1}}]}$, by (\ref{basicconstat}). 
Since $\smash{\zeta'_{u_{l+1}}\!\!\! \eqo \mathscr H_{L_T(l+1)-}(T)} $ $\eqo$ $\smash{ \mathscr C_{\Lambda_T (1+K_l)}}$ and since 
$\Lambda_T (1+K_l)$ $ \leqo $ $2L_T(l+1)$, we get 
\begin{equation}
\label{Kvs2id}
\forall l\ino \bbN^*, \quad \max_{s\in [0, L_T(l)]} \big| 2s \! -\! \mathcal K_T(s) \big| \leq \max_{s\in [0, L_T(l)]} \mathscr H_s(T) \leq \max_{s\in [0, 2L_T(l)]}  \mathscr C_s(T)  \, , 
\end{equation} 
which provides a sufficient control of $\mathcal K_T$. 

\smallskip

\noi
$(d)$ We next check that 
\begin{equation}
\label{capKlowupbounds}
\forall s\ino \bbR_+, \quad s\leqo \mathcal K_T (s) \leqo 2s \; .
\end{equation} 
\emph{Indeed}, let $0\leqo l\leko \mathtt r_N$ and $s\ino [L_T(l), L_T(l+1))$. By (\ref{basicconstat}) and (\ref{defcapK}), 
$\smash{\mathcal K_T(s)\! -s\! \eqo } $ $\smash{\Lambda_T(K_l)\! -\! L_T(l)\eqo }$ $\smash{ L_T(l) \! -\! \zeta'^{\ast}_{u_{l+1}}}$. Then, 
we use (\ref{LTident}) to see that $\smash{L_T(l) \eqo \sum_{u\leq_{\mathtt{lex}} u_l} \ell'_u}$. 
By definition $\smash{\zeta^{'\ast}_{u_{l+1}}\eqo \zeta'_{\overleftarrow{u}_{\! l+1}}}$. Since 
$\smash{\overleftarrow{u}_{\! l+1} \! \leq_{\mathtt{lex}}\!  u_l}$, 
we get $\smash{\sum_{u\leq_{\mathtt{lex}} u_l} \ell'_u \! \geq \!  \zeta'_{\overleftarrow{u}_{\! l+1}}}$. 
This proves $s\leqo  \mathcal K_T (s)$. Right after (\ref{ccCvsccH}), we show that 
$\smash{2s\! -\! \mathcal K_T(s) \ino [0,  \zeta'_{u_{l+1}}]}$, which proves  $\mathcal K_T (s) \leqo 2s $. 
If $s\geqo \mathtt r_N$, $\mathcal K_T(s)\eqo 2s$ and (\ref{capKlowupbounds}) is trivially satisfied.  \cq 
\end{remark}

\noi
\textbf{GW-trees and BRWs.} 
\begin{definition}
\label{randmarktreedef} $\!\! (a)$ We denote by $\ccF(\bbT)$ the sigma-field on $\bbT$ generated by the subsets 
$\{ t\ino \bbT\!  : \! v\ino t\}$, where $v$ ranges in $\bbU$. A random tree is a $(\mathscr F, \mathscr F(\bbT))$-measurable function $\tau\! : \! \Omega \! \to \! \bbT$. 

\smallskip

\noi
$(b)$ Let $(E,d)$ be a metric space equipped with its Borel sigma field $\mathscr B(E)$. 
On the space  $\bbT_{\! E}\! :=\!  \{ (t, (s_u)_{u\in t})\, ; \, t \ino \bbT, \, s_u\ino E, u\ino t \}$ of \emph{$E$-marked trees}, we denote by 
$\ccF(\bbT_{\! E})$ the sigma field generated by the sets $\{ (t, (s_u)_{u\in t}) \ino \bbT_{\! E}\! : \! v\ino t \; \textrm{and} \; s_v\ino B \}$, where $v$ ranges in $\bbU$ and $B$ ranges in $\mathscr B(E)$. A random $E$-marked tree is a $(\mathscr F, \ccF(\bbT_{\! E}))$-measurable function 
$\mathcal T\eqo (\tau, (s_u)_{u\in \tau} ) \! : \! \Omega \! \to \! \bbT_{\! E}$. 
Although a $d$-dimensional BRW $(\tau, (S_u)_{u\in \tau})$ is formally a 
$(\ccF, \ccF(\bbT_{\bbR^d}))$-measurable r.v.~we simply denote it 
by $S\eqo (S_u)_{u\in \tau}$. \cq 
\end{definition}
We next briefly recall the definition of GW-trees and forests. 
\begin{definition}
\label{GWfordef} Let $\mu\eqo (\mu(k))_{k\in \bbN}$ be a probability measure on $\bbN$ and let $q\ino \bbR_+^*$. 

\smallskip

\noi
$(a)$ A \textit{Galton--Watson tree with offspring distribution $\mu$} (a \textit{GW($\mu$)-tree}, for short) is a $(\ccF, \ccF(\bbT))$-measurable r.v.~$\tau\! : \! \Omega \!\to \! \bbT$ which satisfies the following. 
\begin{compactenum}
\item[$(i)$] $k_\varnothing (\tau)$ has law $\mu$. 
\item[$(ii)$] For all $k\ino \bbN^*$ such that $\mu(k)\! >\! 0$, the subtrees $\theta_{[1]} \tau, \ldots , \theta_{[k]} \tau$ under $\bP (\, \cdot \, | k_\varnothing (\tau)\! = \! k)$ are independent with the same law as $\tau$ under $\bP$. 
\end{compactenum}

\smallskip

\noi
$(b)$ A \emph{GW($\mu$)-tree with independent and exponentially distributed lifespans with parameter $q$} (a \emph{GW($\mu, q$)-tree} for short) is a $(\ccF, \ccF(\bbT_{\bbR_+}))$-measurable 
r.v.~$\mathcal T\eqo (\tau, (\ell_u)_{u\in \tau}) \! : \! \Omega \!\to \! \bbT_{\bbR_+}$ such that $\tau$ is a GW($\mu$)-tree and conditionally given $\tau$, the r.v.s $\ell_u$, $u\ino \tau$, are independent and exponentially distributed with parameter $q$. 

\smallskip

\noi
$(c)$  A \emph{GW($\mu$)-forest} (resp.~a \emph{GW($\mu,q$)-forest}) is a (possibly infinite) sequence $\tau\eqo (\tau (p))_{1\leq p\leq N}$ (resp.~$\mathcal T\eqo (\mathcal T(p)\eqo (\tau (p), (\ell_u(p))_{u\in \tau(p)} ))_{1\leq p\leq N}$) such that conditionally given $N$, the $\tau(p)$ (resp.~the $\mathcal T(p)$), $1\leqo p\leqo N$ are i.i.d.~GW($\mu$)-trees (resp.~ i.i.d.~GW($\mu,q$)-trees). \cq
\end{definition}

\section{Brownians snakes}
\label{Brosnasec}
\noi
\textbf{Definitions and general results.} 
Recall that $\bC^{_0}_{^{\! d}}$ is the space $\bC (\bbR_+, \bbR^d)$, equipped with the topology of uniform convergence on every compact intervals. 
\begin{definition}  
\label{snadef} Let $h\! : \! \bbR_+ \! \to \! \bbR_+$ be lower semicontinous, and for all $s\ino \bbR_+$, let
$w_s(\cdot) \ino \bC^{_0}_{^{\! d}}$. Then $(h,w)$ is a $\bbR^d$-valued \emph{snake} if the following holds. 
\begin{compactenum}

\smallskip

\item[$(a)$] For all $s\ino \bbR_+$ and for all $r\ino [h(s), \infty)$, $w_s (r)\! = \! w_s (h(s))\! =: \! \widehat{w}_s $. 

\smallskip

\item[$(b)$] For all $s_1, s_2 \ino \bbR_+$, 
$w_{s_1} (r)\! = \! w_{s_2} (r)$ for all $r\ino [0,m_h (s_1, s_2)]$ where we have set $m_h (s_1, s_2)\!   =\!  \inf_{s\in [s_1\wedge s_2, s_1\vee s_2]} h(s)$. 

\smallskip

\end{compactenum}
\noi
The snake $(h,w)$ is said continuous if $h$ is continuous and if the function $s\ino \bbR_+ \! \mapsto \! w_s (\cdot) \ino \bC^{_0}_{^{\! d}}$ is continuous. We denote by $\fSigma_d$ the \emph{space of $d$-dimensional continuous snakes}. \cq 
\end{definition}

\begin{remark}
\label{easy} 
$(a)$ For all $f,f'\ino \bC^{_0}_{^{\! d}}$, we set  
$\delta^{{\mathtt{unif}}}_{^{\! d}} (f,f')\eqo \sum_{p\in \bbN^*} 2^{-p} (1\wedge \max_{s\in [0, p]} |f(s)\! -\! f'(s)|)$. Then, the space 
$(\bC^{_0}_{^{\! d}} , \delta^{\mathtt{unif}}_{^{\! d}} )$ is Polish and $\delta^{\mathtt{unif}}_{^{\! d}} $ metrizes the topology of uniform convergence on every compact intervals. For all $(h,w)$, $(h',w') \ino \bC^{_0}_{^1} \! \times \!  \bC(\bbR_+, \bC^{_0}_{^{\! d}})$, we next set 
\begin{equation}
\label{snaspacedist}
D_d \big(  (h,w), (h',w')\big)\! := \! \delta^{\mathtt{unif}}_{{1}} (h,h') + \sum_{p\in \bbN^*} 2^{-p}\!\! \! \sup_{s\in [0, p]} 
\delta^{\mathtt{unif}}_{{\! d}} \big(w_s (\cdot), w_s'(\cdot)\big).
\end{equation}
Then $(\bC^{_0}_{^1} \! \times \!  \bC(\bbR_+, \bC^{_0}_{^{\! d}}), D_d)$ is a Polish metric space yielding the product topology of 
$(\bC^{_0}_{^{\! 1}}, \delta^{\mathtt{unif}}_1)$ with $\bC(\bbR_+, \bC^{_0}_{^{\! d}})$ equipped 
with the topology of uniform convergence on every compact intervals.  
We easily check that $\fSigma_d$ is closed in $\bC^{_0}_{^1} \! \times \!  \bC(\bbR_+, \bC^{_0}_{^{\! d}})$, equipped with the product topology, and $(\fSigma_d, D_d)$ is thus Polish.  

\smallskip

\noi
$(b)$ If $(h,w)\ino \fSigma_d$, then note that $\widehat{w}\ino \bC^{_0}_{^d}$ and that $(h,w)\ino \fSigma_d \! \mapsto \! \widehat{w}\ino \bC^{_0}_{^d}$ is continuous.

\smallskip

\noi
$(c)$ Let $(h,w)\ino \fSigma_d$ and let $s, r\ino [0, \infty)$. There is a measurable function 
$G_{r,s}$ from $\bC^{_0}_{^1}\! \times \! (\bC^{_0}_{^d})^2$ to $\bbR^d$ such that 
$\smash{G_{r,s} (h,w_0, \widehat{w})\eqo w_{s}(r)}$. \emph{Indeed}, for any function $h'\ino \bC^{_0}_{^1}$ we set 
$\alpha_{r,s}(h')\eqo \sup \{ s^\prime\ino [0, s]: h'(s^\prime)\leko r\wedge h'(s)\}$ with the convention that $\sup \emptyset\eqo 0$. Observe that for all $s'\ino (0, s)$, $\alpha_{r,s}(h')\geko s'$ iff $m_{h'}(s',s)\leko  r\wedge h'(s)$, which implies that $h' \ino \bC^{_0}_{^1}\mapsto \alpha_{r,s} (h')\ino [0, s]$ is measurable. Then we take 
$\smash{G_{r,s} ( h,w_0, \widehat{w})\eqo w_0(r)}$ if $r \leko m_h(0, s)$ and 
$\smash{G_{r,s} ( h,w_0, \widehat{w})\eqo \widehat{w}_{\alpha_{r,s} (h)}}$ otherwise; $G_{r,s}$ is clearly measurable 
and the snake property then entails the desired result. 

Now recall that the Borel sigma-field of $\bC^{_0}_{^d}$ is generated by the functions $f\ino \bC^{_0}_{^d} \!  \mapsto\!  f(r)\ino \bbR^d$, with $r$ ranging in $\bbR_+$. Similarly the  Borel sigma-field of $\bC(\bbR_+, \bC^{_0}_{^{\! d}})$ is generated by the functions $w\ino \bC(\bbR_+, \bC^{_0}_{^{\! d}}) \! \mapsto \! w_s (\cdot) \ino \bC^{_0}_{^d}$, with $s$ ranging in $\bbR_+$. 
This combined with the previous result implies that there is a measurable 
function $G$ from $\bC^{_0}_{^1}\! \times \! (\bC^{_0}_{^d})^2$ to $\fSigma_d$ such that $\smash{G(h,w_0, \widehat{w})\eqo (h, w)}$ for all $(h,w)\ino \fSigma_d$.

\smallskip

\noi
$(d)$ Let $(S_u)_{u\in \fftree}$ be a $\bbR^d$-valued BRW and let $W_\cdot (S, \cdot)$ stand for its discrete snake as in Definition \ref{contsnadef}. Then, we observe that $(C_\cdot (t), W_\cdot (S, \cdot))\ino \fSigma_d$. A similar result holds for the height snake. \cq 
\end{remark}

We next prove that the law of a continuous snake is characterised by the finite dimensional marginal laws of its lifetime process and 
its endpoint process.

\begin{lemma}
\label{endmargcar} Let $(C, W)$ and $(C ', W')$ be two $\fSigma_d$-valued r.v.s such that a.s.~$C_0\eqo C_0'\eqo 0$. 
We assume for any $p\ino \bbN^*$ and any real numbers $0\leqo s_1\leqo \ldots \leqo s_p$ that $\smash{(C_{\!s_j}; \widehat{W}_{\! s_j})_{1\leq j\leq p}}$ has the same law as 
$\smash{(C'_{\!s_j}; \widehat{W}'_{\! s_j})_{1\leq j\leq p}}$. Then $(C,W)$ has the same law as $(C',W')$. 
\end{lemma}

\noi
\textbf{Proof.} Since a.s.~$C_0\eqo C_0'\eqo 0$, $W_0$ and $W'_0$ are a.s.~constant. 
We denote by $h\ino \bC^{_0}_{^1}\! \mapsto \! P_h$ (resp.~$P'_h$) a regular version of the conditional law of $W$ 
given $C$ (resp.~of  $W'$ given $C'$). Namely $h\! \mapsto \! P_h$ and $h\! \mapsto \! P'_h$ are Borel measurable functions from $\bC^{_0}_{^1}$ 
to the Polish space $\mathcal M_1(\fSigma_d)$ of the Borel probability 
measures on $\fSigma_d$ equipped with the topology of 
weak convergence, such that for all Borel subsets $B$ of $\fSigma_d$ a.s.~$\smash{\bE [\un_{B}(W)| C]\eqo P_C (B)}$ 
and $\smash{\bE [\un_B(W')| C']\eqo P_{C'}' (B)}$. 
Let $(C,W'')$ be a $\fSigma_d$-valued r.v.~such that the conditional law of $W''$ given $C$ is $\smash{P'_C }$. 
Since $\smash{(C,\widehat{W})}$ has the same law as $\smash{(C',\widehat{W}')}$, by our assumptions, $C$ and $C'$ have the same 
law and thus $(C',W')$ and $(C, W'')$ have the same law. By Remark \ref{easy} $(b)$, 
$\smash{(C',\widehat{W}')}$ and $\smash{(C,\widehat{W}'')}$ must have the same law. Therefore $\smash{(C,\widehat{W})}$ 
and $\smash{(C,\widehat{W}'')}$ have the same law too. Let $\{ B_n; n \ino \bbN\}$ be a pi-system generating the Borel sigma field of $\bC^{_0}_{^d}$. Thus a.s.~for all $n\ino \bbN$, $\smash{P_C( \{ w\ino \bC(\bbR_+, \bC^{_0}_{^{\! d}})  \! : \! \widehat{w} \ino B_n\} )\eqo P_{C}'(\{ w\ino \bC(\bbR_+, \bC^{_0}_{^{\! d}}) \!  :\!  \widehat{w} \ino B_n\})}$. 
Consequently there is a Borel subset $B$ of $\bC^{_0}_{^1}$ such that $\smash{\bP (C\ino B)\eqo 1}$ 
and such that for all $h\ino B$, $\widehat{w}$ under $P_h(dw)$ has the same law as $\smash{\widehat{w}}$ under 
$\smash{P_h'(dw)}$. By using the function $G$ in Remark \ref{easy} $(c)$, we thus get that $(h,w)$ under $P_h(dw)$ 
has the same law as $(h,w)$ under $\smash{P_h'(dw)}$, which implies the desired result. \cqfd 

\smallskip

In what follows, $\bC(\bbR_+, E)$ stands for 
the space of continuous functions from $\bbR_+$ to a Polish space $(E, d_E)$, which is 
is equipped with the Polish topology of uniform convergence on every compact intervals. 
We generically denote the \emph{$\eta$-modulus of continuity of $f\ino \bC(\bbR_+, E)$ on $[a, b]$} by  
\begin{equation}
\label{etamodu} 
\omega (f , \eta, [a, b]) =\sup \big\{ d_E(f(s), f(s^\prime))\, ; \;  s,s^\prime\ino [a, b]:  \! |s\! -\! s^\prime| \leqo \eta  \big\}. 
\end{equation}
We recall the following 
from D., Khanfir, Lin \& Torri \cite{DuKhLiTo22} and Marckert \& Mokkadem \cite{MaMo03}.
\begin{proposition}
\label{endsnake1} Let $(h,w)$ be a $d$-dimensional snake. Then, the following holds true.
\begin{compactenum}

\smallskip

\item[$(i)$] We suppose that $h\ino \bC^{_0}_{^1}$ and that 
$\widehat{w}\ino \bC^{_0}_{^d}$. Then, $(h,w)\ino \fSigma_d$. 
Moreover, $\smash{\omega (w, \eta, [0, N])} $ $\leqo $ $\smash{2 \omega (\widehat{w}, \eta, [0, N])}$ for all $\eta\ino \bbR_+^*$ and all $N\ino \bbN^*$.

\smallskip

\item[$(ii)$] Let $(h,w)$ and $\smash{(h^{(n)}\! , w^{(n)})\ino \fSigma_d}$, $n\ino \bbN$, be such that 
$\smash{\lim_{n\rightarrow \infty} h^{(n)}\eqo h}$ in $\bC^{_0}_{^1}$, 
$\smash{\lim_{n\rightarrow \infty} w_{^0}^{_{(n)}}}$ $\eqo$ $\smash{ w_0}$ and $\smash{\lim_{n\rightarrow \infty} \widehat{w}^{{(n)}}\eqo \widehat{w}}$ in $\bC^{_0}_{^d}$. Then $\smash{\lim_{n\rightarrow \infty} w^{(n)}\eqo w}$ in $\bC(\bbR_+, \bC^{_0}_{^d})$. 
\end{compactenum}
\end{proposition}
\noi
\textbf{Proof.} For $(i)$, see Lemma 4.18 \cite{DuKhLiTo22}. For $(ii)$, see Theorem 2.1 \cite{MaMo03}. \cqfd 

\smallskip

We next give a specific representation of finite dimensional marginal laws of random snakes which 
is used in the proof of our results. Let $w_0, w_1\ino \bCqq$ and $m\ino \bbR_+$. We define $w_2\! :=\! w_0 \! \oplus_m \! w_1\ino \bCqq$ by setting $w_2(s)\eqo w_0(s)$ if $s\ino [0, m]$ and $w_2(s)\eqo w_0(m)+ w_1(s\! -\! m)\!-\! w_1(0)$ if $s\ino [m, \infty)$ we easily check 
\begin{equation}
\label{contoplus_prel}
(w_0, w_1, m, r)\ino (\bCqq)^2\! \times \! \bbR_+^2\longmapsto  (w_0 \! \oplus_m\!  w_1) (\cdot \wedge r) \ino \bCqq \quad \textrm{is continuous.}
\end{equation}
Let $p\ino \bbN^*$. We denote by $\mathtt{Dom}_p$ 
the space of $(h, (w_j)_{ 0\leq j\leq p}\, ,  (s_j)_{1\leq j\leq p})$ such that $h\ino \bCpl$ is nonnegative, 
the $s_j$ are real numbers such that $ 0\leqo s_1\leko \ldots \leko s_p$ and  
$w_0, \ldots, w_p\ino \bCqq$. We set $s_0\!:=\! 0$ and we define $\mathtt{Marg}_p (h, (w_j)_{ 0\leq j\leq p}\, ,  (s_j)_{1\leq j\leq p})\eqo (w^\prime_j)_{0\leq j\leq p} \ino (\bCqq)^{p+1}$ by setting 
\begin{equation}
\label{defFpmarg}
w^\prime_0\eqo w_0(\cdot\wedge h(0)) \quad \textrm{and} \quad w^\prime_{j+1} \eqo \big( w^\prime_j \oplus_{m_h(s_j, s_{j+1})} w_{j+1} \big) (\cdot \wedge h(s_{j+1})) , \quad 0\leqo j\leko p\; ,
\end{equation}
where we recall that $m_h (s_j, s_{j+1})\!  =\!  \inf_{s\in [s_j, s_{j+1}]} h(s)$. By (\ref{contoplus_prel}) 
and the continuity of the function $(h, (s_j)_{0\leq j\leq p})\! \mapsto \! (h(0), (h(s_{j+1}), m_h (s_j, s_{j+1}))_{0\leq j<p})$, we see that 
\begin{equation}
\label{contoplus}
\mathtt{Marg}_p : \mathtt{Dom}_p
 \longmapsto (\bCqq)^{p+1}  \; \textrm{is continuous.}
\end{equation}

The functions $\smash{(\mathtt{Marg}_p)_{p\in \bbN^*}}$ provide a convenient way of expressing  finite dimensional marginals of snakes. Indeed, let $\smash{\bS\eqo (S_u)_{u\in \btt}}$ be a $\smash{\bbR^d}$-valued BRW starting at the origin and whose jumps 
$\smash{(\xi_u)_{u\in\btt  \backslash \{ \varnothing\}}}$ are conditionally independent given $\btt$, with the same deterministic law $\bgam$. Then for all $\bbN$-\emph{valued times} $0\eqo s_0 \leko s_1\leko \ldots \leko s_p$, we easily check that 
\begin{equation}
\label{margdissna}
\big(W_{\! s_0} (\bS, \cdot ), W_{\! s_1} (\bS, \cdot ), \ldots , W_{\! s_p} (\bS, \cdot)  \big) \; \overset{\textrm{(law)}}{=} \;  \mathtt{Marg}_p \big( C_{ \cdot} (\btt)\,  , 
(U^{_{(j)}}_{\cdot}) _{0\leq j\leq p}  , (s_j)_{1\leq j\leq p} \big)\; , 
\end{equation}
where the $U^{_{(j)}}_{\cdot}$ are independent copies of the continuous affine interpolation of a RW
whose jumps are independent and distributed according to $\bgam$. Note that in (\ref{margdissna}) $W_{\! s_0} (\bS, \cdot )$ is the null process.

We now define Brownian snakes as follows. 
\begin{definition}
\label{Brosnadef} 
A $d$-dimensional Brownian snake is a pair of processes $(C, W)$ such that $C$ is $\bbR_+$-valued and lower semicontinuous, $W$ is a snake with lifetime $C$ as in Definition \ref{snadef} and for all $p\ino \bbN^*$, and all real numbers 
$0\leqo s_1\leko \ldots \leko s_p$, 
\begin{equation}
\label{margBrosna}
\big(W_0,W_{\! s_1}, \ldots , W_{\! s_p} \big) \; \overset{\textrm{(law)}}{=} \;  \mathtt{Marg}_p \big( C \, , (B^{_{(j)}}_{\cdot}) _{0\leq j\leq p}  , (s_j)_{1\leq j\leq p} \big)\; , 
\end{equation}
where $\smash{B^{_{(0)}}_{\cdot} \!\! \eqo W_{\! 0} }$ and where the $\smash{(B^{_{(j)}}_{\cdot})_{1\leq j\leq p}}$ are independent standard $d$-dimensional Brownian motions with initial value $0$ that are independent from $W_0$ and $C$. We say that $(C,W)$ is continuous if, in addition, a.s.~$(C,W)\ino \bSigma_d$. \cq 
\end{definition}
\begin{remark}
\label{Bsnarem} The family of laws provided by the right hand side 
of (\ref{margBrosna}) actually forms 
a consistent family of finite dimensional marginal laws (so that the previous definition makes sense). \cq 
\end{remark}

 As mentioned in the introduction, even if $C$ is a continuous process there is not necessarily a continuous Brownian snake with lifetime process $C$. However we recall from  D.\& Le Gall \cite{DuLG02} the following results on Brownian snakes whose lifetime process is a $\psi$-height process $(H_s)_{s\in \bbR_+}$.
\begin{proposition}
\label{BrosnapsiH} Let $\psi\ino \mathscr L$ satisfy \emph{\texttt{Grey}$(\psi)$}. Let $(H_s)_{s\in \bbR_+}$ be a continuous $\psi$-height process as defined in (\ref{defH}). We set $C_s \eqo H_{s/2} $, $s\ino \bbR_+$. Then the following holds true. 
\begin{compactenum}

\smallskip

\item[$\!\! \!\! (i)\! $] $\! C$ is the lifetime-process of a continuous $d$-dimensional Brownian snake iff $\psi$  satisfies 
\emph{\texttt{Sheu}$_{\,}$($\psi$)}.

\smallskip

\item[$(ii)\!\!$] We assume \emph{\texttt{Sheu}$_{\,}$($\psi$)} and we denote by $W$ a continuous $1$-dimensional Brownian snake with lifetime process $C$. Let $x\ino \bbR_+^*$. We recall from (\ref{infireach}) the definition of $\varsigma_{-x}$. Then for all $z\ino \bbR_+^*$, 
\begin{equation}
\label{exitsys}
\bP \Big( \! \! \! \! \! \! \max_{\quad s\in [0, 2\varsigma_{-x} ]} \! \! \! \! \! \!  \widehat{W}_{\! s} \leqo z  \Big)= e^{-xw(z)} \quad \textrm{where} \quad \int_{w(z)}^\infty \frac{ds}{2\sqrt{\int_0^s \psi (r) dr}}= z \; . 
\end{equation}
\end{compactenum} 
\end{proposition}
\noi
\textbf{Proof.} For $(i)$ see Theorem 4.5.2 \cite{DuLG02}, p.~121. Let us explain $(ii)$. As discussed in Section \ref{4casessec}, the Brownian snake is defined under the excursion measure $\bN$. An elementary argument on Poisson point processes 
implies that $\smash{\bP (\max_{s\in [0, 2\varsigma_{-x} ]} \widehat{W}_{\! s} \leqo z )\eqo e^{-xw(z)}}$ where $\smash{w(z)\eqo \bN (\max \widehat{W} \geko z)}$. As proved by Dynkin \cite{Dyn91} (see also  \cite{DuLG02} p.~131), 
$w$ is $C^2$ on $\bbR_+^*$ and satisfies 
$\frac{_1}{^2} w''(z)\eqo \psi(w(z))$ with boundary conditions $\lim_{z\to 0^+} w(z) \eqo \infty$ and $\lim_{z\to \infty} w(z)\eqo 0$, which entails (\ref{exitsys}). \cqfd

\smallskip

\noi
\textbf{Oscillations of 1-dimensional snakes.} We state and prove a lemma on 
oscillations of snakes, which is used in the proof of Theorem \ref{Sheuexplain}. More precisely let $h\! : \! \bbR_+ \!\! \to \! \bbR_+ $ be lower semicontinous and such that $h(0)\eqo 0$, and let $w$ be a snake with lifetime process $h$. We suppose that $w_0(\cdot )$ is the null function.   
Recall that $m_h(s,s')$ stands for $\inf_{r\in [s\wedge s', s\vee s']} h(r)$. For all 
$z\ino \bbR_+^*$ we recursively define the sequence $(\bsigma_{\! q} (w,z))_{q\in \bbN}$ of successive 
\emph{$z$-oscillation times of $w$} as follows: we set $\bsigma_{0} (w,z)\eqo 0 $ and for all $q\ino \bbN$, 
\begin{equation}
\label{osci1def} \bsigma_{\! q+1} (w, z) \eqo \inf 
\big\{ s\ino (\bsigma_{\! q} (w,z), \infty): \big|\widehat{w}_s\! -\! w_{s} (m_h (\bsigma_{\! q} (w,z),s ) )  \big| \geko z \big\} 
\end{equation}
with the convention that $\inf \emptyset \eqo \infty$. 
\begin{lemma}
\label{osc1lemma} We keep the above notations. We assume that $(h,w)$ is a $\bbR$-valued \emph{continous} snake. Let $s_0, s_1 \ino \bbR_+$ be such that $s_0 \leko s_1$ and let $z\ino \bbR_+^*$. We suppose the following. 
\begin{compactenum}

\smallskip

\item[$(a)$] For all $r\ino [ m_h (s_0,s_1), h(s_0)]$, $ \big| \widehat{w}_{s_0} \! -\! w_{s_0} (r)  \big| \leqo z$. 

\smallskip

\item[$(b)$] $\max_{s\in [s_0, s_1]} \big| \widehat{w}_{s} \! -\! \widehat{w}_{s_0}\big| \geko 7z$. 

\smallskip

\end{compactenum}
Then there exists $q\ino \bbN^*$ such that $s_0 \leko \bsigma_{\! q} (w,z) \leko  \bsigma_{\! q+1}(w,z) \leko  s_1$. 
\end{lemma}
\noi
\textbf{Proof.} To simplify notations we set $\sigma_{\! p} \eqo \bsigma_{\! p} (w,z)$, $p\ino \bbN$. We first prove the following.

\smallskip

\begin{compactenum}
\item[$\mathbf{(A)}$] \emph{We fix $s'_0\ino [0, s_1]$ and set $q'\eqo \max\{ p\ino \bbN^*\! : \! \sigma_{\! p-1} \leqo s'_0\}$, which is well-defined and finite since $\sigma_{ 0}\eqo 0$ and $\lim_{p\to \infty} \sigma_p\eqo \infty$, by continuity of the snake (note that $\sigma_{q'} \geko s'_0$ and that $\sigma_{q'}$ is possibly infinite). Let $z_0\ino \bbR_+^*$. We assume that}
\begin{equation}
\label{endpointcontrol}
\forall r\ino [m_h(s'_0, s_1) , h(s_0')], \quad \big|  \widehat{w}_{s'_0} \! -\! w_{s_0'} (r)\big| \leq z_0 \; .
\end{equation}
\emph{Then for all $s\ino [s_0', s_1  \wedge \sigma_{\! q'} ]$, $ \big|  \widehat{w}_{s'_0} \! -\!\widehat{w}_{s} \big| \leqo z+ (z_0\vee z)$. Here $\sigma_{\! q'}$ is possibly infinite but if we furthermore assume that $\max_{s\in [s'_0, s_1]} \big|  \widehat{w}_{s'_0} \! -\!\widehat{w}_{s} \big| \geko z+ (z_0\vee z)$, then $s_0'\leko \sigma_{\! q'} \leko s_1$. }
\end{compactenum}

\smallskip

\noi
\emph{Proof of $\mathbf{(A)}$.} {\smash Let $s\ino [s'_0,  s_1  \wedge \sigma_{\! q'} )$. Since $s,s'_0\ino [\sigma_{\! q'-1}, \sigma_{\! q'} )$, we get 
\begin{equation}
\label{nonosci}
 \big|  \widehat{w}_{s'_0} \! -\! w_{s_0'} ( m_h(\sigma_{q'-1}, s'_0))\big| \leqo z  \quad \textrm{and} \quad  \big|  \widehat{w}_{s} \! -\! w_{s} ( m_h(\sigma_{q'-1}, s))\big| \leqo z .
\end{equation}
We observe that $m_h(\sigma_{\! q'-1}, s)\eqo  m_h(\sigma_{\! q'-1}, s'_0)\wedge  m_h(s'_0,s)$ and by the snake property, we see in particular that $w_{s'_0} (r)\eqo w_s(r)$ for all $r\ino [0, m_h(\sigma_{\! q'-1}, s)]$. This entails the following.

\noi
$-$ If $m_h(\sigma_{\! q'-1}, s)\eqo  m_h(\sigma_{\! q'-1}, s'_0)$, then  
$\big|  \widehat{w}_{s} \! -\! w_{s_0'} ( m_h(\sigma_{\! q'-1}, s'_0))\big|\eqo \big|  \widehat{w}_{s} \! -\! w_{s} ( m_h(\sigma_{\! q'-1}, s)) \big| \leq z$ by the second inequality in (\ref{nonosci}). Combined with the first one, it entails $  \big|  \widehat{w}_{s'_0} \! -\! \widehat{w}_{s} \big| \leqo 2z$. 

\noi
$-$ If $m_h(\sigma_{\! q'-1}, s)\eqo  m_h(s'_0,s)$, then  
$\big|  \widehat{w}_{s} \! -\! w_{s_0'} ( m_h(s'_0,s))\big|\eqo \big|  \widehat{w}_{s} \! -\! w_{s} ( m_h(\sigma_{\! q'-1}, s)) \big| \leq z$ by the second inequality in (\ref{nonosci}). Since $m_h(s_0',s) \ino [m_h(s_0',s_1), h(s'_0)] $, our assumption implies 
that $\big|  \widehat{w}_{s'_0} \! -\! w_{s_0'} ( m_h(s'_0,s))\big| \leqo z_0$ and we get 
$  \big|  \widehat{w}_{s'_0} \! -\!  \widehat{w}_{s} \big| \leqo z+z_0$, which completes the proof of $\mathbf{(A)}$. \cq

\medskip

We next set $q\eqo \max\{ p\ino \bbN^*: \sigma_{\! p-1} \leqo s_0\}$, which is well-defined and finite.  
We first apply $\mathbf{(A)}$ with $s_0'\eqo s_0$ (which implies that $q\eqo q'$) and with $z_0 \eqo z$: it shows that $s_0\leko \sigma_{\! q} \leko s_1$ and 
\begin{equation}
\label{nouvosci}
\forall s\ino [s_0, \sigma_{\! q}] ,\quad 
\big|  \widehat{w}_{s_0} \! -\! \widehat{w}_{s} \big| \leqo 2z \quad \textrm{and thus} \quad 
\max_{s\in [\sigma_{\! q}, s_1]} \big| \widehat{w}_{s} \! -\! \widehat{w}_{\sigma_{\! q}}\big| \geko 5z.
\end{equation} 
To complete the proof of Lemma \ref{osc1lemma}, we want to apply $\mathbf{(A)}$ again with $s'_0\eqo \sigma_{\! q}$ (which implies that $q'\eqo q+1$) and $z_0\eqo 4z$. Namely, to that end, we have to prove first that 
\begin{equation}
\label{endpointsigma}
\forall r\ino [m_h (\sigma_{\! q}, s_1), h(\sigma_{\! q})], \quad  \big|  \widehat{w}_{\sigma_{\! q}} \! -\! w_{\sigma_{\! q}} (r)\big| \leq 4z .
\end{equation}

\smallskip

\noi
\emph{Proof of (\ref{endpointsigma}).} 
We first fix $r$ in the interval $[m_h(s_0 , s_1), m_h(s_0, \sigma_{\! q})]$ that possibly reduces to a point. Then the snake property implies that $w_{s_0}(r)\eqo w_{\sigma_{\! q}} (r)$. This implies $\big| \widehat{w}_{\sigma_{\! q}} \! -\! w_{\sigma_{\! q}} (r) \big| \eqo  \big|\widehat{w}_{\sigma_{\! q}} \! -\! w_{s_0} (r) \big|$ and thus 
\begin{equation}
\label{ineq3z}
\forall r\ino [m_h(s_0 , s_1), m_h(s_0, \sigma_{\! q})] , \quad
 \big| \widehat{w}_{\sigma_{\! q}} \! -\! w_{\sigma_{\! q}} (r) \big|  \leqo  \big|\widehat{w}_{\sigma_{\! q}} \! -\! \widehat{w}_{s_0} \big|+ \big| \widehat{w}_{s_0} \! -\! w_{s_0} (r) \big| \leqo 3z
\end{equation} 
by (\ref{nouvosci}) and by our assumption. 
We next fix $r\ino [m_h (s_0, \sigma_{\! q}), h(\sigma_{\! q})]$. There exists $s\ino [s_0, \sigma_{\! q}]$ such that $h(s)\eqo r$ and 
$\widehat{w}_s \eqo w_{\sigma_{\! q}} (r)$ which implies $|\widehat{w}_{\sigma_{\! q}} \! -\! w_{\sigma_{\! q}} (r)| \eqo |\widehat{w}_{\sigma_{\! q}} \! -\! \widehat{w}_s \big|$ and thus 
\begin{equation}
\label{ineq4z}
\forall r\ino [m_h(s_0 , \sigma_{\! q}), h(\sigma_{\! q})] , \quad
 \big| \widehat{w}_{\sigma_{\! q}} \! -\! w_{\sigma_{\! q}} (r) \big|  \leqo  \big|\widehat{w}_{\sigma_{\! q}} \! -\! \widehat{w}_{s_0} \big|+ \big| \widehat{w}_{s_0} \! -\! \widehat{w}_{s} \big| \leqo 4z
\end{equation} 
by (\ref{nouvosci}). We finally fix $r\ino [m_h(\sigma_{\! q}, s_1), h(\sigma_{\! q})]$. Since
$m_h(s_0 , s_1)\eqo m_h(s_0, \sigma_{\! q})\wedge m_h (\sigma_{\! q}, s_1)$, there are two cases to consider. 

\smallskip

\noi
$-$ \emph{Case $1$:} $m_h(s_0 , s_1)\eqo m_h(s_0, \sigma_{\! q})$. Then 
$m_h(s_0 , \sigma_{\! q})\leqo m_h(\sigma_{\! q}, s_1)$ and (\ref{ineq4z}) implies (\ref{endpointsigma}).

\smallskip

\noi
$-$ \emph{Case $2$:} $m_h(s_0 , s_1)\eqo m_h(\sigma_{\! q}, s_1)$. Then 
$m_h(\sigma_{\! q}, s_1) \leqo m_h(s_0 , \sigma_{\! q})$. If $r\ino [m_h(\sigma_{\! q}, s_1), m_h(s_0 , \sigma_{\! q})]$ $=$  $[m_h(s_0, s_1), m_h(s_0 , \sigma_{\! q})]$, then (\ref{ineq3z}) implies (\ref{endpointsigma}) and if $r\ino [ m_h(s_0 , \sigma_{\! q}), h(\sigma_{\! q})]$, (\ref{ineq4z}) implies (\ref{endpointsigma}). This completes the proof of  (\ref{endpointsigma}). \cq 

\smallskip

Let us complete the proof of Lemma \ref{osc1lemma}: (\ref{endpointsigma}) allows to apply $(\mathbf{A})$ with 
$s'_0\eqo \sigma_{\! q}$, which implies that $q'\eqo q+1$, and with $z_0\eqo 4z$. Namely, for all $s\ino [\sigma_{\! q}, s_1\wedge \sigma_{\! q+1} ]$, we get $\big|\widehat{w}_{\sigma_{\! q}} \! -\! \widehat{w}_{s} \big| \leqo 5z$ and the second inequality in (\ref{nouvosci}) implies that $\sigma_{\! q+1} \leko s_1$.   \cqfd

}

\section{Preliminary results on BRWs}
\label{prelimBRWsec}

\subsection{Brownian snakes indexed by GW-trees with exponential lifespans} 
\label{RealBrosnaGWsec}

\noi
\textbf{Total height of GW-trees with exponential lifespans.}
In this section, we compare the total height of a discrete GW$(\mu$)-tree and the total height of a GW($\mu, q$)-tree (i.e., GW($\mu$)-trees with i.i.d.~$\mathtt{expo}(q)$ lifespans: see Definition \ref{GWfordef} $(c)$), which behaves in a more regular way. We shall always assume that 
the offspring distribution $\mu$ satisfies (\ref{nontricri}) (namely, it is nontrivial and critical: $\sum_{k\in \bbN} k\mu (k) \eqo 1$ and $\, \mu (0)\! +\!  \mu (1)\leko 1$). 
The estimates that we provide involve functions that are derived from $\mu$ as follows:
\begin{equation}
\label{gfpsiR}
\forall r\ino [0, 1], \qquad g_\mu (r)\eqo \sum_{k\in \bbN} r^k\mu (k), \quad f_\mu (r)\eqo 1\! -\! g_\mu (1\! -\! r) \quad \textrm{and} \quad  \Psimu (r)\eqo  r\! -\!  f_{\mu} (r) 
\end{equation}
and whose basic properties are stated in the following elementary lemma. 
\begin{lemma}
\label{fungeneprop} $f_\mu \! : \! [0, 1] \! \rightarrow \!  [0, 1\! -\! \mu (0)]$ is a strictly concave increasing bijection and  $g_\mu\!  : \! [0, 1]\!  \rightarrow \! [\mu (0), 1]$ and $ \Psimu \! : \! [0, 1] \! \rightarrow \! [0, \mu (0)]$ are 
strictly convex increasing bijections. 
\end{lemma}
\noi
\textbf{Proof.} It is left to the reader. \cqfd

\smallskip

Let $T\eqo (t, (\ell_u)_{u\in t})$ be a single $\bbR_+$-marked finite tree. 
Total heights of $T$ and $t$ are defined by 
$$ \Gamma (t) \eqo \max_{u\in t} |u| \quad \textrm{and} \quad \Gamma (T) \eqo \max_{u\in t}\sum_{v\in \lgeo \varnothing , u \rgeo} \ell_v \; .$$
If $\mathcal T\eqo (\tau, (\ell_u)_{u\in \tau})$ is a GW($\mu, q$)-tree, then 
standard results on continuous-times Markovian Galton-Watson processes imply 
\begin{equation}
\label{semigrey}
\int_{\bP( \Gamma (\mathcal T) >s)}^1 \!\!  \frac{\mathrm dr}{\, \Psimu (r)} = qs \; .
\end{equation}
See e.g.~Athreya \& Ney \cite{AtNe72} Chapter III, Section 3, Equation (7) p.~106 and Section 4, 
Equation (1) p.~107. We first prove the following. 
\begin{lemma}
\label{hghtestimm} Let $\mu$ satisfy (\ref{nontricri}). Let $\tau$ be a GW($\mu$)-tree. 
There is a universal constant $c_0\ino \bbR_+^*$ such that for all $s\ino\bbR_+$
\begin{equation}
\label{hghtestii}
\int_{\! c_0\bP(\Gamma(\tau) >s)}^1\!\!  \frac{\mathrm dr}{\, \Psimu (r)} \, \geq  \tfrac{1}{2}s \quad \textrm{and if $\mu(1) \geko 0$, \quad } \int_{\bP (\Gamma (\tau) > s)}^1 \!\!  \frac{\mathrm dr}{\, \Psimu (r)} \, \leq (s + 1) \frac{\log \frac{1}{\mu(1)}}{1\! -\! \mu(1)} \; .
 \end{equation}
\end{lemma}
\noi
\textbf{Proof.} The first inequality in (\ref{hghtestii}) is proved in Broutin, D.~\& Wang \cite{BrDuWa21}, Lemma 7.4, but somehow implicitly. Let us provide here a brief and self contained argument. 
Let $\mathcal T\eqo (\tau, (\ell^1_u)_{u\in \tau})$, where 
$\tau$ is a GW($\mu$)-tree and conditionally given $\tau$, the $(\ell^1_u)_{u\in \tau}$ are independent 
$\mathtt{expo} (1)$ r.v.s. Let $(\mathtt e_l)_{l\in \bbN^*}$ be independent $\mathtt{expo} (1)$ r.v.s. We assume that 
$(\mathtt e_l)_{n\in \bbN^*}$ and $\mathcal T$ are independent and we set 
$S_l\eqo \sum_{1\leq j\leq l}\mathtt e_j$, $l\ino \bbN^*$. By the law of large numbers a.s.~$S_l/l \! \to \! 1$, 
which implies that $c_0\! :=\!  \inf_{l\in \bbN^*} \bP (S_l\geko \frac{_1}{^2}l) \geko 0$. Let $u_*$ be the $<_{\mathtt{lex}}$-least 
vertex of $\tau$ such that $|u_*|\eqo \Gamma (\tau)$. Then, 
$\Gamma (\mathcal T) $ $ \geqo$ $ \zeta_{u_*}$ $\! :=\! $ $\sum_{u\in \lgeo \varnothing , u_*\rgeo} \ell^1_u$, and conditionally given $\tau$, $\zeta_{u_*}$ has the same law as 
$S_{1+|u_*|}$. Thus, for all $s'\ino [ \frac{_1}{^2}, \infty)$, we a.s.~get $\bP (\Gamma (\mathcal T) \geko s' | \tau) 
\geqo \bE [\un_{\{ \zeta_{u_*}> s'\}}| \tau]\eqo 
\bE [\un_{\{ S_{1+|u_*|} > s'\}}| \tau] \geqo \bE [\un_{\{ S_{1+|u_*|} > (1+|u_*|)/2\}}| \tau] \un_{\{ 1+ |u_*| > 2s'\}} $. 
Therefore,  
$\bP ( \Gamma (\mathcal T) \geko s') \geqo c_0\bP (\Gamma (\tau) \geko 2s'\! -\! 1)$. Namely, for all $s\ino \bbR_+$, 
$\bP (\Gamma (\mathcal T) \geko  \frac{_1}{^2} (s+1))\geqo c_0 \bP (\Gamma (\tau)\geko s)$ and (\ref{semigrey}) 
implies the first inequality of (\ref{hghtestii}).

Let us prove the second one. 
We assume $\mu(1)\geko 0$ and we denote by $\mu_{\mathrm{pr}}$ the \emph{proper} offspring 
distribution associated with $\mu$: namely $\mu_{\mathrm{pr}}(k)\eqo \mu (k)/ (1\! -\! \mu(1))$ if 
$k\ino \bbN \backslash \{ 1\}$ and $\mu_{\mathrm{pr}} (1)\eqo 0$. It is non-trivial, critical and 
$\Psi_{\! \mu_{\mathrm{pr}}} (\cdot) \eqo \Psimu (\cdot) /(1\! -\! \mu(1))$. 
We take $q\eqo \log 1/\mu (1)$ and we consider 
$\mathcal T_{\! \mathrm{pr}} \eqo(\tau_{\mathrm{pr}},(\ell_u)_{ u\in \tau_{\mathrm{pr}}} ) $, 
where $\tau_{\mathrm{pr}}$ is a GW($\mu_{\mathrm{pr}} $)-tree and where conditionally 
given $\tau_{\mathrm{pr}}$, the $\ell_u$ are i.i.d.~expo($q$). We also introduce 
$\mathcal T^*_{\! \mathrm{pr}} \eqo (\tau_{\mathrm{pr}},( \lceil \mathcal \ell_u\rceil)_{ u\in \tau_{\mathrm{pr}}} )$ 
and $\tau^*\ino \bbT$ such that $k_\varnothing (\tau^*)\eqo 1$ and $\theta_{[1]} \tau^*\eqo \tau$. 
Then observe that $\mathcal T^*_{\! \mathrm{pr}}$ and $\tau^*$ are geometrically the same and their total 
heights have the same law. Therefore 
$ \Gamma (\mathcal T_{\! \mathrm{pr}} ) \leqo  \Gamma (\mathcal T^*_{\! \mathrm{pr}} )$ and 
$\Gamma (\mathcal T^*_{\! \mathrm{pr}} )$ 
has the same law as 
$\Gamma (\tau^*)\eqo \Gamma (\tau)+1 $. Then, for all $s\ino [1, \infty)$ we get 
$ \smash{\int_{\bP (\Gamma (\tau) >s)}^1 \!\!  \frac{\mathrm dr}{\, \Psi_{\! \mu_{\mathrm{pr}}}  (r)}
\leqo \int_{\bP (\Gamma (\mathcal T_{\mathrm{pr}}) >s+1)}^1 \!\! 
\frac{\mathrm dr}{\, \Psi_{\! \mu_{\mathrm{pr}}}  (r)} \eqo  (s+1) \log\tfrac{1}{\mu(1)} } $ by (\ref{semigrey}),
 which easily completes the proof of (\ref{hghtestii}). \cqfd

\smallskip

\noi
\textbf{Brownian snakes indexed by discrete trees with variable lifespans.}
Let $T$ $\eqo$  $\big(T(p) $ $ \eqo$  $ (t(p),$ $ (\ell_u(p))_{u\in t(p)})\big)_{1\leq p\leq N}$ be 
a forest of $\bbR_+$-marked finite trees. We recall from 
Definitions \ref{Contlifespan} and \ref{lifespanforest} that $(\mathscr C_s(T))_{s\in \bbR_+}$ 
stands for the contour process of $T$ and that 
$(\mathscr H_s(T))_{s\in \bbR_+}$ stands for its height process, which is c{\`a}dl{\`a}g. 
As explained in Remark \ref{contlifespanrem} $(c)$, there exists an increasing 
time-change $\mathcal K_T \! :\! \bbR_+ \! \to \! \bbR_+$ such that (\ref{ccCvsccH}) holds: $\mathscr C_{\mathcal K_T (s)} (T)\eqo \mathscr H_s(T)$, for all $s\ino \bbR_+$. 
\begin{definition}
\label{1brosnadistree} We keep the above notations and we denote by $(\mathscr W_{s} (T,r))_{r,s\in \bbR_+}$ 
the $1$-dimensional Brownian snake with initial value $0$ and lifetime process $\mathscr C_\cdot(T)$. We also set 
\begin{equation}
\label{defviatimechange}\mathtt W_s(T,r)= \mathscr W_{\mathcal K_T(s)} (T,r), \quad r,s\ino \bbR_+, 
\end{equation}
which we call the  \emph{$1$-dimensional height Brownian snake} with initial value $0$. \cq 
\end{definition}
\begin{remark}
\label{constrheightBrosna}
$(a)$ Note that $s\ino \bbR_+\! \mapsto \!  \mathtt W_s(T, \cdot) \ino \bC^{_0}_{^1}$ is c{\`a}dl{\`a}g and that actually 
it is distributed as a $1$-dimensional Brownian snake with initial value $0$ and lifetime process $\mathscr H_\cdot(T)$. Indeed, $\mathscr H_\cdot (T)$ is c{\`a}dl{\`a}g without positive jump, hence it is lower semicontinuous, as required in Definition \ref{Brosnadef} of Brownian snakes. 

\smallskip

\noi
$(b)$ We shall need a more direct construction of $1$-dimensional height Brownian snake of $T$. To that end, we denote by $T'\eqo (t', (\ell'_u)_{u\in t'})$ the $\bbR_+$-marked tree associated with the forest $T$ as explained in Remark \ref{contlifespanrem} and Definition \ref{Ulamtree} $(c)$. Namely, for all $1\leqo p\leqo N$, $\theta_{[p]} t'\eqo t(p)$, $\ell'_{[p]\ast u} \! :=\!  \ell_u(p)$, for all $u\ino t(p)$ and $\ell'_\varnothing\eqo 0$. We recall that the death-time of $u\ino t'$ is $\zeta'_u\eqo \sum_{v\in \lgeo \varnothing, u\rgeo} \ell'_v$ and its birth-time is denoted by $\zeta'^*_u\eqo \zeta'_u \! -\! \ell'_u$. We denote by $(u_l)_{0\leq l\leq \# t'\! -\! 1}$ the sequence of the vertices of $t'$ listed in increasing depth-first order, i.e., lexicographical order. We recall the process $(L_T (l))_{l\in \bbN}$ 
from Definitions \ref{Contlifespan} and \ref{lifespanforest}. We recall from (\ref{LTident}) 
that $ \# t'\! -\! 1\eqo \mathtt r_N\eqo \sum_{1\leq p\leq N} \# t(p)$ and that for all $0\leqo l\leqo \mathtt r_N$,  
$L_T (t) \eqo \sum_{0\leq k\leq l} \ell'_{u_j}$. 
We recall from (\ref{cHindetof}) that 
$\mathscr H_{L_T(l) +s} (T)\eqo \zeta_{u_{l+1}}'^*\!\!\! +s$, for all $0\leqo l \leko \mathtt r_N$ and $s\ino [0 ,\ell'_{u_{l+1}})$.   

By definition of $\mathtt W_{\cdot} (T,\cdot)$, the paths 
$s\ino [0, \ell'_{u_{l+1}}) \mapsto \widehat{\mathtt W}_{L_T(l) +s}(T) \! -\! \widehat{\mathtt W}_{L_T(l)} (T)$ are independent 
stopped Brownian motions starting at $0$. Note that the left-limit of $\smash{\widehat{\mathtt{W}}_\cdot }$ 
at $\smash{\ell'_{u_{l+1}}}$ exits and is equal to 
$\smash{\widehat{\mathtt W}_{L_T(l+1)-}(T)}$.  
We now extend these Brownian paths thanks to other independent Brownian motions as follows.  
Let $(B'_s(u))_{s\in \bbR_+}$, $u\ino \bbU$, be a family of auxiliary i.i.d.~$1$-dimensional Brownian motions with initial value $0$, which are furthermore supposed to be independent from $\mathtt W_\cdot( T, \cdot)$. For all 
$0\leqo l\leko \mathtt r_N$ and all $s\ino \bbR_+$, we set 
\begin{equation}
\label{pathudef}
B_s(u_{l+1})= \left\{ \begin{array}{ll}
\widehat{\mathtt W}_{ L_T(l)+s} (T) \! -\! \widehat{\mathtt W}_{L_T(l)} (T)  & \textrm{if $s \leko \ell'_{u_{l+1}}$,} \\
\widehat{\mathtt W}_{L_T(l+1)-}(T)  \! -\! \widehat{\mathtt W}_{L_T(l)} (T)  + B'_{s-\ell'_{u_{l+1}}} (u_{l+1})& \textrm{if $s \geqo \ell'_{u_{l+1}}$,}
  \end{array} \right.
\end{equation} 
We also set $B_\cdot (\varnothing)\eqo B'_\cdot (\varnothing)$, which plays no role in what follows.   
This defines the paths $(B_s(u))_{s\in \bbR_+}$ for all $u\ino t'$. 
We recall from Definition \ref{randmarktreedef} that $\smash{ \bbT_{\bbR_+\times \bC^{_0}_{^1}} }$ stands for the space of $\smash{(\bbR_+\times \bC^{_0}_{^1})}$-marked ordered rooted discrete trees  and we observe the following. 
\begin{compactenum}

\smallskip

\item[$(i)$] The $(B_s(u))_{s\in \bbR_+}$, $u\ino t'$, are i.i.d.~$1$-dimensional Brownian motions with initial value $0$. 

\smallskip

\item[$(ii)$] There is a measurable function $F\! : \! \bbT_{\bbR_+\times \bC^{_0}_{^1}} \! \to \! \bD (\bbR_+,  \bC^{_0}_{^1})$ such that 
\begin{equation}
\label{iidpathssnaconstr}
 F \big( \Theta \big)  = \big( \mathtt W_s(T, \cdot) \big)_{s\in \bbR_+} .
\end{equation} 
where we have set $\Theta\eqo  \big( t', (\ell'_u, B_\cdot (u) )_{u\in t'} \big) $.  \cq 
\end{compactenum}
\end{remark}

We next discuss useful properties of the contour-Brownian snake of the forest $T$. For all integers $1\leqo p \leqo N$, 
we recall the notation $\mathtt r_p\eqo \sum_{1\leq q\leq p} \# t(q)$ and we denote by $(v_k)_{1\leq k\leq 2\mathtt r_N}$ 
the contour exploration of the tree $t'$ associated with the forest $(t(p))_{1\leq p\leq N}$. We also recall from 
Definitions \ref{Contlifespan} and \ref{lifespanforest}, $\Lambda_T\! :\! \bbR_+\! \to \! \bbR_+$, which is continuous, 
increasing and such that $\Lambda_T(0)\eqo 0$ and $\Lambda_T(s+k\! -\! 1)\eqo \Lambda_T(k\! -\! 1)+ 
s|\zeta_{v_k} \!\! -\! \zeta_{v_{k-1}}|$, $s\ino [0, 1]$, $1\leqo k\leqo 2\mathtt r_N$. We recall from Remark \ref{constrheightBrosna} $(b)$ the definition of the independent Brownian motions $(B_\cdot (u))_{u\in t'}$. We also recal the notation $t'$, $\ell_u', \zeta'_u$, $u\ino t'$.   
\begin{lemma}
\label{compareBrowsnakeBRW} 
We keep the previous notation and we assume that $((\ell_u(p))_{u\in t(p)})_{1\leq p\leq N}$ are 
independent exponentially distributed r.v.s with mean $1$. Let $\mathbf S\eqo 
((S_u(p) )_{u\in t(p)})_{1\leq p\leq N}$ be the BRW whose jumps are given by 
$\smash{\xi_u(p)\! :=\! B_{\ell'_{{[p]  \ast  u} }}\!  ([p ]  \ast  u)}$, $u\ino t(p)\backslash 
\{ \varnothing\}$, $1\leqo p\leqo N$. We recall that $\smash{(\widehat{W}_{\! s}(\mathbf S))_{s\in \bbR_+}}$ 
is the endpoint process of the snake of $\mathbf S$ as in Definition \ref{contsnadef}, and that $\smash{(\widehat{\mathscr W}_{\! s} (T))_{s\in \bbR_+}}$ is the endpoint process of the contour-Brownian 
snake of $T$ as in Definition \ref{1brosnadistree}. Then, the following holds true.

\begin{compactenum}

\smallskip

\item[$(i)$] The jumps $((\xi_u(p))_{u\in t(p)\backslash \{ \varnothing\}})_{1\leq p\leq N}$ are i.i.d.~with law $\bgam (dx)\eqo \frac{1}{\sqrt{2}} e^{-|x|\sqrt{2}} dx$. 

\smallskip

\item[$(ii)$] Let $p\ino \bbN^*$ such that $p\leqo N$. For all $z\ino \bbR_+$, 
\begin{equation}
\label{closeWWexpli}
\bP \Big( \max_{s\in [0, 2\mathtt r_p]} \big|\widehat{W}_{\! s}(\mathbf S) \! -\! \widehat{\mathscr W}_{\Lambda_T (s)} (T) \big|\geko z\sqrt{2}\Big)\leq 1-\Big( 1-2e^{-z/3}\Big)^{\! \mathtt r_p} . 
\end{equation}

\smallskip

\item[$(iii)$] For all $z\ino \bbR_+$, 
$ 1\!-\! \big( 1\! -\!  e^{-z\sqrt{2}} \big)^{\mathtt r_p -p}  \leq  2\bP \big( \max_{s\in [0, 2\mathtt r_p]} 
\widehat{W}_{\! s} (\mathbf S) \geko \tfrac{1}{2}z\big) \; .$ 

\smallskip

\item[$(iv)$] For all $z\ino \bbR_+$, 
$ 1\!-\! \big( 1\! -\!  e^{-z\sqrt{2}} \big)^{\mathtt r_p}  \leq  2\bP \big( \max_{s\in [0, \Lambda_T (2\mathtt r_p)]} \widehat{\mathscr W}_s (T) \geko \tfrac{1}{2} z\big) \; .$ 
\end{compactenum}
\end{lemma}
\noi
\textbf{Proof.} Since the r.v.s $\smash{(\ell'_u, B_\cdot (u))_{u\in t'\backslash \{ \varnothing\}}}$ are i.i.d., so are the jumps of $\mathbf S$. Moreover, since for all $u\ino t'\backslash \{ \varnothing \}$, $\ell'_u $ and $B_\cdot (u)$ are independent, an elementary computation entails that the law of the jumps are $\bgam$ (see e.g.~p.~153, Formula 1.0.5 in Borodin \& Salminen \cite{BoroSalm}). This proves $(i)$.

To prove $(ii)$, we extend the BRW and its jumps to $t'$ by setting $\smash{\xi'_u \eqo S'_u\eqo  0}$ if $|u| \ino \{ 0, 1\}$, and 
$\smash{\xi'_{[p]\ast u}} $ $\! :=\!$ $ \xi_u(p)$ and $\smash{S'_{[p]\ast u}\! :=\!  S_u(p)}$ for all $u\ino t(p)$, 
$1\leqo p\leqo N$. For all $\smash{u\ino t'\backslash \{ \varnothing \}}$, we set $\smash{J_u}$ $ \eqo$ $\smash{\max_{s\in [0, \ell'_u]} |B_s(u)|}$. 
The previous arguments on independence imply that the $\smash{(J_u)_{u\in t'\backslash \{ \varnothing \}}}$ are independent 
and distributed as $\smash{\max_{s\in [0, \mathcal E]} |B_s|}$, where $\mathcal E$ is an $\mathtt{expo} (1)$ r.v.~and 
$\smash{B_\cdot}$ is an independent $1$-dimensional Brownian motion starting at $0$. Standard computations 
(see e.g.~Formula 1.1.2, p.~153 in Borodin \& Salminen \cite{BoroSalm}) show that $\smash{\max_{s\in [0, \mathcal E]} B_s}$ is a 
$\smash{\mathtt{expo} (\sqrt{2})}$ r.v., which implies that 
$\smash{\bP (\max_{s\in [0, \mathcal E]} |B_s| \geko z) \leqo 2 \exp (-z\sqrt{2})}$ since $\smash{(-B_s)_{s\in [0, \mathcal E]}}$ 
has the same distribution as $\smash{(B_s)_{s\in [0, \mathcal E]}}$. 
We set $\mathcal D\! :=\! \max_{1\leq q \leq p} \max_{u\in t(q)} J_{[q] \ast u} $ and we thus get 
\begin{equation}
\label{Deltadelta}
\bP (\mathcal D \geko z)= 1- \big( 1- \bP \big(\max_{s\in [0, \mathcal E]} |B_s|  \geko z \big) \big)^{\mathtt r_p} \leqo 1-\big( 1-2e^{-z \sqrt{2}} \big)^{\mathtt r_p}. 
\end{equation}

We next fix $1\leqo k\leqo 2\mathtt r_p$ and we first observe that $\smash{\widehat{W}_{\! k-1} (\mathbf S)\eqo S'_{v_k}}$. If 
$v_k\ino t(q)$, we then note that $\smash{\widehat{\mathscr W}_{\Lambda_T(k)}\eqo S'_{v_k}+   B_{\ell_{[q]}}([q]) 
\eqo \widehat{W}_{\! k-1} (\mathbf S) + B_{\ell_{[q]}}([q])}$. Thus 
$$\big|  \widehat{W}_{\! k-1} (\mathbf S) \! -\! \widehat{\mathscr W}_{\Lambda_T(k-1)} (T)\big| \leq
 \big|  \widehat{\mathscr W}_{\Lambda_T(k)} (T) \! -\! \widehat{\mathscr W}_{\Lambda_T(k-1)}(T) \big| + \mathcal D.$$
Moreover, for all $r\ino [0, 1]$, $\smash{\widehat{W}_{k-1+r} (\mathbf S) \! -\!  \widehat{W}_{k-1} (\mathbf S)  \eqo r\xi'_{v_{k+1}}\un_{\{ v_k\eqo \overleftarrow{v}_{\! k+1}\}} \! -  r\xi'_{v_{k+1}}\un_{\{ v_{k+1}\eqo \overleftarrow{v}_{\! k}\}}     }$. 
Thus, we get $\smash{\max_{s\in [k-1,k]}} $ $\smash{ | \widehat{W}_{\! s} (\mathbf S)}$ $\smash{ \! -\!  \widehat{W}_{\! k-1} (\mathbf S) | \leqo \mathcal D}$. Finally, we observe for all $\smash{ r\ino [ 0, |\zeta'_{v_k}\!\! -\! \zeta'_{v_{k-1}} | ]}$ that 
$$ 
\widehat{\mathscr W}_{\Lambda_T(k-1) +r} (T)-\widehat{\mathscr W}_{\Lambda_T(k-1)} = \left\{ \begin{array}{ll}
 B_r(v_k)  & \textrm{if $v_{k-1}\eqo \overleftarrow{v}_{\!\!k}$,}\\
B_{\ell'_{v_{k-1}}\!\! -r \, } (v_{k-1}) - B_{\ell'_{v_{k-1}}} (v_{k-1}) & \textrm{if $v_{k}\eqo \overleftarrow{v}_{\!\! k-1}$.}
\end{array} \right.
 $$
Thus, $\smash[t]{\max_{s\in [k-1, k]} |\widehat{\mathscr W}_{\! \Lambda_T(s)} (T) \! -\! \widehat{\mathscr W}_{\! \Lambda_T(k-1)} (T)| 
\leqo J_{v_k} \! \un_{\{ v_{k-1}\eqo \overleftarrow{v}_{\! k} \} }+ 2J_{v_{k-1}} \un_{\{v_{k}\eqo \overleftarrow{v}_{\! k-1} \} }\leqo 
2\mathcal D}$. 
By the previous inequalities, we get
$ \smash{\max_{s\in [0, 2\mathtt r_p]}}$ $\smash{\big|\widehat{W}_{\! s}(\mathbf S)}$ 
$ \! -\!$ $\smash{ \widehat{\mathscr W}_{\Lambda_T (s)} (T) \big| \leqo 6 \mathcal D}$ and thus (\ref{closeWWexpli}) by (\ref{Deltadelta}). 

To prove $(iii)$, we first note that 
$$\max_{1\leq q\leq p\; } \max_{u\in t(q)\backslash \{ \varnothing \}} |\xi_u(q)| \eqo \max_{1\leq k\leq 2\mathtt r_p}
|\widehat{W}_{\! k} (\mathbf S)\! -\! \widehat{W}_{\! k-1} (\mathbf S)| \leq 2\!\! \max_{s\in [0, 2\mathtt r_p]}\! 
|\widehat{W}_{\! s} (\mathbf S)| .$$
Since $\bgam$ is symmetric, $\smash{(-\widehat{W}_{\! s} (\mathbf S))_{s\in \bbR_+}}$ has the same law as 
$\smash{(\widehat{W}_{\! s} (\mathbf S))_{s\in \bbR_+}}$. Thus 
\begin{eqnarray*}
1\!-\! \big( 1\! -\!  e^{-z\sqrt{2}} \big)^{\mathtt r_p -p} \!\!\! &= &\!\!\! 
\bP \big( \max_{1\leq q\leq p} \max_{u\in t(q)\backslash \{ \varnothing \}} |\xi_u(q)| \geko z \big) \\
\!\!\! & \leq & \!\!\!  \bP \big( \max_{s\in [0, 2\mathtt r_p]} |\widehat{W}_{\! s} (\mathbf S)| 
\geko  \tfrac{1}{2}z \big) \leqo 2\bP \big( \max_{s\in [0, 2\mathtt r_p]} \widehat{W}_{\! s} (\mathbf S) \geko \tfrac{1}{2}z \big), 
\end{eqnarray*}
which is $(iii)$. Similarly, to prove $(iv)$ we see that 
$$\max_{1\leq q\leq p\; } \max_{u\in t(q)} \big| B_{\ell_{[q]  \ast  u}} ([q] \! \ast \! u) \big| \eqo 
\max_{1\leq k\leq 2\mathtt r_p} |\widehat{\mathscr W}_{\! \Lambda_T(k)} (T)\! -
\! \widehat{\mathscr W}_{\! \Lambda_T(k-1)} (T)| \leq 2\!\! \max_{s\in [0, 2\mathtt r_p]} \!
 |\widehat{\mathscr W}_{\! \Lambda_T(s)} (T)| $$
Since $\smash{(-\widehat{\mathscr W}_{\! s} (T))_{s\in \bbR_+}}$ has the same law as 
$\smash{(\widehat{\mathscr W}_{\! s} (T))_{s\in \bbR_+}}$, we get 
\begin{eqnarray*}
 1\!-\! \big( 1\! -\!  e^{-z\sqrt{2}} \big)^{\mathtt r_p } \!\!\! &= &\!\!\! 
\bP \big(\max_{1\leq q\leq p} \max_{u\in t(q)} \big| B_{\ell_{[q]  \ast  u}} ([q] \! \ast\!  u) \big| \geko z \big) \\
\!\!\! & \leq & \!\!\!  \bP \big(\max_{s\in [0, 2 \mathtt r_p]} |\widehat{\mathscr W}_{\! \Lambda_T(s)} (T)| \geko \tfrac{1}{2} z \big) \leqo 2\bP \big( \max_{s\in [0, 2\mathtt r_p]} \widehat{\mathscr W}_{\! \Lambda_T(s)} (T)  \geko \tfrac{1}{2}z \big), 
\end{eqnarray*}
which implies $(iv)$. \cqfd

\smallskip

\noi
\textbf{Maximal displacement of GW-indexed Brownian snakes.} We keep the previous notation and we introduce the process $(\mathscr Z_a (T))_{a\in \bbR_+}$, which is the \emph{exit branching process} of the snakes $\mathscr W_\cdot (T, \cdot)$ (or of  $\mathtt W_\cdot (T, \cdot)$). Namely, $\mathscr Z_a (T)$ is the number of individuals in $T$ that are the first among their ancestors to reach the value $a$. More precisely, for any function 
$w\ino \bC^{_0}_{^1}$ such that $w(0)\eqo 0$, we set $\varsigma_a (w)\eqo \inf \{ r\ino \bbR_+: w(r) \geqo a \}$, with the convention $\inf \emptyset \eqo \infty$. Then 
$$\mathscr Z_a (T) \eqo 
\# \big\{s\ino \bbR_+: \varsigma_a(\mathtt W_s (T, \cdot) )\eqo \mathscr H_s (T)  \big\} $$
which is a finite quantity. Note that a.s.~$\smash{\mathscr Z_a (T)\eqo \frac{_1}{^2}\# \big\{s\ino \bbR_+: \varsigma_a(\mathscr W_s (T, \cdot) )\eqo \mathscr C_s  (T) \big\}}$ 
because edges are crossed twice in contour order but once only during the depth-first exploration. We then see that 
$ \smash{\bP (\max_{s\in \bbR_+} \! \widehat{\mathscr W}_s(T)  \geqo a) \eqo \bP (\mathscr Z_a (T) \! \neq \! 0)}$. 
We recall from Kaj \& Salminen \cite{KajSal93} the following lemma which computes this probability 
when $\mathcal T$ a GW($\mu, q$)-tree. 
\begin{lemma}
\label{maxBrosnadiscr} Let $\smash{\mathcal T\eqo (\tau, (\ell_u)_{u\in \tau})}$ be a GW$(\mu, q)$-tree as in Definition \ref{GWfordef} 
$(b)$. To simplify notation we set $\smash{\mathscr W\eqo \mathscr W_\cdot (\mathcal T, \cdot) }$, $\smash{\mathtt W\eqo \mathtt W_\cdot (\mathcal T, \cdot) }$ and $\smash{\mathscr Z_\cdot \eqo \mathscr Z_\cdot  (\mathcal T)}$. 
For all $r\! \in [0, 1]$ and all $a\ino \bbR_+$, we also set 
$\smash{w_r(a)\eqo 1\! -\! \bE [r^{\mathscr{Z}_a}]}$ and $\smash{w(a) \eqo \lim_{r\to 0^+} w_r(a)\eqo 
 \bP (\max_{s\in \bbR_+} \!\!  \widehat{\mathscr W}_{ s}  \geqo a)}$. Then, the following holds.  
 \begin{equation}
\label{exitBroSna}
\int_{w_r(a)}^{1-r}\frac{dy}{2\sqrt{\int_0^y \Psi_{\! \mu} (z) dz}}= \int_{w(a)}^{1} \frac{dy}{2\sqrt{\int_0^y \Psi_{\! \mu} (z) dz}}= a\sqrt{q}. 
\end{equation}
\end{lemma}
\noi
\textbf{Proof.} See Theorem 1 in \cite{KajSal93}. \cqfd

\smallskip

\noi
\textbf{Renewal property and oscillations of GW-indexed Brownian snakes.}
We next prove a renewal property of height-Brownian snakes indexed by 
GW forests. We apply this renewal property to control oscillations of such snakes (see (\ref{osci1def})) in the proof of Theorem \ref{Sheuexplain}. 

Let $\smash{\mathcal T\eqo \big(\mathcal T (p)\eqo (\tau(p) , (\ell_u(p))_{u\in \tau(p)})\big)_{p\in \bbN^*}}$ be an infinite 
GW($\mu,q$)-forest as in Definition \ref{GWfordef}.
To simplify notation for all $r,s,s'\ino \bbR_+$, we set $\smash{\mathscr C_s\eqo \mathscr C_s(\mathcal T)}$, 
$\smash{\mathscr H_s\eqo \mathscr H_s(\mathcal T)}$, $\smash{\mathscr W_s(r)\eqo \mathscr W_s(\mathcal T, r)}$, 
\begin{equation}
\label{Wshiftdef}
   \mathtt W_s(r)\eqo \mathtt W_s(\mathcal T, r) \quad \textrm{and} \quad  \mathtt W^{[s]}_{s'} (r) \eqo \mathtt W_{s'+s} \big( r+ m_{\mathscr H  } (s,s'+s) \big) -\mathtt W_{s} \big( m_{\mathscr H  } (s,s'+s) \big).
\end{equation}   
   and we also denote by $\mathscr F_{\! s}$ the sigma-field generated by $\mathscr H_{s''}$ and $\mathtt W_{s''} (\cdot) $, $s''\ino [0, s]$. 
\begin{lemma}
\label{renewsnake} We keep the above notations. Let $\mathbf s$ be a $(\mathscr F_{\! s})_{s\in \bbR_+}$-stopping time. We assume that $\mathbf s$ is a.s.~finite. Then $\mathtt W^{_{[\mathbf s]}}_{\cdot} (\cdot)$ is independent from $\mathscr F_{\! \mathbf s}$ and it has the same distribution as $\mathtt W_\cdot (\cdot)$. 
\end{lemma}
\noi
\textbf{Proof.} We denote by $\mathcal T'\! \eqo (\tau', (\ell'_u)_{u\in \tau'})$ the $\bbR_+$-marked tree associated with the forest $\mathcal T$ as explained in 
Remark \ref{contlifespanrem}. We denote by $(u_l)_{l\in \bbN}$ the vertices of $\tau'$ listed in increasing depth-first order. 
To simplify we set $L _l\! :=\! L_{\mathcal T} (l)  \eqo \sum_{0\leq j\leq l} \ell'_{u_j}$. We recall 
$\zeta'_u\eqo\zeta_u'^*+ \ell'_u\eqo  \sum_{v\in \lgeo \varnothing , u\rgeo} \ell'_v$ and $\mathscr H_{L_l +s} \eqo \zeta'^*_{u_{l+1}}\!\! + s$, $s\ino [0, \ell'_{u_{l+1}})$.   
From (\ref{pathudef}) in Remark \ref{constrheightBrosna} $(b)$, we also recall that $B_\cdot (u)$ is the path of $u$ that is shifted from the death position of $\overleftarrow{u}$ and independently extended after the death of $u$. 
From (\ref{iidpathssnaconstr}) we also recall that there is a measurable function $\smash{F\! : \! \bbT_{\bbR_+\! \times \bC^{_0}_{^1}} \! \to \! \bD (\bbR_+,  \bC^{_0}_{^1})}$ such that $\smash{F \big( \Theta \big)  \eqo  \mathtt W_\cdot (\cdot) }$, where we have set $\smash{\Theta\eqo  \big( t', (\ell'_u, B_\cdot (u) )_{u\in t'} \big) }$. 

  We next observe that the r.v.s $k_{u_l} (\tau')$, $\ell'_{u_l}$, $B_\cdot (u_l)$, $l\ino \bbN^*\! $, are independent: 
 $k_{u_l} (\tau')$ has law $\mu$, $\ell'_{u_l}$ is exponentially distributed with parameter $q$ and $B_\cdot (u_l)$ is a 
 $1$-dimensional Brownian motion with initial value $0$. So we easily see that there is a measurable function 
 $G$ from $\smash{(\bbN\times \bbR_+ \times \bC^{_0}_{^1})^{\bbN^*} } $ to  $\smash{ \mathbf D (\bbR_+,  \bC^{_0}_{^1})}$ such that 
\begin{equation}
\label{Markovsnake}
G \Big( \big(k_{u_l} (\tau'), \ell'_{u_l}, B_\cdot (u_l) \big)_{l\in \bbN^*}\Big)= \big(\mathtt W_s (\cdot))_{s\in \bbR_+}\; .
\end{equation}  

We fix $s\ino \bbR_+$ and we now define a forest 
$\smash{\Theta_s\eqo (\tau'_s , (\ell^{_{[s]}}_u, B^{_{[s]}}_{\cdot} (u))_{u\in \tau'_s} )}$ 
such that $F (\Theta_s)\eqo \mathtt W^{_{[s]}}_{^\cdot} (\cdot)$. To that end, we 
note that there is a unique random $l(s) \ino \bbN^*$ such that $L_{l(s)-1} \leqo s \leko L_{l(s)} $. 
We denote by $(w_p)_{p\in \bbN^*}$ the $<_{ \mathtt{lex}}$-increasing sequence of vertices of $\tau'$ 
which are $\leq_{\mathtt{lex}}$-larger or equal to $u_{l(s)}$ and that are grafted on the ancestral line of 
$u_{l(s)}$. Namely, $w_p\leko_{\mathtt{lex}} w_{p+1}$ for all $p\ino \bbN^*$ 
and $\{ w_p\, ;\,  p\ino \bbN^*\}\eqo \big\{ v\ino \tau'\! :\!  u_{l(s)} \leq_{\mathtt{lex}} v \; \textrm{and} \; \overleftarrow{v} 
\in \lgeo  \varnothing , u_{l(s)}\lgeo \big\} $. Note that $w_1\eqo u_{l(s)}$ and for all $p\ino \bbN^*$, 
we set $\tau_s(p) \eqo \theta_{w_p} \tau'$, which is the subtree of $\tau'$ stemming from $w_p$. 
We denote by $\tau'_s $ the tree associated with the forest $\tau_s \eqo (\tau_s(p))_{p\in \bbN^*}$ 
as explained in Remark \ref{contlifespanrem}. For all $u\ino \tau'_s$, we define the marks 
$\ell^{_{[s]}}_u$ and $B^{_{[s]}}_{^\cdot} (u)$ as follows. 

\begin{compactenum}

\smallskip

\item[$(i)$] We set $\ell^{_{[s]}}_\varnothing\eqo 0$ and we take $B^{_{[s]}}_{^\cdot} (\varnothing)$ as the null function. 

\smallskip

\item[$(ii)$] We set $\ell^{_{[s]}}_{^{[1]}}\eqo \zeta_{u_{l(s)}}\! \! -\! \mathscr H_{s}  \eqo L_{l(s)} \! -\! s$ and for all $r\ino \bbR_+$, 
\begin{equation}
\label{pathcutats} B^{_{[s]}}_{r} ([1]) \eqo B_{r+ s-L_{l(s)-1}} (u_{l(s)}) -B_{s-L_{l(s)-1}} (u_{l(s)}) \; . 
\end{equation}
\item[$(iii)$] Let $u\ino \tau'_s\backslash \{\varnothing, [1] \}$. Then there is $p\ino \bbN^*\! $, $u'\ino \tau_s (p) $ and $v\ino \tau'$ such that  $u\eqo [p] \ast u'$ and $v\eqo w_p \ast u'$. Then we set $\ell^{_{[s]}}_{u} \eqo \ell'_v$ and $B^{_{[s]}}_{\cdot} (u)\eqo B_\cdot (v)$. 
\end{compactenum}
We set $\Theta_s\eqo \big( \tau'_s ,(\ell^{_{[s]}}_{u}, B^{_{[s]}}_{\cdot})_{u\in \tau_s'}\big)$. Then, it is easy to check (deterministically) that $F( \Theta_s) \eqo \mathtt W^{_{[s]}}_{^\cdot} $. 
 
Let us denote by $(u_j(s))_{j\in \bbN}$ the vertices of $\tau'_s$ listed in increasing depth-first order. By definition, 
$$ \forall j\geqo 2, \quad  k_{u_j(s)} (\tau'_s)\eqo k_{u_{l(s)+j -1}} (\tau') , \quad \ell_{u_j(s)}^{[s]}\eqo \ell'_{u_{l(s)+j-1}} \quad \textrm{and} \quad B^{[s]}_\cdot (u_j(s))\eqo B_\cdot (u_{l(s)+j -1}). $$   
We next introduce the sigma field $\mathscr F_{\! s}'$ generated by $l(s)$ and for all $l\leqo l(s)$, by $u_l$, $(s\wedge L_{l}) \! -\! L_{l-1}$, $B_{\cdot \wedge ((s\wedge L_{l})  -\! L_{l-1})} (u_l)$. We see that $(\mathscr F_{\! s}')_{s\in \bbR_+}$ is a filtration such that $\mathscr F_{\! s}  \! \subset \! \mathscr F_{\! s}'$. Let us fix $A\ino \mathscr F_{\! s}'$ and $l \ino \bbN^*$. There is, by definition, a measurable function $\phi_l$ such that 
$$ \un_{A \cap \{ l(s)=l\} }= \phi_l \Big( (u_k, \ell'_{u_k}, B_\cdot (u_k))_{0\leq k\leq l-1} ; u_l; s\! -\! L_{l-1} ; B_{\cdot \wedge (s-L_{l-1} ) }(u_l) \Big)\un_{\{ \ell'_{u_l} > s-L_{l-1} \geq 0 \}} \; .$$
This implies that the event $A \cap \{ l(s)\eqo l\}$ is independent from 
$(k_{u_j(s)} (\tau'_s),\ell_{u_j(s)}^{_{[s]}},B^{_{[s]}}_\cdot (u_j(s)))_{j\geq 2}$. Since exponential r.v.s are memoryless 
and 
since Brownian motion has independent and homogeneous increments, $A \cap \{ l(s)\eqo l\}$ is also 
independent from $k_{u_l } (\tau')$, $\ell^{_{[s]}}_{^{[1]}} \eqo \ell'_{u_l} \! -\! (s\! -\! L_{l-1})$ and $B^{_{[s]}}_{\cdot} ([1])$ 
as in (\ref{pathcutats}) which are distributed resp.~as $\mu$, as an $\mathtt{expo} (q)$ r.v.~and as a 
$1$-dimensional Brownian motion starting at $0$. Namely, $A \cap \{ l(s)\eqo l\}$ is independent 
from $(k_{u_j(s)} (\tau'_s),\ell_{u_j(s)}^{_{[s]}},B^{_{[s]}}_\cdot (u_j(s)))_{j\geq 1}$, which is distributed as 
$(k_{u_{l'}} (\tau'),\ell'_{u_{l'}} 
,B_\cdot (u_{l'}))_{l'\geq 1}$. Thus for all measurable $\Phi\! :\! \bD (\bbR_+, \bC^{_0}_{^1}) \! \rightarrow \! \bbR_+$, we get 
\begin{eqnarray*}
\bE \big[  \un_{A}\Phi \big( \mathtt W^{_{[s]}}_{^\cdot} \big)\big] \!\!\!\!  & =&\!\!\!\!\!  \sum_{l\in \bbN^*} \bE \big[  \un_{A \cap \{ l(s)=l\} }\Phi \big( G \big( (k_{u_j(s)} (\tau'_s),\ell_{u_j(s)}^{_{[s]}},B^{_{[s]}}_\cdot (u_j(s)))_{j\geq 1}\big) \big)  \big]  \\
\!\!\!\!  & =& \!\!\!\!\!  \sum_{l\in \bbN^*} \bP (A \cap \{ l(s)\eqo l \}) \bE \big[ \Phi \big( G \big( (k_{u_{l'}} (\tau'),\ell'_{u_{l'}},B_\cdot (u_{l'}))_{l'\geq 1}\big) \big)  \big] = \bP (A) \bE \big[\Phi (\mathtt W_\cdot) \big] , 
\end{eqnarray*}  
This proves that $\mathtt W^{_{[s]}}_\cdot$ has the same law as $\mathtt W_\cdot$ and that it is independent from $\mathscr F_{\! s}'$ and thus from $\mathscr F_{\! s}$.  

Let us fix $s_1, \ldots, s_p \ino \bbR_+$. By (\ref{Wshiftdef}), observe that 
$s\ino \bbR_+ \! \mapsto \!  \mathtt W^{_{[s]}}_{s_j} (\cdot)$ is right-continuous in 
$\bC^{_0}_{^1}$. Let $f\! :\! (\bC^{_0}_{^1})^p \! \to \! \bbR$ be bounded and continuous and set 
$\Phi (\mathtt W^{_{[s]}}_\cdot)\eqo f \big(  \mathtt W^{_{[s]}}_{s_1} (\cdot), \dots, \mathtt W^{_{[s]}}_{s_p} (\cdot) \big)$. 
Then 
$s \! \mapsto\!  \Phi (\mathtt W^{_{[s]}}_\cdot)$ is right-continuous. We now consider $\mathbf s$, an a.s.~finite 
$(\mathscr F_{\! s})_{s\in \bbR_+}$-stopping time. For all $n\ino \bbN$, we set 
$\mathbf s_n\eqo 2^{-n} \lceil 2^n\mathbf s\rceil$. Let $A \ino \mathscr F_{\! \mathbf s}$. Then 
$\bE \big[\un_A \Phi (\mathtt W^{_{[\mathbf s]}}_\cdot )   \big] \eqo \lim_{n\to \infty} 
\bE \big[\un_A \Phi (\mathtt W^{_{[\mathbf s_n]}}_\cdot )   \big]$. Now observe that 
$\bE \big[\un_A \Phi (\mathtt W^{_{[\mathbf s_n]}}_\cdot )   \big]\eqo \sum_{k\in \bbN} 
\bE \big[ \un_{A\cap \{ k2^{-n} \leq \mathbf s < (k+1) 2^{-n}\} } \Phi (\mathtt W^{_{[(k+1) 2^{-n}]}}_\cdot ) \big]$. 
By the previous arguments, 
$ A\cap \{ k2^{-n} \leqo \mathbf s \leko (k+1) 2^{-n}\} \ino \mathscr F_{\! (k+1)2^{-n}}$, is independent from 
$\Phi (\mathtt W^{_{[(k+1) 2^{-n}]}}_\cdot )$, which has the same law as $\mathtt W$. 
Therefore $\bE \big[\un_A \Phi (\mathtt W^{_{[\mathbf s_n]}}_\cdot )   \big]\eqo \bP(A) \bE \big[\Phi (\mathtt W_\cdot )   \big]$. 
We have proved for all 
$A\ino \mathscr F_{\! \mathbf s}$, for all $s_1, \ldots, s_p\ino \bbR_+$ and for all bounded and continuous 
$f\! :\! (\bC^{_0}_{^1})^p \! \to \! \bbR$, that 
$$ \bE \big[ \un_A f \big(  \mathtt W^{_{[\mathbf s]}}_{s_1} (\cdot), \dots, \mathtt W^{_{[\mathbf s]}}_{s_p} (\cdot) \big)\big] 
= \bP (A)  \bE \big[  f \big(  \mathtt W_{s_1} (\cdot), \dots, \mathtt W_{s_p} (\cdot) \big)\big] $$
which easily entails the desired result.   \cqfd

\medskip

Let $z\ino \bbR_+^*$. We recall from (\ref{osci1def}) the definition of the successive 
$z$-oscillation times of $\mathscr W$ and $\mathtt W$. To simplify notation we set here for all $q\ino \bbN$, 
$$ \bsigma_{\! q} (z) = \bsigma_{\! q} (\mathscr W, z) \quad \mathbf s_q (z)= \bsigma_{\! q} (\mathtt W, z) \; .$$ 
Since the forest $\mathcal T$ is infinite, we see that a.s.~for all $q\ino \bbN$, $\bsigma_{\! q} (z)$ and $\mathbf s_q(z)$ are finite and by (\ref{defviatimechange}) we check that a.s.~for all $q\ino\bbN$,  
\begin{equation}
\label{changetimeosci}
\mathcal K_{\mathcal T} ( \mathbf s_q (z))\eqo \bsigma_{\! q} (z) \; .
\end{equation} 
Since $\mathscr W$ is continuous, Lemma \ref{osc1lemma} applies to the times $\bsigma_{\! q} (z)$, which allows 
to control of the oscillations of $\mathscr W$. However the random times $\bsigma_{\! q} (z)$ have no nice probabilistic properties, to the contrary of the times $\mathbf s_q (z)$, as asserted by the following lemma which is an easy consequence of the previous one.
\begin{lemma}
\label{renewosci} We keep the previous notations and assumptions. Then a.s.~for all $q\ino \bbN$, $\mathbf s_q (z) \leko \infty$ and the r.v.s $(\mathbf s_{q+1} (z) \! -\! \mathbf s_q (z))_{q\in \bbN}$ are i.i.d.
\end{lemma}
\noi
\textbf{Proof.} Observe that $\mathbf s_{q+1} (z) \! -\! \mathbf s_q (z) \eqo \mathbf s_1 (\mathtt W^{_{[\mathbf s_q (z)]}}_{\, \cdot} \! , z)$ and the desired result is an immediate consequence of Lemma \ref{renewsnake} and of a recursive argument. \cqfd

\subsection{Coupling for BRWs, Gaussian estimates}
\label{couplingsec}
In this section we first state a coupling result for BRWs which is used to derive Theorems \ref{extenscv} and \ref{maincvsnake} from Theorem \ref{Sheuexplain}. 
More precisely, our result relies on the following result which is due to Sakhanenko \cite{Sak91} and which extends 
Theorem 4 in Koml{\'o}s, Major \& Tusn{\'a}dy  \cite{KomMajTus76} to RWs with independent jumps but with possibly distinct laws. 
\begin{theorem}
\label{recallSakhanenko}
Let $G$ be a moment gauge function as in (\ref{momgaudef}). Let $ K_0, y \ino \bbR_+^*$. 
Let $n\ino \bbN\backslash\{ 0, 1\}$
and let $(\xi_j)_{1\leq j\leq n}$ be $\bbR$-valued r.v.s. We assume that the probability space $(\Omega, \ccF\! , \bP)$ is sufficiently rich to carry an additional r.v.~$U$ that is uniformly distributed on $[0, 1]$ and that is independent of $(\xi_j)_{1\leq j\leq n}$. We assume the following. 

\smallskip

\begin{compactenum}
\item[$(a)$] The r.v.s $(\xi_j)_{1\leq j\leq n}$ are independent.

\smallskip

\item[$(b)$] For all $j\ino \{ 1, \ldots , n\}$,  $\bE[ G(\xi_j)] \leqo K_0 \sigma_j^2\leko \infty$, where $\sigma^2_j\eqo \bE[ \xi_j^2]$, and $\bE[ \xi_j]\eqo 0$. 

\smallskip

\item[$(c)$] $y^{-2} G(y) \geqo K_0$. 

\smallskip

\end{compactenum}

\noi
Then, there is $K_1\ino \bbR_+^*$ which only depends on $K_0$ and $G$, and there are 
independent standard Gaussian r.v.s $(Y_j)_{1\leq j\leq n}$, which 
depend measurably on $K_0, G, (\xi_j)_{1\leq j\leq n}, U$ and $y$, such that 
$$ \bP \Big(\max_{1\leq k\leq n} 
\Big|\! \sum_{1\leq j\leq k} \!\!  \xi_j  \! -\! \sigma_jY_j \,  \Big| \! >\!  K_1y \Big) \leq \frac{4K_0 (\sigma_1^2+ \ldots + \sigma_n^2)}{G(y)} .$$
\end{theorem}
\noi
\textbf{Proof.} See \cite{Sak91}, Corollary 12, p.~81.  \cqfd 

\medskip

\begin{remark}
\label{kappax0} In the definition (\ref{momgaudef}) of the moment gauge function $G$, two positive parameters $\kappa$ and $x_0$ are involved and are part of the very definition of $G$. Therefore in the previous statement and in the next statements, when we mention that a quantity depends on $G$ we actually mean that it depends on $G$, $\kappa$ and $x_0$. \cq 
\end{remark}
The coupling for BRWs is stated in Theorem \ref{brwcouplingth} below. It involves the \emph{Horton-Strahler number} of the indexing tree that is defined as follows. 
\begin{definition}
\label{HSnumber} ({\small \textbf{Horton-Strahler number}}) 
Let $\fftree\ino \bbT$ be a finite tree. For any function $\varphi \! : \! \fftree \! \rightarrow \! \bbN$ and for all
 $u\ino \fftree \backslash \mathtt{Lf} (\fftree)$, we set $\overline{\varphi} (u)\eqo \max\{ \varphi (v); v\ino \fftree: \overleftarrow{v}\eqo u \}$. 
Then the \emph{Horton-Strahler flow} is the $\bbN$-valued function $\phi$ on $\fftree$ that is characterised by the following. 

\smallskip

\begin{compactenum}
\item[$(a)$] For all $v\ino \mathtt{Lf} (\fftree)$, $\phi (v)\eqo 1$. 

\smallskip

\item[$(b)$] For all $u\ino \fftree \backslash \mathtt{Lf} (\fftree)$, $ \phi (u)= \overline{\phi} (u) + \un_{\{ \exists v, v^\prime \in \fftree\; : \;  v\neq v^\prime, \;  \overleftarrow{v}= \overleftarrow{v}^\prime= u  \; \textrm{and} \;  \phi (v)= \phi (v^\prime) =  \overline{\phi} (u)\}}. $
\end{compactenum} 

\smallskip

\noi
The \emph{Horton-Strahler number} of $\fftree$ is then defined as $\mathtt{HS} (\fftree)\! : = \! \phi (\varnothing)$. \cq 
\end{definition}
Note that $\mathtt{HS} (\theta_u \fftree)\eqo \phi (u)$ for all $u\ino \fftree$ and that $\phi$ is nonincreasing with respect to the genealogical order $\preceq$. 
The Horton-Strahler number measures tree-complexity. It originates in works on hydrology by Horton \cite{Hor45} and Strahler \cite{Str52}. It also appeared independently in computer science under the name of \emph{register numbers} (see Flajolet, Raoult \& Vuillemin \cite{FlaRaoVui79} and independently Kemp \cite{Kem79}), and also in other scientific fields: we refer to the mathematical survey of Kovchegov \& Zaliapin \cite{KovZal20} for a more comprehensive account on that topics. 
The role of the Horton-Strahler number in our coupling is explained at the beginning of the proof of Theorem \ref{brwcouplingth} below. 

The Horton-Strahler number is in 
general not easy to analyse directly. 
Let us mention that in practice we shall only need the upper bound $\mathtt{HS} (\fftree) \leqo 1+ 
\log_2 (\# \mathtt{Lf} (\fftree))$ that is easy to prove
(see Lemma \ref{HSloppconn}). 
For a large class of GW-trees $\log (\# \fftree)$ turns out to be the right order of magnitude (up to a multiplicative factor that is in general difficult to compute or to characterise): 
see for instance Devroye \& Kruszewski \cite{DevKru95} and also the recent works by 
Brandenberger, Devroye \& Reddad \cite{BraDevRed21} and Khanfir \cite{khanfir_2023_horton_strahler, khanfir2024fluctuations} for quite precise results on Horton-Strahler of conditioned GW-trees. 

Our main result on coupling of BRWs can be stated as follows. 
\begin{theorem}
\label{brwcouplingth} Let $G$ be a moment gauge function as in (\ref{momgaudef}). Let $ C_1, C_2, C_3, y \ino \bbR_+^*$. 
Let $\fftree\ino \bbT$ be finite. Its Horton-Strahler number is denoted by $\mathtt{HS} (\fftree)$ and its set of leaves by $\mathtt{Lf} (t)$. Let $(\xi_u)_{u\in \fftree \backslash \{ \varnothing\}}$ be $\bbR$-valued r.v.s. 
We assume that the probability space $(\Omega, \ccF\! , \bP)$ is sufficiently rich to carry an additional r.v.~$U$ that is uniformly distributed on $[0, 1]$ and that is independent of $(\xi_u)_{u\in \fftree \backslash \{ \varnothing\}}$. We assume the following. 

\smallskip

\begin{compactenum}
\item[$(a)$] The r.v.s $\big( (\xi_{[u]\ast i} )_{1\leq i\leq k_u (\fftree)}; u \! \in \fftree \backslash \mathtt{Lf} (\fftree)\big)$ are independent.

\smallskip

\item[$(b)$] For all $u\ino \fftree\backslash \{ \varnothing \} $, $\;  \bE \big[ G(\xi_u) \big] \leqo C_1$, $\; C_2 \leqo\sigma_u^2\! := \!  \bE \big[ \xi_u^2\big] \leqo C_3$ and $\; \bE [ \xi_u] \eqo 0$.

\smallskip

\item[$(c)$]  $ \big( \frac{y}{\mathtt{HS} (\fftree)} \big)^{\! -2} G\big( \frac{y}{\mathtt{HS} (\fftree)}\big)\geqo C_1/C_2$.

\smallskip

\end{compactenum}

\noi
Then, there exists $K\ino \bbR_+^*$ and independent standard Gaussian r.v.s $(Y_u)_{u\in \fftree \backslash \{ \varnothing\}}$ such that 
\begin{equation}
\label{coupling1}
 \bP \Big(\max_{u\in \fftree} \big|S_u \! -\! R_u \big| \! >\!  Ky \Big) \leq 2C_1\!  \Big(1+ \frac{2C_3}{C_2} \Big) \, \frac{\# \fftree}{G  \big(\frac{y }{\mathtt{HS}(\fftree)}  \big)} 
 \end{equation}
 where for all $u\ino \fftree\backslash \{ \varnothing \}$, we have set 
$S_u \eqo \sum_{v\in \rgeo \varnothing , u \rgeo} \xi_v$, $\, R_u \eqo \sum_{v\in \rgeo \varnothing , u \rgeo} \sigma_vY_v$ and $S_\varnothing\eqo  R_\varnothing\eqo  0$. 
Here, $K$ depends on $C_1$, $C_2$, $C_3$ and $G$ only; the r.v.s
$(Y_u)_{u\in \fftree\backslash \{ \varnothing\} }$ depend measurably on $C_1$, $C_2$, $C_3$, $G$, $(\xi_u)_{u\in \fftree \backslash \{ \varnothing\}}$, $U$ and $y$. 
\end{theorem}
\noi
\textbf{Proof.}
The main idea is to apply Theorem \ref{recallSakhanenko} to a decomposition of the tree $\fftree$ into 
pairwise disjoint lines that 
we call 
\emph{lopping}. 
It is defined as follows. 
We first give each leaf $v\ino \mathtt{Lf} (\fftree)$ a \emph{rank} $r (v)\ino \{ 1, \ldots, \# \mathtt{Lf} (\fftree)\}$ 
and we assume that the ranks $\mathbf{r} \eqo (r(v))_{v\in \mathtt{Lf} (\fftree)}$ are distinct. Then, the $\mathbf r$-lopping of $\fftree$ is a finite increasing sequence of distinct graph-subtrees (i.e., subsets $\fftree^\prime \! \subset \! \fftree$ such that $\fftree^\prime \eqo \bigcup_{u\in \fftree^\prime} \lgeo \varnothing , u \rgeo$): 
$$\fftree_0\eqo \{ \varnothing\}  \subsetneq \fftree_1 \subsetneq \ldots \subsetneq \fftree_m \subsetneq \fftree_{m+1} \subsetneq \ldots \subsetneq \fftree_{M(\mathbf{r})}\eqo \fftree $$
that are defined recursively as follows. Suppose that $\fftree_m$ is already defined. 
If $\fftree_m\eqo \fftree$, then $m\eqo M(\mathbf r)$ and the construction is completed. 
Suppose that $\fftree_m \! \subsetneq \!  \fftree$. We denote by $C_{m+1} (j)$, $1\leqo j\leqo p_{m+1}$, the connected components of $\fftree\backslash \fftree_m$ and for all 
$j\ino \{ 1, \ldots, p_{m+1}\}$, let us denote by $v_{m+1}(j)$ the leaf of  $C_{m+1} (j)$ whose rank is minimal: namely $r (v_{m+1} (j))\eqo \min \{ r(v) ; v \ino C_{m+1} (j)\cap \mathtt{Lf} (\fftree) \}$. We also denote by $b_{m+1} (j)$ the unique vertex of $\fftree_m$ such that $\fftree_m\cap \lgeo \varnothing , v_{m+1} (j) \rgeo \eqo  \lgeo \varnothing , b_{m+1} (j) \rgeo$. We then set 
$$ \fftree_{m+1} =  \fftree_{m}\,  \cup \!\!\!\! \!\! \!\! \!\!  \bigcup_{\quad 1\leq j\leq p_{m+1}} \!\! \!\! \!\! \!\!  \lgeo \varnothing , v_{m+1} (j) \rgeo \; .$$
We observe that the lines $\, \rgeo\,  b_{m+1} (j)  , v_{m+1} (j) \rgeo$, $1\leqo j\leqo p_{m+1}$, form a partition of $\fftree_{m+1} \backslash \fftree_m$: they are 
the \emph{Step ($m+1$)-branches of the lopping process}. Thus,   
\begin{equation}
\label{branches}
 \rgeo \, b_{m} (j)  , v_{m} (j)\rgeo  \, , \; 1\leqo j\leqo p_{m}, \; 1\leqo m\leqo M(\mathbf r) ,\; \textrm{form a partition of $\fftree\backslash\{ \varnothing \}$.}
\end{equation}
Therefore,  
\begin{equation}
\label{sumnmj}
 \# \fftree \eqo 1+\!\! \!\! \!\!\!\!\!\!  \sum_{\quad \substack{1\leq m\leq M(\mathbf r) \\ 1\leq j\leq p_{m}}} \!\! \!\! \!\!\!\!
n_{m} (j)
 \quad \textrm{where} \quad n_{m} (j)\eqo |v_{m} (j)| \! -\!  |b_{m} (j)| \; .
 \end{equation}
Then, for all $u\ino \fftree\backslash \{ \varnothing \}$, there exists a unique $(m, j) \ino \{ 1, \ldots, M (\mathbf r) \}\! \times \!  \{ 1, \ldots, p_{m} \}$ such that $u\! \in \,  \rgeo \, b_{m} (j)  , v_{m} (j)\rgeo$ and we set 
\begin{equation}
\label{piudef}\pi (u)\eqo m, \quad \mathtt{lf} (u)\eqo v_{m} (j) \quad \textrm{and} \quad n(u)\eqo n_{m} (j) . 
\end{equation}
Note that $b_m(j) \ino t_{m-1}$ for all $1\leqo j\leqo p_m$ and all $1\leqo m\leqo M(\mathbf r)$. In particular $b_1(j)\eqo \varnothing$, for all $1\leqo j\leqo p_1$.

As measurable deterministic functions of 
$U$, we produce independent r.v.s $U_0$ and $U_{j,m}$, $1\leqo j\leqo p_m, 1\leqo m\leqo M(\mathbf r)$, which are uniformly distributed on $[0, 1]$. They are necessarily independent from $(\xi_u)_{u\in t\backslash \{ \varnothing\}}$.

To apply Theorem \ref{recallSakhanenko} along branches we need independent jumps. However, 
jumps of siblings may be dependent r.v.s. To overcome this problem we 
introduce an auxiliary BRW with independent jumps as follows. 
We first use the additional uniform r.v.~$U_0$ to generate r.v.s $(\widetilde{\xi}_u)_{u\in \fftree \backslash \{ \varnothing\}}$ such that: 

\noi
$-$  $(\widetilde{\xi}_u)_{u\in \fftree \backslash \{ \varnothing\}}$ are mutually independent, 

\noi
$-$ $(\widetilde{\xi}_u)_{u\in \fftree \backslash \{ \varnothing\}}$ and $(\xi_u)_{u\in \fftree\backslash \{ \varnothing\}}$ are independent, 

\noi
$-$ for all $u\ino \fftree \backslash \{ \varnothing\}$, $\widetilde{\xi}_u$ has the same law as $\xi_u$.  

\smallskip

\noi
Thanks to the $\widetilde{\xi}_u$ we define auxilliary 
jumps $\xi^\prime_u$ as follows. Let $m\ino \{ 1, \ldots, M(\mathbf r) \}$ and let 
$j\ino \{ 1, \ldots, p_m\}$. Then, for all $u\! \in \, \rgeo b_m(j), v_m(j)\rgeo$ we set 
$$ \xi^\prime_u= \xi_u \quad \textrm{if $|u| \! >\!  |b_m(j)|\! +\! 1\; $ and} \quad \xi^\prime_u \eqo \widetilde{\xi}_u \quad \textrm{if $|u| \eqo |b_m (j)| \!+ \! 1$.} $$
By Remarks \ref{siblindep} $(b)$ and (\ref{branches}) we check recursively in $m$ for $ t_1, \ldots, t_m, \ldots, t_{M(\mathbf r)}$ that  

\smallskip

\noi
$-$  $(\xi^\prime_u)_{u\in \fftree \backslash \{ \varnothing\}}$ are mutually independent, 

\smallskip

\noi
$-$ $(\widetilde{\xi}_u)_{u\in \fftree \backslash \{ \varnothing\}}$ and $(\xi^\prime_u)_{u\in \fftree\backslash \{ \varnothing\}}$ have the same law. 

\smallskip 

We denote by $B\eqo \{ b_m(j); 1\leqo m\leqo M(\mathbf r) , \, 1\leqo j\leqo p_m \} $, which is the set of branching points where branches are cut. Let $u\ino \fftree \backslash \{ \varnothing\}$. Let 
\begin{equation}
\label{bjwju}
b_1\eqo \varnothing \preceq w_1 \preceq b_2  \preceq w_2  \preceq \ldots    \preceq b_{\pi (u)}  \preceq w_{\pi (u)} \preceq u
\end{equation}
be such that $\{ b_j\, ;\,  1\leqo j\leqo \pi (u) \}\eqo   B \cap \, \lgeo \varnothing, u \lgeo \, $ 
and for all $j\ino \{ 1 , \ldots , \pi (u) \} $, $|w_j|\eqo |b_j| +1$.
 If we denote by $(S^\prime_u)_{u\in \fftree}$ the BRW associated with the jumps $(\xi^\prime_u)_{u\in \fftree \backslash \{ \varnothing\}}$, then $S^\prime_u\! -\! S_u\eqo \sum_{1\leq j\leq \pi (u)} \widetilde{\xi}_{w_j} \! -\! \xi_{w_j} $. Thus 
$$ \max_{u\in \fftree \backslash \{ \varnothing\}} \frac{_1}{^{\pi (u)}} \big| S^\prime_u\! -\! S_u\big| \, \leq \Delta + \widetilde{\Delta} \quad \textrm{where} \quad \Delta = \!\!\!\!\! \!   \max_{\quad u\in \fftree \backslash \{ \varnothing\}}\!\! \!\!     |\xi_u| \quad \textrm{and} \quad \widetilde{\Delta}= 
\!\!\!  \!\! \!   \max_{\quad u\in \fftree \backslash \{ \varnothing\}} \!\! \!\! | \widetilde{\xi}_u|.  $$
Observe that $\bP (\Delta \geko z)$ and $\bP (\widetilde{\Delta} \geko z)$ are both bounded by 
$\sum_{u\in \fftree \backslash \{ \varnothing\}} \bP(|\xi_u| \! >\! z)$ for all  $z\ino \bbR_+^*$. 
Therefore,
\begin{equation}
\label{controlauxiBRW}
\bP \Big( \max_{u\in \fftree \backslash \{ \varnothing\}} \frac{_1}{^{\pi (u)}} \big| S^\prime_u\! -\! S_u\big|  >2z \Big) \, \leq \!\!\!\!\!\!  \sum_{\quad u\in \fftree \backslash \{ \varnothing\}} \!\!\!\!\!\! \!\! \frac{2\bE [G(\xi_u) ]}{G(z)}\,  \leq \frac{2C_1 \# t}{G(z)}  ,
\end{equation}
by a simple Markov inequality since $G$ is nondecreasing on $\bbR_+$. 

We then fix $m\ino \{ 1, \ldots, M(\mathbf r) \}$ and 
$j\ino \{ 1, \ldots, p_m\}$ and we apply Theorem \ref{recallSakhanenko} with the moment gauge function $G$, with $K_0 \eqo C_1/C_2$ and with $\smash{y_1\! :=\! \frac{y}{\mathtt{HS} (\fftree)}}$, which satisfies $y_1^{-2}G(y_1)\geqo K_0$, to the independent r.v.s $U_{j,m}$ and $(\xi^{j,m}_k)_{1\leq k\leq n_m(j)}$: here $\xi^{j,m}_k\eqo \xi^\prime_w$, where $w\! \in \, 
\rgeo b_m(j), v_m (j) \rgeo $ is such that $|w| \eqo k+  |b_m(j)|$.  
Namely, Theorem \ref{recallSakhanenko} 
asserts that there are:

\smallskip

\noi
$-$ $K_1\ino \bbR_+^*$ that only depends on $G$, $C_1$ and $C_2$ ,

\noi
$-$ independent standard Gaussian r.v.s $(Y^{j,m}_k )_{1\leq k\leq n_m(j)}$ that measurably depend on 
$C_1$, $C_2$, $G$, $(\xi^{j,m}_k)_{1\leq k\leq n_m(j)}$, $U_{j,m}$ and $y_1$, such that:  
\begin{equation}\label{Ajm}
\bP (A_{j,m}) \leq \frac{4C_1C_3\,  n_m(j) }{C_2 \, G(y_1)} \; \textrm{where} \; A_{j,m}\eqo  \Big\{\!  \max_{\; 1\leq k\leq n_m(j)} \Big| \! \sum_{1\leq l\leq k} \! \xi^{j,m}_l  \! \! - \sigma^{j,m}_l  Y^{j,m}_l 
\Big|   > K_1y_1 \Big\} 
\end{equation}
where we have set $ \sigma^{j,m}_l \eqo (\bE [ (\xi^{j,m}_l )^2])^{1/2}$. 

Here, we note that \emph{all} standard Gaussian r.v.s~$Y^{j,m}_k$, $m\ino \{ 1, \ldots, M(\mathbf r) \}$, $j\ino \{ 1, \ldots, p_m\}$, $k\ino \{ 1, \ldots, n_m(j)\}$ are mutually independent (of course there are highly dependent on the $\xi^\prime_u$ and on $U$). 
Then, we recall that for all $u\ino \fftree \backslash \{ \varnothing\}$, there exists a unique  
$(m,j)\ino \{ 1, \ldots, M(\mathbf r) \} \times  \{ 1, \ldots, p_m\}$ such that 
$u\! \in\,  \rgeo b_m(j), v_m(j) \rgeo $. We set $Y_u\eqo Y^{j,m}_{|u| -|b_m(j)|}$. Therefore, the  $(Y_u)_{u\in \fftree \backslash \{ \varnothing\}}$ are independent standard Gaussian r.v.s.

We recall the notation 
$R_u \eqo \sum_{v\in \rgeo \varnothing , u \rgeo} \sigma_vY_v$ and let $b_1, \ldots , b_{\pi (u)}$ be as in (\ref{bjwju}). First observe that 
$$ \big| S^\prime_u \! -\! R_u\big| \leq \big| S^\prime_u \!\!  -\! S^\prime_{b_{\pi (u)}} \!\!  - (R_u\! -\! R_{b_{\pi (u)}} )\big| \; + \!\! \!\! \!\! \!\! \sum_{\quad 1\leq j < \pi (u)} \!\! \!\! \!\!  \big|S^\prime_{b_{j+1}} \!\!  - S^\prime_{b_{j}} \!\!  -(R_{b_{j+1}} \!  \! - R_{b_{j}} ) \big| $$
with the convention that the sum in the right hand side is null if $\pi(u)\eqo 1$. 
Then on $\Omega_0:= \bigcap_{1\leq m\leq M (\mathbf r ), 1\leq j\leq p_m} (\Omega \backslash A_{j,m} ) $, we get 
$\big| S^\prime_u\! -\! R_u\big| \leq \pi (u) K_1 y_1 $ for all $u\ino t\backslash \{ \varnothing\}$. Therefore, by (\ref{Ajm}) and (\ref{sumnmj})
\begin{equation}
\label{couplauxiBRW}
\bP \Big( \max_{u\in \fftree \backslash \{ \varnothing\}} \frac{_1}{^{\pi (u)}} \big| S^\prime_u\! -\! R_u\big|  > K_1 y_1 \Big) \leq  \bP (\Omega \backslash  \Omega_0) \leq \sum_{m,j} \bP (A_{j,m}) \leq 
\frac{4C_1C_3 \# t }{C_2 \, G(y_1)}, 
\end{equation}
where in the sum, $m$ and $j$ range resp.~in $\{ 1, \ldots, M(\mathbf r) \}$ and $\{ 1, \ldots, p_m\}$.

Observe that the previous couplings heavily rely on the ranking $\mathbf r$ that is initally chosen. The following lemma shows that there is a ranking $\mathbf r $ where $M(\mathbf r )\eqo \max_{u\in \fftree} \pi (u)$ is minimal  and equal to the Horton-Strahler number $\mathtt{HS} (\fftree)$ of $\fftree$. 
\begin{lemma}
\label{HSloppconn}  We keep the previous notations. Then 
\begin{equation}
\label{connectHS}
\mathtt{HS} (\fftree)= \min M (\mathbf{r})  \leq 1+ \log_2 \big(\# \mathtt{Lf} (\fftree) \big) \; , 
\end{equation}
where the minimum is taken over all the possible rankings $\mathbf{r}$ of the leaves of $\fftree$. 
\end{lemma}

We postpone the proof of the lemma to complete the proof of Theorem \ref{brwcouplingth} first. 
We suppose that $\mathbf r$ is such that $M(\mathbf r)\eqo \mathtt{HS} (\fftree)$. We now recall that $y_1\eqo y/ \mathtt{HS} (\fftree)$. 
By (\ref{controlauxiBRW}) with $z\eqo y_1$ and (\ref{couplauxiBRW}), we easily get 
\begin{eqnarray*}
\bP \Big(\! \max_{u\in \fftree } \big|S_u \! -\! R_u \big| \!\!  \!\!     & > & \!\!  \!\!    (K_1+2)y \Big) \;  \leq \; 
\bP \Big(\max_{u\in \fftree  \backslash \{ \varnothing\} } \frac{_1}{^{\pi (u)}}\big|S_u \! -\! R_u \big| \! >\!  (K_1+2)y_1 \Big) \\
& \leq &  \bP \Big( \max_{u\in \fftree \backslash \{ \varnothing\}} \frac{_1}{^{\pi (u)}} \big| S^\prime_u\! -\! S_u\big|  >2y_1 \Big) + \bP \Big( \max_{u\in \fftree \backslash \{ \varnothing\}} \frac{_1}{^{\pi (u)}} \big| S^\prime_u\! -\! R_u\big|  > K_1 y_1 \Big) \\
& \leq & \frac{2C_1 \# t}{G(y_1)} + \frac{4C_1C_3 \# t }{C_2 \, G(y_1)}= 2C_1\!  \Big(1+ \frac{2C_3}{C_2} \Big) \, \frac{\# \fftree}{G  \big(\frac{y }{\mathtt{HS}(\fftree)}  \big)}, 
\end{eqnarray*}
and we get (\ref{coupling1}) by taking $K$ equal to $K_1+2$. This completes the proof of Theorem \ref{brwcouplingth}. \cqfd 

\medskip

\noi
\textbf{Proof of Lemma \ref{HSloppconn}.} 
First let us recursively construct a ranking $\mathbf r$ that satisfies 
$M (\mathbf{r})  \leqo 1+ \log_2 \big(\# \mathtt{Lf} (\fftree) \big)$. To that end, we introduce the following notations. Let $S$ be any non empty subset of leaves of $\fftree$; we set $n\eqo \#S$ and we call $v$ the \emph{midpoint} of $S$ if 
$\# \{w\ino S:  w\!  \leq_{\mathtt{lex}} \! v \}\eqo \lceil \frac{_1}{^2}n \rceil$ 
(recall that $\lceil \cdot \rceil $ means $\lfloor \cdot \rfloor +1$ for us, so 
$\lceil \frac{_1}{^2}n \rceil\eqo  \frac{_1}{^2}n+1$ if $n$ is even). 
Note that 
$\max \big( \# \{w\ino S:  w\!  <_{\mathtt{lex}} \! v \} , \# \{w\ino S:  w\!  >_{\mathtt{lex}} \! v \} \big) \eqo \lfloor n/2 \rfloor$. Then in the previous construction, for all $m\ino \{ 1, \ldots, M(\mathbf r)\}$ and for all 
$j\ino \{ 1, \ldots, p_m\}$, we can always choose the ranking such that $v_m(j)$ is the midpoint of $\mathtt{Lf} (\fftree) \cap C_m(j)$ and we easily check that for all $m\ino \{ 1, \ldots, M(\mathbf r)\}$
$$\max \big\{ \#  (\mathtt{Lf} (\fftree) \cap C_{m} (j) ) \, ; \, j\ino \{ 1, \ldots, p_{m}\} \big\} \leq 
2^{-(m-1)}\# \mathtt{Lf} (\fftree). $$
Thus $1\leqo 2^{-M(\mathbf r)+ 1}\# \mathtt{Lf} (\fftree)$ which implies the desired upper bound. 

\smallskip 

  Let $\mathbf{r}$ be any ranking of the leaves of $\fftree$: we next 
prove that $\mathtt{HS} (\fftree) \leqo M(\mathbf{r})$. To that end, we recall 
from Definition \ref{HSnumber} that $\phi$ stands for the Horton-Strahler flow and we 
derive another flow $\varphi\! : \! \fftree \! \rightarrow \! \bbN$ from $\mathbf r$ as follows: for all $u\ino \fftree\backslash \{ \varnothing\}$, we recall from  
(\ref{piudef}) the definition of $\pi (u)$. We agree on the convention $\pi (\varnothing)\eqo 1$ and we set $\varphi (u)\eqo \max_{v\in \, \theta_{ u}  \fftree} \pi (u\ast v) \! -\! \pi (u)+1$. We claim that  $\varphi \geqo \phi $ which implies that 
$M(\mathbf{r}) \eqo \varphi (\varnothing) \geqo \phi (\varnothing) \eqo \mathtt{HS} (\fftree)$. 
Indeed, first observe that $\phi (u) \eqo \varphi (u)\eqo 1$ for all $u\ino \mathtt{Lf} (\fftree)$. Then, it is sufficient to prove the following implication for all 
$u \ino \fftree \backslash \mathtt{Lf} (\fftree)$: 
\begin{equation}
\label{descendingarg}
\textrm{\emph{If $\varphi (u\ast [j]) \geqo \phi (u\ast [j])$ for all $j\ino \{ 1, \ldots, k_u (\fftree) \}$, then $\varphi (u) \geqo \phi (u)$.}}
\end{equation}
\emph{Indeed,} suppose that there is $u_0\ino \fftree  \backslash \mathtt{Lf} (\fftree)$ such that $\varphi(u_0)\leko \phi(u_0)$. W.l.o.g.~we can assume that $|u_0|\eqo \max \{ |u| \ino \fftree \! : \! \varphi(u)\leko \phi(u)\}$. 
Therefore $\varphi (u_0\ast [j]) \geqo \phi (u_0\ast [j])$ for all $j\ino \{ 1, \ldots, k_{u_0} (\fftree) \}$ and (\ref{descendingarg}) entails a contradiction. Thus $\varphi(u)\geqo \phi(u)$, for all $u\ino \fftree$. 

\smallskip

\noi
\emph{Proof of (\ref{descendingarg})}. We fix $u\ino \fftree  \backslash \mathtt{Lf} (\fftree)$ such that 
$\varphi (u\ast [j]) \geqo \phi (u\ast [j])$ for all $j\ino \{ 1, \ldots, k_u (\fftree) \}$. We first observe that 
\begin{equation}
\label{maxequ} \varphi (u) \eqo \max_{1\leq j\leq  k_u (\fftree) } ( \varphi (u\ast [j]) + \pi (u\ast [j]) -  \pi (u) ) \; ,
\end{equation}
By definition, there is a unique $j_0\ino \{ 1, \ldots, k_u (\fftree) \}$ such that 
$\pi (u\ast [j])\eqo \pi (u)+ \un_{\{ j\neq  j_0\}} $ for all $j\ino \{ 1, \ldots, k_u (\fftree) \}$. Therefore $\varphi (u)-\overline{\varphi} (u)\eqo 0$ or $1$, where we use the notation 
$\overline{\varphi} (u)\eqo \max\{ \varphi (v); v\ino \fftree: \overleftarrow{v}\eqo u \}$ 
(and a similar notion concerning $\phi$ as in Definition \ref{HSnumber}). There are several cases to consider. 

\noi
$-$ If $\varphi (u)\eqo \overline{\varphi} (u) +1$, then $\varphi (u)\eqo  \overline{\varphi} (u) +1 \geqo \overline{\phi} (u) +1\geqo \phi (u)$ by Definition \ref{HSnumber} $(b)$.

\noi
$-$ Let us suppose that $\varphi (u)\! =\! \overline{\varphi} (u)$. 
Then (\ref{maxequ}) entails that $\varphi (u)\eqo \varphi (u\ast [j_0]) \geko \varphi (u\ast [j])$ for all $j\ino \{ 1, \ldots, k_u (\fftree) \}\backslash \{ j_0\}$ (here, we use the fact that for all $m, n\ino \bbN$, $(1+m)\vee n \eqo m\vee n $ implies $n\geko m$). Then we have the following cases. 
\begin{compactenum}
\item[$-$] If $\overline{\varphi} (u) \eqo \overline{\phi} (u)$ 
then $ \overline{\phi} (u)\eqo \overline{\varphi} (u)\eqo \varphi (u)
\eqo \varphi (u\ast [j_0]) \geko \varphi (u\ast [j]) 
 \geqo \phi (u\ast [j])$ for all $j\ino \{ 1, \ldots, k_u (\fftree) \}\backslash \{ j_0\}$, which implies that $\phi (u)\eqo  \overline{\phi} (u)$ (by Definition \ref{HSnumber} $(b)$) and thus $\varphi (u)\! =\!\phi (u)$. 
\item[$-$] If $ \overline{\varphi} (u) \geko \overline{\phi} (u)$, then $\varphi (u) =\overline{\varphi} (u)\geqo  \overline{\phi} (u)+ 1 \geqo \phi (u)$ by Definition \ref{HSnumber} $(b)$. 
\end{compactenum}
This completes the proof of the inequality $\mathtt{HS} (\fftree) \leqo M(\mathbf{r})$.

\smallskip

  We next find a ranking $\mathbf{r}$ such that $M (\mathbf r)\eqo \mathtt{HS} (\fftree) $ thanks to the Horton-Strahler flow $\phi$. To that end, for all 
$u\ino \fftree \backslash \mathtt{Lf} (\fftree)$ we define $u^\circ$ (resp.~$u^\bullet$) as 
the $<_{\mathtt{lex}}$-first (resp.~last) vertex $v$ such that $\overleftarrow{v}\eqo u$ and 
$\phi (v)\eqo \overline{\phi} (u)$. By Definition \ref{HSnumber} $(b)$, $\phi (u)\eqo \overline{\phi} (u)+ \un_{\{ u^\circ \neq u^\bullet \}}$. 
 
We next define $\pi \! : \! \fftree \! \rightarrow \! \bbN$ by specifying the following: $\pi (\varnothing) \eqo 1$ and  
for all $u \ino \fftree \backslash\mathtt{Lf} (\fftree)$ and for all $j\ino \{ 1, \ldots, k_u (\fftree)\}$, 
$\pi (u \ast  [j])\eqo \pi (u)+ \un_{\{ u\ast [j] \neq u^\circ \}}$. 
We then define a linear order $<^\circ$ 
on $\mathtt{Lf} (\fftree)$ by specifying that $v\! <^\circ \! v^\prime$ if $\pi (v) \leko \pi (v^\prime)$ or if 
$\pi (v) \eqo \pi (v^\prime)$ and $v <_{\mathtt{lex}} v^\prime$. Then, for all $v\ino \mathtt{Lf} (\fftree)$, we define $r(v)$ as the rank of $v$ with respect to $<^\circ $. It is easy to check that 
$\pi$ is associated with $\mathbf r$ as in (\ref{piudef}). 
Since, by definition, $\max_{v\in \mathtt{Lf} (\fftree)} \pi (v) \eqo M(\mathbf r)$, 
it remains to prove that $\max_{v\in \mathtt{Lf} (\fftree)} \pi (v)\eqo \mathtt{HS} (\fftree)$. 

To that end, we set $h\eqo \pi+ \phi $ and for all $u\ino \fftree \backslash \mathtt{Lf} (\fftree)$ we prove 
\begin{equation}
\label{maxpath} 
  h(u^\bullet)\eqo h(u) \geqo h(v) , \quad v\ino \fftree\backslash \{ u^\bullet\} : \overleftarrow{v} \eqo u . 
\end{equation} 
Before proving (\ref{maxpath}), let us explain why it implies 
$M(\mathbf r) \eqo \max_{v\in \mathtt{Lf} (\fftree)} \pi (v)\eqo \mathtt{HS} (\fftree)$. \emph{Indeed}, first 
note that (\ref{maxpath}) implies that $h(\varnothing)\eqo \max_{u\in t} h (u) $ 
((\ref{maxpath}) entails that $h(u) \eqo \overline{h} (u)$ for all $u\ino t\backslash \mathtt{Lf} (t)$, thus $h(\overleftarrow{u}) \geqo h(u)$ for all $u\ino t\backslash \{ \varnothing\}$: namely, $h$ is nonincreasing with respect to $\preceq$). Thus for all $v\ino \mathtt{Lf} (\fftree)$, 
$\pi (v)+1\eqo \pi (v)+ \phi (v)\eqo  h(v) \leqo h(\varnothing)\eqo \pi (\varnothing )+ \phi (\varnothing) \eqo 1 + \mathtt{HS} (\fftree)$, which implies that $\max_{v\in \mathtt{Lf} (\fftree)} \pi (v)\leqo \mathtt{HS} (\fftree)$ and thus $M(\mathbf r) \eqo \mathtt{HS} (\fftree)$ by the previous result.

\smallskip

It remain to prove  (\ref{maxpath}). We consider two cases. 

\noi
$-$ Suppose first that 
$u^\circ \eqo u^\bullet$, then $\pi (u^\bullet) \eqo \pi (u^\circ)\eqo \pi (u)$ and $ \phi (u)\eqo \overline{\phi} (u)+0\eqo \phi (u^\bullet)$, which implies first $h(u^\bullet)\eqo h(u)$. If $v\! \neq \! u^\bullet$ is such that $\overleftarrow{v}\eqo u$, then $\phi (v) \leko \phi (u^\bullet)\eqo \phi (u)$ but $\pi (v)\eqo \pi (u)+ 1$ and thus $h(v) \leqo h(u)$.  

\noi
$-$ Next suppose that $u^\circ \! \neq \! u^\bullet$. Then $\phi (u)\eqo \phi(u^\bullet) +1 $ and $\pi (u^\bullet) \eqo \pi (u)+1$, which first implies that $h(u^\bullet)\eqo h(u)$. Let 
$v\! \neq \! u^\bullet$ such that $\overleftarrow{v}\eqo u$. Then $\phi (v) \leq \phi (u^\bullet) \leko \phi (u)$ and $\pi (v) \leqo \pi (u)+1$ and therefore $h(v) \leqo h(u)$. This proves (\ref{maxpath}) and it also completes the proof of the lemma. \cqfd 

\medskip

\noi
\textbf{Gaussian estimates for BRWs.}  We next want to remove the factors $\sigma_u$ in (\ref{coupling1}). To that end, we rely on the following bound for BRWs with Gaussian jumps. 
\begin{lemma} 
\label{Gaussjumps}
Let $\fftree\ino \bbT$ be finite.
We denote by $\Gamma (\fftree)\eqo \max_{u\in \fftree} |u|$ its \emph{total height}. 
Let $(Y_u)_{u\in \fftree \backslash \{ \varnothing\}}$ be independent $\bbR$-valued standard Gaussian r.v.s. Let 
$\sigma$ and $(\sigma_u)_{u\in \fftree \backslash \{ \varnothing\}}$ belong to $\bbR_+^*$. We set 
$$ \delta \eqo \max_{u\in \fftree \backslash \{ \varnothing\}} |\sigma^2_u \! -\! \sigma^2| , \quad 
\, R_u =\!  \sum_{v\in \rgeo \varnothing , u \rgeo} \sigma_vY_v \, , \quad    
R^\bullet_u  = \! \sum_{v\in \rgeo \varnothing , u \rgeo} Y_v, \quad u\ino \fftree \backslash \{ \varnothing\} $$
and $R_\varnothing\eqo R^\bullet_\varnothing \eqo 0$. To avoid trivialities, we assume that $\delta\! >\! 0$. Then, for all $z_1, z_2\ino \bbR_+^*$,
\begin{equation}
\label{Gausscontrol}
\bP \Big( \! \max_{u\in \fftree} \big| R_u\! -\! \sigma R^\bullet_u \big| \! >\!  z_1+z_2 \Big)\leq 4 \, \bP \Big( \! \max_{u\in \fftree} \big| R^\bullet_u \big| \! >\!  \sigma z_1/\delta\Big)+ 4 \exp\!  \left( \! \! -\frac{z_2^2\sigma^2}{ 2\, \delta^2  \Gamma (\fftree)}  
 \right) .
\end{equation}
If we furthermore assume that $\sigma^2  \geko \delta$, then 
\begin{equation}
\label{Gausscontrol2}
\bP \Big( \! \max_{u\in \fftree} \big| R^\bullet_u \big| \! >\!  z_1+z_2 \Big)\leq 4 \, \bP \Big( \! \max_{u\in \fftree} \big| R_u \big| \! >\! z_1\sqrt{\sigma^2 \! -\! \delta} \Big)+ 4 \exp\!  \left( \! \! -\frac{z_2^2 (\sigma^2 \! -\! \delta)}{ 2\, (\sigma^2 \! +\! \delta)  \Gamma (\fftree)}  
 \right). 
\end{equation}

\end{lemma}
\noi
\textbf{Proof.} We use 
the following comparison result for Gaussian processes (see e.g.~Ledoux and Talagrand \cite{LedTal91}, Corollary 3.12, p.75): 

\smallskip

\begin{compactenum}
\item[(\textbf{Comparison 1})] \emph{Let $n\ino \bbN^*$ and let $\mathbf X\eqo (X_1, \ldots, X_n)$ and  $\mathbf Y\eqo (Y_1, \ldots, Y_n)$ be centered Gaussian random vectors such that $\bE [X_i^2]\eqo \bE [Y^2_i]$ and 
$\bE \big[ (Y_i  - Y_j )^2 ]\leqo  \bE \big[ (X_i  - X_j )^2 ]$ for all $1\leqo i,j\leqo n$. Then for all $z \ino \bbR$, we get 
$ \bP \big( \max_{1\leq i\leq n} Y_i  \geko z \big) \leqo  \bP \big( \max_{1\leq i\leq n} X_i  \geko z \big) $.}
\end{compactenum}

\smallskip

\noi
Thanks to (\textbf{Comparison 1}) we prove the following.  

\smallskip

\begin{compactenum}
\item[(\textbf{Comparison 2})] \emph{Let $n\ino \bbN^*$, let $\mathbf X\eqo (X_1, \ldots, X_n)$ and  $\mathbf Y\eqo (Y_1, \ldots, Y_n)$ be centered Gaussian random vectors such that 
$\bE \big[ (Y_i \! -\! Y_j )^2 ]\leqo  \bE \big[ (X_i \! -\! X_j )^2 ]$ 
and $\bE\big[Y_i^2]\leqo \bE\big[X_i^2]$ for all $1\leqo i,j\leqo n$. We set 
$\overline{\sigma} \eqo   \max_{1\leq i\leq n} (\bE [X_i^2] )^{1/2}$. 
Then for all $z_1, z_2 \ino [0, \infty)$, we get }
\end{compactenum}
\begin{equation}
\label{comp2}
 \bP \big( \max_{1\leq i\leq n} |Y_i | \geko z_1+z_2 \big) \leqo  4\, \bP \big( \max_{1\leq i\leq n} |X_i|   \geko z_1 \big) +4 \exp\!  \left( \! \! -\frac{z_2^2}{ 2\, \overline{\sigma}^2 } \right) 
 \end{equation}
\noi 
\emph{Proof of (\ref{comp2}).} We follow the idea of the proof of Corollary 3.14 
Ledoux and Talagrand \cite{LedTal91}: w.l.o.g., we assume that $\mathbf X$ and $\mathbf Y$ 
are independent. Let $Z$ be a standard $\bbR$-valued Gaussian r.v.~that is independent of $\mathbf X$ and $\mathbf Y$. 
For all $1\leqo i\leqo n$, we let   $X^\prime_i\eqo X_i + (\overline{\sigma}^2 + \bE [Y_i^2] \! -\! \bE [X_i^2])^{1/2} Z$ 
and  $Y^\prime_i\eqo Y_i + \overline{\sigma}Z$ so that $(X^\prime_i)_{1\leq i\leq n}$ and $(Y^\prime_i)_{1\leq i\leq n}$ 
satisfy the assumption of (\textbf{Comparison 1}). Thus 
$ \bP \big( \max_{1\leq i\leq n} Y^\prime_i  \geko z \big)$ $\leqo$  $\bP \big( \max_{1\leq i\leq n} X^\prime_i  \geko z \big) $ 
for all $z\ino \bbR$. Since $X^\prime_i \leqo X_i + \overline{\sigma} (Z)_+$, we get for all $z_1, z_2\ino [0, \infty)$, 
\begin{eqnarray*}
\bP \big( \max_{1\leq i\leq n} Y_i  \geko z_1+z_2 \big) &=&  
2\bP \big( \overline{\sigma} Z \geqo 0 \, ; \,  \max_{1\leq i\leq n} Y_i  \geko z_1+z_2 \big) 
\leq 2\bP \big( \max_{1\leq i\leq n} Y^\prime_i  \geko z_1+z_2 \big) \\
& \leq & 2 \bP \big( \max_{1\leq i\leq n} X^\prime_i  \geko  z_1+z_2 \big) \leqo   
2 \bP \big( \max_{1\leq i\leq n} X_i + \overline{\sigma} (Z)_+  \geko  z_1+z_2 \big) \\
& \leq & 2\bP \big( \max_{1\leq i\leq n} X_i  \geko  z_1 \big)+  2\bP \big(  (Z)_+ \geko  z_2/\overline{\sigma} \big) \\
& \leq & 2\bP \big( \max_{1\leq i\leq n} |X_i|  \geko  z_1 \big) +2\exp\!  \left( \!  -\tfrac{1}{2} (z_2/\overline{\sigma})^2\right) 
\end{eqnarray*}
since $\bP (  (Z)_+\! >\!   z_2/\overline{\sigma}) \leqo \exp  ( -\frac{_1}{^2} (z_2/\overline{\sigma})^2 )  $.
We easily get a similar inequality for $-\mathbf X$ and $-\mathbf Y$. Since $\max_{1\leq i\leq n} |Y_i|$ is 
either equal to $\max_{1\leq i\leq n} Y_i $ or to $\max_{1\leq i\leq n} (-Y_i)$, it implies the desired inequality. \cq

\smallskip

We then use (\textbf{Comparison 2}) to prove (\ref{Gausscontrol}). To that end, for all $u\ino \fftree$, we set $R^\prime_u\eqo R_u\! - \! \sigma R^\bullet_u$ and $\mathtt D\eqo \max_{u\in \fftree \backslash \{ \varnothing\}} |\sigma_u \! -\! \sigma| $. 
Let $u^*$ such that $\mathtt D\eqo |\sigma_{u^*} \! -\! \sigma|$. Then $\delta \geqo |\sigma_{u^*}^2\! -\! \sigma^2|\eqo \mathtt D(\sigma_{u^*} + \sigma)\geqo \mathtt D \sigma$. 
Then observe that 
$(R^\prime _u)_{u\in \fftree}$ and  $(\frac{\delta}{\sigma} R^\bullet _u)_{u\in \fftree}$ are two centered Gaussian random vectors 
such that, for all $u,v\ino \fftree$, 
$$ \bE \big[ |R^\prime_u \! -\! R^\prime_v|^2 \big]=  \!\!   \!\!   \!\!   \!\!  \!\!   \!\! \!\!  \sum_{\qquad w\in \lgeo u,v\rgeo\backslash \{ u\wedge v\}}  \!\!   \!\!   \!\!   \!\!  \!\!   \!\!   (\sigma_w\! -\! \sigma)^2 \leq \, \mathtt D^2 d_{\mathtt{gr}} (u,v) \leq  \tfrac{\delta^2}{\sigma^2} d_{\mathtt{gr}} (u,v)=
\bE \big[ |\tfrac{\delta}{\sigma}  R^\bullet_u \! -\! \tfrac{\delta}{\sigma}  R^\bullet_v|^2 \big] , $$
where $d_{\mathtt{gr}}$ stands for the graph-distance on $\fftree$. Thus (\ref{Gausscontrol}) is a consequence of (\textbf{Comparison 2}) since here $\overline{\sigma}^2\! :=\!  \max_{u\in \fftree} \frac{\delta^2}{\sigma^2}\bE [|R^\bullet_u|^2 ]\eqo \frac{\delta^2}{\sigma^2} \Gamma (\fftree) $. 

  We next assume that $\sigma^2 \geko \delta$. Then 
$$ (\sigma^2 \! -\! \delta)^{-1}\bE \big[ |R_u \! -\! R_v|^2 \big]=  \!\!   \!\!   \!\!   \!\!  \!\!   \!\! \!\!  \sum_{\qquad w\in \lgeo u,v\rgeo\backslash \{ u\wedge v\}}  \!\!   \!\!   \!\!   \!\!  \!\!   \!\!   \frac{\sigma^2_w}{\sigma^2 \! -\! \delta} \geq \, d_{\mathtt{gr}} (u,v) = \bE \big[ | R^\bullet_u \! -\! R^\bullet_v|^2 \big] .$$
Then (\ref{Gausscontrol2}) follows again from (\textbf{Comparison 2}) since 
$$\overline{\sigma }^2\! :=\!  \, \max_{u\in t} \frac{\bE  [R^2_u  ]}{\sigma^2 \! -\! \delta}  \, = \,  \max_{u\in t}\sum_{v\in \rgeo \varnothing, u \rgeo} \frac{\sigma_v^2}{\sigma^2 \! -\! \delta}  \, \leq \frac{\sigma^2\!  + \delta}{\sigma^2 \! -\! \delta} \Gamma (t).$$  
This completes the proof of Lemma \ref{Gaussjumps}. \cqfd

\smallskip

 Recall that $(\mathtt{e}_1, \ldots, \mathtt{e_d})$ stands for the canonical basis of $\bbR^d$ 
that is equipped with the canonical scalar product $\langle \cdot, \cdot \rangle $ and with 
the associated Euclidean norm $\lvert \,  \cdot \, \rvert$. We also recall the notation 
$\sqrt{\beta} \ino \mathtt{Sym}_{^d}^{_+}$ for the square root of $\beta\ino  \mathtt{Sym}_{^d}^{_+}$. 
In this section, we provide a key estimate aimed to control the finite dimensional marginal distributions of 
$\bbR^d$-valued BRWs whose jumps are independent by group of siblings and have uniformly bounded $G$-moments. 

\begin{proposition}
\label{fdestimate} 
Let $\fftree\ino \bbT$ be finite. Let $\kappa \ino (0, 1)$ and let $G$ be a moment gauge function such that 
$G(x) \geqo |x|^{2+ \kappa}$ for all $x\ino \bbR$. 
Let $(\xi_u)_{u\in \fftree \backslash \{ \varnothing\}}$ be $\bbR^d$-valued r.v.s~that satisfy the following.

\smallskip

\begin{compactenum}
\item[$(a)$] The r.v.s~$\big( (\xi_{[u]\ast i} )_{1\leq i\leq k_u (\fftree)}; u \! \in \fftree \backslash \mathtt{Lf} (\fftree)\big)$ are independent.

\smallskip

\item[$(b)$] $M \! :=\! \max_{u\in \fftree \backslash \{ \varnothing\}} \bE \big[ G(|\xi_u|)  \big] \leko \infty$ 
and $\bE [ \xi_u] \eqo 0\, $ for all $u\ino \fftree \backslash \{ \varnothing\}$. 

\smallskip

\end{compactenum}
Let $(Y_u)_{u\in \fftree \backslash \{ \varnothing\}}$  be independent $d$-dimensional 
Gaussian random vectors whose covariance matrix is the identity. Let $\beta\eqo (\beta_{i,j})_{1\leq i,j \leq d} \ino  \mathrm{Sym}^{_+}_{^d}$.
We denote by $(S_u)_{u\in \fftree}$ and 
$(R_u)_{u\in \fftree}$ the BRWs whose jumps are resp.~$(\xi_u)_{u\in \fftree \backslash \{ \varnothing\}}$ and $(Y_u)_{u\in \fftree \backslash \{ \varnothing\}}$ and whose initial value is $S_\varnothing \eqo R_\varnothing \eqo0$. We set $\lVert \beta \rVert$ $ \!= \!$ $ \max_{1\leq i,j\leq d} |\beta_{i,j}|$, $ m_1 \eqo \max_{u\in \fftree \backslash \{ \varnothing\}} \bE \big[ |\xi_u| \big]$, 
and  
$$ \delta \eqo \max_{u\in \fftree \backslash \{ \varnothing\}} \max_{1\leq i,j\leq d}\big| \bE \big[ \langle \mathtt{e}_i , \xi_u \rangle  \langle \mathtt{e}_j , \xi_u \rangle \big] \! -\! \beta_{i,j}\big| \; .$$
Let $p\ino \bbN^*$, $u_1, \ldots, u_p \ino \fftree \backslash \{ \varnothing \}$ and $y_1, \ldots, y_p\ino \bbR^d$. We set $r\! = \!  \max_{1\leq j\leq p} |y_j| $. Then, 
\begin{eqnarray}
\label{fdcontrol}
\Big| \bE \Big[ \exp \!\! \!\! \!\! \!\! & & \!\! \!\! \!\! \!\! \Big( \!\! \!\! \sum_{\; \; 1\leq k\leq p}   \mathrm{i}\,  \langle y_k , S_{u_k} \rangle   \Big) \Big] \!
- \! \bE \Big[ \exp \Big( \!\! \!\! \sum_{\; \; 1\leq k\leq p} \!\! \!\!  \mathrm{i}\,   \langle y_k , \sqrt{\beta}.R_{u_k} \rangle   \Big) \Big]  \Big| \nonumber  \\
\!\! \!\!  & \leq & \!\! \!\!  p^3 \big( m_1\! +\!  \sqrt{d \lVert \beta \rVert} \big)r + d^2p^4 \big( |u_1| +\ldots + |u_p| \big) r^2 \big( \delta +
\lVert \beta \rVert^2 r^2+ M r^\kappa \big)\, .
\end{eqnarray}
\end{proposition}
\noi
\textbf{Proof.} 
Let $u\ino \fftree$ and $y\eqo (y^{(1)}\! , \ldots, y^{(d)}) \ino \bbR^d$. 
Since $\xi_u$ is centered, we first observe that  
$$ \big| \bE \big[ e^{\mathrm{i}  \langle y , \xi_u \rangle } \big] \! -\! 1 + \tfrac{1}{2} \bE [\langle y, \xi_u \rangle^2]  \big| \leqo 
\bE\big[ \langle y, \xi_u \rangle^2 \min (1,  \tfrac{1}{6}  |\langle y, \xi_u \rangle | )  \big] \leqo 
\bE\big[ | \langle y, \xi_u \rangle |^{2+\kappa} \big] \leqo M |y|^{2+\kappa} . $$
We next note that $| \bE [\langle y, \xi_u \rangle^2] \! -\!  \langle y, \beta.y \rangle| \eqo  | \sum_{1\leq i, j\leq d} y^{(i)}y^{(j)}\bE [\langle\mathtt e_i ,\xi_u  \rangle \langle\mathtt e_j ,\xi_u  \rangle  \! -\! \beta_{i,j}] | \leqo d |y|^2 \delta $. 
Since $0\leqo e^{-\frac{1}{2}  \langle y, \beta.y \rangle } \! -\! 1+ \frac{1}{2}  \langle y, \beta.y \rangle  \leqo \frac{1}{8}  \langle y, \beta.y  \rangle^2 \leq \frac{1}{8} d^2\lVert \beta \rVert^2 |y|^4$, we get 
 \begin{equation}
\label{caraxiu}
\forall u\ino \fftree\backslash \{ \varnothing \} , \; \forall y \ino \bbR^d, \quad 
\Big| \bE \big[ e^{\mathrm{i}  \langle y , \xi_u \rangle } \big] \! -\! e^{-\tfrac{_1}{^2}  \langle y, \beta.y \rangle} \Big| \leq \tfrac{1}{8}d^2 \lVert \beta \rVert^2 |y|^4+ \tfrac{1}{2}d |y|^2 \delta+M|y|^{2+ \kappa}\; .
\end{equation}

We next fix $p\ino \bbN^*\!$, $u_1, \ldots, u_p \ino \fftree \backslash \{ \varnothing \}$ and $y_1, \ldots, y_p\ino \bbR^d$. 
We define $T\eqo \bigcup_{1\leq j\leq p} \,  \rgeo \varnothing , u_j \rgeo$. Note that $\varnothing \! \notin \! T$. For all 
$u\ino T  $ we also set $I (u)\eqo \{ j\ino \{ 1, \ldots, p\} \! : \! u \! \in\,  \rgeo \varnothing , u_j \rgeo \}$ and  $y (u) \eqo \sum_{j\in I(u)} y_j$. 
First observe that 
$$ Z_1\! := \!  \sum_{u\in T } \langle y(u), \xi_u \rangle \eqo \!\! \sum_{1\leq j\leq p} \!\! \langle y_j, S_{u_j} \rangle
\quad \textrm{and} \quad   Z_2\! := \!  \sum_{u\in T } \langle y(u), \sqrt{\beta}.Y_u \rangle \eqo \!\! \sum_{1\leq j\leq p} \!\! \langle y_j, \sqrt{\beta}. R_{u_j} \rangle . $$
Note that $|y(u)| \leqo r \# I(u) \leqo p r$. 

We next introduce $B = \{ u_j \! \wedge u_k ; 1\leqo j\leko k \leqo p\}$ and $C\eqo \{ u\ino T: \overleftarrow{u} \ino B \}$. Consequently, for all distinct $u, v\ino T\backslash C$ we get $\overleftarrow{u}\! \neq \! \overleftarrow{v}$. Thus by Assumption $(a)$, the r.v.s~$(\xi_u)_{u\in T\backslash C}$, are independent. If we set $Z_1^*\eqo Z_1\! -\! \sum_{u\in C} \langle y(u), \xi_u\rangle$ and 
$Z_2^*\eqo Z_2\! -\! \sum_{u\in C} \langle y(u), \sqrt{\beta}.Y_u\rangle$, we therefore get 
$$ \bE \big[ e^{\mathrm{i} Z_1^*}\big]\eqo\!\!\!\! \!\!  \prod_{\quad u\in T\backslash C} \!\!\!\! \!\! \bE \big[ e^{\mathrm{i} \langle y(u), \xi_u \rangle}\big] \quad \textrm{and} \quad  \bE \big[ e^{\mathrm{i} Z_2^*}\big]\eqo \!\!\!\! \!\!  \prod_{\quad u\in T\backslash C} \!\!\!\! \!\!  \bE \big[ e^{\mathrm{i} \langle y(u), \sqrt{\beta }.Y_u \rangle}\big] .$$
Observe next that 
\begin{eqnarray*}
\big|  \bE \big[ e^{\mathrm{i} Z_1}\big] \! -\!  \bE \big[ e^{\mathrm{i} Z_1^*}\big]\big| \!\!\! &\leq & \!\!\! \bE \big[ |Z_1\! -\! Z^*_1|\big] \leq \sum_{u\in C} \bE \big[  |\langle y(u), \xi_u \rangle|  \big] \leq \sum_{u\in C} |y(u)| \bE \big[  |\xi_u |  \big] \\
\!\!\! & \leq & \!\!\! pr (\#C) m_1 \leq  p^2 r (\#B) m_1 \leq p^3 m_1r .
\end{eqnarray*}
Similarly, 
\begin{eqnarray*}
\big|  \bE \big[ e^{\mathrm{i} Z_2}\big] \! -\!  \bE \big[ e^{\mathrm{i} Z_2^*}\big] \big|  \!\!\! &\leq & \!\!\! \bE \big[ |Z_2\! -\! Z^*_2|\big] 
\leq \sum_{u\in C} \bE \big[  |\langle y(u), \sqrt{\beta}.Y_u \rangle|  \big] \leq \sum_{u\in C} \sqrt{\bE \big[ \langle y(u), \sqrt{\beta}.Y_u\rangle^2  \big] } \\
 \!\!\! &\leq & \!\!\!  \sum_{u\in C} \sqrt{\langle y(u), \beta. y(u) \rangle} \leq \sqrt{d\lVert \beta \rVert}  \sum_{u\in C} |y(u)| \leq \sqrt{d\lVert \beta \rVert} p^3 r .
\end{eqnarray*}
Consequently, 
\begin{equation}
\label{stepinde}
\big|  \bE \big[ e^{\mathrm{i} Z_1}\big] \! -\!  \bE \big[ e^{\mathrm{i} Z_2}\big]\big|  \leq p^3  (m_1\! +\!  \sqrt{d \lVert \beta \rVert}  )r + \Big|\!\!\!\! \!\!  \prod_{\quad u\in T\backslash C} \!\!\!\! \!\! \bE \big[ e^{\mathrm{i} \langle y(u), \xi_u \rangle}\big]   -  \!\!\!\! \!\!  \prod_{\quad u\in T\backslash C} \!\!\!\! \!\!  \bE \big[ e^{\mathrm{i}\langle y(u), \sqrt{\beta}.Y_u \rangle}\big]  \, \Big| . 
\end{equation}
We use next the following inequality: \emph{let $a_1, \ldots, a_n, b_1, \ldots, b_n\ino \mathbb C$ such that $|a_k|$ and $|b_k| \leqo 1$, $k\ino \{ 1, \ldots, n\}$; Then $\big| a_1\!  \ldots a_n   -  b_1\!  \ldots b_n \big| \leq \sum_{1\leq k \leq n} |a_k  - b_k | $.}
(\emph{Indeed,} set $p_k\eqo a_1 \! \ldots a_{k} b_{k+1} \! \ldots b_n$, with the convention that $p_0\eqo b_1\! \ldots b_n$ and $p_n \eqo a_1\ldots a_n$ and observe that $|p_{k} \! -\! p_{k-1}| \leqo  |a_k \! -\! b_k |$; since $|p_n \! -\! p_0| \leqo \sum_{1\leq k\leq n} |p_{k} \! -\! p_{k-1}|$, we get the desired inequality.) Therefore, 
\begin{eqnarray*} 
\Big|\!\!\!\! \!\!  \prod_{\quad u\in T\backslash C} \!\!\!\! \!\! \bE \big[ e^{ \mathrm{i} \langle y(u), \xi_u \rangle}\big]   -  \!\!\!\! \!\!  \prod_{\quad u\in T\backslash C} \!\!\!\! \!\!  \bE \big[ e^{ \mathrm{i} \langle y(u), \sqrt{\beta}.Y_u \rangle}\big]  \, \Big|  
\!\!\!\!  & \leq & \!\!\!\!\!\!   \sum_{u\in T\backslash C} \big| \bE \big[ e^{\mathrm{i}\langle y(u), \xi_u \rangle}\big]  \! -\! \bE \big[ e^{ \mathrm{i} \langle y(u), \sqrt{\beta}. Y_u \rangle}\big]  \big| \nonumber  \\
\textrm{by (\ref{caraxiu})} \qquad & \leq &\!\!\!\!\!\! \sum_{u\in T\backslash C} \big( \tfrac{1}{8} d^2 \lVert \beta \rVert^2|y(u)|^4+ \tfrac{1}{2}d |y(u)|^2 \delta+M |y(u)|^{2+\kappa}\big) \nonumber \\
 \leq  \tfrac{1}{{8}}\big(|u_1| +\!\!\!  & \ldots&\!\!\!   + |u_p| \big) p^2 r^2 \big( d^2 \lVert \beta \rVert^2 p^2r^2 \! +\!  4 d \delta \! +\!  8 p^\kappa Mr^\kappa \big) , \\
  \leq  d^2p^4 \big(|u_1| +\!\!\!  & \ldots&\!\!\!   + |u_p| \big) r^2 \big( \lVert \beta \rVert^2 r^2  \! +\!   \delta \! +\!  Mr^\kappa  \big) 
\end{eqnarray*}
since 
$\# T\leqo |u_1| + \ldots  + |u_p|$. This implies (\ref{fdcontrol}) by (\ref{stepinde}). \cqfd 

\smallskip
  
We next prove the following estimates on moment gauge funtions, and two consequences on \textbf{Sib-Ind} Assumptions that are used in the proofs of Theorems \ref{extenscv} and \ref{maincvsnake}. 
\begin{lemma}
\label{parasubaddgauge}  Let $\bS_n\eqo (S_{n,u})_{u\in \btt_n}$, $n\ino \bbN$, a sequence of $\bbR^d$-valued BRWs such that the assumptions \emph{\textbf{Sib-Ind}$_{_{\, }}$($\baa, \bbb, \bS_\cdot, G, C_\cdot, \fbeta$)} hold.  
Then, the following holds true. 
\begin{compactenum}

\smallskip

\item[$(i)$] There are $b\ino \bbR^*_+$ and $c\ino [1, \infty)$ such that for all $x,y\ino \bbR$, $G(x+y)\leqo cG(x)e^{b |y|} $.

\smallskip

\item[$(ii)$] We set $G_* (x) \! :=\! cG(0)e^{b|x|} $, $x\ino \bbR$, and $\beta_*\! :=\!  2(1+b)^{-2}$. Let $(Z_u)_{u\in \bbU}$, be i.i.d.~$\bbR$-valued centered 
r.v.s whose density is $\bgam (dx) \! :=\!  (2\beta_*)^{-1/2} \exp (-|x| (2/\beta_*)^{1/2}) dx$. 
For all $n\ino \bbN$ and all $u\ino \btt_n\backslash \{ \varnothing\}$, we set $S^*_{n,u} \eqo \sum_{v\in \, \rgeo \varnothing , u\rgeo} Z_v$. Then $G_*$ is a moment gauge function and the assumptions \emph{\textbf{Sib-Ind}$_{_{\, }}$($\baa, \bbb, \bS^*_\cdot, G_*, C_\cdot, \beta_*$)} hold.

\smallskip

\item[$(iii)$] We assume that the BRWs $(\mathbf S_n)_{n\in \bbN}$ are $\bbR$-valued. 
Let $(Y_u)_{u\in \bbU}$, be i.i.d.~standard $\bbR$-valued Gaussian r.v.s, which are furthermore assumed to be independent of the $(\mathbf S_n)_{n\in \bbN}$. For all $n\ino \bbN$ and all $u\ino \btt_n\backslash \{ \varnothing\}$, we set $S'_{n,u} \eqo S_{n,u}+  \sum_{v\in \, \rgeo \varnothing , u\rgeo} Y_v $. 
 Then, the assumptions \emph{\textbf{Sib-Ind}$_{_{\, }}$($\baa, \bbb, \bS'_\cdot, G, C_\cdot, 1+\fbeta$)} hold. 
\end{compactenum}
\end{lemma}
\noi
\textbf{Proof.} We recall that $G\! : \! \bbR \! \to \! \bbR_+^*$ is even, continuous and that 
there are $\kappa , x_0\ino \bbR_+^*$ such that 
$x\ino [x_0, \infty)\! \mapsto\! x^{-2 -\kappa} G(x)$ is nondecreasing and $x\ino [x_0, \infty)\! \mapsto \! x^{-1} \log G(x)$ is nonincreasing. 
Let $a\ino [1, \infty)$ such that $a^{-1} G(x_0)\leqo G(z) \leqo aG(x_0)$ for all $z\ino [-x_0, x_0]$. Since $G$ is nondecreasing on $[x_0, \infty)$, we get $G(x_0)\leqo a G(z)$ for all $z\ino \bbR_+$. This inequality also holds for all $z\ino \bbR$ since $G$ is even. Namely, 
\begin{equation}
\label{Gx0piv}
\forall z\ino \bbR, \quad G(x_0)\leq a G(z) . 
\end{equation}

Next, let $z\ino [x_0, \infty)$ and $y\ino \bbR_+$. 
We set $b\! :=\! x_0^{-1}\log G(x_0)$. Since $G(x) \! \to \! \infty$ as $x\! \to \! \infty$, there is $x\geko x_0$ such that $G(x) \geko 1$, which implies $0\leko x^{-1}\log G(x) \leqo b$ since $x\ino [x_0, \infty)\! \mapsto \! x^{-1}\log G(x) $ is nonincreasing. Thus $b\ino \bbR_+^*$. Similarly, $ \frac{\log G(z+y)}{z+y} \leqo \frac{\log G (z)}{z}\leqo b$. Thus $\log G(z+y) \leqo \log G(z) + by$, i.e., 
\begin{equation}
\label{xx0ypos}
\forall z\ino [x_0, \infty), \forall y\ino \bbR_+ , \quad G(z+y)\leqo G(z)e^{b y}.
\end{equation}

Since $G$ is even, we only need to prove $(i)$ for $x\ino \bbR_+$ and $y\ino \bbR$. 
If $x+y \geqo x_0$, and $x \geqo x_0$, then $G(x+y) \leqo G( x)\un_{\{ y\leq 0\} } +  G(x)e^{b y}\un_{\{ y>0\}}\leqo G(x)e^{b |y|}$ by (\ref{xx0ypos}) and since $G$ is nondecreasing on $[x_0, \infty)$. 
If $x+y\geqo x_0 $ and $x\ino [0, x_0]$, then $G(x+y)\leqo G(x_0) e^{b(x+y-x_0)}\leqo aG(x)e^{b|y|}$, by (\ref{xx0ypos}) 
and (\ref{Gx0piv}). If $x+y\ino [-x_0, x_0]$, then $G(x+y) \leqo aG(x_0)\leq a^2 G(x)\leqo a^2 G(x)e^{b|y|}$ by (\ref{Gx0piv}). 
Let suppose that $x+y\leqo -x_0$. Then $x+y\eqo -x_0\! -\! y'$, where $0\leqo y'\eqo |y| -x-x_0 \leqo |y|$. Consequently, 
$G(x+y)\eqo G(x_0+y')$ since $G$ is even, $G(x_0+y')\leqo G(x_0)e^{by'}$ by (\ref{xx0ypos}) and $G(x_0)  e^{by'} \leqo 
aG(x) e^{b|y|}$ by (\ref{Gx0piv}). Namely, $G(x+y)\leqo aG(x)e^{b|y|}$. This proves $(i)$ with $c\eqo a^2$. 

We next prove $(ii)$. Clearly $G_*$ is a moment gauge function such that $G (x)\leqo G_*(x) $ for all $x\ino \bbR$, 
by $(i)$. Moreover 
$\int_{\bbR} x^2 \bgam (dx) \eqo \beta_*$. 
To simplify notation, we refer to \textbf{Sib-Ind}$_{_{\, }}$($\baa, \bbb, \bS_\cdot, G, C_\cdot, \fbeta$) as to  \textrm{Sib-Ind}, and to \textbf{Sib-Ind}$_{_{\, }}$($\baa, \bbb, \bS^*_\cdot, G_*, C_\cdot, \beta_*$) as to \textrm{Sib-Ind}$^*$. 
We first note that \textrm{Sib-Ind} (6) $=$ \textrm{Sib-Ind}$^*$ (6). Since the jumps of $\mathbf S^*_n$ are i.i.d.~with a fixed density $\bgam$, which admits exponential moments and which is centered, \textrm{Sib-Ind}$^*$ (2) and (4) hold true. 
Then observe that \textrm{Sib-Ind} (5) immediately implies \textrm{Sib-Ind}$^*$ (5). 
Since $G \leqo G_*$, \textrm{Sib-Ind} (1) implies \textrm{Sib-Ind}$^*$ (1). Finally, we observe that 
a.s.~$M (\mathbf S_n^*, G_*)\eqo \int_{\bbR} G_*(x) \bgam(dx) \leqo cG(0)(2\beta_*)^{-1/2}\int_{\bbR} \exp (-|x| ) dx= cG(0) (b+1)$ and  \textrm{Sib-Ind}$^*$ (3) holds. This completes the proof of $(ii)$.

 Let us prove $(iii)$. We denote by 
$(\xi_{n,u})_{u\in \btt_n \backslash \{ \varnothing \}}$ the jumps of $\mathbf S_n$. Then, the jumps of  $\mathbf S'_n$ are 
$(\xi_{n,u}+Y_u)_{u\in \btt_n \backslash \{ \varnothing \}}$.  
 To simplify, we refer to \textbf{Sib-Ind}$_{_{\, }}$($\baa, \bbb, \bS'_\cdot, G, C_\cdot, 1+\fbeta$) as to \textrm{Sib-Ind}$'$. 
 First note that \textrm{Sib-Ind} (1)$=$\textrm{Sib-Ind}$'$ (1) and \textrm{Sib-Ind} (6) $=$ \textrm{Sib-Ind}$'$ (6). Then \textrm{Sib-Ind} (2) easily entails \textrm{Sib-Ind}$'$ (2). Let $\fbeta_n$ and $\fdelta_n$ be as in \textrm{Sib-Ind} (4) and recall for all $u\ino \btt_n$ that $\beta(\mathbf S_n, u)\eqo \bE [ \xi_{n,u}^2 | \btt_n ]$. Thus $\beta(\mathbf S'_n, u)\eqo \bE [(\xi_{n,u}+Y_u)^2|\btt_n]\eqo 1+ \beta(\mathbf S_n, u)$. Then \textrm{Sib-Ind}$'$ (4) and \textrm{Sib-Ind}$'$ (5) hold with 
 $ \fbeta_n'\! :=\! 1+ \fbeta_n$. Finally, by $(i)$, we get $\bE [G( \xi_{n,u} +Y_u) | \btt_n ] \leq c'\bE [G(\xi_{n,u} )|\btt_n ]$, where 
 $c'\eqo c\bE[e^{b|Y|}]$ and where $Y$ is a $\bbR$-valued standard Gaussian r.v. This a.s.~implies that $M(\mathbf S_n', G)\leq c'M(\mathbf S_n, G)$ and \textrm{Sib-Ind} (3) entails \textrm{Sib-Ind}$'$ (3), which completes the proof of the lemma.\cqfd 
 
\section{Proofs of Theorem \ref{Sheuexplain}, Theorem \ref{extenscv} and Proposition \ref{excntrex}}
 \label{Thm1pfsec}

\subsection{Preliminary estimates} 
\label{Thm1prelsec}

We discuss estimates related to \texttt{\L{}uka}$_{^{\,}}$($\baa, \bbb, \bmu, \psi$), \texttt{Grey}$_{^{\,}}$($\baa, \bbb, \bmu$) and \texttt{Sheu}$_{^{\,}}$($\baa, \bbb, \bmu$). 
Here, $(a_n)_{n\in \bbN}$ and $(b_n)_{n\in \bbN}$ stand for two sequences of positive real numbers that only satisfy $\lim_{n\to \infty} a_n \eqo \lim_{n\to \infty} b_n\eqo \infty$.  
We also fix a sequence $(\mu_n)_{n\in \bbN}$ of non-trivial, critical offspring distributions, i.e., $\mu_n$ satisfies (\ref{nontricri}).   
We recall from (\ref{genemuenn}) and  (\ref{psidisdef}) that for all $n\ino \bbN$, $r\ino [0, 1]$ and $\lambda \ino [0, a_n]$, 

\vspace{-4mm}

$$ g_{\mu_n} (r) \eqo \sum_{k\in \bbN} r^k \mu_n (k) \quad \textrm{and} \quad  \psi_n (\lambda)\eqo b_n \big( g_{\mu_n} \! \big(1\! -\! \tfrac{\lambda}{a_n} \big) \! -\! 1+ 
\tfrac{\lambda}{a_n} \big)= b_n \Psi_{\! \mu_n} (\lambda/a_n) \; .$$

\vspace{-2mm}

\noi
where $\Psi_{\! \mu_n} (r)\eqo g_{\mu_n} (1\! -\! r)\! -\! 1+r$. 

Let $\psi \ino \mathscr L$ be of the L\'evy-Khintchine form (\ref{LK}), with Brownian component 
$\beta_\psi$ and L\'evy measure $\pi_\psi$ and let $(X_s)_{s\in \bbR_+}$ be a Lévy process without 
negative jump whose Laplace exponent is $\psi$. 
We  denote by $(\xi_{n,k})_{k\in \bbN}$ a sequence of i.i.d.~r.v.s such that $\bP(\xi_{n,k} \eqo j)\eqo \mu_n (j+1)$ 
for all $j\ino \{-1\} \cup \bbN$. We recall that 
$V^n_k\eqo  \xi_{n, 1}+ \ldots + \xi_{n, k}$ and that Assumption 
\texttt{\L{}uka}$_{^{\,}}$($\baa, \bbb, \bmu, \psi$) means the convergence 
$V^n_{\lfloor b_n\rfloor }/a_n \! \to \! X_1$ weakly on $\bbR$. By Theorem 2.14 in Jacod \& Shiryaev \cite{JaSh02} 
(Chapter VII, Section 2.a, p.~398) \texttt{\L{}uka}$_{^{\,}}$($\baa, \bbb, \bmu, \psi$) holds true \emph{iff} the following is satisfied. 
\begin{compactenum}

\medskip

\item[$(\alpha)$] \emph{There is a bounded continuous $f_0\! :\! \bbR \! \to \! \bbR$ such that $f_0 (x)\eqo x$ in a neighbourhood of $0$ and}  
\end{compactenum}
$$  \lim_{n\to \infty} b_n  \bE \big[ f_0 \big( \tfrac{\xi_{n,1}}{a_n} \big) \big]\eqo \! \int_{\bbR_+^*} \!\! \! \! ( f_0(x) \!  -\! x ) \pi_\psi (dx) \; \textrm{\emph{and}} \,   \lim_{n\to \infty} b_n  \bE \big[ f^2_0  \big( \tfrac{\xi_{n,1}}{a_n} \big)  \big]\eqo  2\beta_\psi \! + \! \!  
\int_{\bbR_+^*} \!\! \!\! \! f^2_0(x)\pi_\psi (dx) .$$
\begin{compactenum}
\item[$(\beta)$] \emph{$\lim_{n\to \infty} b_n \bE \big[ g\big( \tfrac{\xi_{n,1}}{a_n} \big)\big]  \eqo  \int_{\bbR_+^*}  \!\! g(x) \pi_\psi (dx) $, for all bounded continuous $g \! : \! \bbR\! \rightarrow \! \bbR$ vanishing on a neighbourhood of $0$.} 

\medskip

\end{compactenum}
The following lemma discusses consequences of \texttt{Var}$_{\infty}$($\psi$) and \texttt{\L{}uka}$_{^{\,}}$($\baa, \bbb, \bmu, \psi$) concerning the relative growth of $b_n$ in terms of $a_n$ and the convergence of $\psi_n$.
\begin{lemma}
\label{Lukagrowth} We keep the above notations and assumptions. 
We assume \emph{\texttt{\L{}uka}$_{^{\,}}$($\baa, \bbb, \bmu, \psi$)}. Then the following holds true.
\begin{compactenum}

\smallskip

\item[$(i)$] Let $c_n \ino \big(1 \! -\!  \mu_n (1), \infty \big),$ $n\ino \bbN$, be such that $\lim_{n\to \infty}c_n \eqo \infty$. 
For all $n\ino \bbN$, we set $a^*_n\eqo a_n$, $b^*_n \eqo c_nb_n$, $\mu^*_n  (k) \eqo c_n^{-1}\mu_n (k)$ for all $k \ino \bbN\backslash \{ 1\}$ and $\mu^*_n (1)\eqo 1\! -\!  c^{-1}_n (1\! -\! \mu_n (1))$. Then, 
the $(\mu^*_n)_{n\in \bbN}$ are non-trivial, critical, they satisfy \emph{\texttt{\L{}uka}$_{^{\,}}$($\baa^*\! , \bbb^*\! , \bmu^*\! , \psi$)} and for all $z \in [0, a_n]$, 
\begin{equation}
\label{psipsietoi} \psi_n (z) \! :=\!  b_n \big( g_{\mu_n} \big( 1\! -\! \tfrac{z}{a_n}\big)\! -\! 1 + \tfrac{z}{a_n}\big)\eqo  b^*_n \big( g_{\mu^*_n} \big( 1\! -\! \tfrac{z}{a^*_n}\big)\! -\! 1 + \tfrac{z}{a^*_n}\big)\! =:\! \psi^*_n (z) .
\end{equation}

\smallskip

\item[$(ii)$] $\limsup_{n\to \infty} (1\! -\! \mu_n (1)) b_n a_n^{-2} \leq 4 \beta_\psi$.

\smallskip

\item[$(iii)$] If $\psi$ satisfies \emph{\texttt{Var}$_{\infty}$($\psi$)}, then $\lambda_n\!  := \! b_n / a_n \to \infty$. 

\smallskip

\item[$(iv)$]  For all $y\ino \bbR_+$, $\lim_{n\to \infty} \max_{\lambda \in [0, y]}|\psi (\lambda) \! -\! \psi_n (\lambda)| \eqo 0$. 

\end{compactenum}
\end{lemma}
\noi
\textbf{Proof.} We first prove $(i)$. The fact that $\mu^*_n$ is a non-trivial critical offspring distribution satisfying (\ref{psipsietoi}) is an immediate consequence of the definition. We then denote by $\xi^*_{n,1}$ a r.v.~such that $\bP(\xi^*_{n,1} \eqo k)\eqo \mu^*_n (k+1)$ for all $k\ino \{-1\} \cup \bbN$. 
We observe that $b_n^*\bE [h( \xi^*_{n,1} / a_n^*) ]\eqo b_n\bE [h( \xi_{n,1} / a_n) ]$ for any measurable and bounded $h\! : \! \bbR \to \bbR$ such that $h(0)= 0$. Thus, if $\baa$, $\bbb$ and $\bmu$ satisfy $(\alpha)$ and $(\beta)$, then 
$\baa^*$, $\bbb^*$ and $\bmu^*$ meet these conditions too, which implies \texttt{\L{}uka}$_{^{\,}}$($\baa^*\! , \bbb^*\! , \bmu^*\! , \psi$). 

Let us prove $(ii)$. Let $r\geko \epp\geko 0 $ and let 
$f_{r, \epp}\! :\! \bbR\! \to \! \bbR$ be an odd function such that $\smash{f_{r, \epp} (x)} $ $\eqo$ $\smash{ \big( x  -  (r+\epp)( \varepsilon^{-1}(x\! -\! r))_+ \big)_+ }$, $x\ino \bbR_+$. 
We first observe that $b_n a_{n}^{-2}\mu_n (0) \leqo  b_n\bE [ f^2_{r, \epp} ( \xi_{n,1} / a_n) ]$ for all $n\in \bbN$ such that $a_n r \geqo 1$. Since $\mu_n$ is critical, we then get $1\geqo 2 \mu_n (\bbN \backslash \{ 0, 1\} )+ \mu_n (1)\eqo 2-\mu_n (1) -2\mu_n (0)$, which entails $\mu_n (0) \geqo \frac{_1}{^2} (1\! -\! \mu_n (1))$. By $(\beta)$, $(\alpha)$ holds with $f_0\eqo f_{r, \epp}$. Then, the previous inequality and the second limit in $(\alpha)$ imply $\frac{_1}{^2}\limsup_{n\to \infty} (1\! -\! \mu_n (1))b_n a_n^{-2} \leqo 2\beta_\psi + \! \! 
\int_{\bbR_+^*} f^2_{r, \epp}(x)\pi_\psi (dx)$. Since it holds for arbitrarily small $r, \epp$, we get the desired inequality.

Let us prove $(iii)$. We assume \texttt{Var}$_{\infty}$($\psi$). 
Since $\mu_n$ is critical, $\bE [\xi_{n, 1}]\eqo 0$. Thus 
for all $n\ino \bbN$ such that $ra_n \geqo 1$, 
$$ -a_n  \bE \big[ f_{r, \epp} \big( \tfrac{\xi_{n,1}}{a_n} \big) \big]\leqo  -\bE \big[ \xi_{n,1}\un_{\{ \xi_{n,1 } \leq a_n r\}}]
=\bE \big[ \xi_{n,1}\un_{\{ \xi_{n,1 } > a_n  r\}}] \leqo \!\!\! \!\!  \sum_{\quad k > a_nr}\!\!\! \!\! k\mu_n (k)\leqo 1 .  $$
By $(\alpha)$ with $f_0\eqo f_{r, \epp}$, 
we get $\smash{\liminf_{n} \lambda_n \geqo \int_{\bbR_+^*} (x- f_{r, \epp}(x) ) \pi_\psi (dx) 
\geqo \int_{(r+\epp, \infty)}  x\pi_\psi (dx)}$. Since $r, \epp$ can be arbitrarily small and $\psi$ satisfies \texttt{Var}$_\infty(\psi)$,  if $\beta_\psi \eqo 0$, it implies $\lambda_n \! \to \! \infty$. 

  We now assume that $\beta_\psi \geko 0$. For all $n$ such that $a_nr \geqo 1$ we get the following. 
\begin{eqnarray*}
b_n \bE \big[ f^2_{r, \epp} \big( \tfrac{\xi_{n,1}}{a_n}\big)\big]  &\leq & b_n a_n^{-2}\mu_n (0)+ b_n a_n^{-2} \bE \big[\xi^2_{n,1}\un_{\{ 0 \leq \xi_{n,1}  \leq a_n (r+\epp)\} }\big] \\
&\leq & b_n a_n^{-2}\mu_n (0)+ (r+\epp) b_n a_n^{-1}\bE \big[\xi_{n,1}\un_{\{ \xi_{n,1} \geq 0 \} }\big] \leq \lambda_n \big( a_n^{-1} + r+\epp \big).
\end{eqnarray*}
The second limit in $(\alpha)$ then implies that $2\beta_\psi /(r+ \epp) \leqo \liminf_{n\to \infty} \lambda_n$. Since it holds for arbitrarily small $r, \epp$ and since we assume $\beta_\psi \geko0$, it implies $\lambda_n \! \to \! \infty$. 

We finally prove $(iv)$. By a result due to Grimvall on Laplace transforms (see Grimvall \cite{Gr74}, Theorem 2.1, p.~1029) we get for all $\lambda \ino \bbR_+$, that 
$\lim_{n\to \infty} \bE [\exp (-\lambda  V^{n}_{^{\lfloor b_n\rfloor}}/a_n)]\eqo \exp (\psi (\lambda))$. To simplify we set $h_n (\lambda) \eqo a_n (1\! -\! \exp (-\lambda/a_n))$, which tends to $\lambda$. Then we observe that 
$$ \bE \Big[ e^{- \frac{\lambda}{a_n} V^{n}_{\! {\lfloor b_n\rfloor}}} \Big] \eqo \Big(e^{\frac{\lambda}{a_n}} g_{\mu_n} \big(e^{-\frac{\lambda}{a_n}} \big) \Big)^{\lfloor b_n\rfloor} = \Big( 1+ \tfrac{b_n^{-1} \psi_n (h_n(\lambda))}{1- a^{-1}_n h_n(\lambda)}   \Big) ^{\lfloor b_n\rfloor}, $$
which implies that $\lim_{n\to \infty}\psi_n (h_n(\lambda))\eqo \psi (\lambda)$. Since $\psi_n$ and $\psi$ are continous and increasing, Dini's theorem implies the desired result. \cqfd

 \begin{remark}
 \label{Lukagrowthrem} By Lemma \ref{Lukagrowth} $(i)$, we see that \texttt{\L{}uka}$_{^{\,}}$($\baa , \bbb , \bmu , \psi$) imposes no limitation on how big $b_n$ can be with respect to $a_n$ but if $b_n a_n^{-2}$ is unbounded, then Lemma \ref{Lukagrowth} $(ii)$ entails that $\limsup_{n\to \infty} 
 \mu_n (1)\eqo 1$ and $(b_n)_{n\in \bbN}$ can not be seen as the genuine time-renormalization sequence for the RWs $(V^{_{(n)}}_{^k})_{k\in \bbN}$. However, if $(b_n)_{n\in \bbN}$ is the true time-renormalization, namely if $\limsup_{n\to \infty} 
 \mu_n (1) \leko 1$,  
then \texttt{\L{}uka}$_{^{\,}}$($\baa , \bbb , \bmu , \psi$) implies that $b_n \eqo \mathcal O (a_n^2)$ and furthermore we get $b_na_n^{-2} \! \to \! 0$ if $\beta_\psi \eqo 0$. We also observe that \texttt{Var}$_{\infty}$($\psi$) and \texttt{\L{}uka}$_{^{\,}}$($\baa , \bbb , \bmu , \psi$) actually impose the \emph{minimal} growth condition $b_n / a_n \! \to \! \infty$. \cq   
 \end{remark}

We next set 
\begin{equation}
\label{tauxenndef}
\tau^{_\infty}_n= (\tau_n (p))_{p\in \bbN^*} \quad \textrm{and} \quad \tau^{_x}_{n}= (\tau_n (p))_{1\leq p \leq \lfloor xa_n \rfloor } 
\end{equation}
where the $\tau_n(p)$ are i.i.d.~GW($\mu_n$)-trees and where $x\ino \bbR_+^*$. 
To simplify notation we set the following. For all $n\ino \bbN$ and for all $s\ino \bbR_+$, 
\begin{equation}
\label{rerescale}
C^{_{(n)}}_{s} \! \eqo \tfrac{1}{\lambda_n} C_{b_ns} (\tau^{_\infty}_n), \; H^{_{(n)}}_{s} \eqo  \tfrac{1}{\lambda_n} H_{b_ns} (\tau^{_\infty}_n), \;  \; \textrm{and} \;  \; V^{_{(n)}}_{s} \eqo \tfrac{1}{a_n} V_{b_ns} (\tau^{_\infty}_n) 
\end{equation}
where we recall that $\lambda_n \eqo b_n / a_n$, that $C_\cdot (\tau^{_\infty}_n)$ is the contour process of $\tau^{_\infty}_n$, $H_\cdot (\tau^{_\infty}_n)$ is its height process and $V_\cdot (\tau^{_\infty}_n)$ is its \L{}ukasiewicz path. 

In the following, we assume that $\psi $ satisfies at least \texttt{Var}$_{\infty}$($\psi$). 
We denote by $(H_s)_{s\in \bbR_+}$ the height process associated with the L\'evy process $X$ as defined by (\ref{defH}). The height process $H$ is at least lower semicontinuous and we recall from the introduction that $H$ can be chosen continuous 
\emph{iff} $\psi$ satisfies \texttt{Grey}$_{^{\,}}$($\psi$).
As already reminded in the introduction, under the assumptions 
\texttt{Norm}$_{^{\,}}$($\baa, \bbb$), \texttt{\L{}uka}$_{\,}$($\baa, \bbb, \bmu,\psi$) and \texttt{Grey}$_{^{\,}}$($\psi$),  Assumption \texttt{Hght}$_{^{\,}}$($\baa, \bbb, \bmu, \psi$) is equivalent to the following convergence  in law in $(\bC^{_0}_{^1})^2 \! \times \! \bD (\bbR_+, \bbR)$. 
\begin{equation}
\label{basiccvtreee} 
\big(C^{_{(n)}}_{\cdot}\! ,H^{_{(n)}}_{\cdot}\! ,V^{_{(n)}}_{\cdot} \big) \!  \longrightarrow \! \big( (H_{s/2})_{s\in \bbR_+} ,  (H_{s})_{s\in \bbR_+} ,  (X_{s})_{s\in \bbR_+} \big)
\end{equation}

 Let us discuss several basic consequences of the scaling limit (\ref{basiccvtreee}) that are used next. 
We  first set $\# \tau^{_x}_n\eqo \sum_{1\leq p\leq \lfloor a_n x\rfloor} \# \tau_n(p)$ 
and we denote by $\Gamma (\tau^{_x}_n)$ 
the total height of $\tau^{_x}_n$. Namely, 
$\Gamma (\tau^{_x}_n )\eqo \max \big\{ H_s (\tau^{_x}_n)  ; s\ino [0 , \# \tau^{_x}_n)]\big\}$. 
We also observe that $\# \tau^{_x}_n \eqo \inf \{ k \ino \bbN: V_k(\tau^{_\infty}_n) \eqo -\lfloor a_n x\rfloor\}$, 
by definition of the \L{}ukasiewicz path of $\tau^{_x}_n$. We set $\varsigma_{-x}\eqo \inf \{ s\ino \bbR_+ : 
X_s \leko -x\}$. Standard arguments combined with (\ref{basiccvtreee}) 
imply the following convergences in law in $\bbR_+$, \emph{jointly with (\ref{basiccvtreee})}
\begin{equation}
\label{szhghtGW}
\varsigma_{n,x} \! : =\! \tfrac{1}{b_n} \# \tau^{_x}_n  \xrightarrow[n\to \infty]{\; }\varsigma_{-x} \quad \textrm{and} \quad
 \tfrac{1}{\lambda_n} \Gamma( \tau^{_x}_n)  \xrightarrow[n\to \infty]{\; }\Gamma_{\! x}:=\max_{s\in [0, \varsigma_{-x}]} H_s .
\end{equation}
\emph{Indeed,} for the first limit, see Jacod \& Shiryaev \cite{JaSh02} Proposition 2.11, Chapter VI, Section 2a p.~341 or more precisely Broutin, D.~\& Wang \cite{BrDuWa21} Lemma B.3 $(iv)$; the second limit uses the first one, (\ref{basiccvtreee}) and quite standard arguments.

We then recall from Bertoin \cite{Be} Chapter VII, Theorem 1,  that $(\varsigma_{-x})_{x\in \bbR_+}$ is a subordinator with Laplace exponent $\psi^{-1}$, the inverse of $\psi$, and we recall from Theorem 1.4.1 and Corollary 1.4.2 in D.~\& Le Gall \cite{DuLG02}, that $\Gamma_{\! x}$ is distributed as the extinction time of a CSBP($\psi, x$).
Under Assumption \texttt{Grey}$_{^{\,}}$($\psi$), for all $\lambda, x, z\ino \bbR_+^*$, 
\begin{equation}
\label{szhghtLAW}
\bE \big[ e^{-\lambda \varsigma_{-x}} \big]\eqo e^{-x\psi^{-1} (\lambda)} \quad \textrm{and} \quad \bP (\Gamma_{\! x} \leqo z) \eqo e^{-xv(z)} \quad \textrm{where} \quad \int_{v(z)}^{\infty} \frac{\mathrm d r}{\, \psi (r)} = z\; .
\end{equation}

Let us assume that $\psi $ satisfies \texttt{Grey}$_{^{\,}}$($\psi$). We define the functions $G_{ n} \! :\! (0, a_n] \! \to \! \bbR_+$ and $G\! : \! \bbR_+^*\! \to \! \bbR_+^*$ by
\begin{equation}
\label{defGnG}
G_{ n} (y) = \int_{y}^{a_n}\!\!\!\! 
 \frac{ds}{\psi_n (s) } \quad \textrm{and} \quad G(y)= \int_{y}^\infty\!\!\!\! 
 \frac{ds}{ \psi (s)} \; . 
\end{equation} 
Note that (\ref{szhghtLAW}) implies $\lim_{y\to 0^+} G(y)\eqo \infty$. Similarly, we also get 
$\lim_{y\to 0^+} G_n(y)\eqo \infty$. \emph{Indeed}, we fix $y_0 \ino (0, a_n)$ and for all $s\ino (0, y_0]$, we get $\psi_n (s) \leqo s\psi_n (y_0) / y_0$ since $\psi_n$ is convex by Lemma \ref{fungeneprop}. 
Thus, for all $y\ino [0, y_0]$, we get  $G_n(y)\! -\! G_n(y_0) \geqo \frac{y_0}{\psi_n (y_0)} \log (y_0/y)$, which implies the desired result. \cq 

\smallskip

Therefore $G_n $ (resp.~$G$) is a decreasing one-to-one function from $(0, a_n]$ onto $\bbR_+$ (resp.~from $\bbR_+^*$ onto $\bbR_+^*$) and we denote by $v_n \! : \! \bbR_+ \! \to \! (0, a_n]$ its inverse (and $v$ in (\ref{szhghtLAW}) is the inverse of $G$). 
Observe that \texttt{Grey}$_{\,}$($\baa, \bbb, \bmu $) actually means that 
$\lim_{y\to \infty}\limsup_{n\to \infty} G_n(y)\! =\!  0$. 
The following lemma shows that \texttt{Grey}$_{\,}$($\baa, \bbb, \bmu $) implies (\texttt{Hght}$_{^{\,}}$($\baa, \bbb, \bmu $) $+$ \texttt{Grey}$_{^{\,}}$($\psi$)) and that in many cases these assumptions are equivalent. The proof relies on the inequalities (\ref{hghtestii}) in Lemma \ref{hghtestimm}. 
\begin{lemma}
\label{controlgrey}
 Let $\psi \ino \mathscr L$ satisfy \emph{\texttt{Var}$_{\infty}$($\psi$)}. Let 
$(a_n)_{n\in \bbN}$ and $(b_n)_{n\in \bbN}$ satisfy \emph{\texttt{Norm} $_{^{\,}}\! $($\baa, \bbb$)}. Let $(\mu_n)_{n\in \bbN}$ satisfy (\ref{nontricri}). We assume \emph{\texttt{\L{}uka}$_{\,}$($\baa, \bbb, \bmu,\psi$)}. Then the following holds true. 
\begin{compactenum}

\smallskip

\item[$(i)$]  \emph{\texttt{Grey}$_{^{\,}}$($\baa, \bbb, \bmu$)} implies  \emph{\texttt{Grey}$_{^{\,}}$($\psi$)} and  \emph{\texttt{Hght}$_{^{\,}}$($\baa, \bbb, \bmu$)}. 

\smallskip

\item[$(ii)$] Let us assume 
\emph{\texttt{Grey}$_{^{\,}}$($\psi$)}
and $\liminf_{n\to \infty} \mu_n(1) \geko 0$. Then \emph{\texttt{Grey}$_{^{\,}}$($\baa, \bbb, \bmu$)} is equivalent to  \emph{\texttt{Hght}$_{^{\,}}$($\baa, \bbb, \bmu$)}. 

\item[$(iii)$] Let $G_n$ and $G$ be as in (\ref{defGnG}) and let $v_n$ and $v$ stand for their inverses. We assume  \emph{\texttt{Grey}$_{^{\,}}$($\baa, \bbb, \bmu$)}. Then $\lim_{n\to \infty} v_n (z) \eqo v(z)$ for all $z\ino \bbR_+^*$. 
\end{compactenum}
\end{lemma}

\noi
\textbf{Proof.} Let us prove $(i)$. First note that \texttt{Grey}$_{\,}$($\baa, \bbb, \bmu$) implies that there exists $y_0 \ino \bbR_+^*$ such that $\limsup_{n\to \infty} \int_{y_0}^{a_n} \frac{\mathrm ds}{\psi_n (s)} \leqo 1$. By Lemma \ref{Lukagrowth} $(iv)$ and Fatou's Lemma, we get 
$$ \int_{y_0}^\infty\!\!\!\! 
 \frac{\mathrm d s}{ \psi (s)} \, \leq \, \liminf_{n\to \infty}  \int_{y_0}^{a_n} \frac{\mathrm ds}{\psi_n (s)}  \leq  \limsup_{n\to \infty}  \int_{y_0}^{a_n} \frac{\mathrm ds}{\psi_n (s)}  \leq 1 $$
and therefore $\psi$ satisfies \texttt{Grey}$_{^{\,}}$($\psi$). 

Let $x, \delta\ino \bbR_+^*$. Since $\bP (\Gamma(\tau^{x}_{n}) \leqo \lfloor \lambda_n \delta \rfloor )\eqo \exp \big(\lfloor a_n x\rfloor \log \bP (\Gamma (\tau_n(1)) \leqo \lfloor \lambda_n \delta \rfloor  ) \big)$, we get 
\begin{equation}
\label{equivheight}
\textrm{\texttt{Hght}$_{^{\,}}$($\baa, \bbb, \bmu$)}\; \; \Longleftrightarrow \; \forall \delta\ino \bbR_+^*\; : \; \limsup_{n\to \infty} 
a_n \bP \big( \Gamma (\tau_n(1)) \geko \lambda_n \delta   \big) \leko \infty\; .
\end{equation}
We now assume \texttt{Grey}$_{^{\,}}$($\baa, \bbb, \bmu$). Then for all $\delta \ino \bbR_+^*$, there necessarily 
exists $y_\delta \ino \bbR_+^*$ such that $\smash{\limsup_{n}\int_{y_\delta}^{a_n} \mathrm d r/\psi_n (r) \leko \frac{_1}{^2} \delta}$. By  (\ref{hghtestii}), we get $\smash{\int_{c_0 a_n  \bP ( \Gamma (\tau_n(1)) > \lambda_n \delta  )  }^{a_n} \mathrm d r/\psi_n (r) \geqo \frac{_1}{^2} \delta}$ for all $n\ino \bbN$, which implies $c_0 a_n  \bP ( \Gamma (\tau_n(1)) \geko \lambda_n \delta  ) \leqo y_\delta$ for all sufficiently large $n$. It entails \texttt{Hght}$_{^{\,}}$($\baa, \bbb, \bmu$) by (\ref{equivheight}), which completes the proof of $(i)$. 

Let us prove $(ii)$. We first assume \texttt{Hght}$_{^{\,}}$($\baa, \bbb, \bmu$), \texttt{Grey}$_{^{\,}}$($\psi$) and $2a\! :=\! \liminf_{n\to \infty} \mu_n (1) \geko 0$. 
By (\ref{szhghtGW}) and (\ref{szhghtLAW}), we get $\lim_{n\to \infty}\bP (\Gamma(\tau^{x}_{n}) \leqo \lfloor \lambda_n \delta \rfloor )\eqo \exp(-xv(\delta))$, which implies $v(\delta) \eqo \lim_{n\to \infty} a_n \bP \big( \Gamma (\tau_n(1)) \geko \lambda_n \delta   \big) $ 
for all $\delta \ino \bbR_+^*$. Consequently, 
$a_n \bP \big( \Gamma (\tau_n(1)) \geko \lambda_n \delta   \big) \leqo 2 v(\delta)$, for all sufficiently large $n$,  
and by the second inequality in (\ref{hghtestii}), we get 
$$ \int_{2 v(\delta)}^{a_n} \frac{\mathrm d r}{\psi_n (r)} \leq  
\int_{a_n \bP \big( \Gamma (\tau_n(1)) > \lambda_n \delta \big)}^{a_n} \frac{\mathrm d r}{\psi_n (r)} \leq 2 \delta 
\frac{\log \frac{1}{a}}{1\! -\! a} \; , $$
since $x\ino (0, 1] \mapsto \! (1\! -\! x)^{-1}\log 1/x$ decreases. This entails 
\texttt{Grey}$_{^{\,}}$($\baa, \bbb, \bmu$) since $\lim_{\delta\to 0^+} v(\delta)\eqo \infty$. This completes the proof of $(ii)$. 

Let us prove $(iii)$. By Lemma \ref{Lukagrowth} $(iv)$ and by Dini's result on uniform convergence of continuous monotone functions, for all $y_1, y_2 \ino \bbR_+^*$ such that $y_1\leko y_2$ we get 
$$ \mathtt{Err}_n (y_1,y_2)\! :=\! \max_{y\in [y_1, y_2]} \Big| G_n(y\! \wedge \! a_n) \! -\! G_n (y_2 \! \wedge \! a_n)  \! -\! \big( G (y) \! -\! G(y_2)\big) \Big| \xrightarrow[n\to \infty]{\; } 0 . $$
We next oberse that  
$ | G_n (y \wedge  a_n)\! -\! G(y)| \leqo \un_{[y_1, y_2]} (y)  \mathtt{Err}_n (y_1,y_2)+ G_n (y_2 \! \wedge \! a_n) + G(y_2) $, for all $y\ino [y_1, \infty)$. This implies for all $y_2 \ino (y_1, \infty)$, 
$$ \limsup_{n\to \infty}\max_{y\in [y_1, \infty)} | G_n (y)\! -\! G(y)| \leq G(y_2) + \limsup_{n\to \infty} G_n (y_2) \underset{y_2\to \infty}{-\!\!\! -\!\!\! -\!\!\! \longrightarrow} 0, $$ 
by \texttt{Grey}$_{\,}$($\baa, \bbb, \bmu$), which implies \texttt{Grey}$_{^{\,}}$($\psi$) by $(i)$. Elementary arguments then entail $(iii)$. \cqfd 

\medskip

Let $\psi \ino \mathscr L$ satisfy \texttt{Sheu}$_{^{\,}}$($\psi$). We define $F_{\! n} \! :\! (0, a_n] \! \to \! \bbR_+$ and $F\! : \! \bbR_+^*\! \to \! \bbR_+^*$ by
\begin{equation}
\label{defFnF}
F_{\! n} (y) = \int_{y}^{a_n}\!\!\!\! 
 \frac{ds}{2\sqrt{\! \int_0^s\! \psi_n (r) dr\, }} \quad \textrm{and} \quad F(y)= \int_{y}^\infty\!\!\!\! 
 \frac{ds}{2\sqrt{\! \int_0^s\! \psi (r) dr\, }} \; . 
\end{equation} 
Let us prove that $\lim_{y\to 0^+} F(y)\eqo \lim_{y\to 0^+} F_n(y) \eqo \infty$. \emph{Indeed}, we fix $y_0 \ino \bbR_+^*$. Since $\psi$ is convex, for all $s\ino (0, y_0]$, we get $\psi (s) \leqo s\psi (y_0) / y_0$ and thus $\smash{\int_0^s\! \psi (r) dr \leqo \frac{_1}{^2} s^2 \psi (y_0) / y_0}$. Therefore, for all $y\ino [0, y_0]$, we get $F(y)\! -\! F(y_0) \geqo (2y_0/\psi (y_0))^{1/2} \log (y_0/y)$, which implies the desired result for $F$. Since $\psi_n$ is also convex, similar arguments entail the same result for $F_n$. \cq 

\smallskip

Hence, $F_n $ (resp.~$F$) is a decreasing one-to-one function from $(0, a_n]$ onto $\bbR_+$ (resp.~from $\smash{\bbR_+^*}$ onto $\smash{\bbR_+^*}$) and we denote by $w_n\! : \! \bbR_+ \! \to \! (0, a_n]$  (resp.~$\smash{w\! : \! \bbR_+^*\! \to \! \bbR_+^*}$) its bijective inverse. 
Observe that \texttt{Sheu}$_{\,}$($\baa, \bbb, \bmu$) actually means that $\lim_{y\to \infty} \limsup_{n\to \infty} F_n(y)\eqo 0$. 
\begin{lemma}
\label{controlpsi} Let $\psi \ino \mathscr L$ satisfy \emph{\texttt{Var}$_{\infty}$($\psi$)}. Let 
$(a_n)_{n\in \bbN}$ and $(b_n)_{n\in \bbN}$ satisfy \emph{\texttt{Norm} $_{^{\,}\! }$($\baa, \bbb$)}. Let $(\mu_n)_{n\in \bbN}$ satisfy (\ref{nontricri}). Let $F_n$ and $F$ be as in (\ref{defFnF}) and let $w_n$ and $w$ be their respective inverses. We assume \emph{\texttt{\L{}uka}$_{\,}$($\baa, \bbb, \bmu,\psi$)} and \emph{\texttt{Sheu}$_{\,}$($\baa, \bbb, \bmu$)}.  
Then the following holds true. 
\begin{compactenum}

\smallskip

\item[$(i)$] $\psi$ satisfies \emph{\texttt{Sheu}$_{^{\,}}$($\psi$)}.

\smallskip

\item[$(ii)$] \emph{\texttt{Grey}$_{\,}$($\baa, \bbb, \bmu$)} holds true (and so does 
\emph{\texttt{Hght}$_{\,}$($\baa, \bbb, \bmu $)} by Lemma \ref{controlgrey}). 

\smallskip

\item[$(iii)$]  $\lim_{n\to \infty} w_n (z)$ $\eqo$  $w (z)$ for all $z\ino \bbR_+^*$. 
\end{compactenum}
\end{lemma}
\noi
\textbf{Proof.} $(i)$ is a consequence of Lemma \ref{Lukagrowth} $(iv)$ and Fatou's lemma: we argue as in the proof of Lemma \ref{controlgrey} $(i)$. 
Let us prove $(ii)$. Since $\psi_n$ is convex, for all $s\ino [1, a_n]$, we get $\psi_n (s)/s \geqo \psi_n (1)$. By Lemma \ref{Lukagrowth} $(iv)$, there is $n_0\ino \bbN$ such that for all 
integers $n\geqo n_0$, $\psi_n(1) \geqo \frac{_1}{^2} \psi (1) \geko 0$. Thus, for all $n\geqo n_0$ and $s\ino [1, a_n]$,
$ \smash{\int_0^s\! \psi_n (r) dr \leqo s\psi_n (s) \leqo \frac{2}{{\psi (1)}} \psi^2_n (s)}$. Consequently, 
$$ \forall n\geqo n_0, \; \forall y\ino [1, a_n], \quad \big(\frac{_{_1}}{^{^2}} \psi(1)\big)^{\frac{1}{2}} \int_{y}^{a_n}\!\!\!\! 
 \frac{ds}{\psi_n (s)} \, \leq  \int_{y}^{a_n}\!\!\! 
 \frac{ds}{\sqrt{\! \int_0^s\! \psi_n (r) dr\, }} , $$
which implies $(ii)$. The proof of $(iii)$ is quite similar to the proof of Lemma \ref{controlgrey} $(iii)$: we leave the details to the reader.  \cqfd  

\begin{lemma}
\label{sheustable} We assume the following: for all $n\ino \bbN$, $\mu_n\eqo \mu$, $b_n\eqo n$ and 
for all $\lambda\ino \bbR_+$, $\psi(\lambda) \eqo \lambda^\alpha$, with $\alpha \ino (1, 2]$. We also assume that there 
exists a sequence of positive real numbers $(a_n)_{n\in \bbN}$ tending to $\infty$ such that \emph{\texttt{\L{}uka}$_{\,}$($\baa, \bbb, \bmu,\psi$)} holds true. Then \emph{\texttt{Sheu}$_{\,}$($\baa, \bbb, \bmu$)} holds true too as well as 
\emph{\texttt{Grey}$_{\,}$($\baa, \bbb, \bmu$)} and \emph{\texttt{Hght}$_{\,}$($\baa, \bbb, \bmu$)}.  
\end{lemma}
\noi
\textbf{Proof.} Recall that $\psi_n (\lambda) \eqo n\Psi_{\! \mu} (\lambda/a_n)$, $\lambda\ino [0, a_n]$. Thus  
$\smash{2F_n (y)\eqo (\frac{a_n}{n})^{\frac{1}{2}} \! \int_{y/a_n}^1\!\!  \mathrm ds\,  (\int_0^s \Psi_{\! \mu} (r) \mathrm d r)^{-\frac{1}{2}} }$. Under our assumptions, (\ref{genedomain}) asserts that $ \smash{\Psi_{\! \mu} (r) \sim_{0^+} C_\alpha r^\alpha L(1/r)}$, where $L$ varies slowly at $\infty$ and 
$C_\alpha$ is a positive constant which only depends on $\alpha$. By the Karamata Theorems (after the change of variable $x\eqo 1/r$, see e.g.~Proposition 1.5.10 in Bingham, Goldies \& Teugels \cite{BiGoTe} p.~27) we see that 
$\smash{\int_0^s \Psi_{\! \mu} (r) \mathrm d r\sim_{0^+} C_\alpha(\alpha+1)^{-1} s^{\alpha+1} L(1/s)}$. Again by the 
Karamata Theorems (after the change of variable $x\eqo 1/s$, see e.g.~Proposition 1.5.8 in Bingham, Goldies \& Teugels \cite{BiGoTe} p.~26) we see that 
$$ F_n (y) \sim_{n\to \infty} C'_\alpha y^{-\frac{1}{2} (\alpha -1)} \sqrt{\frac{a_n^\alpha}{nL(a_n/y)}} \sim_{n\to \infty}  C'_\alpha y^{-\frac{1}{2} (\alpha -1)} $$
since $\smash{a_n^\alpha \sim_{n\to \infty} nL(a_n)\sim_{n\to \infty}  nL(a_n/y)}$. Here we have set $\smash{C'_\alpha \eqo ((\alpha \! -\! 1)^2 (\alpha+1)^{-1}C_\alpha)^{-1/2}}$. Thus $\smash{\limsup_{n\to \infty} F_n (y) \eqo C'_\alpha  y^{-\frac{1}{2} (\alpha -1)}}$, which implies \texttt{Sheu}$_{\,}$($\baa, \bbb, \bmu$). Lemma \ref{controlpsi} then implies that \texttt{Grey}$_{\,}$($\baa, \bbb, \bmu$) and \texttt{Hght}$_{\,}$($\baa, \bbb, \bmu$) hold true. \cqfd

 \subsection{Asymptotics of the contour of GW-forests with exponential lifespans}
 \label{PfsecThmGreyexplain}

To prove Theorem \ref{Sheuexplain}, we first need to prove a limit theorem for the contour process of GW-forests with exponential lifespans. This result also explains precisely the role of the assumption \texttt{Grey}$_{^{\,}}$($\baa, \bbb, \bmu$). 
Here we keep the basic assumptions and notations of the previous section. Namely, $\psi \ino \mathscr L$ satisfies \texttt{Var}$_{\infty}$($\psi$), the renormalization sequences $(a_n)_{n\in \bbN}$ and $(b_n)_{n\in \bbN}$ tend to 
$\infty$ and $(\mu_n)_{n\in \bbN}$ satisfy (\ref{nontricri}), i.e., they are non-trivial and critical. For all $n\in \bbN$, we set 
\begin{equation}
\label{defcaltaun}
\mathcal T_{ n} = \big( \mathcal T_n(p)\eqo \big(\tau_n(p), (\ell_{n,u} (p))_{u\in \tau_n(p)}  \big) \big)_{p\in \bbN^*} 
\end{equation}
where $\tau^{_\infty}_n\eqo (\tau_n(p))_{p\in \bbN^*}$ is as in (\ref{tauxenndef}), i.e., sequence of i.i.d.~GW($\mu_n$)-trees, and where, conditionally given $\tau^{_\infty}_n$, the $\ell_{n,u}(p)$ are i.i.d.~exponentially distributed r.v.s with parameter $1$.  We set 

\vspace{-4mm}

\begin{equation}
\label{rererescale}
\forall s\ino \bbR_+, \quad \mathscr C_{s}^{(n)} =  \tfrac{1}{\lambda_n}\mathscr C_{b_ns} (\mathcal T_n)
\end{equation}
where $(\mathscr C_s(\mathcal T_n))_{s\in \bbR_+}$ stands for the contour process of $\mathcal T_n$ (see Definitions \ref{Contlifespan} and \ref{lifespanforest}). For all $x\ino \bbR_+^*$, we set 
\begin{equation}
\label{forestxx}
\tau^{x}_n\eqo (\tau_n(p))_{1\leq p\leq \lfloor a_n x\rfloor} \quad \textrm{and} \quad  \mathcal T^x_n\eqo \big( \mathcal T_n (p)\big)_{1\leq p\leq \lfloor a_n x\rfloor} 
\end{equation}
and we denote by $\Gamma (\mathcal T^{x}_n)$ the total height of $\mathcal T^x_n$ and by $\Gamma (\tau^{_x}_n)$  the total height of  $\tau^{_x}_n$. Namely 
\begin{equation}
\label{heightforestxx}\Gamma (\mathcal T^{x}_n)\eqo \max_{1\leq p\leq \lfloor a_n x\rfloor} \max_{u\in \tau_n (p)} \zeta_{[p]\ast u}\quad \textrm{and} \quad \Gamma (\tau^{_x}_n)\eqo \max_{1\leq p\leq \lfloor a_n x\rfloor} \max_{u\in \tau_n (p)} |u|. 
\end{equation}
\begin{theorem}
\label{Greyexplain} Let $\psi \ino \mathscr L$ satisfy \emph{\texttt{Var}$_\infty$($\psi$)}. Let $a_n$, $b_n$, $\mu_n$, $\btt_n\eqo \tau^{_\infty}_n$, $n\ino \bbN$ and $(C,Y)$ be as in $\emph{\textbf{Case (0)}}$
and let $(\mathcal T_n)_{n\in \bbN}$ and $(\mathscr C^{_{(n)}}_s)_{s\in \bbR_+}$ be as above. 
We assume \emph{\texttt{\L{}uka}$_{^{\,}}$($\baa, \bbb, \bmu, \psi$)}. Then the following assertions are equivalent.

\begin{compactenum}

\smallskip

\item[$(a)$]  \emph{\texttt{Grey}$_{^{\,}}$($\baa, \bbb, \bmu$)}. 

\smallskip

\item[$(b)$] \emph{\textrm{\texttt{Grey}$_{^{\,}}$($\psi$)}} and $\smash{\tfrac{1}{{\lambda_n}}\Gamma (\mathcal T_n^x)\! \to \! \Gamma_{\! x}}$ in law on $\smash{\bbR_+}$ (where $\smash{\Gamma_{\! x}}$ is as in (\ref{szhghtLAW})).
%The laws of $\smash{(\tfrac{1}{\lambda_n}\Gamma (\mathcal T_n^x))_{n\in \bbN}}$ 
%are tight on $\bbR_+$. 

\smallskip

\item[$(c)$]$\lim_{n\to \infty} ( \mathscr C^{_{(n)}}_\cdot \! , C^{_{(n)}}_{\cdot}\! ,V^{_{(n)}}_{\cdot} )\eqo ( C_\cdot , C_\cdot, Y_\cdot)$ in law  in $(\bC^{_0}_{^1})^2 \! \times \! \bD (\bbR_+, \bbR)$. 

\smallskip

\item[$(d)$] \emph{\textrm{\texttt{Grey}$_{^{\,}}$($\psi$)}}, \emph{\texttt{Hght}$_{^{\,}}$($\baa, \bbb, \bmu $)} $\; $  and $\; $ $b_n / (a_n\log a_n)  \! \to \! \infty $. 

\smallskip

\end{compactenum}
Moreover if  \emph{\texttt{Grey}$_{^{\,}}$($\baa, \bbb, \bmu$)} holds true, then the processes $(\mathscr C^{_{(n)}}_{s}\!\! - C^{_{(n)}}_{s})_{s\in \bbR_+}$ tend to the null function in probability in $\bC^{_0}_{^1}$. 
\end{theorem}
\noi
\textbf{Proof.} We fix $n\ino \bbN$ and we denote by $\mathcal T'_n \eqo (\tau'_n, (\ell_{u})_{u\in \tau_n'})$ the $\bbR_+$-marked tree associated with $\mathcal T_n$ as defined in Remark \ref{contlifespanrem}. We denote by $(v_k)_{k\in \bbN}$ the contour exploration and
by $(u_l)_{l\in \bbN}$ the depth-first exploration of the vertices of $\tau'_n$. For all $u\ino \tau_n'$, we denote by $\zeta_u \eqo \sum_{v\in \lgeo \varnothing , u\rgeo} \ell_v$ and $\zeta^*_u\eqo \zeta_u \! -\! \ell_u$ the death- and birth-times of $u$ (we also recall the convention $\ell_\varnothing\eqo 0$). 
We recall $\Lambda_{\mathcal T_n}$ and $L_{\mathcal T_n}$ from Definitions \ref{Contlifespan} and \ref{lifespanforest}, and  Remark \ref{contlifespanrem}: they are continuous increasing functions from $\bbR_+$ onto $\bbR_+$, affine on each interval 
$[k\! -\! 1, k]$ and such that 
$\Lambda_{\mathcal T_n} (k)\eqo \sum_{1\leq i\leq k} |\zeta_{v_{i}} \!\!  -\! \zeta_{v_{i-1}} |$ and $L_{\mathcal T_n} (k)\eqo \sum_{0\leq j\leq k} \ell_{u_j}$, for all $k\ino \bbN^*$. 
We also recall from 
Remark \ref{contord} $(b)$ the definition of the continuous increasing time-change $\phi_{\tau^{\infty}_n}\! : \! \bbR_+ \! \to \! \bbR_+$ such that $H_s(\tau^{_\infty}_n) \eqo C_{\phi_{\tau^\infty_n} (s)} (\tau^{_\infty}_n)$, $s\ino \bbR_+$. 
For all $x\ino \bbR_+^*$, we recall from (\ref{forestxx}) and (\ref{heightforestxx}) the definition of the forests $ \tau^{_x}_n$ and $ \mathcal T^x_n$ and their respective total heights $\Gamma (\tau^{_x}_n)$ and $\Gamma (\mathcal T^{x}_n)$. 
We introduce the following additional notations. 
\begin{compactenum}

\smallskip

\item[$(i)$] $\Upsilon^x_n\eqo  2+\max_{1\leq p\leq \lfloor a_n x\rfloor} \max_{u\in \tau_n (p)} \big|\zeta_{[p] \ast u} \! -\! |u| \! -\! 1 \big|$. 

\smallskip

\item[$(ii)$]$\varsigma_{n, x}\eqo \frac{1}{{b_n}}\# \tau^{_x}_n$, where $\# \tau^{_x}_n \eqo 
\sum_{1\leq p \leq \lfloor a_n x\rfloor} \# \tau_n (p) \eqo \inf \big\{ k\ino \bbN: V_{k} (\tau^{_\infty}_n )\eqo -  \lfloor a_n x\rfloor \big\}$. 

\smallskip

\item[$(iii)$] $\mathbf r_{n,x}\eqo \frac{2}{b_n}L_{\mathcal T_n} (\# \tau^{_x}_n \! -\! 1)\eqo \frac{1}{b_n} 
\Lambda_{\mathcal T_n} (2(\# \tau^{_x}_n \! -\! 1))$. 

\smallskip

\item[$(iv)$] For all $s\ino \bbR_+$, $\phi_n (s)\eqo  \tfrac{1}{b_n} \phi_{\tau^{\infty}_n} (b_ns)$, $L_n (s)\eqo \tfrac{1}{b_n} L_{\mathcal T_n} (b_ns)  $ and $\Lambda_n (s)\eqo \tfrac{1}{b_n} \Lambda_{\mathcal T_n} (b_ns) $.

\smallskip

\end{compactenum}

\noi
Note that $\Gamma (\mathcal T^{x}_n)\eqo\lambda_n  \max_{s\in [0, \mathbf r_{n,x}]} \! \mathscr C_s^{_{(n)}}$ and $\Gamma (\tau^{_x}_n)\eqo \lambda_n  \max_{s\in [0, 2\varsigma_{n,x}]} C_s^{_{(n)}}$. Then by (\ref{compCTCt1}) in Remark \ref{contlifespanrem} 
$(a)$, we can compare the normalized contour processes of $\mathcal T_n$ and of $\tau^{_\infty}_n$ and 
we get for all $x\ino \bbR_+^*$ and all $n\ino \bbN$, 
\begin{equation}
\label{Cccontrol}
\max_{s\in [0, 2\varsigma_{n,x}]} \big|\mathscr C^{_{(n)}}_{\! \Lambda_n(s)}\! \! -\! C^{_{(n)}}_s   \big| \leq \tfrac{1}{\lambda_n}\Upsilon^x_n. 
\end{equation}
Moreover by (\ref{controphit}) in Remark \ref{contord} $(b)$ and by (\ref{controlLambda}) in Remark \ref{contlifespanrem} $(b)$ we get 
\begin{align}
\!\!\!\! \! \sup_{\quad s\in [0,  \varsigma_{n,x}]}\! \!\! \!\! \!\!\!  \big|\tfrac{1}{2} \phi_n(s)  -s \big| \leqo \tfrac{1}{2b_n}\Gamma (\tau^{_x}_n)\, & \qquad  \textrm{and}  \nonumber
 \\
   \sup_{\quad s\in [0,  \varsigma_{n,x}]} \!\!\!\! \! \!\!\!  \big|\Lambda_n(\phi_n(s))  & - 2 L_n(s) \big| \leqo  
\tfrac{2}{b_n} \Gamma (\mathcal T^x_n) \overset{\textrm{by (\ref{Cccontrol})}}{\leqo}   \tfrac{2}{b_n} \Gamma (\tau^{_x}_n)+ \tfrac{2}{b_n}\Upsilon^x_n. \label{timechangecontr}
\end{align}
To prove Theorem \ref{Greyexplain} we proceed in several steps. 

\smallskip

\noi
$\textbf{ Step (A)}.$  \emph{We first prove that the $(L_n(\cdot))_{n\in \bbN}$ tend to the identity, in probability in $\bC^{_0}_{^1}$, and that}
\begin{equation}
\label{rnxcv}
\mathbf r_{n,x} \xrightarrow[n\to \infty]{\textrm{(law)}}
 2\varsigma_{-x}\, , 
\end{equation}
\noi
\emph{Proof.} Observe that 
$(L_{\mathcal T_n}(l))_{l\in \bbN}$ is a RW whose jump law is $\mathtt{expo} (1)$. Since $L_{\mathcal T_n} (\cdot)$ is the continuous affine interpolation of this RW,  we get $\max_{r\in [0, s]} |L_n (r)\! -\! r| \leqo \frac{1}{b_n} \max_{1\leq l\leq \lceil b_ns \rceil } |L_{\mathcal T_n}(l) \! -\! l|$. 
The $L^2$-maximal inequality and Tchebychev inequality imply $\bP (\max_{r\in [0, s]} |L_n (r)\! -\! r| \geko \epp)\leqo  4\lceil b_ns \rceil /(b_n\epp)^2 $, which tends to $0$ as $n$ goes to $\infty$. 
Then (\ref{rnxcv}) follows from $\varsigma_{n,x} \! \to \! \varsigma_{-x}$ in (\ref{szhghtGW}). \cq

\smallskip

\noi
$\textbf{Step (B)}.$ \emph{We next prove that for all $x\ino \bbR_+^*$, }
\begin{equation}
\label{Greyandheight}
\textrm{ \texttt{Grey}$_{^{\,}}$($\baa, \bbb, \bmu$)} \quad \Longleftrightarrow \quad 
\textrm{\texttt{Grey}$_{^{\,}}$($\psi$) and $\smash{\tfrac{1}{{\lambda_n}}\Gamma (\mathcal T_n^x)\! \to \! \Gamma_{\! x}}$ in law on $\smash{\bbR_+}$,}
%The laws of $\smash{(\tfrac{1}{\lambda_n}\Gamma (\mathcal T_n^x))_{n\in \bbN}}$ 
%are tight on $\bbR_+$}.}
\end{equation}
%\emph{Moreover if one of the two previous conditions holds, then 
%$\smash{\tfrac{1}{\lambda_n}\Gamma (\mathcal T_n^x)\! \to \! \Gamma_{\! x}}$ in law on $\smash{\bbR_+}$, 
\emph{where $\smash{\Gamma_{\! x}}$ is distributed as specified in (\ref{szhghtLAW}).}

\noi
\emph{Proof.} We first observe that 
$\smash{\bP (\Gamma (\mathcal T_n^x) \leqo \lambda_n z)\eqo (1\! -\! \frac{1}{a_n} v_n(z))^{\lfloor a_n x \rfloor}}$ 
where we have set $\smash{v_n (z)}$ $\eqo$ $\smash{a_n\bP (\Gamma (\mathcal T_n (1)) \geko \lambda_n z)}$ 
and where we recall that $\smash{\mathcal T_n (1)}$ is a single GW($\mu_n,\! 1$)-tree. 
By (\ref{semigrey}), $v_n$ is the inverse of $G_n$ as defined by (\ref{defGnG}). 
Let us assume \texttt{Grey}$_{^{\,}}$($\baa, \bbb, \bmu$). By Lemma \ref{controlgrey}, \texttt{Grey} ($\psi$) holds and 
$\smash{\lim_{n\to \infty} v_n (z)\eqo v(z)}$ for all $\smash{z\ino \bbR_+^*}$. Here $v$ is as in (\ref{szhghtLAW}): 
the inverse of $G$, which is defined in (\ref{defGnG}). The previous arguments imply 
$\smash{\bP (\Gamma (\mathcal T_n^x) \leqo \lambda_n z)\! \to \! \bP (\Gamma_{\! x} \leqo z)}$, 
for all $\smash{z\ino \bbR_+^*}$. Namely, $\smash{\frac{_1}{^{\lambda_n}}\Gamma (\mathcal T_n^x)\! \to \! \Gamma_{\! x}}$. 
in law on $\smash{\bbR_+}$. 
Conversely, let us assume \texttt{Grey}$_{^{\,}}$($\psi$) and $\smash{\tfrac{1}{{\lambda_n}}\Gamma (\mathcal T_n^x)\! \to \! \Gamma_{\! x}}$ in law on $\smash{\bbR_+}$. Then, for all $z\ino \bbR^*_+$, 
$\lim_n \bP (\Gamma (\mathcal T_n^x) \leqo \lambda_n z)\eqo \bP (\Gamma_{\! x} \leqo z)$, which implies that 
$\lim_{n} v_n (z)\eqo v(z)$, where $v$ is as in (\ref{szhghtLAW}). Elementary arguments entail for all $y\ino \bbR_+^*$ 
that $\lim_n G(y)\eqo G(y)$, which entails \texttt{Grey}$_{^{\,}}$($\baa, \bbb, \bmu$). \cq

\smallskip

\noi
$\textbf{Step (C)}.$ \emph{We next prove that}
\begin{equation}
\label{Greybngrowth}
\textrm{ \texttt{Grey}$_{^{\,}}$($\baa, \bbb, \bmu$)} \quad \Longrightarrow \quad \frac{b_n}{a_n \log a_n} \xrightarrow[n\to \infty]{\;}\infty \; .
\end{equation}
\emph{Proof.} We set $\mathbf m_{n,x} \! : =\!  \max_{1\leq p\leq \lfloor a_n x\rfloor} \max_{u\in \tau_n (p)} \ell_{n,u}(p) $ and we observe that $\mathbf m_{n,x} \leqo \Gamma (\mathcal T^x_n) $. 
Then for all $\delta, z \ino \bbR_+^*$ we get 
\begin{eqnarray}
\bP \big( \Gamma (\mathcal T^x_n) \geqo \lambda_n z \big)  \!\! \!\! &\geq &  \!\! \!\!  \bP \big(\mathbf m_{n,x}  \geqo   \lambda_n z  \big) =  \bE \big[1\! -\!  \big(1\! -\! e^{-\lambda_n z} \big)^{\# \tau^x_n} \big]  \nonumber  \\
 \!\! \!\! & \geq &  \!\! \!\! \big( 1\! -\! e^{b_n \delta  \log ( 1 - e^{-\lambda_n z} ) }\big)  \bP \big( \# \tau_n^{_x} \geqo b_n \delta \big) \nonumber \\
   \!\! \!\! \!& \geq &  \!\! \!\! \!\big( 1\! -\! e^{- \delta e^{\log b_n -\lambda_n z}  }\big) \bP \big( \# \tau_n^{_x} \geqo b_n \delta \big) \label{extremearg}
\end{eqnarray}
We then recall from (\ref{szhghtGW}) that $\frac{1}{b_n}\# \tau_n^{_x}$ tends to $\varsigma_{-x}$ in law and we also 
recall from $\textbf{Step (B)}$ that \texttt{Grey}$_{^{\,}}$($\baa, \bbb, \bmu$) implies \texttt{Grey}$_{^{\,}}$($\psi$) and  $\frac{1}{\lambda_n} \Gamma (\mathcal T^x_n)\! \to \! \Gamma_{\! x}$ in law on $\bbR_+$. Let us suppose that there is an increasing sequence of integers $(n_k)_{k\in \bbN}$ such that $\lim_{k\to \infty} 
\log b_{n_k} \! -\! z\lambda_{n_k}\eqo \infty$. Then the previous convergences combined with (\ref{extremearg}) imply 
$ \bP (\varsigma_{-x} \geqo \delta) \leqo \bP (\Gamma_x \geqo z)\eqo 1\! -\! e^{-xv(z)} $. Since it holds for arbitrarily small $\delta$ and 
since $\lim_{\delta\to 0^+} \bP (\varsigma_{-x} \geqo \delta)\eqo \bP (\varsigma_{-x} \geko 0)\eqo 1$, we get a contradiction. Thus $\limsup_n (\log b_n \! -\! z\lambda_n) \leko \infty$. This implies that $\limsup_n (\log a_n \! -\! z\lambda_n) \leko \infty$ and 
$\limsup_n (a_n \log a_n)/b_n \leqo z$. Since it holds for arbitrarily small $z$, we get the desired result. \cq

\smallskip

\noi
$\textbf{Step (D)}.$ \emph{For all $\epp, x, y\ino \bbR_+^*$ and for all $n\ino \bbN$ such that $\lambda_n \geqo 1\vee y^{-1}$, we get} 
\begin{equation}
\label{controlUpsi}
\bP \big( \Upsilon^x_n \geq 2\! + \! \lambda_n \epp   \big) \leq  \bP \big( \Gamma (\tau^{_x}_n) \geqo\lambda_n y \big) + 2b_nx(y+1) \exp \big( \! -\! \tfrac{\lambda_n\epp^2}{2( 2y+ \epp)}\big)  
\end{equation}

\noi
\emph{Proof.} Let $(\mathcal E_l)_{l\in \bbN^*}$ be i.i.d.~$\mathtt{expo} (1)$ r.v.s. We set $S_l \eqo \sum_{1\leq k\leq l} \mathcal E_k $ and we use the following elementary deviation inequality whose proof is left to the reader. 
\begin{equation}
\label{devRWexpo} 
\forall l\ino \bbN^*, \, \forall y\ino \bbR_+,  \quad \bP \big( |S_l \! -\! l | \geqo y \big) \leq 2e^{-\frac{y^2}{2(y+l)}}\; .
\end{equation}

\noi
We first observe for all $\epp, y \ino \bbR_+^*$ that
$$ \bP \big( \Upsilon^x_n \geqo 2+ \lambda_n \epp\,  \big| \, \tau^{_x}_n  \big) \leq \un_{\{ \Gamma (\tau^{_x}_n) \geq \lambda_n y \}}+ \!\!\! \sum_{1\leq p\leq \lfloor a_n x \rfloor} \sum_{u\in \tau_n (p)} \un_{\{ 0\leq   |u| < \lambda_n y \}} \bE \big[ \un_{\{ |\zeta_{[p] \ast u} -|u|-1 | \geq \lambda_n \epp \}} \big|   \tau^{_x}_n  \big] \; .$$
Then, we note that if $u\ino \tau_n(p)$ is such that $l\! :=\! |u| +1\leko \lambda_ny+1\leqo 2\lambda_ny$, then  
$$ \bE \big[ \un_{\{ |\zeta_{[p]\ast u} -|u| -1| \geq \lambda_n \epp \}} \big|   \tau^{_x}_n  \big]\eqo \bP (|S_{l}\! -\! l| \geq \lambda_n \epp ) \leq 2\exp \big(\! -\! \tfrac{\lambda^2_n\epp^2}{2(l+ \lambda_n \epp)}\big)  \leq  2\exp \big( \! -\! \tfrac{\lambda_n\epp^2}{2( 2y+ \epp)}\big)$$
Then, we recall that $Z^{_x}_{^{n,l}}\eqo \sum_{1\leq p\leq \lfloor a_n x \rfloor} \sum_{u\in \tau_n (p)} \# \{ u\ino \tau^{_x}_n\!  :\!  |u|\eqo l\}$, $l\ino \bbN^*$, is a GW($\mu_n$) critical Markov chain. In particular $\bE [Z^{_x}_{^{n,l}}]\eqo Z^{_x}_{^{n,0}}\eqo  
\lfloor a_n x \rfloor$ and we get 
$$\bP \big( \Upsilon^x_n \geqo 2\! + \! \lambda_n \epp   \big) \leq   \bP \big( \Gamma (\tau^{_x}_n) \geqo\lambda_n y \big) + 2\exp \big( \! -\! \tfrac{\lambda_n\epp^2}{2( 2y+ \epp)}\big)  \!\!\!\!\!\! \!\! \sum_{\quad 0\leq l< \lambda_n y}\!\!\!\!\!\!  \bE \big[ Z^{_x}_{^{n,l}} ] $$
which easily implies (\ref{controlUpsi}) since $\lceil \lambda_n y \rceil \lfloor a_n x \rfloor \leqo b_n x (y+ \lambda_n^{-1})$.   \cq

\smallskip

\noi
$\textbf{Step (E)}.$ \emph{Assume \emph{\texttt{\L{}uka}$_{^{\,}}$($\baa, \bbb, \bmu, \psi$)}, \emph{\texttt{Grey}$_{^{\,}}$($\psi$)}, \emph{\texttt{Height}$_{^{\,}}$($\baa, \bbb, \bmu$)} and $\smash{\lim_{n} \frac{b_n }{a_n\log a_n}\eqo \infty}$.
Then $\lim_{n\to \infty}\frac{1}{\lambda_n}\Upsilon^x_{\! n}\eqo 0$ in probability for all $x\ino \bbR_+^*$ and }
\begin{equation}
\label{Lambdaphicv}
\textrm{\emph{$\frac{_{_1}}{^{^2}} \phi_n$ and $\Lambda_n$ both converge to the identity in probability in $\bC^{_0}_{^1}$.}}
\end{equation}

\noi
\emph{Proof.} As already mentioned, under \texttt{\L{}uka}$_{^{\,}}$($\baa, \bbb, \bmu, \psi$), \texttt{Grey}$_{^{\,}}$($\psi$) 
and \texttt{Height}$_{^{\,}}$($\baa, \bbb, \bmu$), (\ref{basiccvtreee}) and (\ref{szhghtGW}) hold true. 
This, combined with (\ref{controlUpsi}) in $\textbf{Step (D)}$ and the assumption $\lim_{n\to \infty} b_n /(a_n\log a_n)\eqo \infty$, implies that 
$\limsup_{n\to \infty} \bP (\Upsilon^x_{\! n}\geqo 2+ \lambda_n \epp) \leqo \bP (\Gamma_{\! x} \geqo y)$ for all $x,y, \epp\ino \bbR_+^*$. Therefore, it shows that  
$\lim_{n\to \infty}\frac{1}{\lambda_n}\Upsilon^x_{\! n}\eqo 0$ in probability for all $x\ino \bbR_+^*$. 
By (\ref{timechangecontr}), so does $\sup_{s\in [0,  \varsigma_{n,x}]}|\frac{_{1}}{^2} \phi_n(s) \! -\! s |  $. Then, for all $s, \epp \ino \bbR_+^*$, 
$$ \bP \Big( \sup_{r\in [0,  s]}\big|\tfrac{1}{2} \phi_n(r) \! -\! r \big| \geko \epp \Big) \leq  \bP \Big( \sup_{r\in [0,  \varsigma_{n,x}]}\big|\tfrac{1}{2} \phi_n(r) \! -\! r \big| \geko \epp \Big) + \bP \big(  \varsigma_{n,x} \leqo s \big) $$
Thus by (\ref{szhghtGW}), $\limsup_{n\to \infty}  \bP \big( \sup_{r\in [0,  s]}\big|\tfrac{1}{2} \phi_n(r) \! -\! r \big| \geko \epp) \leqo \bP (\varsigma_{-x} \leqo s)$, which tends to $0$ as $x$ goes to $\infty$ since a.s.~$\lim_{x\to \infty} \varsigma_{-x}\eqo \infty$. Thus 
$\frac{_1}{^2}\phi_n$ tends to the identity in probability, in $\bC^{_0}_{^1}$
and so does $s \mapsto 2\phi_n^{-1}(s)$, by an elementary argument.   
Similarly, by the second bound in (\ref{timechangecontr}) and \textbf{Step (A)}, $\frac{_1}{^2}\Lambda_n \circ \phi_n$ tends to the identity in probability in $\bC^{_0}_{^1}$ and by easy arguments $\Lambda_n$ too. \cq

\smallskip

Let us denote by $\Lambda_n^{-1}$ the inverse of the continuous increasing process $\Lambda_n$ from $\bbR_+$ onto $\bbR_+$.

\smallskip

\noi
$\textbf{Step (F)}.$ \emph{Assume \emph{\texttt{\L{}uka}$_{^{\,}}$($\baa, \bbb, \bmu, \psi$)}, \emph{\texttt{Grey}$_{^{\,}}$($\psi$)}, \emph{\texttt{Height}$_{^{\,}}$($\baa, \bbb, \bmu$)} and  $\smash{\lim_{n} \frac{b_n }{a_n\log a_n}\eqo \infty}$.
Then }
\begin{equation}
\label{scrCCcontigu}
\textrm{\emph{$ (\mathscr C^{_{(n)}}_s \! -\! C_{{\!\Lambda_n^{\! {-1}} \! (s)}}^{_{(n)}} )_{s\in \bbR_+}$ converges to the null process 
in probability in $\bC^{_0}_{^1}$.}}
\end{equation}
\noi
\emph{Proof.} By $\textbf{Step (E)}$, 
we get that $\smash{\lim_{n\to \infty}\frac{1}{\lambda_n}\Upsilon^x_{\! n}\eqo 0}$ in probability for all $\smash{x\ino \bbR_+^*}$. By (\ref{Cccontrol}), so does  
$\smash{\max_{s\in [0, 2\varsigma_{n,x}]} |\mathscr C^{_{(n)}}_{\! \Lambda_n(s)}\! \! -\! C^{_{(n)}}_s  |}$. 
Then, for all $s, \epp \ino \bbR_+^*$, 
$$ \bP \Big(\max_{r\in [0, s]} \big|\mathscr C^{_{(n)}}_{\! \Lambda_n(r)}\! \! -\! C^{_{(n)}}_r   \big|\geko  \epp \Big) \leq  \bP \Big(\max_{r\in [0, 2\varsigma_{n,x}]} \big|\mathscr C^{_{(n)}}_{\! \Lambda_n(r)}\! \! -\! C^{_{(n)}}_r   \big| \geko \epp \Big) + \bP \big(  2\varsigma_{n,x} \leqo s \big) $$
Thus by (\ref{szhghtGW}), $\smash{\limsup_{n\to \infty}  \bP (\max_{r\in [0, s]} |\mathscr C^{{(n)}} (\Lambda_n(r)) \! -\! C^{_{(n)}}_r  |\geko  \epp ) \leqo \bP (2\varsigma_{-x} \leqo s)}$, which tends to $0$ as $x$ goes to $\infty$ since a.s.~$\lim_{x\to \infty} \varsigma_{-x}\eqo \infty$. Thus 
$\smash{(\mathscr C^{{(n)}} ( \Lambda_n(s))\! -\! C^{_{(n)}}_s)_{s\in\bbR_+}}$ 
tends to the null process in probability in $\bC^{_0}_{^1}$. 
By (\ref{Lambdaphicv}) in $\textbf{Step (E)}$,  
$\Lambda_n$ converges to the identity in probability in $\bC^{_0}_{^1}$. So does its inverse $\smash{\Lambda_n^{-1}}$ and easy arguments entail (\ref{scrCCcontigu}).  \cq

\smallskip

\noi
\textbf{Proof of Theorem \ref{Greyexplain} $(a) \Longleftrightarrow (b)$}. This is $\textbf{Step (B)}$. \cq

\smallskip

\noi
\textbf{Proof of Theorem \ref{Greyexplain} $(a) \Longleftrightarrow (d)$}. We overall assume \texttt{\L{}uka}$_{^{\,}}$($\baa, \bbb, \bmu, \psi$) and \texttt{Var}$_\infty$ ($\psi$). 
By Lemma \ref{Lukagrowth} $(iii)$, $\lim_{n\to \infty}\lambda_n\eqo \infty$, so \texttt{Norm}$_{^{\,}}$($\baa, \bbb$) holds true too. Then Lemma \ref{controlgrey} $(i)$ and 
(\ref{Greybngrowth}) in $\textbf{Step (C)}$ show that $(a) \Longrightarrow (d)$. 

Conversely, let us assume $(d)$. By $\textbf{Step (E)}$, $\smash{\lim_{n\to \infty}\frac{1}{\lambda_n}\Upsilon^x_{\! n}\eqo 0}$ 
in probability for all $\smash{x\ino \bbR_+^*}$. Since $\smash{|\Gamma (\mathcal T^x_n) \! -\! \Gamma (\tau^{_x}_n) |
\leqo \Upsilon^x_{\! n}}$, (\ref{szhghtGW}) entails that $\smash{\frac{1}{\lambda_n} \Gamma (\mathcal T^x_n)\!
\to \! \Gamma_{\! x}}$ weakly on $\smash{\bbR_+}$ and $\textbf{Step (B)}$ entails \texttt{Grey}$_{^{\,}}$($\baa, \bbb, \bmu$). \cq  

\smallskip

\noi
\textbf{Proof of Theorem \ref{Greyexplain} $(a) \Longleftrightarrow (c)$}. We overall 
assume \texttt{Var}$_{\infty}$ ($\psi$) and \texttt{\L{}uka}$_{^{\,}}$($\baa, \bbb, \bmu, \psi$), which imply \texttt{Norm}$_{^{\,}}$($\baa, \bbb$). Let us first assume $(c)$, which entails implicitly \texttt{Grey}$_{^{\,}}$($\psi$). 
By (\ref{rnxcv}) in \textbf{Step (A)} and by an elementary argument, we get 
that $\smash{\max_{s\in [0, \mathbf{r}_{n,x}]} \mathscr C^{_{(n)}}_{s}\! \to \! \max_{s\in [0, 2\varsigma_{-x}]} H_{s/2}}$, i.e., that $\smash{\frac{1}{\lambda_n} \Gamma (\mathcal T^x_n) \! \to \! \Gamma_{\! x}}$ weakly on $\smash{\bbR_+}$, which entails $(b)$ and thus $(a)$ as already proved.

Conversely, let us assume $(a)$. We already proved that it implies $(d)$. Thus (\ref{basiccvtreee}) 
and (\ref{szhghtGW}) hold true and, by $\textbf{Step (E)}$, $\Lambda_n$ converges to the identity in 
probability in $\bC^{_0}_{^1}$, and thus its inverse $\smash{\Lambda_n^{-1}}$ too.  
This, combined with the weak convergence $\smash{C^{_{(n)}}_\cdot\! \to \! C_\cdot}$ in $\bC^{_0}_{^1}$, 
implies that $\smash{(C^{{(n)}} (\Lambda_n^{-1} (s)) \! -\! C^{_{(n)}}_s)_{s\in \bbR_+}}$ converges to the null function in probability in $\bC^{_0}_{^1}$. By (\ref{scrCCcontigu}), $\smash{(\mathscr C^{_{(n)}}_{s}\!\! - \! C^{_{(n)}}_{s})_{s\in \bbR_+}}$ tends to the null function in probability in $\bC^{_0}_{^1}$. By (\ref{basiccvtreee}), it entails $(c)$. This completes the proof of Theorem 
\ref{Greyexplain}. \cqfd

\subsection{Proof of Theorem \ref{Sheuexplain}}
 \label{ThmSheuexplainPfsec}

We first prove the convergence of finite dimensional marginal laws of discrete snakes with i.i.d.~jumps and of Brownian snakes indexed by GW-trees with exponential lifespans thanks to the two following general lemmas, which are a simple consequence of the continuity property (\ref{contoplus}). 
\begin{lemma}
\label{fdcvsna} Let $\smash{\beta \ino \mathrm{Sym}^{+}_{d}}$ and let $\smash{\sqrt{\beta}\!  \in \! \mathrm{Sym}^{+}_{d}}$ be a 
square root of $\beta$, asdiscussed before. Let $C_\cdot$ be a continous nonnegative process such that $C_0\eqo 0$. We denote by $W$ the (not necessarily continuous) $d$-dimensional Brownian snake with lifetime process 
$C_\cdot$ as in Definition \ref{Brosnadef}. Let $\smash{(a_n)_{n\in \bbN}}$ and $\smash{(b_n)_{n\in \bbN}}$ satisfy \emph{\texttt{Norm} ($\mathbf a, \mathbf b$)}. For all $n\ino \bbN$, let $\smash{\bS_n \eqo (S_{n, u})_{u\in \btt_n}}$ be a $\smash{\bbR^d}$-valued BRW indexed by a forest of finite trees $\btt_n$. We denote by $\smash{C^{_{(n)}}_{\cdot}}$ and $\smash{W^{_{(n)}}_\cdot}$ respectively the contour process of $\btt_n$ and the snake associated with $\bS_n$, which are normalized as in (\ref{renormCW}).

We assume that $C^{_{(n)}}_{\cdot}\! \to \! C_\cdot$ in law in $\bC^{_0}_{^1} $ and that conditionally given $\btt_n$, the jumps of $\bS_n$ are i.i.d.~r.v.s whose 
deterministic law $\bgam_n$ satisfies 
\begin{equation}
\label{bgamdeux}
\int_{\bbR^d}\!\! y \bgam_n(dy)\eqo 0 \, , \; \int_{\bbR^d}\!\! |y|^2 \bgam_n(dy) \leko \infty  \; \textrm{and} \; 
\beta_n \! : =\!  \Big( \int_{\bbR^d}\!\! y_iy_j \bgam_n(dy) \Big)_{\! \! 1\leq i,j \leq d} \!\!\! \underset{\, n\to\infty}{-\!\!\!-\!\!\!-\!\!\!\longrightarrow}\,  \beta. 
\end{equation}
Then, for all $p\ino \bbN$ and all nonnegative real numbers $s_p \geqo \ldots \geqo s_1 \geqo 0$, the following convergence holds weakly in $\bC^{_0}_{^1} \! \times \!(\bC^{_0}_{^d})^{p}$: 
\begin{equation}
\label{cvsnake11bis}
\big( C^{_{(n)}}_\cdot;    W^{_{(n)}}_{\! s_1}, \ldots , W^{_{(n)}}_{\! s_p}\big)  \xrightarrow[n\to \infty]{\;}\big(C_\cdot \, ; \, \sqrt{ \beta} .W_{\! s_1}, \ldots,\sqrt{\beta} .W_{\! s_p} \big)  .
\end{equation}
\end{lemma}
\noi
\textbf{Proof.} We set $s_0\eqo 0$ and for all $n\ino \bbN$ and all $j\ino \{ 0, \ldots, p\}$, we denote by $(U^{_{\circ , j }}_{ n,r})_{r\in \bbR_+}$  independent copies of a continuous affine interpolation of a RW that starts at $0$ and 
whose i.i.d.~jumps are distributed according to $\bgam_{\! n}$. We also assume the RWs $U^{_{\circ , j }}_{n,\cdot}$ to be independent from $C^{_{(n)}}_\cdot$. 
We next set $U^{_{(j)}}_{n,r}\eqo \lambda_n^{_{-1/2}} 
U^{_{\circ, j}}_{^{\! n,\lambda_n r}}$ and we set $s_{n,j} \eqo \lfloor b_n s_j \rfloor/ b_n$. We recall from (\ref{defFpmarg}) the definition of the function $\mathtt{Marg}_p$. Then (\ref{margdissna}) implies that 
\begin{equation}\label{snjtimes}
 \big( C^{_{(n)}}_\cdot; 
   \big( W^{_{(n)}}_{{\! s_{n,j}}} \big)_{0\leq j\leq p}  
  \big)  \overset{\textrm{(law)}}{=}  \Big( C^{_{(n)}}_\cdot; 
   \mathtt{Marg}_p \big( C^{_{(n)}}_{\cdot}   ; 
(U^{_{(j)}}_{n, \, \cdot}) _{0\leq j\leq p}  ,  (s_{n,j} )_{1\leq j\leq p} \big)  \Big). 
\end{equation}
By (\ref{bgamdeux}), Donsker's invariance principle applies and we thus get the joint convergence 
$$   \big( C^{_{(n)}}_\cdot;  (U^{_{(j)}}_{n, \, \cdot}) _{0\leq j\leq p} \big) \xrightarrow[n\to \infty]{\; } \big( C_\cdot; (\sqrt{\beta}.B^{_{(j)}}_{ \cdot}) _{0\leq j\leq p} \big)$$
weakly in $\bC^{_0}_{^1} \! \times \! (\bC^{_0}_{^d})^{p+1}$. Here, the 
$B^{_{(j)}}_{ \cdot}$ are independent standard $\bbR^d$-valued Brownian motions starting at the origin, which are furthermore independent of $C$. As a consequence of the continuity of $\mathtt{Marg}_p$ stated in (\ref{contoplus}) and of (\ref{margBrosna}) in Definition \ref{Brosnadef} of finite dimensional marginal laws of Brownian snakes, (\ref{snjtimes}) converges to the right hand side of (\ref{cvsnake11bis}). To complete the proof, we observe that for all $j\ino \{ 1, \ldots, p\}$, there exists a $\bgam_n$-distributed r.v.~$Z_{n,j}$ such that 
$\max_{r\in \bbR_+} |W^{_{(n)}}_{\! s_{n,j}} (r) \! -\! W^{_{(n)}}_{\! s_j}\! (r) |\leqo \lambda_n^{_{-1/2}}Z_{n,j} \! \to \! 0$ 
in probability as $n\! \to \! \infty$, since $\int_{\bbR^d}\! |y|^2 \bgam_n(dy) \! \to \! \sum_{1\leq i\leq d} \beta_{i,i}$. \cqfd 

\smallskip

\begin{lemma}
\label{fdcvBrosna} Let $C_\cdot$ be a continous nonnegative process such that $C_0\eqo 0$. We denote by $W$ the (not necessarily continuous) $\bbR^d$-valued Brownian snake with lifetime process 
$C_\cdot$ as in Definition \ref{Brosnadef}. Let $\mathscr C^{_{(n)}}_{\cdot}$ be a sequence of $\bbR_+$-valued continuous processes such that $\mathscr C^{_{(n)}}_{0}\eqo 0$. Conditionally given $\mathscr C^{_{(n)}}_{\cdot}\! $, let $\mathscr W^{_{(n)}}_\cdot (\cdot)$ be distributed as the $d$-dimensional Brownian snake with lifetime process 
$\mathscr C^{_{(n)}}_{\cdot}$ and initial value the null function. 
We assume that $\mathscr C^{_{(n)}}_{\cdot}\! \to \! C_\cdot$ in law in $\bC^{_0}_{^1} $. 
Then, for all $p\ino \bbN$ and all nonnegative real numbers $s_p \geqo \ldots \geqo s_1\geqo 0$, the following convergence holds weakly on $\bC^{_0}_{^1} \! \times \!(\bC^{_0}_{^d})^{p}$: 
\begin{equation}
\label{cvbrosnake}
\big(\mathscr C^{_{(n)}}_\cdot;   \mathscr W^{_{(n)}}_{\! s_1}, \ldots ,  \mathscr W^{_{(n)}}_{\! s_p}\big)  \xrightarrow[n\to \infty]{\;}\big(C_\cdot \, ; \,  W_{\! s_1}, \ldots,W_{\! s_p} \big)  
\end{equation}
\end{lemma}
\noi
\textbf{Proof.} We set $s_0\eqo 0$. Then Definition \ref{Brosnadef} implies that 
$$ \big( \mathscr C^{_{(n)}}_\cdot;  \mathscr W^{_{(n)}}_{\! s_0},  \mathscr W^{_{(n)}}_{\! s_1}, \ldots , \mathscr W^{_{(n)}}_{\! s_p}\big)  \overset{\textrm{(law)}}{=}  \Big(\mathscr C^{_{(n)}}_\cdot ; \mathtt{Marg}_p \big( \mathscr C^{_{(n)}}_\cdot \, ; (B^{_{(j)}}_{\cdot}) _{0\leq j\leq p}  ;  (s_j)_{1\leq j\leq p} \big) \Big), $$
and (\ref{cvbrosnake}) is a consequence of the continuity of $\mathtt{Marg}_p$ stated in (\ref{contoplus}). \cqfd

\vspace{-3mm}

\paragraph*{Proof of Theorem \ref{Sheuexplain} in $\textbf{Case (0)}$.}  

In this case, for all $n\ino \bbN$ we set $\smash{\btt_n\eqo \tau^{_\infty}_n\eqo (\tau_n (p))_{p\in \bbN^*}}$, 
which is an infinite GW($\mu_n$)-forest and 
$\smash{\mathcal T_n\eqo \big( \mathcal T_n(p)\eqo (\tau_n (p), (\ell_{n,u} (p))_{u\in \tau_n (p)} ) \big)_{\! p\in \bbN^*}}$, 
which is an infinite GW($\mu_n, 1$)-forest, i.e., conditionally given $\smash{\tau^{_\infty}_n}$, the $\smash{\ell_{n,u} (p)}$ 
are independent 
$\mathtt{expo} (1)$ r.v.s. We recall here that $\smash{C^{_{(n)}}_{\cdot}\! }$, $\smash{H^{_{(n)}}_{\cdot}\! }$ and $\smash{V^{_{(n)}}_{\cdot}\! }$ are resp.~the contour process, the 
height process and the \L{}ukasiewicz path of $\smash{\tau^{_\infty}_n}$, which are normalized as in (\ref{rerescale}) 
(or (\ref{renormCW})). 
$$\emph{\textrm{We overall assume }} \textrm{\texttt{\L{}uka}$_{^{\,}}$($\baa, \bbb, \bmu, \psi$)} \emph{\textrm{ and }} 
\textrm{\texttt{Grey}$_{^{\,}}$($\baa, \bbb, \bmu$).} $$
 By Lemma \ref{Lukagrowth} $(iii)$, $\lim_{n\to \infty}\lambda_n\eqo \infty$, so \texttt{Norm}$_{^{\,}}$($\baa, \bbb$) holds true too. Moreover, \texttt{Hght}$_{^{\,}}$($\baa, \bbb, \bmu$) holds true, by Lemma \ref{controlgrey} $(i)$, as well as the convergences (\ref{basiccvtreee}) and (\ref{szhghtGW}). 

We denote by $\smash{\mathscr C^{_{(n)}}_{\cdot}}$ the contour process  $\smash{(\mathscr C_{s} (\mathcal T_n))_{s\in \bbR_+}}$ of $\mathcal T_n$, which is normalized as in (\ref{rererescale}). Under the previous assumption, Theorem \ref{Greyexplain} asserts that 
\begin{equation} 
\label{CCproba0}
\textrm{ $(\mathscr C^{_{(n)}}_{s}\!\! - C^{_{(n)}}_{s})_{s\in \bbR_+}$ tends to the null function in probability in $\bC^{_0}_{^1}$.}
\end{equation}

We denote by $(\mathscr W_s(\mathcal T_n, r))_{s, r\in \bbR_+}$ the contour-Brownian snake of $\mathcal T_n$, i.e., the $1$-dimensional Brownian snake with lifetime process $\mathscr C_{\cdot} (\mathcal T_n)$ as introduced in Definition \ref{1brosnadistree}. 
We recall that $\bS_n$ is a BRW indexed by the forest $\tau^{_\infty}_n$ such that, conditionally given $\tau^{_\infty}_n$, the jumps of $\bS_n$ are $\bbR$-valued i.i.d.~r.v.s with law $\smash{\bgam (dx)\eqo (2\beta)^{-1/2}e^{-(|x|\sqrt{2})/\sqrt{\beta}} dx}$. We first observe that \emph{we only need to prove the $\beta\eqo 1$ case}, for the general case follows by multiplying the $\beta\eqo 1$ snakes by $\sqrt{\beta}$. W.l.o.g., we can also assume that {\emph{the BRW $\mathbf S_n$ is coupled with the contour-Brownian snake $(\mathscr W_s(\mathcal T_n, r))_{s, r\in \bbR_+}$ as in Lemma \ref{compareBrowsnakeBRW}}, which therefore applies (conditionally given $\tau^{_\infty}_n$).  
To simplify, we also set
\begin{equation}
\label{rerererescale}
\forall r,s\ino \bbR_+, \quad W_{s}^{_{(n)}} (r) = \tfrac{1}{\sqrt{\lambda_n}}W_{b_ns} (\bS_n , \lambda_n r) \quad \textrm{and} \quad \mathscr W_{s}^{_{(n)}} (r)=  \tfrac{1}{\sqrt{\lambda_n}}\mathscr W_{b_ns} (\mathcal T_n, \lambda_n r).
\end{equation}
Note that $\smash{(\mathscr W_{\cdot}^{_{(n)}}\! (r))_{r,s\in \bbR_+}}$ is the $1$-dimensional Brownian snake with lifetime process 
$\smash{\mathscr C^{_{(n)}}_{\cdot}}$. 
We recall that $\smash{\widehat{W}_{s}^{_{(n)}}\! \eqo W_{s}^{_{(n)}} (C^{_{(n)}}_{s})}$ and 
$\smash{\widehat{\mathscr W}_{s}^{_{(n)}}\! \eqo  \mathscr W_{s}^{_{(n)}} (\mathscr C^{_{(n)}}_{s})}$ are the endpoint processes 
of $\smash{W_{\cdot}^{_{(n)}}\! }$ and $\smash{\mathscr W_{\cdot}^{_{(n)}}\! }$. 
We also denote by $W$ the $1$-dimensional Brownian snake with lifetime $C$ which is not assumed, a priori, to be continuous. Let $x\ino \bbR_+^*$. We recall that 
$$\varsigma_{n,x} = \frac{1}{b_n} \!\!\! \sum_{\; 1\leq p\leq \lfloor a_n x\rfloor} \!\!\! \!\!\! \#\tau_n(p)\, =\inf \big\{ s\ino \bbR_+ : V^{_{(n)}}_s\eqo -\lfloor a_n x\rfloor /a_n \big\}\; .$$
We also recall from the proof of Theorem \ref{Greyexplain} that  
$ \mathbf r_{n,x}$ stands for $\frac{1}{b_n}$ times the amount of time needed to complete the contour exploration of 
$\smash{\mathcal T^x_n\! :=\!  (\mathcal T_n(p))_{1\leq p\leq \lfloor a_n x\rfloor}}$. Namely 
\begin{equation}
\label{recallrnx} \mathbf r_{n,x}\eqo \frac{1}{b_n}\sum_{1\leq p\leq \lfloor a_n x\rfloor} \sum_{u\in \tau_n(p)} 2\ell_{n,u} (p)\; .
\end{equation} 

The proof proceeds in several steps. 

\smallskip

\noi
$\textbf{Step (1)}$. \emph{If the laws of the endpoint processes $\smash{(\widehat{\mathscr W}^{_{(n)}}_\cdot)_{n\in \bbN}}$ are tight in $\bC^{_0}_{^1}$, then \emph{\texttt{Sheu} $(\psi)$} holds true, which allows to define a continuous version $W$ of the $1$-dimensional Brownian snake with lifetime process $C$. Moreover the following limit }
\begin{equation}
\label{limiBrosna1}
( \mathscr C^{_{(n)}}_\cdot \! , C^{_{(n)}}_{\cdot}\! ,V^{_{(n)}}_{\cdot} , \mathscr W^{_{(n)}}_{\cdot} , \varsigma_{n,x}, \mathbf r_{n,x} )\xrightarrow[n\to \infty]{\;}( C_\cdot , C_\cdot, Y_\cdot, W_\cdot , \varsigma_{-x}, 2 \varsigma_{-x})
\end{equation}
\emph{holds weakly in $(\bC^{_0}_{^1})^2 \! \times \! \bD (\bbR_+, \bbR) \! \times \! \mathbf C (\bbR_+, \bC^{_0}_{^1})\! \times \! \bbR_+^2$, for all $x\ino \bbR_+^*$. }

\noi
\emph{Proof.} By Proposition \ref{endsnake1} $(i)$, the laws of the whole snakes 
$(\mathscr W^{_{(n)}}_{\cdot})_{n\in \bbN}$, are tight in $\bC(\bbR_+, \bC^{_0}_{^1})$. 
Thus the laws of $\smash{( \mathscr C^{_{(n)}}_\cdot \! , C^{_{(n)}}_{\cdot}\! ,V^{_{(n)}}_{\cdot} , 
\mathscr W^{_{(n)}}_{\cdot}  )_{n\in \bbN}}$ are tight in $(\bC^{_0}_{^1})^2 \! \times \! \bD (\bbR_+, \bbR) 
\! \times\!  \mathbf C (\bbR_+, \bC^{_0}_{^1})$. Let us consider a limiting process along a subsequence 
$(n_k)_{k\in \bbN}$: w.l.o.g.~it can be written $(C_\cdot, C_\cdot, Y_\cdot, W')$, by 
Theorem \ref{Greyexplain}. Namely, 
$\smash{\lim_{k\to \infty} (\mathscr C^{_{(n_k)}}_\cdot \! , C^{_{(n_k)}}_{\cdot}\! ,V^{_{(n_k)}}_{\cdot} , \mathscr W^{_{(n_k)}}_{\cdot})}$ 
$\eqo $  $(C_\cdot, C_\cdot, Y_\cdot, W')$ weakly in the appropriate space. 
Then Lemma \ref{fdcvBrosna} implies that $\smash{W'}$ and 
$W$ have the same finite dimensional marginal laws. Since $W'$ is $\bC (\bbR_+, \bC^{_0}_{^1})$-valued, 
this shows that there is a continuous version of the $1$-dimensional Brownian snake with lifetime process 
$C$. Thus \texttt{Sheu} $(\psi)$ holds true by Proposition \ref{BrosnapsiH} $(i)$. By taking $W$ 
continuous, we easily see that 
$( C_\cdot, Y_\cdot, W')$ and $( C_\cdot, Y_\cdot, W)$ have the same law. This proves 
existence and uniqueness of weak limits of the laws of the processes 
$( \mathscr C^{_{(n)}}_\cdot \! , C^{_{(n)}}_{\cdot}\! ,V^{_{(n)}}_{\cdot} , 
\mathscr W^{_{(n)}}_{\cdot}  )$, $n\ino \bbN$. Namely, it proves that 
\begin{equation}
\label{limiBrosna1bbis}
( \mathscr C^{_{(n)}}_\cdot \! , C^{_{(n)}}_{\cdot}\! ,V^{_{(n)}}_{\cdot} , \mathscr W^{_{(n)}}_{\cdot} ) \xrightarrow[n\to \infty]{\;} ( C_\cdot , C_\cdot, Y_\cdot, W_\cdot )
\end{equation}
weakly in $(\bC^{_0}_{^1})^2 \! \times \! \bD (\bbR_+, \bbR) \! \times \! \mathbf C (\bbR_+, \bC^{_0}_{^1})$. 
The same arguments (see Jacod \& Shiryaev \cite{JaSh02} Proposition 2.11, Chapter VI, Section 2a p.~341) 
used to prove that the weak limit 
$\lim_{n\to \infty} \varsigma_{n,x} \eqo \varsigma_{-x}$ holds jointly with (\ref{basiccvtreee}) can be used to prove that it holds jointly with (\ref{limiBrosna1bbis}). Next we recall from $\textbf{Step (A)}$ in the proof of Theorem \ref{Greyexplain} that 
$\mathbf r_{n,x}\eqo 2L_n (\varsigma_{n,x})$ where the process $L_n (\cdot)$  tends to the identity in probability in 
$\bC^{_0}_{^1}$. Combined with the previous arguments, it implies (\ref{limiBrosna1}). \cq

\smallskip

Quite similar arguments relying on Lemma \ref{fdcvsna} instead of Lemma \ref{fdcvBrosna} entail the following.

\smallskip

\noi
$\textbf{Step (2)}$. \emph{If the laws of the endpoint processes $\smash{(\widehat{W}^{_{(n)}}_s)_{s\in \bbR_+}}$, $n\ino \bbN$, are tight in $\bC^{_0}_{^1}$, then \emph{\texttt{Sheu} $(\psi)$} holds true, which allows to define a continuous version $W$ of the $1$-dimensional Brownian snake with lifetime process $C$. Moreover the following limit}
\begin{equation}
\label{limiBrosna2}
(C^{_{(n)}}_{\cdot}\! ,V^{_{(n)}}_{\cdot}\!\! ,W^{_{(n)}}_{\cdot}\!  ,\varsigma_{n,x} ) \xrightarrow[n\to \infty]{\;}
 ( C_\cdot, Y_\cdot, W_\cdot  , \varsigma_{-x})
\end{equation}
\emph{holds weakly in $\bC^{_0}_{^1} \! \times \! \bD (\bbR_+, \bbR) \! \times \! \mathbf C (\bbR_+, \bC^{_0}_{^1})\! \times \! \bbR_+$. }

\smallskip

\noi
$\textbf{Step (3)}$. \emph{We prove that \emph{\texttt{Sheu}$_{\,}(\mathbf a, \mathbf b, \bmu)$} implies for all $x\ino \bbR_+^*$ that}
\begin{equation}
\label{limmaxdispl}
\mathscr M_{n, x}\! :=\! \max_{s\in [0, \mathbf r_{n,x}] } \widehat{\mathscr W}^{_{(n)}}_{\! s} \,\xrightarrow[n\to \infty]{\emph{\textrm{$\, $ in law on $\bbR_+\, $}}} \,\max_{s\in [0, 2\varsigma_{-x}] } \widehat{W}_{s} \; , 
\end{equation}
\emph{where $W$ is the continous version of the $1$-dimensional Brownian snake with lifetime $C$, which is well-defined  since \emph{\texttt{Sheu}$_{\,}(\mathbf a, \mathbf b, \bmu)$} implies \emph{\texttt{Sheu}$_{\,}(\psi)$} by Lemma \ref{controlpsi} $(i)$, and by Proposition \ref{BrosnapsiH} $(i)$.}

\smallskip

\noi
\emph{Proof.} Proposition \ref{BrosnapsiH} $(ii)$ asserts that $\smash{\bP (\max_{s\in [0, 2\varsigma_{-x}] } \widehat{W}_{s} \leqo z) \eqo e^{-xw(z)}}$ where $w$ stands for the inverse of $F$ as defined in (\ref{defFnF}). We next observe that 
$$\mathscr M_{n,x} = \tfrac{1}{\sqrt{\lambda_n}} \max_{1\leq p\leq \lfloor a_nx \rfloor} \max_{s\in \bbR_+} \widehat{\mathscr W}_{\! s} (\mathcal T_n(p)) \; .$$
Since the r.v.s $\smash{\max_{s\in \bbR_+} \widehat{\mathscr W}_{\! s} (\mathcal T_n(p))}$, $p\ino \bbN^*$, are i.i.d., we get $\smash{\bP (\mathscr M_{n,x} \leqo z)\eqo (1-\frac{1}{a_n} w_n (z))^{\lfloor a_n x \rfloor} }$ where we have set $\smash{w_n (z)\eqo a_n \bP (\max_{s\in \bbR_+} \widehat{\mathscr W}_{\! s} (\mathcal T_n(1)) \geko \sqrt{\lambda_n} z)}$. By (\ref{exitBroSna}) in Lemma \ref{maxBrosnadiscr}, the function $w_n$ is the inverse of $F_n$ as defined in (\ref{defFnF}) and Lemma \ref{controlpsi} $(iii)$ implies that $\lim_{n\infty} w_n\eqo w$ pointwise on $\bbR_+^*$. This entails $\smash{\lim_{n\to \infty}\bP (\mathscr M_{n,x} \leqo z)\eqo \bP (\max_{s\in [0, 2\varsigma_{-x}] } \widehat{W}_{s} \leqo z)}$ for all $\smash{z\ino \bbR_+^*}$, which completes the proof of (\ref{limmaxdispl}). \cq 

\smallskip

\noi
$\textbf{Step (4)}$. \emph{We next assume (\ref{limiBrosna1}) for a fixed $x\ino \bbR_+^*$,  
and we prove that it implies \emph{\texttt{Sheu}$_{\,}(\mathbf a, \mathbf b, \bmu)$}. } 

\smallskip

\noi
\emph{Proof.} Let $x\ino \bbR_+^*$. Standard arguments show that (\ref{limiBrosna1}) implies (\ref{limmaxdispl}). 
Then, the arguments used in the proof of $\textbf{Step (3)}$ easily imply that $\lim_{n\to \infty} w_n \eqo w$ pointwise on $\bbR_+^*$, where $w_n$ (resp.~$w$) is the inverse of $F_n$ (resp.~$F$) as defined in (\ref{defFnF}). 
Therefore, for all sufficiently large $n$, we get $w_n (z) \leko 2w(z)$ and thus $F_n (2w(z)) \leqo F_n (w_n(z))\eqo z$. Consequently, 
$\lim_{z\to 0^+}\limsup_{n\to \infty} F_n (2w(z))\eqo 0$, which implies \texttt{Sheu}$_{\,}(\mathbf a, \mathbf b, \bmu)$ since $\lim_{z\to 0^+} w(z) \eqo \infty$. \cq 

\smallskip

\smallskip

\noi
$\textbf{Step (5)}.$ \emph{We suppose} \texttt{Sheu} ($\psi$) \emph{and either (\ref{limmaxdispl}) for all $x\ino \bbR^*_+$ or}
\begin{equation}
\label{limmaxdispl2}
M_{n, x} :=\!\!\!\! \!\!\!  \max_{\quad s\in [0,  2\varsigma_{n,x}] } \!\!\!\!\!\! \! \widehat{ W}^{_{(n)}}_{\! s} \,\xrightarrow[n\to \infty]{\emph{\textrm{$\, $ in law on $\bbR_+\, $}}} \, \!\!\!\!\!\! \!  \max_{\quad s\in [0, 2\varsigma_{-x}] } \!\!\!\!\!\! \! \widehat{W}_{s} \; ,
\end{equation}
\emph{for all $x\ino \bbR^*_+$}. Then $\lim_{n\to \infty} b_n / (a_n (\log a_n)^2)\eqo \infty$.}

\smallskip

\noi
\emph{Proof.} Let us first assume \texttt{Sheu} ($\psi$) and (\ref{limmaxdispl2}) for all sufficiently small $x\ino \bbR^*_+$. 
We fix $z \ino \bbR_+^*$ (small). 
By \texttt{Sheu} ($\psi$), (\ref{exitsys}) in Proposition \ref{BrosnapsiH} applies, and there exists $x\ino \bbR_+^*$, which depends on $z$, such that 
\begin{equation}
\label{concontra}
  \bP \big(\!\! \!\!\!\! \max_{\quad s \in [0, 2\varsigma_{-x }]} \!\!\!\!\!\! \widehat{W}_s \geq z \big) \eqo 1\! -\! e^{-xw(z)} \leko \tfrac{1}{2} \; .
\end{equation}
%We fix $z\ino (0, z_0]$. 
By Lemma \ref{compareBrowsnakeBRW} $(iii)$, conditionally given $\tau^{_\infty}_n$, we almost surely get 
$\smash{2 \bP \big( M_{n, x}\geko  z\sqrt{\lambda_n} \, \big| \, \tau^{_\infty}_n \big)} $ $\geqo$ $\smash{ 1\! -\! \big( 1\! -\! e^{-2z\sqrt{2\lambda_n}} \big)^{b_n\varsigma_{n,x} -\lfloor a_nx\rfloor}  } $. Then, for any $\delta \ino \bbR_+^*$, we get  
\begin{eqnarray}
2 \bP \big( M_{n, x}\geko  z\sqrt{\lambda_n} \big) \!\!\!\! &\geq & \!\!\!\! \Big( 1 \! -\!  \exp
\Big( b_n \delta \log \big( 1\! -\! e^{-2z\sqrt{2\lambda_n} }\big) \Big) \Big) \bP \big( \varsigma_{n,x}-\tfrac{x}{\lambda_n}  \geko \delta \big) \nonumber \\
 \!\!\!\! &\geq & \!\!\!\! \Big(1\! -\! \exp\Big(\!\!   -\! \delta e^{\log b_n -2z\sqrt{2\lambda_n}} \Big) \Big)   \bP \big( \varsigma_{n,x}-\tfrac{x}{\lambda_n}  \geko \delta \big) \label{minomaxmov}
\end{eqnarray}
We first observe that $\smash{\liminf_{n\to \infty} \bP ( \varsigma_{n,x}\! -\! \tfrac{x}{\lambda_n}  \geko \delta )
 \geqo \bP (\varsigma_{-x} \geko \delta )}$ since $\varsigma_{n,x}$ converges in law to $\varsigma_{-x}$ by (\ref{szhghtGW}). 
Suppose that $\smash{\limsup_{n\to \infty} 
\log b_{n} \! -\! 2z\sqrt{2\lambda_{n}}\eqo \infty}$. Then (\ref{limmaxdispl2}) and (\ref{minomaxmov}) 
entail 
$\smash{\tfrac{1}{2}\bP (\varsigma_{-x} \geqo \delta) }$ $\leqo$ $\smash{ \bP (\max_{s \in [0, 2\varsigma_{-x }]} \widehat{W}_s \geqo z) }$. 
Since it holds for arbitrary small $\delta$ and since $\smash{\lim_{\delta\to 0^+} \bP (\varsigma_{-x} \geqo \delta) } $ $\eqo$  
$\smash{\bP (\varsigma_{-x} \geko 0)\eqo 1}$, we get $\smash{ \bP (\max_{s \in [0, 2\varsigma_{-x }]} \widehat{W}_s \geqo z) \geqo \frac{1}{2}}$, which contradicts (\ref{concontra}). 
Thus, $\smash{\limsup_{n}}$ $\smash{ \log b_n \! -\! 2z\sqrt{2\lambda_n}} $ $ \leko$ $\infty$. Consequently, $\smash{\limsup_{n\to \infty} \lambda_{n}^{_{-1/2}} \log b_n \leqo 2z\sqrt{2}}$. 
Since it holds for arbitrary small $z\ino \bbR_+^*$, it implies $\smash{\lim_{n\to \infty} \lambda_{n}^{_{-1/2}} \log b_n\eqo 0}$, which implies  the desired result  by an elementary argument. 

Under the assumptions \texttt{Sheu} ($\psi$) and (\ref{limmaxdispl}) for all $x\ino \bbR_+^*$, we argue in a quite similar way: we use Lemma \ref{compareBrowsnakeBRW} $(iv)$ instead of $(iii)$. We leave the details to the reader. \cq

\smallskip

We denote by $\smash{(v_k)_{k\in \bbN}}$ the contour exploration of the tree corresponding to the GW forest 
$\smash{(\tau_n(p))_{p\in \bbN^*}}$ and we recall from 
Definitions \ref{Contlifespan} and \ref{lifespanforest}, $\smash{\Lambda_{\mathcal T_n}\! :\! \bbR_+\! \to \! \bbR_+}$, which is continuous, 
increasing and such that $\smash{\Lambda_{\mathcal T_n} (0)\eqo 0}$ and $\smash{\Lambda_{\mathcal T_n}(s+k\! -\! 1)\eqo \Lambda_{\mathcal T_n}(k\! -\! 1)+ 
s|\zeta_{v_k} \!\! -\! \zeta_{v_{k-1}}|}$, $s\ino [0, 1]$, $k\ino \bbN^*$. We recall from the proof of Theorem \ref{Greyexplain} that we have set 
$$\Lambda_n (s)\! :=\! \tfrac{1}{b_n} \Lambda_{\mathcal T_n} (b_ns) , \quad s\ino \bbR_+.$$ 

\noi
$\textbf{Step (6)}$. \emph{We assume 
$\lim_{n\to \infty} b_n/(a_n (\log a_n)^2)\eqo \infty$. }
\begin{equation}
\label{scrWWcontig}
\textrm{\emph{Then $( \widehat{\mathscr W}^{_{(n)}}_{\! \Lambda_n (s)} \!\!  -\! \widehat{W}^{_{(n)}}_s)_{s\in \bbR_+}$ converges to the null function in probability in $\bC^{_0}_{^1}$. }}
\end{equation}
\noi
\emph{Proof.} Let $\smash{x\ino \bbR_+^*}$. To simplify, we set $\smash{D_{n,x}\! :=\!  \max_{ s\in [0, 2\varsigma_{n,x}] }\big| \widehat{\mathscr W}^{_{(n)}}_{\! \Lambda_n (s)} \!\!  -\! \widehat{W}^{_{(n)}}_s\big| }$. 
By (\ref{closeWWexpli}) in Lemma \ref{compareBrowsnakeBRW} $(ii)$, 
conditionally given $\smash{\tau^{_\infty}_n}$, for all $\smash{\epp\ino \bbR_+^*}$, we a.s.~get 
$$ \bP \big(D_{n,x}\geko  3\epp \sqrt{2}  \, \big| \, \tau^{_\infty}_n  \big) \leq 1\! -\! \big(1\! -\! 2e^{-\epp \sqrt{\lambda_n}} \big)^{b_n \varsigma_{n,x}}= 1\! -\! e^{-\varsigma_{n,x} f_n(\epp) }   , $$
where we have set $\smash{f_n(\varepsilon)\eqo -b_n \log (1\! -\! 2e^{-\epp \sqrt{\lambda_n}} )}$. Consequently, for all $z\ino \bbR_+^*$, 
$\smash{\bP (D_{n,x}\geko  3\epp \sqrt{2} )}$ $\leqo$ $\smash{ 1\! -\! e^{-zf_n(\epp)}}$ $ +$ $\smash{ \bP (\varsigma_{n,x} \geqo z)}$. Our assumption implies that 
$\lim_{n\to \infty} f_n(\epp)\eqo 0$. Since  $\varsigma_{n,x}$ converges in law to $\varsigma_{-x}$ by (\ref{szhghtGW}), we get 
$\smash{\limsup_{n\to \infty} \bP (D_{n,x}\geko  3\epp \sqrt{2} )\leqo \bP (\varsigma_{-x} \geqo z)\to 0}$ as $z\! \to \! \infty$. This proves that $\lim_{n\to \infty}D_{n,x}\eqo 0$ in probability, for all $x\ino \bbR_+^*$. Then, for all $s, \epp \ino \bbR_+^*$, 
$$ \bP \Big(\max_{r\in [0, s]} \big| \widehat{\mathscr W}^{_{(n)}}_{\! \Lambda_n(r)}\! \! -\! \widehat{W}^{_{(n)}}_r   \big|\geko  \epp \Big) \leq  \bP \big( D_{n,x} \geko \epp \Big) + \bP \big(  2\varsigma_{n,x} \leqo s \big) $$
Thus $\smash{\limsup_{n\to \infty}  \bP (\max_{r\in [0, s]} | \widehat{\mathscr W}^{{(n)}} (\Lambda_n(r)) \! -\! \widehat{W}^{_{(n)}}_r  |\geko  \epp ) \leqo \bP (2\varsigma_{-x} \leqo s)}$, which tends to $0$ as $x\! \to \! \infty$ since a.s.~$\lim_{x\to \infty} \varsigma_{-x}\eqo \infty$. This entails (\ref{scrWWcontig}).\cq

\smallskip

Here is the technical step of the proof of Theorem \ref{Sheuexplain}.

\smallskip 

\noi
$\textbf{Step (7)}.$ \emph{We assume 
\emph{\texttt{Sheu}$_{^{\,}}$($\baa, \bbb, \bmu$)}. Then the laws of $\smash{(\widehat{\mathscr W}^{_{(n)}}_s)_{s\in \bbR_+}}$ are 
tight.}

\smallskip

\noi
\emph{Proof.} By Lemma \ref{controlpsi} $(i)$, \texttt{Sheu} ($\psi$) holds true and by Proposition \ref{BrosnapsiH} $(i)$, $W_\cdot$ can be taken continuous. Since $\smash{\widehat{\mathscr W}^{_{(n)}}_0\eqo 0}$, standard results (see e.g.~Billingsley \cite{Bil68} Thm 7.3 p.~82), imply that we only need to prove  
\begin{equation}
\label{scrWendtight}
\forall N\ino \bbN^*, \; \forall \epp\ino  (0, 1), \qquad \lim_{p\to \infty} \limsup_{n\to \infty} \bP \big( \omega_{N , p} \big( \widehat{\mathscr W}^{_{(n)}}_{\cdot} \big) \geqo \epp\big) = 0,
\end{equation}
where for all $N,p\ino \bbN^*$ and for all $f\ino \bC^{_0}_{^1}$ we have set 
$$ \omega_{N ,p} (f) \eqo \max_{0\leq j< N2^p} \max \big\{\,  | f(s)\! -\! f(j2^{-p}) |\, ; \, s\ino [j2^{-p} , (j+1)2^{-p}] \big\}  \; .$$ 
To that end, we first prove that 
\begin{equation}
\label{scrtheWmgtight}
\forall N\ino \bbN^*, \; \forall \epp\ino  (0, 1), \qquad \lim_{p\to \infty} \limsup_{n\to \infty} \bP \big( \vartheta_{N,p} \big( \mathscr W^{_{(n)}}_{\cdot} \big) \geqo \epp\big) = 0,
\end{equation}
where for all $N,p\ino \bbN^*$ and for all  snakes $(h, w)$ we have set 
$$ \vartheta_{N\! ,\, p} (w) \eqo \!\! \max_{0\leq j< N2^p} \! \! \max \big\{\,  \big| \widehat{w}_{j2^{-p}} - w_{j2^{-p}} \big( h(j2^{-p})  - r\big) \big| \, ; \, r\ino \big[ 0\, , h(j2^{-p}) - m_h (j2^{-p} ,(j+1)2^{-p}) \big] \big\}  \; .$$ 
We recall here that  $m_h(r_1,r_2)\eqo \min_{r\in [r_1, r_2]} h(r)$, for all real numbers $r_2\geqo r_1\geqo 0$. 

\smallskip

\noi
\emph{Proof of (\ref{scrtheWmgtight})}. By Lemma \ref{fdcvBrosna} and basic results on weak convergence in $\bC_{^0}^{_1}$, we get $\vartheta_{N\! ,\, p} ( \mathscr W^{_{(n)}}_{\cdot} )\! \to \! 
 \vartheta_{N\! ,\, p} ( W_{\cdot} )$ in law in $\bbR_+$. Since $(C,W)$ is a continuous snake, we also get 
$\lim_{p\rightarrow \infty}  \vartheta_{N\! ,\, p} ( W_{\cdot} )\eqo 0$ a.s.~for all $N\ino \bbN^*$. This easily entails  (\ref{scrtheWmgtight}).   \cq

\smallskip

Let $z\ino \bbR_+^*$. Recall from (\ref{osci1def}) the definition of the oscillation times $(\bsigma_q (w,z))_{q\in \bbN}$ 
of the snake $(h,w)$. Namely, they are recursively defined as follows: $\bsigma_{\! 0} (w,z)\eqo 0$ and 
$$\bsigma_{\! q+1}(w,z) \eqo \inf\!  \big\{ s\ino (\bsigma_{\! q} (w,z), \infty) :  \big| \widehat{w}_s\! -\! 
w_s(m_h (\bsigma_q(w,z), s)  \big|  \geko  z\big\} \; .$$ 
To simplify we set $\smash{\bsigma^{_{(n)}}_{\! q} (z) \eqo \bsigma_{\! q} (\mathscr W^{_{(n)}}_\cdot \!\!  , z)}$ 
and we define the events $\smash{A_{n}\eqo \big\{ \omega_{N,p} \big(\widehat{\mathscr W}^{_{(n)}}_\cdot \big)
 \geko 7\epp \big\}}$ and 
$\smash{B_{n}  \eqo \big\{  \vartheta_{N,p} \big( \mathscr W^{_{(n)}}_\cdot \big) \leqo \epp\big\}}$. 
Observe that 
\begin{eqnarray*}
A_n\cap B_n =\Big\{ \, \forall j\ino \{ 0, \ldots,  N2^{p} \! -\! 1\} \!\! \!\!\!\! & , &\!\!\!\!  
\forall r \ino [0, \mathscr C^{_{(n)}}_{j2^{-p}} \!- \! m_{\mathscr C^{_{(n)}} } (j2^{-p}, (j+1)2^{-p} )]   \, : \, \\
\big| \widehat{\mathscr W}^{_{(n)}}_{j2^{-p}} \!\!\!\! &- & \!\!\!\!  \mathscr W^{_{(n)}}_{j2^{-p}} 
 \big(  \mathscr C^{_{(n)}}_{j2^{-p}} \!  -\! r\big) \big|  \leqo \epp  \quad  \textrm{and} \quad  \\
\exists i\ino \{ 0, \ldots,  N2^{p} \! -\! 1\}  \!\!\!\!\!\! & , &\!\!\!\!  \exists s\ino [i2^{-p}, (i+1)2^{-p}]  \, : \, 
\big|  \widehat{\mathscr W}^{_{(n)}}_{\! s} \! \! -\!  \widehat{\mathscr W}^{_{(n)}}_{\! i2^{-p}} \big|\geko 
 7 \varepsilon  \Big\}. 
\end{eqnarray*}
By Lemma \ref{osc1lemma} with $(h,w)\eqo (\mathscr C^{_{(n)}}_\cdot\! , \mathscr W^{_{(n)}}_\cdot )$, 
$s_0\eqo i2^{-p} $, $s_1\eqo (i+1)2^{-p}$ and $z\eqo \epp$, we get 
\begin{eqnarray}
A_n \cap B_n  \!\!\! & \subset & \!\!\! \big\{ \exists q\ino \bbN^*, \, \exists i\ino \{ 0, \ldots,  N2^{p} \! -\! 1\} \, : \, i2^{-p} \leqo 
\bsigma_{\! q}^{_{(n)}} (\epp) \leko \bsigma_{\! q+1}^{_{(n)}} (\epp) \leqo (i+1)2^{-p} \big\} \nonumber \\
\!\!\! & \subset & \!\!\! \Big\{\!  
\min \big\{\bsigma^{_{(n)}}_{\! {q+1}} ( \epp) \! -\!  \bsigma^{_{(n)}}_{\! q} ( \epp) \,  ; \, q\ino \bbN^* :  
\bsigma^{_{(n)}}_{q}\!  (\epp)\leqo N  \big\}  \leq 2^{-p}  \Big\}.\label{trueuseosci}
\end{eqnarray}

We now want to replace the oscillation times of 
$\mathscr W^{_{(n)}}_{\! \cdot} $ 
by the oscillation times of the height-Brownian 
snake which enjoy the renewal property stated in Lemma \ref{renewosci}. 
To that end, we recall from (\ref{ccCvsccH}) in Remark \ref{contlifespanrem} $(c)$ 
that there is an increasing c{\`a}gl{\`a}d time-change $\mathcal K_{\mathcal T_n} : \bbR_+ \! \to \! \bbR_+$ such that 
$\mathscr C_{\mathcal K_{\mathcal T_n} (s)} (\mathcal T_n) \eqo \mathscr H_s (\mathcal T_n)$. 
We denote by $(\mathtt W_s (\mathcal T_{\! n}, r))_{s,r\in \bbR_+}$  the height-Brownian snake 
as in (\ref{defviatimechange}), Definition \ref{1brosnadistree}. Namely, $\mathtt W_s(\mathcal T_n ,r)\eqo 
\mathscr W_{\mathcal K_{\mathcal T_n} (s)} (\mathcal T_n, r)$, $r,s\ino \bbR_+$. 
We denote $(u_l)_{l\in \bbN}$ the depth-first indexed vertices of the tree corresponding to the infinite GW forest 
$(\tau_n(p))_{p\in \bbN^*}$. We recall from Definition \ref{lifespanforest} and (\ref{LTident}) in Remark \ref{contlifespanrem} 
that 
$L_{\mathcal T_n} \! : \! \bbR_+\! \to \! \bbR_+$ is the increasing continuous function which is affine on each interval 
$[l, l+1]$ and which satisfies $L_{\mathcal T_n} (l)\eqo \sum_{0\leq j\leq l} \ell_{n,u_j}$, $l\ino \bbN$ (we recall that 
$\ell_{n, \varnothing}\eqo 0$). 
We recall from (\ref{Kvs2id}) that 
$\max_{s\in [0, L_{\mathcal T_n} (l)]} |2s \! -\! \mathcal K_{\mathcal T_n} (s) | \leqo 
\max_{s\in [0, 2L_{\mathcal T_n}(l) ]} \mathscr C_s (\mathcal T_n)$ and 
to simplify notation, we set 
$$\mathcal K_n(s) \eqo \tfrac{1}{b_n}\mathcal K_{\mathcal T_n}  (b_ns), \, L_n(s) \eqo \tfrac{1}{b_n}L_{\mathcal T_n}  (b_ns), \, \mathtt W^{_{(n)}}_s(r)\eqo \tfrac{1}{\sqrt{\lambda_n }}  \mathtt W_{b_ns}(\mathcal T_n ,\lambda_n r)\; \textrm{and} \;  \mathbf s_{q}^{_{(n)}} (\epp)\eqo \bsigma_{\! q} (\mathtt W^{_{(n)}}_\cdot \!\! , \epp).$$  We recall the following results. 
\begin{compactenum}

\smallskip

\item[$(i^*)$] By (\ref{changetimeosci}), $\mathcal K_n(\mathbf s_{q}^{_{(n)}} \! (\epp))\eqo \bsigma_{\! q}^{_{(n)}} \! (\epp)$, for all $q\ino \bbN$. 

\smallskip

\item[$(ii^*)$] For all $s\ino \bbR_+$, $\max_{s'\in [0, L_n (s)] } |2s'\! -\! \mathcal K_n (s')| \leqo \frac{1}{a_n} \max_{s'\in [0, 2L_n (b_n^{-1} \lceil b_ns \rceil )]} \mathscr C^{_{(n)}}_{\! s'}$. 

\smallskip

\item[$(iii^*)$] By $\textbf{Step (A)}$ of the proof of Theorem \ref{Greyexplain}, 
$L_n$ tends to the identity map in probability in $\bC^{_0}_{^1}$.  

\smallskip

\item[$(iv^*)$] By Lemma \ref{renewosci}, The r.v.s $(\mathbf s_{q+1}^{_{(n)}} (z)\! -\! \mathbf s_{q}^{_{(n)}} \! (z))_{q\in \bbN}$ are i.i.d.

\smallskip

\item[$(v^*)$] By (\ref{capKlowupbounds}) in Remark \ref{contlifespanrem} $(d)$, $s\leqo \mathcal K_n (s) \leqo 2s$, for all $s\ino \bbR_+$. 

\smallskip

\end{compactenum}

To simplify notation, we set 
$N_{n}\eqo  b_n^{-1} \lfloor  b_n(N +2)\rfloor$. By $(i^*)$ and $(ii^*)$, 
for all $q\ino \bbN$ such that $\mathbf s_{q+1}^{_{(n)}} (\epp) \leqo L_n (N_n)$, we first get 
\begin{equation}
\label{deltasigandS}
\big| \bsigma_{\! q+1}^{_{(n)}} (\epp)\! -\! \bsigma_{\! q}^{_{(n)}}\! (\epp)-2\big(\mathbf s_{q+1}^{_{(n)}} (\epp)-\mathbf s_{q}^{_{(n)}} \! (\epp) \big)  \big|   \leq \tfrac{2}{a_n}\!\!\! \!\!\! \!  \max_{\quad s\in [0, 2L_n (N_n)]}\!\!\!\!\!\! \!\!\! \!  \mathscr C^{_{(n)}}_{{ s}}. 
\end{equation}
We then fix $z\ino \bbR_+^*$ and we introduce the event 
$$\smash{C_n \eqo \big\{ N+1 \leqo L_n(N_n) \leqo  N+3 \;\,  ; \max_{^{s\in [0,2( N+3)]}} \!\! \mathscr C^{_{(n)}}_{s} \leqo  z \big\}}$$ and the integer 
$n_0\eqo n_0(p,z)$ such that $2z/a_n \leqo 2^{-p}$ for all integers $n\geqo n_0$. 
Let $n\geqo n_0$. On $C_n$, if $\bsigma_{q}^{_{(n)}} \! (\epp)  \leqo N$, $(i^*)$ and $(v^*)$ imply that 
$\mathbf s^{_{(n)}}_q \! (\epp) \leqo N$; if furthermore $\bsigma_{q+1}^{_{(n)}} (\epp) \! -\! \bsigma_{\! q}^{_{(n)}} \! (\epp)   \leqo 2^{-p}$, then $\bsigma_{q+1}^{_{(n)}} (\epp) \leqo N+2^{-p} \leqo N+1$ and by $(i^*)$ and $(v^*)$, we get 
$\mathbf s^{_{(n)}}_{q+1} \! (\epp) \leqo N+1\leqo L_n (N_n)$:  (\ref{deltasigandS}) applies and yields 
$\mathbf s_{q+1}^{_{(n)}} (\epp) \! -\! \mathbf s_{q}^{_{(n)}} \! (\epp)   \leqo 2^{-p}$. Therefore by (\ref{trueuseosci}) 
\begin{equation}
\label{ABCsub} \forall n\geqo n_0, \quad A_n \cap B_n \cap C_n \,  \subset  \Big\{\!  
\min \big\{\mathbf s^{_{(n)}}_{\! {q+1}} ( \epp) \! -\!  \mathbf s^{_{(n)}}_{q} \! ( \epp)  ; \, q\ino \bbN :  \mathbf s^{_{(n)}}_{ q}\!  (\epp)\leqo N \big\}  \leq 2^{-p}  \Big\}.
\end{equation}
Thus, for all $n\geqo n_0$, we get 
\begin{eqnarray}
\bP \big( \omega_{N,p} \big( \widehat{\mathscr W}^{_{(n)}}_{\cdot} \big) \geko 7 \epp\big)  
\!\!\!\! & \leq & \!\!\!\! \bP \big( \vartheta_{N,p} \big( \mathscr W^{_{(n)}}_{\cdot} \big) \geko \epp\big) +
 \bP \big( \Omega\backslash C_n  \big)
   \nonumber  \\
\label{scrdecompeven1} &+ &   \bP  \Big( \min \big\{\mathbf s^{_{(n)}}_{ {q+1}} ( \epp) \! -\!  \mathbf s^{_{(n)}}_{q} \! ( \epp) 
 ; \, q\ino \bbN :  \mathbf s^{_{(n)}}_{ q}\!  (\epp)\leqo N  \big\}  \leq 2^{-p} \Big)
\end{eqnarray}
Since we assume \texttt{\L{}uka}$_{^{\,}}$($\baa, \bbb, \bmu, \psi$) and 
\texttt{Sheu}$_{^{\,}}$($\baa, \bbb, \bmu$), and thus \texttt{Grey}$_{^{\,}}$($\baa, \bbb, \bmu$), Theorem 
\ref{Greyexplain} applies and $\mathscr C^{_{(n)}}_\cdot \! \to \! C$ weakly in $\bC^{_0}_{^1}$. By $(iii^*)$ 
above and since $\lim_{n\to \infty}N_n \eqo N+2$, we get 
\begin{equation}
\label{scrcontrHetvarth1}
\lim_{z\to \infty} \limsup_{n\to \infty} \bP \big(   \Omega\backslash C_n \big) =0
\end{equation}
Then, by a general result stated in Ethier \& Kurtz \cite{EtKu86} Chapter 3, Lemma 8.2 p.~134, we get that 
\begin{eqnarray}
\lim_{p\to \infty} \limsup_{n\to \infty}   \bP  \Big( \!\!\! \!\!\! &  &  \!\!\!  \!\!\!  \min \big\{\mathbf s^{_{(n)}}_{ {q+1}} ( \epp) 
\! -\!  \mathbf s^{_{(n)}}_{q} \! ( \epp)  ; \, q\ino \bbN :  \mathbf s^{_{(n)}}_{ q}\!  (\epp)\leqo N  \big\}  \leq 2^{-p} \Big)\eqo 0 
\quad \Longleftrightarrow \;  \nonumber \\
\label{EthKurlemscr}   \lim_{p\to \infty}    \!\!\!\!  \!  & &\!\!\!\!  \!   \limsup_{n\to \infty} \, \sup_{q\in \bbN} \, 
 \bP \big(\,  \bs^{_{(n)}}_{^{q+1}} ( \epp) \! -\!  \bs^{_{(n)}}_{^q} ( \epp) \, 
 \leq 2^{-p} \, ; \, \bs^{_{(n)}}_{q}\!  (\epp)\leqo N \big) \eqo 0 .
\end{eqnarray}
We use $(iv^*)$ above, combined with (\ref{scrtheWmgtight}), (\ref{scrdecompeven1}), 
(\ref{scrcontrHetvarth1}) and (\ref{EthKurlemscr}), to see that (\ref{scrWendtight}) is implied by 
$\lim_{p}$ $ \limsup_{n}$ $\bP \big(  \bs^{_{(n)}}_{1}\!  (\epp)  $ $\leqo 2^{-p} \big) \eqo 0$, 
which is equivalent to $\lim_{p} \limsup_{n}   \bP \big(  \bsigma^{_{(n)}}_{1}\!  (\epp) \leqo 2^{-p} \big) \eqo 0$, by $(v^*)$. Namely,  
\begin{equation}
\label{itboilsdownto}
\lim_{p\to \infty} \limsup_{n\to \infty} \,   \bP \Big( \!\!\!\!\!\!\!\!\!\!\! \sup_{\qquad s\in 
[0, 2^{-p}] } \!\!\!\!\!\!\!\!\!\! \widehat{\mathscr W}^{_{(n)}}_{\! s} \geko \epp  \Big)\eqo 0\; . 
\end{equation}
Now observe for all $x\ino \bbR_+^*$ that 
$$  \bP \big( \!\!\!\!\!\!\!\!\!\!\!\!  \sup_{\qquad s\in [0, 2^{-p}] } \!\!\!\!\!\!\!\!\!\! \widehat{\mathscr W}^{_{(n)}}_{\! s} \geko 
\epp  \big) \leq \bP (\mathbf r_{n,x} \leqo 2^{-p} ) + \bP \big( \mathscr M_{n,x} \geko \epp  \big) $$
where we recall from (\ref{recallrnx}) the notation $\mathbf r_{n,x} $ and from (\ref{limmaxdispl}) the notation 
$\mathscr M_{n,x}$. 
We also recall that under our assumption, (\ref{rnxcv}) holds true, i.e., $\mathbf r_{n,x} \! \to \!  
2\varsigma_{-x}$ in law on $\bbR_+$. By $\textbf{Step (3)}$, \texttt{Sheu}$_{^{\,}}$($\baa, \bbb, \bmu$) 
implies (\ref{limmaxdispl}) for all $x\ino \bbR_+^*$. 
By (\ref{exitsys}) in Proposition \ref{BrosnapsiH} $(ii)$, $\smash{\bP (\max_{s \in [0, 2\varsigma_{-x}] } \widehat{W}_{\! s} 
\geko \epp )\eqo 1\! -\! e^{-xw(\epp)}}$. Thus,
$$ \lim_{p\to \infty} \limsup_{n\to \infty}   \bP \Big( \!\!\!\!\!\!\!\!\!\!\!  \sup_{\qquad s\in [0, 2^{-p}] } 
\!\!\!\!\!\!\!\!\!\! \widehat{\mathscr W}^{_{(n)}}_{\! s} \geko \epp  \Big) \leq  \lim_{p\to \infty}\bP \big(
 2\varsigma_{-x} \leqo  2^{-p} \big)+ 1\! -\! e^{-xw(\epp)} =1\! -\! e^{-xw(\epp)} 
 \underset{x\to 0^+}{-\!\!\!-\!\!\!\longrightarrow }0 .$$
This completes the proof of \textbf{Step (7)}. \cq

\smallskip

We recall that we overall assume \texttt{\L{}uka}$_{^{\,}}$($\baa, \bbb, \bmu, \psi$) and \texttt{Grey}$_{^{\,}}$($\baa, \bbb, \bmu$) and that we only need to prove the theorem when $\beta\eqo 1$.

\smallskip

We first prove Theorem \ref{Sheuexplain} $(i)$ in $\textbf{Case (0)}$. We first note that 
\texttt{Sheu}$_{^{\,}}$($\baa, \bbb, \bmu$) implies  \texttt{Sheu}$_{^{\,}}$($\psi$) by Lemma \ref{controlpsi} $(i)$. 
By (\textbf{Step 3}), \texttt{Sheu}$_{^{\,}}$($\baa, \bbb, \bmu$) also implies (\ref{limmaxdispl}) and 
by (\textbf{Step 5}),    \texttt{Sheu}$_{^{\,}}$($\psi$) and (\ref{limmaxdispl}) imply 
$b_n /(a_n (\log a_n)^2)\! \to \! \infty$, which completes the proof of Theorem \ref{Sheuexplain} $(i)$ in $\textbf{Case (0)}$. 

Let us prove Theorem \ref{Sheuexplain} $(ii)$ in $\textbf{Case (0)}$. We first assume 
\texttt{Sheu}$_{^{\,}}$($\baa, \bbb, \bmu$). Then $\textbf{Step (7)}$ implies that the laws of the endpoint processes 
$\smash{(\widehat{\mathscr W}^{_{(n)}}_s)_{s\in \bbR_+}}$ are tight in $\bC^{_0}_{^1}$ and by $\textbf{Step (1)}$, 
(\ref{limiBrosna1}) holds true. Conversely, by (\textbf{Step 4}), (\ref{limiBrosna1}) implies 
\texttt{Sheu}$_{^{\,}}$($\baa, \bbb, \bmu$),  which completes the proof of Theorem \ref{Sheuexplain} $(ii)$ in 
$\textbf{Case (0)}$.

To prove Theorem \ref{Sheuexplain} $(iii)$ in $\textbf{Case (0)}$, we only need to prove that  
(\ref{limiBrosna1}) is equivalent to (\ref{limiBrosna2}). 

\emph{Indeed}, let us suppose that (\ref{limiBrosna1}) $\Leftrightarrow$ (\ref{limiBrosna2}).
If the laws of $\smash{\widehat{W}^{_{(n)}}_\cdot}$ are tight, (\textbf{Step 2}) implies \texttt{Sheu}$_{^{\,}}$($\psi$) and 
(\ref{limiBrosna2}), and thus \texttt{Sheu}$_{^{\,}}$($\psi$) and (\ref{limiBrosna1}); then (\textbf{Step 4}) entails 
\texttt{Sheu}$_{^{\,}}$($\baa, \bbb, \bmu$). Conversely, let us assume \texttt{Sheu}$_{^{\,}}$($\baa, \bbb, \bmu$). 
We have proved Theorem \ref{Sheuexplain} $(ii)$  in $\textbf{Case (0)}$, which asserts that $(\ref{limiBrosna1})$ 
holds and thus (\ref{limiBrosna2}) too: it implies that  the laws of $\smash{\widehat{W}^{_{(n)}}_\cdot}$ are tight. 
This proves that \texttt{Sheu}$_{^{\,}}$($\baa, \bbb, \bmu$) is equivalent to the tightness of the laws of 
$\smash{\widehat{W}^{_{(n)}}_\cdot}$. Moreover, if the laws of  $\smash{\widehat{W}^{_{(n)}}_\cdot}$ are tight, (\textbf{Step 2}) 
entails (\ref{limiBrosna2}), which completes the proof of Theorem \ref{Sheuexplain} $(iii)$ in $\textbf{Case (0)}$, 
provided that we have proved the equivalence (\ref{limiBrosna1}) $\Leftrightarrow$ (\ref{limiBrosna2}). 

\smallskip

We now prove that  (\ref{limiBrosna1}) $\Leftrightarrow$ (\ref{limiBrosna2}). 
Let us assume (\ref{limiBrosna1}) (resp.~let us assume (\ref{limiBrosna2})). In particular, by (\textbf{Step 1}) 
(resp.~by (\textbf{Step 2})), it implies  \texttt{Sheu}$_{^{\,}}$($\psi$).  
Under \texttt{\L{}uka}$_{^{\,}}$($\baa, \bbb, \bmu, \psi$) and \texttt{Grey}$_{^{\,}}$($\baa, \bbb, \bmu$), 
we have proved (\ref{rnxcv}), i.e., $\mathbf r_{n,x} \! \to \!  2\varsigma_{-x}$ in law on $\bbR_+$ (resp.~we recall 
from (\ref{szhghtGW}) that $\varsigma_{n,x}\! \to \! \varsigma_{-x}$ in law on $\bbR_+$). 
Combined with (\ref{limiBrosna1}) (resp.~combined with (\ref{limiBrosna2}) it implies (\ref{limmaxdispl}) (resp.~it 
implies (\ref{limmaxdispl2}))) for all $x\ino \bbR_+^*$, and by $\textbf{Step (5)}$ we get $\lim_{n\to \infty} b_n / (a_n (\log a_n)^2)\eqo \infty$.

We next recall (\ref{Lambdaphicv}) in $\textbf{Step (E)}$ of the proof of Theorem \ref{Greyexplain} (which holds under 
the assumptions \texttt{\L{}uka}$_{^{\,}}$($\baa, \bbb, \bmu, \psi$) and \texttt{Grey}$_{^{\,}}$($\baa, \bbb, \bmu$)). 
Namely, 
$\Lambda_n $ converges to the identity in probability in $\bC^{_0}_{^1}$. 
By standard arguments, the inverse functions $\smash{\Lambda_n^{-1}}$ converge to the identity in probability in $\bC^{_0}_{^1}$ too. 
Thus (\ref{limiBrosna1}) (resp.~(\ref{limiBrosna2})) combined with an elementary argument entails that the processes 
$\smash{(\widehat{\mathscr W}^{_{(n)}}_{^{\! \Lambda^{-1}_n (s)}} \!\!   -\! \widehat{\mathscr W}^{_{(n)}}_s)_{s\in \bbR_+}}$ 
(resp.~the processes $\smash{(\widehat{ W}^{_{(n)}}_{^{\! \!\Lambda^{-1}_n (s)}}\!\!   -\! 
\widehat{W}^{_{ (n)}}_s)_{s\in \bbR_+}}$) converge to the null function in probability in $\bC^{_0}_{^1}$. 
We observe that $\textbf{Step (6)}$ applies and (\ref{scrWWcontig}) holds. By an elementary argument, 
$\smash{( \widehat{\mathscr W}^{_{(n)}}_{{\!s}}\! \!\!  -\! \widehat{W}^{_{(n)}}_{^{\!\! \Lambda^{-1}_n (s)}})_{s\in \bbR_+}}$ 
also converge to the null function in probability in $\bC^{_0}_{^1}$. 
As a consequence, the processes $\smash{(\widehat{\mathscr W}^{_{(n)}}_{\! s} \!\!  -\! \widehat{W}^{_{(n)}}_s)_{s\in \bbR_+}}$ 
converge to the null function in probability in $\bC^{_0}_{^1}$, which implies (\ref{limiBrosna2}) (resp.~(\ref{limiBrosna1}), with the help of Theorem \ref{Greyexplain}). This shows that  
(\ref{limiBrosna1}) and  (\ref{limiBrosna2}) are equivalent and it completes the proof of Theorem \ref{Sheuexplain} in 
$\textbf{Case (0)}$. \cq

\paragraph*{Proof of Theorem \ref{Sheuexplain} $(iii)$ in $\textbf{Case (1)}$ and $\textbf{Case (2)}$.}
Standard arguments  entail that (\ref{limiBrosna2}) implies $\textbf{Case (1)}$ for all $x\ino (0,\infty)$. 
It remains to prove that (\ref{limiBrosna2}) in $\textbf{Case (0)}$ implies (\ref{limiBrosna2}) in $\textbf{Case (2)}$
as a consquence of a result of D.~\& Le Gall \cite{DuLG02} Proposition 2.5.2. 

\smallskip

\noi
\emph{Proof.} We recall that $\tau^{_\infty}_{n}\eqo (\tau_n(p))_{p\in \bbN}$ is a forest of independent GW($\mu_n$)-trees and we fix $c\ino \bbR_+^*$. We introduce $\rho_{n,c} \eqo \inf\{ s\ino \bbR_+\!  : \! C_s (\tau^{_\infty}_{n})\geqo c\lambda_n \}$, 
$$\ell_{n,c}\eqo \sup \big\{ s\ino [0,  \rho_{n,c}]\! : \! C_s (\tau^{_\infty}_{n})\eqo 0 \big\} \quad \textrm{and} \quad r_{n,c}\eqo \inf \big\{ s\ino [ \rho_{n,c}, \infty )\! : \! C_s (\tau^{_\infty}_{n})\eqo 0 \big\}$$
and we also set $\overline{V}_{\! k}( \tau^{_\infty}_{n}) )\eqo V_k (\tau^{_\infty}_{n}))-\inf_{l\in [0, k]}V_l ( \tau^{_\infty}_{n})$, $k\ino \bbN$. 
Then,  
\begin{eqnarray*}
( \overline{V}_{\!\! (\frac{_{1}}{^{2}}\ell_{n,c} + \cdot) \wedge \frac{_1}{^2}r_{n,c}}( \tau^{_\infty}_{n}) \!\!\!\!\!\!\!\! & ,& \!\!\!\!  \!\!\!\!C_{(\ell_{n,c} + \cdot) \wedge r_{n,c}} (\tau^{_\infty}_{n})  ,  W_{(\ell_{n,c} + \cdot) \wedge r_{n,c}} (\tau^{_\infty}_{n}, \cdot) ) 
\; \textrm{under $\bP$} \\
\!\!\!\! &\overset{\textrm{(law)}}{=} &  (V_\cdot (\tau_n(1)), C_\cdot (\tau_n(1)), W_\cdot (\tau_n(1), \cdot ))  \; \textrm{under $\bP(\, \cdot \,| \,  \Gamma (\tau_n (1))\geqo c\lambda_n )$,}
\end{eqnarray*}
which is the discrete version of (\ref{selectexcu}).  
Then, we recall from (\ref{sigellerrc}) the defintion of  $(\ell_c, r_c)$.  
The arguments of the proof of Prop.~2.5.2 in D.~\& Le Gall \cite{DuLG02} immediately show that the weak limit 
$(b_n^{-1} \ell_{n,c}, b_n^{-1} r_{n,c}) \! \to \! (\ell_c, r_c)$ holds in law jointly with the convergence (\ref{limiBrosna2}) in $\textbf{Case (0)}$: thanks to the previous arguments and (\ref{selectexcu}), it 
entails (\ref{limiBrosna2}) in $\textbf{Case (2)}$. This completes the proof of Theorem \ref{Sheuexplain}. \cqfd

\subsection{Proof of Theorem \ref{extenscv}}
\label{Thms23pfsec}

In this section, $(a_n)_{n\in \bbN}$ and $(b_n)_{n\in \bbN}$ are two renormalization sequences that satisfy 
\texttt{Norm}$_{^{\,}}$($\baa, \bbb$), $(\btt_n)_{n\in \bbN}$ is a sequence of a.s.~finite random trees whose rescaled 
contour processes $C^{_{(n)}}_{\cdot}$ are as in (\ref{renormCW}), $(C_s)_{s\in \bbR_+}$ is a continuous 
nonnegative process and $W_\cdot$ stands for a $\bbR^d$-valued Brownian snake with lifetime process 
$C_\cdot$ as in Definition \ref{Brosnadef} (here $W_\cdot$ is not necessarily continuous).

We recall that $\langle  \cdot \, ,  \cdot \rangle$ stands for the canonical scalar product on $\smash{\bbR^d}$, $\smash{\lvert \, \cdot \, \rvert}$ 
is the associated Euclidean norm and the canonical basis is denoted by $\smash{(\mathtt e_i)_{1\leq i\leq d}}$. We denote by 
$\smash{\mathrm{Sym}^{_+}_{^d}}$ the space of nonnegative $d\! \times \! d$ symmetric real matrices, which is equipped with the max-norm $\smash{\lVert \cdot \rVert}$. For all $\smash{\beta\ino  \mathrm{Sym}^{_+}_{^d}}$, $\smash{\sqrt{ \beta\, } \ino \mathrm{Sym}^{_+}_{^d}}$ a square root of $\beta$ such that $\smash{\beta \! \mapsto \! \sqrt{\beta\, }}$ is measurable, as specified in Section \ref{statementssec}.  

We first discuss convergence of finite dimensional marginal laws of snakes whose jumps are conditionally 
independent by group of siblings. Since we apply the following lemma to several processes, 
we state it for auxiliary BRWs. 
More precisely, we fix $\kappa \ino (0, 1) $ and w.l.o.g.~we 
assume that there is a moment gauge function $\smash{G_{\! \mathtt a}}$ such that $\smash{G_{\! \mathtt a} (x)\geqo |x|^{2+\kappa}} $ for all $\smash{x\ino \bbR}$. Let $\smash{\fbeta_{\! \mathtt a}\! : \! \Omega \! \to \! \mathrm{Sym}^{_+}_{^d}}$ be measurable and for all $n\ino \bbN$, let $\smash{\bS^{\mathtt a}_n\eqo (S^{\mathtt a} _{n,u})_{u\in \btt_n}}$ be a $\smash{\bbR^d}$-valued BRW whose snake $\smash{W^{_{\mathtt a , (n)}}_{\cdot}}$ is normalized as in (\ref{renormCW}). 

\begin{lemma}
\label{fdcvsnakeSibInd} We keep the above notations and suppose that the BRWs $\smash{(\bS^{ \mathtt a}_{n})_{n\in \bbN}}$
 satisfy Assumptions $\mathbf{(2}$\textbf{-}$\mathbf{5)}$ in \emph{\textbf{Sib-Ind}$_{_{\, }}$($\smash{\baa, \bbb, \bS^{\mathtt a}_\cdot, G_{\! \mathtt a}, C_\cdot, \fbeta_{\! \mathtt a}}$)}. Then, the finite dimensional marginal laws of the 
 $\smash{\big(C^{_{(n)}}_{\cdot}\! , \widehat{W}^{_{\mathtt a , (n)}}_{\cdot}\big)}$ weakly converge to those of $\smash{\big( C_\cdot\, , \sqrt{\fbeta_{\! \mathtt a}} . \widehat{W}\big)}$. 
 \end{lemma}
\noi
\textbf{Proof.} To simplify, we refer to \textbf{Sib-Ind}$_{_{\, }}$($\baa, \bbb, \bS^{\mathtt a}_\cdot, G_{\! \mathtt a}, C_\cdot, \fbeta_{\! \mathtt a}$) as to \textbf{Sib-Ind}$_{_{\, }}$. 
Let $(\xi_{{n,u}}^{{\mathtt a}})_{u\in \btt_n \backslash \{ \varnothing \}}$ be the jumps of $\bS_{n}^{{\mathtt a}}$.
By \textbf{Sib-Ind} (\textbf{4}) and (\textbf{5}) 
there are $\btt_n$-measurable r.v.s $\fbeta_{\! n} \! : \! \Omega \! \to \! \mathrm{Sym}^+_d$ such that 
$(C^{_{(n)}}_{\cdot}\! , \fbeta_{\! n}, \bdelta_n) \! \to \! (C_\cdot, \fbeta_{\! \mathtt a} , 0)$ 
in law in $\bC^{_0}_{^1} \! \times \! \mathrm{Sym}^+_d\! \times \! \bbR_+$, where we have set 
$\bdelta_n\!  := \! \max_{u\in \btt_n\backslash \{ \varnothing\}} \lVert \fbeta_{\! n} \! -\! \beta (\bS^{\mathtt a}_n, u) \rVert $ 
and where for all $1\leqo i,j\leqo d$ and all $u\in \btt_n \backslash \{ \varnothing\}$, 
$\beta_{i,j} (\bS^{\mathtt a}_n, u)\eqo \bE [\langle \mathtt{e}_i, \xi_{n,u}^{\mathtt a} \rangle \langle \mathtt{e}_j,
 \xi_{n,u}^{\mathtt a} \rangle| \btt_n]$. To simplify notation we also set 
 $M_n \eqo \max_{u\in \btt_n\backslash \{ \varnothing\}}\bE [ G_{\! \mathtt a} (|\xi_{n,u}^{\mathtt a}|) | \btt_n]$. 
 By \textbf{Sib-Ind}$_{_{\, }}$ (\textbf{3}), $M_n \leko \infty$ a.s.~and $\lim_{y\to \infty}\limsup_{n\to \infty}
 \bP (M_n \geko y)\eqo 0$. 

  Let $p\ino \bbN^*$, $s_p \geko \ldots \geko s_1 \geqo 0$, $z_1, \ldots, z_p \ino \bbR$ and 
  $y_1, \ldots, y_p \ino \bbR^d$. We have to prove 
\begin{equation} 
\label{cvcharacCW}
\lim_{n\to \infty} \bE \Big[\!\! \! \!  \prod_{\;\; 1\leq k\leq p} \!\!\!\! e^{\mathtt i z_kC^{(n)}_{s_k} +
\mathtt i \langle y_k ,  \widehat{W}^{_{\mathtt a , (n)}}_{\! s_k} \rangle } \Big] =  \bE \Big[\!\! \! \! 
 \prod_{\;\; 1\leq k\leq p} \!\!\!\! e^{\mathtt i z_k C_{s_k} + \mathtt i 
\langle y_k ,  \sqrt{\fbeta_{\! \mathtt a}}.\widehat{W}_{\! s_k} \rangle } \Big] .
\end{equation}
To that end, we use Proposition \ref{fdestimate} to compare $\bS^{\mathtt a}_n$ to a BRW with i.i.d.~Gaussian jumps. 
More precisely, let $Y_u$, $u\ino \bbU$, be independent $\bbR^d$-valued Gaussian r.v.s whose covariance 
matrix is the identity. It is also convenient to assume that $(Y_u)_{u\in \bbU}$ and $(\btt_n, \bS^{\mathtt a}_n)_{n\in \bbN}$ are 
independent. 
For all $u\ino \btt_n \backslash \{ \varnothing\}$, we set $R^\bullet_u \eqo \sum_{v\in \rgeo \varnothing, u \rgeo} Y_v$ 
and $R^\bullet_\varnothing \eqo 0$. We also denote by $W^{\bullet, (n)}$ the snake associated with the BRW $(R^\bullet_u)_{u\in \btt_n}$, which is normalized 
as in  (\ref{renormCW}). 

We first observe that conditionally given $\btt_n$, 
$\smash{\big(\sqrt{\fbeta_n}.\widehat{W}^{_{\bullet, (n)}}_{^{\! k/b_n}} \big)_{0\leq k \leq 2(\#\btt_n-1)}}$ is a $\bbR^d$-valued Gaussian process whose conditional covariance is 
$$\big(  \bE \big[ \langle \mathtt e_i , \sqrt{\fbeta_n}.\widehat{W}^{{\bullet, (n)}}_{{\! k/b_n}}\rangle \langle \mathtt e_j, \sqrt{\fbeta_n}.\widehat{W}^{\bullet, (n)}_{l/b_n}\rangle  \big | \btt_n \big] \big)_{1\leq i,j\leq d}= m_{C^{(n)}}\!  \big(\tfrac{k}{b_n}, \tfrac{l}{b_n}\big) \fbeta_n \; , $$
where we recall that $ \smash{m_{C^{(n)}} (s,s')\eqo \min_{r\in [s\wedge s', s\vee s']} C^{(n)}_r}$. We set $s_{n,k} \eqo \lfloor b_n s_k \rfloor/ b_n$, $1\leqo k\leqo p$. Thus, we a.s.~get 
\begin{equation}
\label{explictGausssnake}
 \bE \Big[ \prod_{1\leq k\leq p} e^{\mathtt i  \langle y_k,\sqrt{\fbeta_n}. \widehat{W}^{{\bullet, (n)}}_{^{\! s_{n,k}}} \rangle} \Big| \btt_n \Big] 
= \exp \Big(\!\! -\! \tfrac{1}{2}\!\! \sum_{1\leq k, l\leq p } \!\!\! m_{C^{(n)} } (s_{n,k}, s_{n,l})
 \langle y_k , \fbeta_n. y_l \rangle \Big). 
 \end{equation}
We next observe that there is a r.v.~$Z_{n,k}$ whose conditional law given $\btt_n$ is $\mathcal N(0, \fbeta_n)$ and such that 
$\smash{|\widehat{W}^{_{\bullet, (n)}}_{{\! s_{n,k}}} \!\! -\!  \widehat{W}^{_{\bullet, (n)}}_{{\! s_{k}}} |\leqo \lambda_n^{_{\! -1/2}}Z_{{n,k}} \! \to \! 0}$ in probability as $n\to \infty$. Similarly, $|C^{_{(n)}}_{{s_{n,k}}} \! -\! C^{_{(n)}}_{{s_k}}| \leqo 1/\lambda_n$. 
Since $\smash{(C^{(n)}, \fbeta_n)\to (C, \fbeta_{\mathtt a})}$ in law on $\bC^{_0}_{^1} \! \times \! \mathrm{Sym}^+_d$, we get
\begin{eqnarray} 
\lim_{n\to \infty} \bE \Big[\!\! \! \!  \prod_{\;\; 1\leq k\leq p} \!\!\!\! e^{\mathtt i z_kC^{(n)}_{s_k} + \mathtt i \langle y_k , \sqrt{\fbeta_n}. \widehat{W}^{_{\bullet , (n)}}_{\! s_k} \rangle } \Big] \!\!\!\! &=  &  \!\!\!\! 
\bE \Big[ e^{\,  \sum_{1\leq k\leq p} \mathtt i z_k C_{s_k}  -\tfrac{_1}{^2} \sum_{1\leq k,l\leq p} m_C (s_k, s_l)  
\langle y_k ,  \fbeta_{\! \mathtt a} .y_l \rangle }  \Big] \nonumber \\
\!\!\!\! &=  &  \!\!\!\!  \bE \Big[\!\! \! \!  \prod_{\;\; 1\leq k\leq p} \!\!\!\! e^{\mathtt i z_k C_{s_k} + \mathtt i \langle y_k ,  \sqrt{\fbeta_{\! \mathtt a}}.\widehat{W}_{\! s_k} \rangle } \Big] . \label{cvcharacCbulletW}
\end{eqnarray}
Namely, $\smash{(C^{_{(n)}}_{\cdot}\! ,  \sqrt{\fbeta_{\! n}} .\widehat{W}^{\bullet, (n)})}$ $\! \to \! $ $\smash{(C_\cdot, \sqrt{\fbeta_{\! \mathtt a}} . 
\widehat{W})}$ in 
the sense of the convergence in distribution of their finite dimensional marginal laws. 

Next, we set $u_{n,k} \eqo v (b_ns_{n,k})$ where 
$(v(l))_{0\leq l \leq 2(\# \btt_n -1) }$ stands for the contour exploration of $\btt_n$. Thus for all $k\ino \{ 1, \ldots, p\}$, 
$$\lambda_n C^{{(n)}}_{{s_{n,k}}}\eqo |u_{n,k}| , \quad \lambda^{1/2}_n \, 
\widehat{W}^{{\mathtt a , (n)}}_{{s_{n,k}}}\eqo \bS^{{\mathtt a}}_{{n, u_{n,k}}} \quad \textrm{and} \quad  
\lambda_n^{1/2} \, \widehat{W}^{{\bullet , (n)}}_{{s_{n,k}}}\eqo R^{{\bullet}}_{{n,u_{n,k}}} .$$
Let $c\ino \bbR_+^*$ be such that $cG_{\! \mathtt a}(x) \geqo |x|$ for all $x\ino \bbR$.  
We observe that  
\begin{equation}
\label{approxdiscr}
\bE \big[ \big| \widehat{W}^{{\mathtt a , (n)}}_{{s_{n,k}}}\! -\! \widehat{W}^{{\mathtt a , (n)}}_{{s_{k}}}\big|\,\big| \, \btt_n  \big] \leq 2\lambda_{^n}^{_{\! -1/2}} \!\!  \!\! \!\! \!\! \! \max_{\quad u\in \btt_n\backslash \{ \varnothing \}} \!\! \!\! \!\! \! \bE \big[ \big|\xi_{{n,u}}^{{\mathtt a}} \big| \, \big|\,  \btt_n\big]  \leq 2c\lambda_{^n}^{_{\! -1/2}}\! M_n .
\end{equation}
We next apply Proposition \ref{fdestimate}: we set $r_0\! := \!  \max_{1\leq j\leq p} |y_j| $ and (\ref{fdcontrol}) a.s.~implies 
\begin{eqnarray}
\label{fdcontrolbis}
\Big| \bE \Big[\!\! \! \!  \prod_{\;\; 1\leq k\leq p} \!\!\!\! \!\!\!\! & & \!\!\!\! \!\!\!\! 
e^{ \mathtt i \langle y_k ,  \widehat{W}^{_{\mathtt a , (n)}}_{\! s_{n,k}} \rangle } \Big| \btt_n \Big] -
 \bE \Big[\!\! \! \!  \prod_{\;\; 1\leq k\leq p} \!\!\!\! e^{ \mathtt i \langle y_k ,  
 \sqrt{\fbeta_{\! n}}.\widehat{W}^{_{\bullet , (n)}}_{\! s_{n,k}} \rangle } \Big| \btt_n \Big] \Big| \leq p^3 
 \big( cM_n\! +\!  \sqrt{d\lVert \fbeta_n \rVert }\big)\lambda_{^n}^{_{\! -1/2}} r_0 \nonumber  \\
\!\!  \!\! \!\! \!\!  & +& \!\! \!\! d^2p^4 \big(C^{{(n)}}_{{s_{n,1}}} + \ldots +C^{{(n)}}_{{s_{n,p}}} \big) r_0^2 \big( \bdelta_n +
\lVert \fbeta_n \rVert^2 \lambda_{^n}^{_{\! -1}}r_0^2+ M_n \lambda_{^n}^{_{\! -\kappa /2}} r_0^\kappa \big). 
\end{eqnarray}
Combined with (\ref{approxdiscr}) and the convergence $\smash{(C^{_{(n)}}_{\cdot}\! , \fbeta_{\! n}, \bdelta_n) \! \to \! (C_\cdot, \fbeta_{\! \mathtt a} , 0)}$ 
in law in $\bC^{_0}_{^1} \! \times \! \mathrm{Sym}^+_d\! \times \! \bbR_+$, it implies that 
$$ \lim_{n\to \infty} \Big| \bE \Big[\!\! \! \!  \prod_{\;\; 1\leq k\leq p} \!\!\!\! e^{ \mathtt i \langle y_k ,  \widehat{W}^{_{\mathtt a , (n)}}_{\! s_{k}} \rangle } \Big| \btt_n \Big]- \bE \Big[\!\! \! \!  \prod_{\;\; 1\leq k\leq p} \!\!\!\! e^{ \mathtt i \langle y_k , \sqrt{ \fbeta_{\! n}}.\widehat{W}^{_{\bullet , (n)}}_{\! s_{k}} \rangle } \Big| \btt_n \Big] \Big| =0$$
in probability, which easily entails (\ref{cvcharacCW}) by (\ref{cvcharacCbulletW}).  \cqfd

\medskip

To prove tightness in Theorem \ref{extenscv}, we next use the coupling in Theorem \ref{brwcouplingth} and Lemma \ref{Gaussjumps}. We assume that 
$a_n$, $b_n$, $\btt_n$, $\smash{C^{(n)}_\cdot}$, $C_\cdot$, $\bS_n^*$ and $\bS_n$  are as in Theorem \ref{extenscv}. It 
is convenient to state an intermediate result for other auxilliary processes. More precisely, we fix a moment gauge function $G_0$, a 
r.v.~$\smash{\fbeta^0 \!:\!  \Omega \! \to \! \bbR_+}$ and a sequence of BRWs $\smash{\bS^0_n\eqo (S_{n,u}^0)_{u\in \btt_n}}$ whose jumps are denoted by $\smash{(\xi_{n,u}^0)_{u\in \btt_n \backslash \{ \varnothing\}}}$ and that satisfy the following. 
\begin{compactenum}

\smallskip

\item[$(a_0)$] The BRWs $(\bS^{0}_{n})_{n\in \bbN}$ are $\bbR$-valued and satisfy \textbf{Sib-Ind}$_{_{\, }}$($\baa, \bbb, \bS^0_\cdot, G_{\! 0}, C_\cdot, \fbeta^{0}$). 
\smallskip

\item[$(b_0)$] $\bP$-.a.s.~$\fbeta^0 \in \bbR_+^*$.   

\smallskip

\end{compactenum}
We also introduce i.i.d.~standard $\bbR$-valued Gaussian r.v.s $(Y_u)_{u\in \bbU}$ that are supposed to be independent of 
$\btt_n$, $\bS^{0}_{n}$, $\bS^{*}_{n}$, $\bS_{n}$. For all $n\in \bbN$ and $u\ino \btt_n \backslash \{ \varnothing \}$, we recall from (\ref{concovdef}) that $\beta (\bS_n^0, u) $ $\eqo$  $\bE [(\xi^0_{n,u})^2 | \btt_n]$ and we set 
\begin{equation}
\label{RRbulletdef} 
R_{n,u} \eqo \sum_{v\in \rgeo \varnothing, u\rgeo} \sqrt{\! \beta (\bS^0_n, v)\, } Y_v \quad \textrm{and} \quad  R^\bullet_{n,u} \eqo \sum_{v\in \rgeo \varnothing, u\rgeo}Y_v
\end{equation}
and $\smash{R_{n, \varnothing} \eqo R^\bullet_{n, \varnothing}\eqo 0}$. We recall that $\smash{W^{*,(n)}_\cdot}$ and $\smash{W^{(n)}}$ are the snakes associated with $\bS_n^*$ and $\bS_n$ which are normalized as in (\ref{renormCW}). We denote by $\smash{W^{0, (n)}}$, $\smash{W^{1, (n)}}$ and $\smash{W^{\bullet, (n)}}$ the snakes associated with resp.~$\bS_n^0$, $\smash{\mathbf R_n\eqo (R_{n,u})_{u\in \btt_n} }$ and $\smash{\bR^\bullet_n\eqo (R^\bullet_{n,u})_{u\in \btt_n} }$, which are normalized as in (\ref{renormCW}). For all $N,p\ino \bbN$ and for all continuous snakes $(w,h)$ we also recall the notation 
$$ \omega_{N\! ,p} (w) \eqo \max_{0\leq j< N2^p} \max \big\{\,  | \widehat{w} (s)\! -\! \widehat{w}(j2^{-p}) |\, ; \, s\ino [j2^{-p} , (j+1)2^{-p}] \big\}  \; .$$ 

The following lemma is the key-point of the proof of Theorem \ref{extenscv}. The main result is Lemma \ref{coupleffec} $(v)$ but for sake of clarity the intermediate steps are stated below as Lemma \ref{coupleffec} $(i$-$iv)$.  
\begin{lemma}
\label{coupleffec} We keep the previous notations and assumptions. Then for all $\epp,\eta \ino (0, \frac{_1}{^4})$, there are $\theta_\epp \ino \bbR_+^*$ and $n_{\epp, \eta} \ino \bbN^*$ such that $\lim_{\epp \to 0^+} \theta_\epp\eqo 0$ and such that for all integers $n\geqo n_{\epp, \eta}$, all $N,p\ino \bbN$, all $z\ino \bbR_+^*$, the following holds true.

\smallskip

\noi
$\! (i)$ $\bP \big( \omega_{N,p} (W^{1, (n)} ) \geko 5\eta \big) \leqo \theta_\epp +  \bP \big( \omega_{N,p} (W^{0, (n)} ) \geko 3\eta \big) \leqo 2\theta_\epp+\bP \big( \omega_{N,p} (W^{1, (n)} ) \geko \eta \big) $.

\smallskip

\noi
$\! (ii)$ $\bP \big( \sup_{ s\in \bbR_+}  | \widehat{W}^{1, (n)}_s | \geko z+ \eta \big) \leqo  \theta_\epp + 
\bP \big( \sup_{ s\in \bbR_+}  | \widehat{W}^{0, (n)}_s |  \geko z \big) $.

\noi
We set $q_{n,\eta, \epp}\eqo \theta_\epp + 4 \bP \big( \sup_{s\in \bbR_+} | \widehat{W}^{\bullet, (n)}_s|  \geko \tfrac{\eta}{2\epp^{7/2}} \big) + 4e^{-\eta^2/(16\epp^4)}$ and we get the following.

\noi
$\! (iii)$ $ q_{n,\eta, \epp} \leq 5\theta_\epp +  16\, \bP \big(\sup_{s\in \bbR_+} | \widehat{W}^{1, (n)}_s | \geko \tfrac{\eta}{4\epp^3} \big) + 20\, e^{-\frac{\eta^2}{ 16\epp^4}}  $, 

\smallskip

\noi
$\! (iv)$ $\! \bP \big( \omega_{N,p} (W^{\bullet, (n)} )\!  \geko\!  \tfrac{5\eta}{\sqrt{2\epp}} \big) \leqo 2q_{n,\eta, \epp}\! +  \bP \big( \omega_{N,p} (W^{1, (n)} ) \! \geko\!  3\eta \big) \leqo  4q_{n,\eta, \epp}\! + $ $ \bP \big( \omega_{N,p} (W^{\bullet, (n)} )\!  \geko  \eta\sqrt{2\epp} \big)  $. 

\smallskip

\noi
$(v)$ The laws of $(\widehat{W}^{0, (n)}_\cdot )_{n\in \bbN}$ are tight in $\bC^{_0}_{^1}$ iff the same holds true for the laws of $(\widehat{W}^{\bullet, (n)}_\cdot)_{n\in \bbN}$. 
\end{lemma}
\noi
\textbf{Proof.} We fix $\epp , \eta \ino (0, \frac{_1}{^4})$ and we recall from Theorem \ref{brwcouplingth} the notation 
$K(C_1, C_2, C_3, G)$ for the coupling constant $K$ in (\ref{coupling1}). W.l.o.g., we suppose that $(\Omega, \ccF, \bP)$ is rich enough to carry a uniform r.v.~$U\! : \! \Omega \! \to \! [0, 1]$ that is independent of all the above mentioned r.v.s. By $\mathbf{ (6)}$ and $\mathbf{ (3)}$ in \textbf{Sib-Ind}$_{_{\, }}$($\baa, \bbb, \bS^0_\cdot, G_{\! 0}, C_\cdot, \fbeta^{0}$), there is $c_\epp \ino \bbR_+^*$ such that 
\begin{equation}
\label{tnMsncontrol}
 \forall n\in \bbN, \quad \bP \big( \# \btt_n \geqo \tfrac{1}{2}b_n c_\epp\big) \leqo \epp \quad \textrm{and} \quad \bP \big( M(\bS_n^0, G_0 ) \geqo c_\epp\big)\leqo \epp \; .
\end{equation} 
We shall apply Theorem  \ref{brwcouplingth} with $C_1\eqo c_\epp$, $C_2\eqo \epp$, $C_3 \eqo (2\epp)^{-1}$ and $G\eqo G_0$ and we take $x_0$ and $\kappa$ as in (\ref{momgaudef}) with $G\eqo G_0$.  To simplify notation, we then set $K_\epp\eqo K(c_\epp, \epp, \frac{1}{2\epp}, G_0)$  
and 
$$ \theta_\epp= 4\epp + \sup_{n\in \bbN} \bP \big( \Gamma (\btt_n) \geqo \lambda_n/\epp\big) + 
 \sup_{n\in \bbN} \bP \big( \fbeta_{\! n} \notin (2\epp, \tfrac{1}{2\epp} \! -\! 1 ) \big)\;. $$
Here we recall that $\Gamma (\btt_n)\eqo \max_{u\in \btt_n} |u|$ is the total height of $\btt_n$ and we recall from $\mathbf{(4)}$ in \textbf{Sib-Ind}$_{_{\, }}$($\baa, \bbb, \bS^0_\cdot, G_{\! 0}, C_\cdot, \fbeta^{0}$) that $\fbeta_{\! n}$ is a $\bbR$-valued $\btt_n$-measurable r.v.~such that $\bdelta_n \! :=\! \max_{u\in \btt_n \backslash \{ \varnothing \}} |\fbeta_{\! n} \! -\! \beta (\bS_n^0, u)| $ $\! \to \! 0$ in law. We also observe that $2\varepsilon\leko (2\varepsilon)^{-1} \! -\! 1$ since $\epp \ino (0, \frac{_1}{^4})$. We first prove that 
\begin{equation}
\label{thetaepp} \lim_{\epp \to 0^+} \theta_\epp = 0 \; .
\end{equation}
\emph{Indeed,} by $\mathbf{(5)}$ in \textbf{Sib-Ind}$_{_{\, }}$($\baa, \bbb, \bS^0_\cdot, G_{\! 0}, C_\cdot, \fbeta^{0}$), 
and since we assume that $\bP$-a.s.~$\fbeta^0 \ino \bbR_+^*$, we get 
$\limsup_{n\to \infty} \bP ( \fbeta_{\! n} \notin (2\epp, \tfrac{1}{2\epp} \! -\! 1 ))\leqo 
\bP ( \fbeta^0 \notin (2\epp, \tfrac{1}{2\epp}  \! -\! 1 )) $ $\!\to \! 0$ as $\epp \! \to \! 0^+$. 
Moreover, for all $s\ino \bbR_+^*$, we get $\limsup_{n\to \infty} \bP ( \Gamma (\btt_n) 
\geqo \lambda_n/\epp)\leqo \limsup_{n\to \infty}\bP (\#\btt_n \geqo \frac{_1}{^2}b_ns)+ \bP (\max_{[0, s]} C \geqo 1/\epp)$. 
Thus 
$   \limsup_{\epp \to 0^+} \limsup_{n\to \infty} \bP ( \Gamma (\btt_n) \geqo \lambda_n/\epp)\leqo 
 \limsup_{n\to \infty}\bP (\#\btt_n \geqo \frac{_1}{^2} b_ns) \! \to \! 0$ as $s\! \to \! \infty$ by $\mathbf{(6)}$ in 
\textbf{Sib-Ind}$_{_{\, }}$($\baa, \bbb, \bS^0_\cdot, G_{\! 0}, C_\cdot, \fbeta^{0}$), which easily completes the proof 
of (\ref{thetaepp}).  \cq 

\medskip

We next set $y_n \eqo \eta \sqrt{\lambda_n}/ K_\epp$, for all $n\ino \bbN$. 
By $\mathbf{(1)}$ and $\mathbf{(4)}$ in \textbf{Sib-Ind}$_{_{\, }}$($\baa, \bbb, \bS^0_\cdot, G_{\! 0}, C_\cdot, \fbeta^{0}$), 
there is $n_{\epp, \eta} \ino \bbN$ such that for all integers $n\geqo n_{\epp, \eta}$, $y_n/\log_2 (c_\epp b_n) \geqo x_0$, 
\begin{equation}
\label{reglageducirque} 
\frac{c_\epp^2(1+ \epp^{-2}) b_n}{G_0( y_n / \log_2 (c_\epp b_n)) } \leqo \epp  , \;  \bP \big( \bdelta_n \geqo \epp^4\big) 
\leqo \epp \quad \textrm{and} \; \, \Big( \frac{y_n}{\log_2 (c_\epp b_n)} \Big)^{\! -2} \! 
G_0  \Big( \frac{y_n}{\log_2 (c_\epp b_n)}  \Big) \geqo \frac{c_\epp}{\epp} .
\end{equation}
We recall from Lemma \ref{HSloppconn} that $\mathtt{HS} (\btt_n) \leqo 1+ \log_2 (\# \mathtt{Lf} (\btt_n))\leqo 
\log_2(2\# \btt_n)$. Since $y\ino [x_0, \infty)\! \mapsto \! y^{-2} G_0 (y)$ is nondecreasing, by definition of 
moment gauge functions, (\ref{reglageducirque}) implies 
\begin{equation}
\label{detailtail}  
\forall n\geqo n_{\epp, \eta}, \quad \textrm{if $\# \btt_n \leqo \tfrac{1}{2}b_n c_\epp$, $\,$ then} \; \, 
 \Big( \frac{y_n}{\mathtt{HS} (\btt_n)}\Big)^{\! -2} \! G_0 \Big( \frac{y_n}{\mathtt{HS} (\btt_n)}\Big) \geqo \frac{c_\epp}{\epp} .
\end{equation}

We denote by $\bbT_{\! f} $ the set of finite ordered rooted trees. For all $n\ino \bbN$ and all $t\ino \bbT_{\! f}$, if 
$\bP (\btt_n \eqo t)\geko 0$, there exist r.v.s $(\xi_{n,u} (t))_{u\in t\backslash \{ \varnothing \}}$ such that 
$(\xi_{n,u}(t))_{u\in t\backslash \{ \varnothing \}}$ is distributed as $ (\xi^0_{n,u})_{u\in \btt_n\backslash \{ \varnothing \}} $ 
under $\bP (\, \cdot \, | \, \btt_n \eqo t)$. 
If $\bP (\btt_n \eqo t)\eqo  0$, we simply set $\xi_{n,u}(t)\eqo 0$, $u\ino t\backslash \{ \varnothing\}$. Then, we set 
\begin{equation}
\label{Snudef} \forall u\in t\backslash \{ \varnothing \}, \quad  S_{n, u} (t)\eqo 
\sum_{v\in \, \rgeo \varnothing, u\rgeo} \xi_{n,v} (t)  \quad \textrm{and} \quad S_{n,\varnothing} (t)\eqo 0. 
\end{equation}

 We next introduce the set of trees for which our coupling applies. To that end, we observe that for all 
 $n\ino \bbN$, there are functions 
 $s_n(\cdot )$ and $\beta_{n,u} (\cdot)  \! :\! \bbT_{\! f}\! \to \! \bbR_+$, 
$u\ino \bbU$, such that $\bP$-a.s.~$ s_n (\btt_n) \eqo \fbeta_{\! n} $ and $\beta_{n,u} (\btt_n)\eqo \beta(\bS^0_n, u)$ 
$=$ $\bE \big[ (\xi_{n,u})^2 \big| \btt_n \big]$, for all  $u\ino \btt_n  \backslash \{ \varnothing\} $. For all $t\ino \bbT_{\! f}$, we also set $M_n (t) \eqo \max_{u\in t\backslash \{ \varnothing\}} \bE [G_0 (\xi_{n,u} (t))]$ and 
$D_n (t)\eqo \max_{u\in t\backslash \{ \varnothing\}} |s_n(t) - \beta_{n,u} (t) |$. We also a.s.~get $M_n (\btt_n)\eqo 
M(\bS^0_n, G_0)$ and $D_n (\btt_n)\eqo \bdelta_n$. Then for all $\epp \ino (0, \frac{_1}{^4})$ and for all $n\ino \bbN$ 
we introduce the following set of coupling trees. 
$$ \bbT_{\! \epp, n}\! \eqo \Big\{  t\ino \bbT_{\! f} \! : \bP (\btt_n \eqo t) \geko 0, \, M_n (t)\leko c_\epp, \, D_n (t) \leko \epp^4, \, \# t \leko \frac{_1}{^2}b_n c_\epp , \, s_n (t) \ino \big( 2\epp, \tfrac{1}{2\epp} -\! 1 \big) , \, \Gamma (t) \leko \frac{_{\lambda_n}}{^\epp} \Big\}.$$ 
We fix $n\geqo n_{\epp, \eta}$ and we first observe that $\bP (\btt_n \! \notin \!  \bbT_{\! \epp, n}) \leqo \theta_\epp  - \epp$, by (\ref{tnMsncontrol}) and (\ref{reglageducirque}). Moreover, for any $t\ino \bbT_{\! \epp, n}$, we get the following. 
\begin{compactenum}
\item[$-$] The $(\xi_{n, u} (t))_{u\in t\backslash \{ \varnothing\}}$ are independent by group of siblings and centered. 
\item[$-$] For all $u\in t\backslash \{ \varnothing\}$, $\beta_{n,u} (t) \ino [s_n (t) \! -\! D_n (t) , s_n (t) + D_n (t)] \subset [\epp, \frac{1}{2\epp}]$. 
\item[$-$] By (\ref{reglageducirque}) and (\ref{detailtail}), Theorem \ref{brwcouplingth} applies with $(C_1,C_2,C_3,G,y)\eqo (c_\epp,\epp, \frac{1}{2\epp}, G_0, y_n)$.
Therefore, there are i.i.d.~real valued standard Gaussian r.v.s $(\widetilde{Y}_{n,u} (t))_{u\in t \backslash \{ \varnothing \}}$ that  measurably depend on $\epp$, 
$\eta$, $n$, $t$, $(\xi_{n,u} (t))_{u\in t \backslash \{ \varnothing \}}$ and $U$, and that satisfy the following: for all $u\ino 
t\backslash \{ \varnothing\}$, we set 
$\smash{ \widetilde{R}_{n,u} (t) \! :=\! \sum_{v\in \rgeo \varnothing, u\rgeo} \sqrt{\beta_{n,v} (t)} \, \widetilde{Y}_{n,v} (t)}$ and $\smash{\widetilde{R}_{n,\varnothing} (t)\! := } 0$, and we define the event 
\begin{equation}
\label{anticouplingevent}
A_{n,\epp, \eta, t}\eqo  \Big\{\!  \max_{u\in t} \big| S_{n,u} (t) \! -\! \widetilde{R}_{n,u} (t)\big| \geko \eta\sqrt{\! \lambda_n}  \Big\}.
\end{equation}
Then (\ref{coupling1}) in Theorem \ref{brwcouplingth} and (\ref{reglageducirque}) imply 
\begin{equation}
\label{coucouplpl} 
\bP (A_{n,\epp, \eta, t} ) \leqo 2C_1\!  \Big(1\! + \! \tfrac{2C_3}{C_2} \Big)  \frac{\# \fftree}{G  \big(\frac{y }{\mathtt{HS}(\fftree)}  \big)} \leqo \frac{c_\epp^2(1+ \epp^{-2}) b_n}{G_0( y_n / \log_2 (c_\epp b_n)) } \leqo \epp .
\end{equation} 
Here, note that we use again the bound $\mathtt{HS}(\fftree) \leqo \log_2 (2\# t)$ which follows from (\ref{connectHS}) in Lemma \ref{HSloppconn} and we recall that $K_\epp y_n \eqo \eta \sqrt{\lambda_n}$, where $K_\epp\eqo 
K(c_\epp, \epp, \frac{1}{2\epp}, G_0)$, which is the constant $K\eqo  K(C_1,C_2,C_3,G)$ in Theorem \ref{brwcouplingth}. 
\end{compactenum}

\smallskip

We first prove Lemma \ref{coupleffec} $(i)$. We recall from (\ref{Snudef}) the definition of $(S_{n,u} (t))_{u\in t}$ 
and from (\ref{RRbulletdef}) the definition of $(R_{n,u} )_{u\in \btt_n}$. 
For all $\epp \ino (0, \frac{_1}{^4})$, all $n\geqo n_{\epp, \eta}$ and all $t\ino \bbT_{\! \epp, n}$, 
we denote by $\smash{W^{0, (n)} (t)}$ (resp.~$\smash{W^{1, (n)} (t)}$) the snake associated with $\smash{(S_{n,u} (t))_{u\in t}}$ 
(resp.~$\smash{(\widetilde{R}_{n,u} (t))_{u\in t}}$)
which is normalized as in (\ref{renormCW}). Note that $\smash{(S_{n,u} (t))_{u\in t}}$ (resp.~$\smash{(\widetilde{R}_{n,u} (t))_{u\in t}}$) has the same law as $\smash{(S_{n,u}^0)_{u\in \btt_n}}$ (resp.~as $\smash{(R_{n,u} )_{u\in \btt_n}}$) under $\smash{\bP (\, \cdot \, | \, \btt_n \eqo t)}$. 
For $i= 0$ or $1$ and for all integer $k \geqo 3$, we first a.s.~get the following 
$$ \bP \big(\omega_{N,p} \big(W^{i,(n)}\big) \geko k\eta \, \big| \, \btt_n\big) \leq \un_{\{ \btt_n \notin \bbT_{\! \epp, n}\}}+ \sum_{t \in \bbT_{\! \epp, n}} \un_{\{\btt_n= t \}} \bP \big(\omega_{N,p} \big(W^{i,(n)} (t)\big) \geko k\eta\big). $$
Let $t\ino \bbT_{\! \epp, n}$. We recall from (\ref{anticouplingevent}) the definition of the event  
$A_{n,\epp, \eta, t}$. By (\ref{coucouplpl}) we then get the following inequalities: 
$\smash{ \bP (\omega_{N,p} (W^{i,(n)} (t)) \geko k\eta)} $  $\leqo$  $\epp $ $+ $  $\smash{\bP (\omega_{N,p} (W^{i,(n)} (t)) \geko k\eta\,  ; \,  A_{n, \epp, \eta, t}^c )} $ $\leqo$ $\epp$ $+$ $\smash{ \bP (\omega_{N,p} (W^{1-i,(n)} (t)) \geko (k\! -\! 2)\eta )}$. Thus 
\begin{eqnarray*}
\bP \big(\omega_{N,p} \big(W^{i,(n)}\big) \geko k\eta\,  \big| \,  \btt_n\big)   \!\!\! &\leq & \!\!\! \un_{\{ \btt_n \notin \bbT_{\! \epp, n}\}} + \epp + \sum_{t \in \bbT_{\! \epp, n}} \un_{\{\btt_n= t \}}   \bP \big(\omega_{N,p} \big(W^{1-i,(n)} (t)\big) \geko (k\! -\! 2)\eta \big) \\
\!\!\! & \leq &  \!\!\!   \un_{\{ \btt_n \notin \bbT_{\! \epp, n}\}} + \epp +\bP \big(\omega_{N,p} \big(W^{1-i,(n)}\big) \geko (k\!-\! 2) \eta \, \big| \, \btt_n\big) .
\end{eqnarray*} 
It implies Lemma \ref{coupleffec} $(i)$ by taking the expectation, since 
$\bP(\btt_n \notin \bbT_{\! \epp, n})\leqo \theta_\epp \! -\! \epp$ for all $n\geqo n_{\epp, \eta}$.  \cq 

\smallskip

Lemma \ref{coupleffec} $(ii)$ is derived from (\ref{coucouplpl}) in a quite similar way: we keep the previous notation 
and we first observe that a.s.
$$ \bP \Big( \!\!\!\!\! \sup_{ \quad s\in \bbR_+} \!\!\!\!   \big| \widehat{W}^{1, (n)}_s \big| \geko z+ \eta  \,  \Big|\,  \btt_n \Big)\, 
 \leq \, \un_{\{ \btt_n \notin \bbT_{\! \epp, n}\}} +
 \sum_{t \in \bbT_{\! \epp, n}} \un_{\{\btt_n= t \}} \,  \bP \Big( \!\!\!\!\! \sup_{ \quad s\in \bbR_+} \!\!\!\!   \big|\widehat{W}^{1,(n)}_s (t)  \big| 
 \geko z+ \eta \Big)\; .$$
Let $t\ino \bbT_{\! \epp, n}$. 
By (\ref{coucouplpl}) we get the following inequalities:  
$\smash{\! \bP( \sup_{ s\in \bbR_+} \! |\widehat{W}^{1,(n)}_s (t) | 
 \geko z\! +\!  \eta )}$ $ \leqo$  $\epp$ $ \! + \!$ $\smash{  \bP ( \sup_{ s\in \bbR_+} \! |\widehat{W}^{1,(n)}_s (t) | 
 \geko z\! + \! \eta ;  A_{n, \epp, \eta, t}^c ) }$ $\leqo$  $\epp$ $+$  $\smash{\bP ( \sup_{ s\in \bbR_+}  |\widehat{W}^{0,(n)}_s (t)  |  \geko z )}$. This a.s.~implies 
\begin{eqnarray*}
\bP \Big( \!\!\!\!\! \sup_{ \quad s\in \bbR_+} \!\!\!\!   \big| \widehat{W}^{1, (n)}_s \big| \geko z+ \eta  \,  \Big|\,  \btt_n \Big)\, 
 \!\!\! &\leq & \!\!\!   \un_{\{ \btt_n \notin \bbT_{\! \epp, n}\}} + \epp + 
 \sum_{t \in \bbT_{\! \epp, n}} \un_{\{\btt_n= t \}} \,  \bP \Big( \!\!\!\!\! \sup_{ \quad s\in \bbR_+} \!\!\!\!   \big|\widehat{W}^{0,(n)}_s (t)  \big| 
 \geko z\Big) \\
  \!\!\! &\leq & \!\!\!   \un_{\{ \btt_n \notin \bbT_{\! \epp, n}\}} + \epp + \bP \Big( \!\!\!\!\! \sup_{ \quad s\in \bbR_+} \!\!\!\!   \big| \widehat{W}^{0, (n)}_s \big| \geko z  \,  \Big|\,  \btt_n \Big)\, ,
 \end{eqnarray*}  
which implies Lemma \ref{coupleffec} $(ii)$ since $\bP(\btt_n \notin \bbT_{\! \epp, n})\leqo \theta_\epp -\epp$ for all $n\geqo n_{\epp, \eta}$. \cq 

\smallskip

To prove $(iii)$ and $(iv)$ we rely on Lemma \ref{Gaussjumps}. More precisely, recall from (\ref{RRbulletdef}) the definition of 
$(R^\bullet_{n,u})_{u\in \btt_n}$. We suppose that $n\geqo n_{\epp, \eta}$
If $\btt_n\ino \bbT_{n, \epp}$, then 
$\fbeta_n \eqo s_n (\btt_n) \geqo 2\epp \geko \epp^4\geqo D_n(\btt_n)\eqo \bdelta_n$ and we can apply 
(\ref{Gausscontrol2}) in Lemma \ref{Gaussjumps} with $\smash{z_1\eqo z_2 \eqo \eta\sqrt{\lambda_n}/ ( 4\epp^{7/2})}$. 
Namely, a.s.~on the event $\{ \btt_n \ino \bbT_{n, \epp} \}$, 
$$ \bP \Big( \max_{u\in \btt_n} |R^\bullet_{n,u}| \geqo \tfrac{\eta \sqrt{\lambda_n}}{2\epp^{7/2}} \, \Big| \, \btt_n \Big) \leqo  
4\bP  \Big( \max_{u\in \btt_n} |R_{n,u}| \geqo \tfrac{\sqrt{\fbeta_{\! n} \! - \bdelta_n}\,  \eta \sqrt{\lambda_n}}{4\epp^{7/2}} \, 
\Big| \, \btt_n \Big) 
+  4\exp \Big(\! -\! \tfrac{\eta^2 \lambda_n (\fbeta_{\! n}  - \bdelta_n)}{{ 32}\epp^7\Gamma (\btt_n) (\fbeta_{\! n}  + 
\bdelta_n)} \Big). $$
If $\btt_n \ino \bbT_{\! \epp, n}$, then 
$\fbeta_{\! n} \! -\! \bdelta_n \geqo \epp$, $\fbeta_n + \bdelta_n \leq (2\epp)^{-1}$ and 
$\lambda_n/ \Gamma (\btt_n) \geqo \epp$. Thus 
$$   \bP \Big( \max_{u\in \btt_n} |R^\bullet_{n,u}| \geqo \tfrac{\eta \sqrt{\lambda_n}}{2\epp^{7/2}} \Big) \leq \,  \theta_\epp + 
4\bP  \Big( \max_{u\in \btt_n} |R_{n,u}| \geqo \tfrac{ \eta \sqrt{\lambda_n}}{4\epp^{3}} \Big) 
+  4e^{-\tfrac{\eta^2 }{16\epp^4} }, $$
since $\bP (\btt_n \! \notin \! \bbT_{\! \epp, n})\leqo \theta_\epp-\epp\leqo \theta_\epp$, for all $n\geqo n_{\epp, \eta}$. We get Lemma \ref{coupleffec} $(iii)$, since 
$\smash{\max_{u\in \btt_n} |R_{n,u}|}$ $ \eqo $ $ \smash{\sqrt{\lambda_n} \sup_{s\in \bbR_+} |  \widehat{W}^{1,(n)}_s |}$ and $\smash{\max_{u\in \btt_n} |R^\bullet_{n,u}| \eqo \sqrt{\lambda_n}  \sup_{s\in \bbR_+} |  \widehat{W}^{\bullet,(n)}_s |}$.  \cq 

\smallskip

Let us prove Lemma \ref{coupleffec} $(iv)$. We suppose that $n\geqo n_{\epp, \eta}$. 
We apply (\ref{Gausscontrol}) in Lemma  \ref{Gaussjumps} with $\smash{z_1\eqo z_2 \eqo \tfrac{1}{2} \eta\sqrt{\lambda_n}}$ 
and we a.s.~get 
 $$ \bP \Big( \max_{u\in \btt_n} |R_{n,u}\! -\! \sqrt{ \fbeta_{\! n}} 
 R^\bullet_{n,u}| \geqo \eta \sqrt{\lambda_n}\, \Big| \, \btt_n \Big) \leqo 4\bP \Big( \max_{u\in \btt_n} 
 |R^\bullet_{n,u}| \geqo \tfrac{\eta \sqrt{\lambda_n \fbeta_n}}{2\bdelta_n} \, \Big| \, \btt_n \Big) 
+  4\exp \Big(\! -\! \tfrac{\eta^2 \lambda_n\fbeta_n }{8\bdelta_n^2\Gamma (\btt_n) } \Big) .$$
If $\btt_n \ino \bbT_{\! \epp, n}$, then $\bdelta_n \eqo D_n (\btt_n) \leqo \epp^4$, $\fbeta_n \geqo \epp$, and 
$\lambda_n / \Gamma (\btt_n) \geqo \epp$. Since $\bP (\btt_n \! \notin \! \bbT_{\! \epp, n})\leqo \theta_\epp-\epp \leqo \theta_\epp$ 
for all $n\geqo n_{\epp, \eta}$, we get 
$$ \bP \Big( \max_{u\in \btt_n} |R_{n,u}\! -\! \sqrt{ \fbeta_{\! n}}
 R^\bullet_{n,u}| \geqo \eta \sqrt{\lambda_n} \Big) \leqo \theta_\epp + 4\bP \Big( \max_{u\in \btt_n} 
 |R^\bullet_{n,u}| \geqo \tfrac{\eta \sqrt{\lambda_n }}{2\epp^{7/2}} \Big) 
+  4\exp \Big(\! -\! \tfrac{\eta^2  }{8\epp^6 } \Big) .$$
Note that 
$\smash{ \max_{u\in \btt_n} |R_{n,u}\! -\! \sqrt{\fbeta_{\! n}} R^\bullet_{n,u}| \eqo \sqrt{\lambda_n} \sup_{s\in \bbR_+}
 | \widehat{W}^{1,(n)}_s \!\! -\! \sqrt{\fbeta_n}\widehat{W}^{\bullet,(n)}_s |}$. Thus 
\begin{equation}
\label{gaussdelta} \bP (B_{n})  \leqo  \theta_\epp + 
4\bP \Big( \sup_{s\in \bbR_+} |\widehat{W}^{\bullet, (n)}_s | \geqo \tfrac{\eta}{2\epp^{7/2}}  \Big) 
+  4\exp \Big(\! -\! \tfrac{\eta^2  }{8\epp^6} \Big) \leqo q_{n, \eta, \epp},  
\end{equation}
where $\smash{B_{n}\eqo \big\{ \! \sup_{s\in \bbR_+} | \widehat{W}^{1,(n)}_s \!\! -\! \sqrt{\fbeta_n}\widehat{W}^{\bullet,(n)}_s | \geqo 
\eta \big\}}$. Then, for all $N, p\ino \bbN$ and all $n\geqo n_{\epp, \eta}$, 
\begin{eqnarray*}\bP \big( \sqrt{\fbeta_{\! n}}\omega_{N,p} \big(W^{\bullet, (n)} \big)\!  \geko 5\eta \big) & \leqo &  \bP (B_{n})   +  \bP \big( \omega_{N,p} \big(W^{1, (n)} \big) \! \geko 3\eta \big) \\
\textrm{and} \; \bP \big( \omega_{N,p} \big(W^{1, (n)} \big) \! \geko 3\eta \big) & \leqo &   \bP (B_{n})  +\bP \big(\sqrt{\fbeta_{\! n}} \omega_{N,p} \big(W^{\bullet, (n)} \big)\!  \geko  \eta \big)  , 
\end{eqnarray*}
which easily implies Lemma \ref{coupleffec} $(iv)$ by (\ref{gaussdelta}) and since $\smash{\bP (\fbeta_{\! n} \! \notin \! [2\epp, \tfrac{1}{2\epp} ]) \leqo \theta_\epp }$.  \cq 

\smallskip

To prove Lemma \ref{coupleffec} $(v)$, we first assume that 
the laws of $\smash{(\widehat{W}^{{0, (n)}}_{\cdot} )_{n\in \bbN}}$ are tight in $\bC^{_0}_{^1}$. By $(i)$ and standard argument (see e.g.~Billingsley \cite{Bil68} Thm 7.3 p.~82), the laws of 
$\smash{\widehat{W}^{1, (n)}_\cdot }$, $n\ino \bbN$, are tight in $\bC^{_0}_{^1}$ too. The tightness of the laws of $\smash{\widehat{W}^{{0, (n)}}_{\cdot} }$, $n\ino \bbN$, then implies for all $s\ino \bbR_+$ that $\smash{\lim_{z\to \infty }\limsup_{n\to \infty} \bP \big( \sup_{r\in [0, s]}  | \widehat{W}^{0, (n)}_r |  \geko z \big) \eqo 0}$. Thus,  
$$\limsup_{z\to \infty} \limsup_{n\to \infty} 
\bP \Big( \!\!\!\! \sup_{\quad  r\in \bbR_+} \!\!\!\!  \big|  \widehat{W}^{0, (n)}_r \big|   \geko z \Big) \leq \limsup_{n\to \infty} \bP \big(  \# \btt_n \geqo \frac{_{_1}}{^{^2}}b_n s \big) \xrightarrow[s\to \infty]{\; } 0 \; .$$

\noi
By Lemma \ref{coupleffec} $(ii)$ we get $\smash{\limsup_{z\to \infty} \limsup_{n\to \infty} 
\bP \big( \sup_{ r\in \bbR_+}  | \widehat{W}^{1, (n)}_r |  \geko z \big) \eqo 0}$. We fix $\eta$. By 
Lemma \ref{coupleffec} $(iii)$ we get $  \limsup_{\epp \to 0} \limsup_{n\to \infty} q_{n,  \eta \sqrt{\epp}, \epp}\eqo 0$. Then the first inequality in Lemma \ref{coupleffec} $(iv)$, with $\eta$ replaced by $ \eta\sqrt{\epp}$, entails 
$$ \limsup_{p\to \infty} \limsup_{n\to \infty} \bP \big(  \omega_{N,p} (W^{\bullet, (n)} )\!  \geko \tfrac{5}{\sqrt{2}}\eta \big)  \leqo     \limsup_{n\to \infty} 2q_{n, \eta \sqrt{\epp}, \epp} \xrightarrow[\epp \to 0]{\! } 0, $$
\noi
which implies that the laws of $\smash{\widehat{W}^{\bullet, (n)}_\cdot}$, $n\ino \bbN$, are tight in $\bC^{_0}_{^1}$. 

We now assume that the laws of $\smash{\widehat{W}^{\bullet, (n)}_\cdot}$, $n\ino \bbN$, are tight in $\bC^{_0}_{^1}$. 
We argue as previously to get $\smash{ \limsup_{z\to \infty} \limsup_{n\to \infty} 
\bP \big( \sup_{ r\in \bbR_+}  | \widehat{W}^{\bullet, (n)}_r |  \geko z \big)\eqo 0}$. Thus, for all $\eta \ino (0, \frac{1}{4})$, 
it entails 
$\lim_{\epp\to 0} \limsup_{n\to \infty} q_{n, \eta, \epp}\eqo 0$. 
The second inequality in Lemma \ref{coupleffec} $(iv)$ implies that the laws of $\smash{\widehat{W}^{1, (n)}_\cdot}$, 
$n\ino \bbN$, are tight in $\bC^{_0}_{^1}$ and the second inequality in Lemma \ref{coupleffec} $(i)$ finally implies that the laws of 
$\smash{\widehat{W}^{0, (n)}_\cdot}$, $n\ino \bbN$, are tight in $\bC^{_0}_{^1}$. This completes the proof of the lemma. \cqfd

\medskip

We proceed to the proof of Theorem \ref{extenscv} in several steps.  

\smallskip

\noi
$(\textbf{I})$ First note that Lemma \ref{coupleffec} applies 
to $\bS^0_n\eqo \bS^*_n$ and therefore implies that  
the laws of $\smash{(\widehat{W}^{\bullet, (n)}_\cdot)_{n\in \bbN}}$ are tight in $\bC^{_0}_{^1}$.

\smallskip

\noi
$(\textbf{II})$ By Proposition \ref{endsnake1} $(i)$ the laws of $(W^{_{\bullet, (n)}}_\cdot)_{n\in \bbN}$ are tight in $\bC (\bbR_+, \bC^{_0}_{^1})$.
Thus, the joint laws of $(C^{_{(n)}}_\cdot  , W^{_{\bullet, (n)}}_\cdot)$ are tight. By 
Lemma \ref{fdcvsna}, the finite dimensional marginals laws of these processes weakly converge to those  
of a $\bbR$-valued Brownian snake with lifetime process $C$. Consequently, there is a continuous version of the $\bbR$-valued Brownian snake with lifetime process $C$. This obviously implies that 
there is a continuous version of the $\bbR^d$-valued Brownian snake with lifetime process $C$ that we keep denoting $W$.

\smallskip

\noi
$(\textbf{III})$ We next suppose that the $\bS_n$ are $\bbR$-valued and that $\bP$-a.s.~$\fbeta \ino \bbR_+^*$. Lemma \ref{coupleffec} $(v)$ applies 
to $\bS^0_n\eqo \bS_n$ and $(\textbf{I})$ implies that the laws of $\smash{(\widehat{W}^{(n)}_\cdot)_{n\in \bbN}}$ are tight in $\bC^{_0}_{^1}$.

\smallskip

\noi
$(\textbf{IV})$ We now only suppose that the $\bS_n$ are $\bbR$-valued. We set $\smash{\bS_n'\eqo \bS_n + \bR^\bullet_n}$ whose snake $\smash{W^{_{'}, (n)}}$ is normalized as in (\ref{renormCW}). By Lemma \ref{parasubaddgauge}, $\smash{\bS_n'}$ satisfies \textbf{Sib-Ind}$_{_{\, }}$($\baa, \bbb, \bS'_\cdot, G, C_\cdot, 1+\fbeta$). Thus, $(\textbf{III})$ applies to $\smash{\bS_n'}$ and the laws of $\smash{(\widehat{W}^{_{'}, (n)}_\cdot)_{n\in \bbN}}$ are tight in $\bC^{_0}_{^1}$. Since the laws of $\smash{(\widehat{W}^{\bullet, (n)}_\cdot)_{n\in \bbN}}$ are tight in $\bC^{_0}_{^1}$ and since $\smash{\widehat{W}^{(n)}_\cdot\eqo \widehat{W}^{_{'}, (n)}_\cdot\! -\! \widehat{W}^{\bullet, (n)}_\cdot }$, the laws of $\smash{(\widehat{W}^{(n)}_\cdot)_{n\in \bbN}}$ are tight in $\bC^{_0}_{^1}$.

\smallskip

\noi
$(\textbf{V})$ We now turn to the general case where the $\bS_n$ are $\bbR^d$-valued. 
We fix $j\in \{ 1, \ldots, d\}$. Observe that the $\langle \mathtt e_j, \bS_n \rangle$ are $\bbR$-valued BRWs satisfying \textbf{Sib-Ind}$_{_{\, }}$($\smash{\baa, \bbb, \bS_\cdot, G, C_\cdot, \langle \mathtt e_j ,\fbeta.\mathtt e_j \rangle}$). Therefore we apply $(\mathbf{IV})$ which proves that the laws of $\smash{(\langle \mathtt e_j , \widehat{W}^{(n)}_\cdot \rangle)_{n\in \bbN}}$ are tight in $\bC^{_0}_{^1}$ and since it holds for all $j\ino \{1, \ldots, d\}$, the laws of $\smash{(\widehat{W}^{(n)}_\cdot)_{n\in \bbN}}$ are tight in $\bC^{_0}_{^d}$ by standard arguments. 

\smallskip

\noi
$(\textbf{VI})$ We now apply Lemma \ref{fdcvsnakeSibInd} to $\bS^{\mathtt a}_n \eqo \bS_n$ which implies that the finite dimensional marginal laws of the $\smash{(C^{(n)}_\cdot, \widehat{W}^{(n)}_\cdot)_{n\in \bbN}}$ converge to those of $\smash{(C_\cdot, \sqrt{\fbeta}.\widehat{W}_\cdot)}$. Combined with $(\mathbf{V})$, it implies the weak convergence $\smash{(C^{(n)}_\cdot\! , \widehat{W}^{(n)}_\cdot) \! \to \! (C_\cdot, \sqrt{\fbeta}.\widehat{W}_\cdot)}$ in $\bC^{_0}_{^1} \times \bC^{_0}_{^d}$.

\smallskip

\noi
$(\textbf{VII})$
By Proposition \ref{endsnake1} $(i)$ the laws of $\smash{(W^{_{(n)}}_\cdot)_{n\in \bbN}}$ are tight in $\bC^{_0}(\bbR_+, \bC^{_0}_{^1})$.
Thus, the joint laws of $\smash{(C^{_{(n)}}_\cdot \!  , W^{_{(n)}}_\cdot)}$ are tight too. Then $(\textbf{VI})$ and Lemma \ref{endmargcar} show that 
$\smash{(C_\cdot ,\sqrt{\fbeta}. W_\cdot)}$ is the only weak limit of the $\smash{(C^{_{(n)}}_\cdot \!  , W^{_{(n)}}_\cdot)}$, which completes the proof of Theorem \ref{extenscv}.  \cqfd

\subsection{Proof of Theorem \ref{maincvsnake}}
\label{Thmmainsnapfsec}

Theorem \ref{maincvsnake} is a consequence of Theorem \ref{Sheuexplain} $(iii)$ and of Theorem \ref{extenscv}. 
Let us fix $x,c\ino \bbR_+^*$ and let $\psi\ino \mathscr L$ satisfy \texttt{Var}$_{\infty}$($\psi$). Let $(a_n)_{n\in \bbN}$, 
$(b_n)_{n\in \bbN}$, $(\mu_n)_{n\in \bbN}$, $\bS_n \eqo (S_{n, u})_{u\in \btt_n}$ and $(C,Y,W)$ be as in 
$\textbf{Case (i)}$, $\mathbf{i}\ino \{ \mathbf{1}, \mathbf{2}, \mathbf{3} \}$. 
We recall that $C^{_{(n)}}_{\cdot}$, $V^{_{(n)}}_{\cdot}$ and $W^{_{(n)}}_{\cdot}$ are rescaled as in (\ref{renormCW}).  
We make the following assumptions. 
\begin{compactenum}

\smallskip

\item[$(a)$] We assume \texttt{\L{}uka}$_{^{\,}}$($\baa, \bbb, \bmu, \psi$) and \texttt{Sheu}$_{^{\,}}$($\baa, \bbb, \bmu $).

\smallskip

\item[$(b)$] We suppose that there is $\fbeta\! : \! \Omega \! \to \! \mathrm{Sym}^{_+}_{^d}$ and a moment 
gauge function $G$ as in (\ref{momgaudef}) such that the $(\bS_n)_{n\in \bbN}$ satisfy Assumptions 
\textbf{Sib-Ind}$_{_{\, }}$($\baa, \bbb, \bS_\cdot, G, C_\cdot, \fbeta$) $(\mathbf{1}$-$\mathbf{5})$, 
which are denoted by \textbf{Sib-Ind}  $(\mathbf{1}$-$\mathbf{5})$ below to simplify notation. 

\smallskip

\end{compactenum}

By Remark \ref{afterthm3} $(b)$, \textbf{Sib-Ind} $(\mathbf{6})$ is satisfied. By Lemma \ref{parasubaddgauge}, 
there are $b, c\ino \bbR_+$ such that $\smash{G(x) \leqo cG(0)e^{b|x|}\! =:\! G_*(x)}$, $x\ino \bbR$. 
We set $\smash{\beta_*\eqo 2(1+b)^{-2}}$ and  
$\smash{\bgam (dx)  \! :=\!  (2\beta_*)^{\!-\frac{1}{2}} \! e^{-|x| (2/\beta_*)^{1/2}} \! \! dx}$, which is a centered density such that 
$\smash{\int_{\bbR} x^2 \bgam (dx) \eqo \beta_*}$. Let $\smash{\bS^*_n \eqo (S^*_{n, u})_{u\in \btt_n}}$, $n\ino \bbN$, $\bbR$-valued BRWs whose jumps, conditionally given 
$\btt_n$, are i.i.d.~with density $\bgam$. Then, Lemma \ref{parasubaddgauge} $(ii)$ asserts that $G_*$ is 
a moment gauge function and that the assumptions 
\textbf{Sib-Ind}$_{_{\, }}$($\baa, \bbb, \bS^*_\cdot, G_*, C_\cdot, \beta_*$) hold. 

Let $W^{*,(n)}$ be the snake associated with $\bS_n^*$ and rescaled as in (\ref{renormCW}). 
Then the laws of $(\widehat{W}^{_{*,{(n)}}}_\cdot)_{n\in \bbN}$ are tight in $\bC^{_0}_{^1}$. 
\emph{Indeed,} in \textbf{Case (1)} and \textbf{Case (2)}, it is stated in Theorem \ref{Sheuexplain} $(iii)$ and in 
\textbf{Case (3)} it is proved in Janson \& Marckert \cite{JanMar05} for $\alpha\eqo 2$ and in Marzouk \cite{Mar20} 
for $\alpha \ino (1, 2)$. 
Therefore Theorem \ref{extenscv} applies and it implies (\ref{jcvbis}). The joint convergence with 
$\smash{Y^{(n)}}$ follows from standard arguments: we leave the details to the reader. \cqfd

\subsection{Proof of Proposition \ref{excntrex}}
\label{Proofsecexcntrex}
We fix two renormalization sequences 
$\smash{(a_n)_{n\in \bbN}}$ and $\smash{(b_n)_{n\in \bbN}}$ that satisfy \texttt{Norm} ($\baa, \bbb$) and we fix $\smash{(\mu_n)_{n\in \bbN}}$ a sequence of offspring distributions satisfying (\ref{nontricri}). In this section we discuss the relative growth of $b_n$ in terms of $a_n$ and it is convenient to introduce a function $\phi \! :\! \bbR_+ \! \to \bbR_+$ such that 
$$\forall n\ino \bbN, \quad b_n \eqo  \phi (a_n)$$ 
and such that it faithfully interpolates the values of $b_n$ is terms of $a_n$. We want 
$(b_n)_{n\in \bbN}$ to be the true time-renormalization sequence of $V^{(n)}$: by Remark \ref{Lukagrowthrem} it amounts to assuming (\ref{righttime}). Namely, $q\! :=\! \limsup_{n\to \infty} \mu_n (1)\leko  1$. Then under \texttt{\L{}uka}$_{^{\,}}$($\baa, \bbb, \bmu, \psi$), Lemma \ref{Lukagrowth} $(ii)$ asserts that 
$\limsup_{n}$ $b_n/a_n^2$ $\leqo$ $4\beta_\psi / (1 \! -\!  q)$ and Lemma \ref{Lukagrowth} $(iii)$ shows that \texttt{Var}$_{\infty}$($\psi$) implies that $b_n /a_n \to \infty$. 
It is therefore natural to restrict ourselves to functions 
$\phi$ such that $\phi(x)/x\! \to \!  \infty$ and $\phi(x)\eqo \mathcal O_{\! \infty} (x^2)$. 
More precisely, to provide examples it is convenient to pick $\phi$ in the class $\mathscr L$ of Laplace exponents of L\'evy processes without negative jump which have infinite variation sample paths (see Definition \ref{defLaplclass}), which is sufficiently rich, as shown by the folllowing lemma.   
Here $u_n \asymp v_n$ (resp.~$f\asymp_\infty g$) means that there is $c\ino (1, \infty)$ such that $c^{-1} u_n \leqo v_n \leq c u_n$ (resp.~$c^{-1} f(x) \leqo g(x) \leq c f(x)$)  for all sufficiently large $n\ino \bbN$ (resp.~$x\in \bbR_+$). 
\begin{lemma} 
\label{examplphi} The following holds true. 
\begin{compactenum}

\smallskip

\item[$(i)$] Let $\phi \ino \mathscr L$ satisfy \emph{$\texttt{Var}_\infty(\phi)$}. Then $\phi$ is strictly convex, $\phi'$ is concave, $\phi (z) \asymp_\infty z\phi'(z)$, 
$\lim_{z\to \infty}z^{-1} \phi (z)  \eqo \infty$, $\lim_{z\to \infty}z^{-2} \phi (z)\eqo \beta_\phi$. 

\smallskip

\item[$(ii)$] For all $(\gamma, c)\ino \big( \{ 1\} \! \times \! \bbR_+^*\big) \cup \big( (1,2) \! \times \! \bbR\big) \cup \big(
 \{2\} \! \times \! (-\infty,0] \big) $ 
there exists $\phi \ino \mathscr L$ such that $\phi (z) \! \asymp_\infty \! z^\gamma (\log z)^c$ 
and there exists $\phi \ino \mathscr L$ such that $\phi (z) \! \asymp_\infty \! z \log \log z $. 
\end{compactenum}
\end{lemma} 
\noi
\textbf{Proof.} $(i)$ is standard: the details are left to the reader. 
Let us prove $(ii)$. For all 
$(\gamma, a)\ino \bbR^2$ we set $\mathtt{m}_{\gamma, a} (dx)=  x^{-\gamma-1} (\log 1/x)^a  \un_{\{ 0< x< 1/2\}} dx$. Thus $\smash{\int \! x^2\mathtt{m}_{\gamma, a} (dx) \leko \infty}$ and $\smash{\int \! x\,  \mathtt{m}_{\gamma, a} (dx) \eqo \infty}$ iff $(\gamma, a)\ino \big( \{ 1\} \! \times \! [-1, \infty)  \big) \cup\big( (1,2) \! \times \! \bbR \big) \cup\big( \{ 2\} \! \times \! (-\infty, -1)  \big)$. In these cases, we thus define $\phi \ino \mathcal L$ by (\ref{LK}) with $\pi_\phi\eqo  \mathtt{m}_{\gamma, a}$ and $\beta_\phi\eqo 0$. 

Let us first suppose that $(\gamma, a)\ino \big( \{ 1\}  \times  [-1, \infty)  \big) \cup\big( (1,2) \times  \bbR \big) $. By standard arguments in asymptotic analysis combined with Karamata's Theorem 1.7.1 in Bingham, Goldies \& Teugels \cite{BiGoTe} p.~38, we get $\phi''(z) \! \sim_\infty \! \Gamma (2\!-\!  \gamma) z^{-(2-\gamma)} (\log z)^a$ (we leave the details to the reader). If $(\gamma, a)\ino (1,2)  \times  \bbR $, then $\phi'(z) \! \asymp_\infty  \! z^{\gamma -1} (\log z)^a$ and thus $\phi(z)\!  \asymp_\infty \!  z^{\gamma } (\log z)^a$. If $\gamma\eqo 1$ and $a\geko -1$, then $\phi'(z) \! \asymp_\infty \!  (\log z)^{a+1}$ and 
$\phi'(z)\!  \asymp_\infty \!  z (\log z)^{a+1}$. If $\gamma \eqo -a\eqo 1$, then $\phi (z)\!  \asymp_\infty \!  \log \log z$ and thus $\phi(z)\!  \asymp_\infty  \! z\log \log z$. 

   We suppose that $(\gamma, a)\ino \{ 2\} \! \times \! (-\infty, -1) $. Standard arguments in asymptotic analysis imply  
$-\phi'''(z) \! \sim_\infty\!  z^{-1} (\log z)^{-|a|}$. Thus, $\phi''(z) \! \asymp_\infty \! (\log z)^{-(|a|-1)}$, $\phi'(z) \! \asymp_\infty \! z (\log z)^{-(|a|-1)}$ and $\phi(z) \! \asymp_\infty \! z^2 (\log z)^{-(|a|-1)}$. To complete the proof of $(ii)$, we observe that the case where $c\eqo 0$ and $\gamma\eqo 2$ corresponds to the Brownian Laplace exponent $\phi (z)\eqo z^2$. \cqfd

\medskip
 
To prove Proposition \ref{excntrex}, we construct examples that are derived from functions of $\ccL$ as explained in 
the following lemma. 
\begin{lemma}
\label{rigimuphi} Let $\phi \ino \ccL$ satisfy \emph{$\texttt{Var}_\infty(\phi)$}, let $a\ino \bbR_+^*$ and let $q\ino [0, 1)$. Then there exists an 
offspring distribution $(\mu(k))_{k\in \bbN}$ such that $ \mu (0)+ \mu (1)\!  \leko 1$, $\, \sum_{k\in \bbN} k\mu(k)\eqo 1$, 
$\, \mu(1)\eqo q$ and for all $r\ino [0, 1]$, 
\begin{equation}
\label{gmudephi}
g_\mu (r) \eqo \sum_{k\in \bbN} r^k \mu(k) = r+ \frac{1\! -\! q}{a\phi'(a)} \, \phi \big( a(1\! -\! r)\big) \; .
\end{equation}
\end{lemma} 
\noi
\textbf{Proof.} For all integers $k\geqo 2$, we set $u_k \eqo  \beta_\phi a^2 \un_{\{ k=2\}}+ \int \! \pi_\phi (dx) (ax)^ke^{-ax}/k! \eqo (-a)^k\phi^{(k)} (a)/k!$ and we observe that $u_k \geqo 0$. By Fubini, we get for all $r\ino [0, 1)$, 
$$\sum_{k\geq 2} u_k r^k = \beta_\phi (ar)^2+ \int_{\bbR_+^*} \!\!\!\!\! \pi_\phi (dx) e^{-ax}  \big(e^{rax}\! -\! 1\! -\! rax \big) = \phi \big( a(1\! -\! r)\big) \! -\! \phi (a) + ra\phi'(a). 
$$
This extends to $r\eqo 1$ by monotone convergence. We then set $\mu (k)\eqo (1\! -\! q) u_k/(a \phi' (a))$ if $k\geqo 2$, $\mu(1)\eqo q$ and $\mu (0)\eqo (1\! -\! q)\phi (a)/ (a\phi'(a))$ so that 
(\ref{gmudephi}) holds for all $r\ino [0, 1]$. Since $\phi$ is strictly convex, $\phi (a) \leko a\phi'(a)$ and we get $\mu (0) + \mu(1) \leko 1$. We next get $g_\mu'(r)\eqo 1\! -\!  (1\! -\! q) \phi'(a(1\!- \! r)) /\phi'(a)$. Letting $r\uparrow 1$ shows that $g'_\mu(1^-)\eqo 1$, which proves that $\mu$ is critical, which finishes the proof. \cqfd

\medskip
 
The following lemma proves that offspring distributions derived from Laplace exponents as in Lemma \ref{rigimuphi} behave regularly. 
 
 \begin{lemma}
\label{muphiregu} Let $\psi \ino \ccL$ satisfy \emph{$\texttt{Var}_\infty(\psi)$}. For all $n\ino \bbN$, let $F_{\! n} \ino \ccL $, 
$a_n \ino \bbR_+^*$, $q_n \ino [0, 1)$ and let $\mu_n$ be the non-trivial critical offspring distribution such that 
\begin{equation}
\label{Fnmun}
g_{\mu_n} (r) \eqo \sum_{k\in \bbN} r^k \mu_n(k) = r+ \frac{1\! -\! q_n}{a_nF_{\! n}'(a_n)} \, F_{\! n} \big( a_n(1\! -\! r)\big) \; .
\end{equation}  
We also set $b_n \eqo a_n F_{\! n}'(a_n) / (1\! -\! q_n)$. Then the following holds true. 
\begin{compactenum}

\smallskip

\item[$(i)$] For all $z\ino [0, a_n]$, $\psi_n (z)\eqo  b_n \big(g_{\mu_n} \big( 1\! -\! \frac{z}{a_n}\big)\! -\! 1+\frac{z}{a_n}\big)\eqo F_{\! n} (z)$. 

\smallskip

\item[$(ii)$] If $a_n \! \to \! \infty$ and if for all $z\ino \bbR_+$, $\lim_{n\to \infty} F_{\! n} (z)\eqo \psi (z)$, then  \emph{\texttt{\L{}uka}$_{^{\,}}$($\baa, \bbb, \bmu, \psi$)} holds true. 
\end{compactenum}
\end{lemma} 
\noi
\textbf{Proof.} Note that $(i)$ is an immediate consequence of the definitions. To prove $(ii)$ we introduce $(V^{_{(n)}}_{^{\! k}})_{k\in \bbN}$ the RW whose jump distribution is $\mu_n (1+ \cdot)$ on $\{ -1\} \cup \bbN$. 
Then for all $z\ino \bbR_+$, we get 
$$ \bE \Big[ e^{-zV^{(n)}_{\lfloor b_n \rfloor }/a_n }\Big]\eqo  \big( e^{z/a_n} g_{\mu_n} \big( e^{-z/a_n}\big)\big)^{\lfloor b_n \rfloor} \eqo \exp \Big(\lfloor b_n \rfloor \log \Big(1+ \tfrac{1}{b_n}e^{z/a_n}F_{\! n} \big( a_n \big( 1\! - \! e^{-z/a_n}\big) \big) \Big)  \Big). $$
By the second Dini Theorem, the $F_{\! n}$ converge to $\psi$  uniformly on every compact subsets of $\bbR_+$. Thus   
$\smash{e^{z/a_n}F_{\! n} \big( a_n \big( 1\! - \! e^{-z/a_n}\big) \big)\! \to \! \psi (z)}$ and  
$\smash{\lim_{n\to \infty} \bE \big[ \exp (-zV^{_{(n)}}_{^{\lfloor b_n \rfloor} }/a_n ) ]\eqo \exp (\psi (z))}$, for all $z\ino \bbR_+$. This entails \texttt{\L{}uka}$_{^{\,}}$($\baa, \bbb, \bmu, \psi$) by Lemma A.3 in Broutin, D.~\& Wang \cite{BrDuWa21}.   
\cqfd

\medskip

We then rewrite Assumptions \texttt{Grey}$_{^{\,}}$($\psi$)  and \texttt{Sheu}$_{^{\,}}$($\psi$) in the following way.   
\begin{lemma}
\label{sheu1integ} Let $h\! : \! \bbR_+ \! \to \! \bbR_+$ be $C^1$, convex, increasing and such that $h(0)\eqo 0$. Then for all positive real numbers $z_1\leko z_2$, 
\begin{eqnarray}
\label{GreySheunew} 
\int_{z_1}^{z_2} \!\!\!\!  \frac{ds}{sh'(s)} \!\!\! &  \leqo &  \!\!\!   \int_{z_1}^{z_2} \!\!\! \! \frac{ds}{h(s)} \leq  2\!  \int_{\frac{_1}{^2}z_1}^{\frac{_1}{^2}z_2} \!\!\! \!\! \frac{ds}{\, sh'(s)} \qquad \textrm{and} \; \\
& & \sqrt{2} \! \int_{z_1}^{z_2} \!\!\!\!  \frac{ds}{s\sqrt{h'(s)\,}} \leq  
\int_{z_1}^{z_2} \!\!\! \! \frac{ds}{\sqrt{\int_0^s h(r) dr}} \leq \frac{4}{\sqrt{3}} \int_{\frac{_1}{^4}z_1}^{\frac{_1}{^4}z_2} \!\!\! \!\frac{ds}{\, s\sqrt{h'(s)\, }} \nonumber
\end{eqnarray}

\end{lemma}
\noi
\textbf{Proof.} Since $h$ is $C^1$ and convex, $\frac{_1}{^2} rh'(r/2) \leq \frac{_1}{^2} r (h(r)\! -\! h(r/2))/(r/2) \leqo h(r) \eqo r(h(r)/r) \leqo rh'(r)$, which entails the first inequalities in (\ref{GreySheunew}). 
This also implies the following: $\frac{3}{{16}} s^2h' ( \frac{_1}{^4}s )\leq   \int_{\frac{_1}{^2}s}^s \frac{1}{2} rh' ( \frac{_1}{^2}r )  dr \leq  \int_0^s h(r) dr \leq  \int_0^s rh'(r) dr \leq \tfrac{1}{2} s^2h'(s) $, which easily entails (\ref{GreySheunew}). \cqfd 

\subsection*{Proof of Proposition \ref{excntrex}} 
For all $y\ino \bbR$, we denote by $\delta_y (dx)$ the Dirac mass at $y$. We fix $\psi, \phi \ino \ccL$ satisfying \texttt{Var}$_{\infty}$($\psi$) and \texttt{Var}$_{\infty}$($\phi$) and  $\beta_\phi \leqo  \beta_\psi$. 
For all $\epp , \eta \ino \bbR_+^*$, we  set 
$\mathtt{m}_{\epp, \eta} (dx)\eqo \un_{(0, \epp)} (x) \pi_\phi (dx)$ $+$ $ \un_{[\eta , \infty)} (x) \pi_\psi (dx) $ $+$  $2\eta^{-2} (\beta_\psi \! -\! \beta_\phi)  \delta_{\eta} (dx) $ and 
\begin{equation}
\label{Fepsi}
\forall z\in \bbR_+,  \quad F_{\! \epp, \eta} (z) \eqo \beta_\phi z^2+ \int_{\bbR_+^*} \!\!   \!\!   \!\!  \mathtt{m}_{\epp, \eta}  (dx)\big( e^{-zx} \! -\! 1+ zx\big).
\end{equation}
Clearly $F_{\! \epp, \eta} \ino \ccL$, it satisfies \texttt{Var}$_\infty(F_{\! \epp, \eta})$ and $F'_{\! \epp, \eta} \eqo g_\epp+ f_\eta$ where 
\begin{eqnarray}
\forall z\ino \bbR_+, \quad g_\epp (z)= 2\beta_\phi z \!\!\!\! &+ &\!\!\!\!\!  \int_{(0, \epp)} \!\!\!\!  \pi_\phi (dx) x \big( 1 \! -\! e^{-zx} \big) \nonumber \\
\label{fgepsiloneta}\textrm{and}  & & \!\!\!\! \!\!\!\! \!\!\!\!  f_\eta (z)= \frac{_2}{^{\eta}} (\beta_\psi \! -\! \beta_\phi)  \big(1\! -\! e^{-\eta z} \big) + 
\int_{[\eta , \infty )}  \!\! \!\! \!\! \!\! \!\!  \pi_\psi (dx) \, x \big( 1 \! -\! e^{-zx} \big). 
\end{eqnarray}
We observe for all $z\ino \bbR_+$ that 
\begin{eqnarray}
\label{boundgf}
0\leq \phi'(z) \! -\!  g_\epp (z) \!\!\!\!\!   &\leq & \!\!\!\!\!     \int_{[\epp, \infty)}  \! \!\!\!\!\!  \!    \pi_\phi (dx) x  \quad \textrm{and that} \quad\\
0  \!\!\!\!  &\leq & \!\!\!\!   \psi'(z) - f_\eta (z)   \leq  
  \beta_\psi \eta z^2+ 2\beta_\phi z +  z\! \int_{(0, \eta)} \!\! \!\! \!\! \!\!  \pi_\psi (dx) \,  x^2  \; .\nonumber
\end{eqnarray}
We also observe for all $z\in \bbR_+$ that 
\begin{equation}
\label{sandwichgf}
\phi'(z) - \! \int_{[\epp, \infty)}  \! \!\!\!\!\!  \!    \pi_\phi (dx) x \, \leq \, F_{\! \epp, \eta}' (z) \, \leq \, \phi'(z) + 
 \frac{_2}{^{\eta}} (\beta_\psi \! -\! \beta_\phi)  +
\int_{[\eta , \infty )}  \!\! \!\! \!\! \!\! \! \pi_\psi (dx) x \; .
\end{equation}

We recall that $(a_n)_{n\in \bbN}$ is a fixed sequence of positive real numbers which converges to $\infty$. 
 Let $(a'_n)_{n\in \bbN}$ be such that $\lim_{n\to \infty} a_n'\eqo \infty$ and $a_n'\leko a_n$. 
We specify $(a'_n)_{n\in \bbN}$ later. Let us mention here that since $\phi$ satisfies \texttt{Var}$_\infty(\phi)$, $\lim_{z\to \infty} \phi'(z)\eqo \infty$. 
Thus, it is always possible to find two sequences of positive real numbers $(\eta_n)_{n\in \bbN}$ and 
$(z_n)_{n\in \bbN}$ such that  $\eta_n \! \to \! 0$, $z_n \! \to \! \infty$, 
 \begin{equation}
\label{constraint1}
a'_n > \frac{16}{\eta_n}> 16z^2_n ,  \quad \phi'(a'_n) \geq  \frac{_2}{^{\eta_n}} (\beta_\psi \! -\! \beta_\phi)  +\! 
\int_{[\eta_n , \infty )}  \!\! \!\! \!\! \!\! \! \! \pi_\psi (dx) x \quad \textrm{and} \quad z_n  \int_{(0, \eta_n)}  \! \!\!\!\!\!  \!  \!\!\!\!    \pi_\psi (dx) x^2 \leq 1 .
\end{equation}
Furthermore, it is also always possible to find $(\epp_n)_{n\in \bbN}$ such that $\epp_n \! \to \! 0$ and 
\begin{equation}
\label{constraint2}
\phi'(z_n) \geq 2\!   \int_{[\epp_n, \infty)}  \! \!\!\!\!\!  \!  \!\!\!\!    \pi_\phi (dx) x . 
\end{equation}
By (\ref{boundgf}) and (\ref{constraint2}), for all  $z\ino [z_n, \infty)$, we get $\smash{\phi'(z) \leqo g_{\epp_n} (z)+  \int_{[\epp_n, \infty)}   \pi_\phi (dx) x \leqo  g_{\epp_n} (z)+\frac{1}{2} \phi'(z_n)}$ and thus 
$ g_{\epp_n} (z) \geqo \tfrac{1}{2} \phi'(z) $. 
By (\ref{boundgf}) and (\ref{constraint1}) for all $z\ino [0, z_n]$ we get $f_{\eta_n } (z) \geq \psi' (z)  \! -\! \beta_\psi \! -\! 1\! -\! 2\beta_\phi z$. Thus, 
\begin{equation}
\label{minoFepet}
 \forall z\ino [z_n, \infty), \quad F'_{\! \epp_n, \eta_n} (z) \geq \tfrac{1}{2} \phi'(z) \quad \textrm{and} \quad \forall z\ino [0, z_n], \quad 
F'_{\!\epp_n,  \eta_n } (z) \geq \psi' (z)  \! -\! \beta_\psi \! -1.
\end{equation}
We note that $a'_n \geq z_n$ for all sufficiently large $n$.  Thus by (\ref{sandwichgf}),  (\ref{constraint1})
and (\ref{minoFepet}) we get for all 
$ z\ino [a'_n, \infty)$, 
\begin{equation}
\label{growthF}
 \tfrac{1}{2} \phi'(z) \leqo F'_{\! \epp_n , \eta_n} (z) \leqo  \phi'(z) +\frac{_2}{^{\eta_n}} (\beta_\psi \! -\! \beta_\phi)  +\! 
\int_{[\eta_n , \infty )}  \!\! \!\! \!\! \!\! \! \! \pi_\psi (dx) x \, \leq  \phi'(z)+ \phi'(a_n')\leqo 2\phi'(z).
 \end{equation}
Therefore, by Lemma \ref{examplphi} $(i)$, we get 
\begin{equation}
\label{growthbnan}
a_n F'_{\! \epp_n, \eta_n } (a_n) \asymp a_n \phi'(a_n) \asymp \phi (a_n) \; .
\end{equation}
We next observe that 
\begin{equation}
\label{limFnn}
\forall z\ino \bbR_+, \quad \lim_{n\to \infty} F_{\epp_n, \eta_n } (z)= \psi(z) \; .
\end{equation}
We fix $q \ino [0, 1)$. By Lemma \ref{rigimuphi}, for all $n\ino \bbN$ there is a non-trivial and critical offspring distribution $\mu_n$ such that $\mu_n (1)\eqo q$ and such that for all $r\ino [0, 1]$, 
 \begin{equation}
\label{Fnmunbis}
g_{\mu_n} (r) \eqo \sum_{k\in \bbN} r^k \mu_n(k) = r+ \frac{1\! -\! q}{a_nF_{\! \epp_n, \eta_n}'(a_n)} \, F_{\! \epp_n, \eta_n} \big( a_n(1\! -\! r)\big) \; .
\end{equation}  
We also set $b_n \eqo a_n F_{\! n}'(a_n) / (1\! -\! q)$. By (\ref{growthbnan}), (\ref{limFnn}) and Lemma \ref{muphiregu}, 
$$\psi_n  \eqo F_{\! \epp_n, \eta_n }, \quad  \mu_n (1)\eqo q, \quad b_n \asymp \phi (a_n) \quad \textrm{and} \; \textrm{\texttt{\L{}uka}$_{^{\,}}$($\baa, \bbb, \bmu, \psi$) holds true.} $$

\smallskip

\noi
\textbf{Proof of Proposition \ref{excntrex} $(i)$.} We assume here that $\psi$ satisfies \texttt{Sheu}$_{^{\,}}$($\psi$) and that $\phi$ satisfies \texttt{Sheu}$_{^{\,}}$($\phi$). 
It remains to prove \texttt{Sheu}$_{^{\,}}$($\baa, \bbb, \bmu, \psi$). To that end we fix $y\ino \bbR_+^*$ and for all sufficiently large $n$, $ \frac{1}{4} a_n \geko \frac{1}{4} a'_n\geko 4z_n^2 \geko  z_{n} \geko y$ and we get by (\ref{GreySheunew}) in Lemma \ref{sheu1integ} 
\begin{eqnarray*}
\int_{y}^{a_n} \!\!\! \! \frac{ds}{\sqrt{\int_0^s \psi_n(r) dr}} \!\! \!\! \!\! \!\! \!\! &=& \!\! \!\! \!\! \!\! \!\! \!
\int_{y}^{a_n} \!\!\! \! \frac{ds}{\sqrt{\int_0^s \! F_{\! \epp_n, \eta_n }(r) dr}}\leq  \frac{4}{\sqrt{3}} \int_{\frac{_1}{^4}y}^{\frac{_1}{^4}a_n} \!\!\! \!\frac{ds}{\, s\sqrt{F_{\! \epp_n, \eta_n }'(s)\, }} \\ 
\!\! \!\! \!\! \!\! \!\!  & \overset{\!\! \textrm{by (\ref{minoFepet}) }}{\leq}  &\!\! \!\! \!\! \!\! \!
\frac{\, 4}{\sqrt{3}} \!  \int_{\frac{_1}{^4}y}^{z_{n}} \!\!\! \frac{ds}{\, \, s \sqrt{\psi'(s)\! -\! \beta_\psi\! -\! 1} }+ \frac{4\sqrt{2}}{\sqrt{3}} \int_{z_{n}}^{\frac{_1}{^4}a_n} \!\!\! \!\frac{ds}{\, s\sqrt{\phi'(s)\, }} 
\end{eqnarray*} 
Therefore 
$$ \limsup_{n\to \infty} \int_{y}^{a_n} \!\!\! \! \frac{ds}{\sqrt{\int_0^s \psi_n(r) dr}} \leq \frac{4}{\sqrt{3}} \int_{\frac{_1}{^4}y}^{\infty} \!\!\! \frac{ds}{\,\,  s \sqrt{\psi'(s)\! -\! \beta_\psi\! -\! 1} } \xrightarrow[y\to \infty]{\; } 0, $$
by Lemma \ref{sheu1integ}. This proves \texttt{Sheu}$_{^{\,}}$($\baa, \bbb, \bmu$) and thus Proposition \ref{excntrex} $(i)$.\cqfd 

\smallskip

\noi
\textbf{Proof of Proposition \ref{excntrex} $(ii)$.} 
Here we assume that $\psi$ satisfies \texttt{Grey}$_{^{\,}}$($\psi$) and that $\phi $ satisfies \texttt{Grey}$_{^{\,}}$($\phi$) but not  \texttt{Sheu}$_{^{\,}}$($\phi$).
 Therefore we can choose $a_n'$ such that 
\begin{equation}
\label{nosheuphi}
\int_{a'_{n}}^{a_n} \!\!\! \!\frac{ds}{\, s\sqrt{\phi'(s)}} \geq 1 \; . 
\end{equation} 
Let us prove that \texttt{Grey}$_{^{\,}}$($\baa, \bbb, \bmu$) is satisfied. To that end, we fix $y\ino \bbR_+^*$ and for all sufficiently large $n$, by (\ref{constraint1}) we get $a_n \geko 2 z_n\geko y$ and we get 
\begin{eqnarray*} 
\int_{y}^{a_n} \!\!\! \! \frac{ds}{\psi_n(s) } \!\!= \!\!  \int_{y}^{a_n} \!\!\! \!\!  \frac{ds}{F_{\epp_n, \eta_n}(s) } 
\!\!\!\! &\! \overset{\textrm{by (\ref{GreySheunew})}}{\leq} &\!\!\!\! 2 \!  \int_{\frac{_1}{^2}y}^{\frac{_1}{^2}a_n} \!\!\!\! \!\frac{ds}{sF_{\! \epp_n, \eta_n}'(s)} \\
%=\!  2 \!  \int_{\frac{_1}{^2}y}^{z_n} \!\! \frac{ds}{\, sF_{\! \epp_n, \eta_n}'(s)}+2 \!  \int_{z_n}^{\frac{_1}{^2}a_n} \!\!\! \!\! \!\! \frac{ds}{sF_{\! \epp_n, \eta_n}'(s)}  \\
\!\!\!\! &\overset{\textrm{by (\ref{minoFepet})}}{\leq}  &  \!\!\!\! 2\! \int_{\frac{_1}{^2}y}^{\infty} \!\!\! \frac{ds}{\,\,  s( \psi'(s)\! -\! \beta_\psi\! -\! 1) }+4 \!  \int_{z_n}^{\frac{_1}{^2}a_n} \!\!\! \!\! \frac{ds}{s\phi'(s)} .
\end{eqnarray*} 
which entails  \texttt{Grey}$_{^{\,}}$($\baa, \bbb, \bmu$) by Lemma \ref{sheu1integ} and since $\phi$ and $\psi$ satisfy  \texttt{Grey}$_{^{\,}}$($\phi$) and  \texttt{Grey}$_{^{\,}}$($\psi$).

 We now prove that \texttt{Sheu}$_{^{\,}}$($\baa, \bbb, \bmu$) is not satisfied. Indeed, we observe 
$$ \int_{y}^{a_n} \!\!\! \! \frac{ds}{\sqrt{\int_0^s \psi_n(r) dr}}  \overset{\textrm{by (\ref{GreySheunew})}}{\geq} \sqrt{2}  \int_{y}^{a_n} \!\!\! \!\frac{ds}{\, s\sqrt{F_{\! \epp_n, \eta_n}'(s)\, }}  \overset{\textrm{by (\ref{growthF})}}{\geq} \!  \int_{a'_{n}}^{a_n} \!\!\! \!\frac{ds}{\, s\sqrt{\phi'(s)\, }} 
 \overset{\textrm{by (\ref{nosheuphi})}}{\geq}  1\, .$$
This completes the proof of Proposition \ref{excntrex} $(ii)$. \cqfd

\paragraph*{Proof of Proposition \ref{excntrex} $(iii)$.}  Here we assume that $\phi$ does not satisfy \texttt{Grey}$_{^{\,}}$($\phi$). 
By Lemma \ref{sheu1integ}, we can choose therefore $a'_n$ such that 
$\smash{\int_{a'_{n}}^{a_n} ds/(s\phi'(s)) \geqo 1}$ and we get 
$$ \int_{y}^{a_n} \!\!\! \! \frac{ds}{ \psi_n(s)}  \overset{\textrm{by (\ref{GreySheunew})}}{\geq}  \int_{y}^{a_n} \!\!\! \!\frac{ds}{\, sF_{\! \epp_n, \eta_n}'(s)}  \overset{\textrm{by (\ref{growthF})}}{\geq}  \int_{a'_{n}}^{a_n} \!\!\! \!\frac{ds}{ s\phi'(s)} 
\geq  1 \, , $$
and  \texttt{Grey}$_{^{\,}}$($\baa, \bbb, \bmu$) is not satisfied. It completes the proof of Proposition \ref{excntrex} $(iii)$. \cqfd

\section{Proof of Theorem \ref{applicactus}}
\label{pfsecapplicactus}
\subsection{Preliminary results}
\label{prelsecapplicactu}
This section brings together the results from D., Khanfir, Lin \& Torri \cite{DuKhLiTo22}, on which the proof of Theorem \ref{applicactus} is based, the proof itself being given in the following section. 

\smallskip

\noi
\textbf{Topological framework, snake metric and cactus.}
\label{prelsecapplicactu}
Let us first recall the definition of the Gromov-Hausdorff-Prokhorov space $(\bM, \dGHP)$. 
For any metric space $(X, \delta)$, we denote by $\cM_f(X)$ the space of its finite Borel positive measures.  
\begin{definition}
\label{GHPpointdef} For $i\ino \{ 1, 2\}$, let $(E_i, d_i, \rho_i, \mu_i )$ be compact metric spaces where 
$\rho_i \! \in \! E_i$ and $\mu_i\ino \cM_f(E_i)$.
We define their 
\textit{Gromov-Hausdorff-Prokhorov pseudo-distance} by 
\begin{equation}
\label{GHP}
\dGHP \big(E_1,   E_2\big) \!  = \!    \inf \! \Big\{  \delta \big( \phi_1 (\rho_1), \phi_2 (\rho_2) \big) \! \vee\!  \deHaus (
 \phi_1 (E_1), \phi_2(E_2) ) \! \vee \!  \dePro(  
\mu_1 \! \circ \! \phi_1^{-1}\! ,  \mu_2\!   \circ \! \phi_2^{-1} ) \Big\}, 
\end{equation}
where the infimum is taken on all Polish metric spaces $(X, \delta)$ and all isometrical embeddings $\phi_i \! : \! 
E_i \! \rightarrow \! X$, $i\ino \{ 1,2\}$; 
here, $ \deHaus$ stands for the Hausdorff distance on the space of compact subsets of $X$ and $ \dePro$ stands for
the Prokhorov distance on $\cM_{\! f} (X)$. \cq 
\end{definition}
We note that $\dGHP (E_1, E_2) \! =\!  0$  iff there is a bijective isometry $\varphi \! : \! 
E_1\rightarrow E_2$ such that $\varphi (\rho_1)\! = \! \rho_2$ and $\mu_1 \! \circ  \! \varphi^{-1}\!\! = \! \mu_2$. 
Let $\bM$ be the \emph{space of 
isometry classes of pointed finitely measured compact spaces}. 
Then $(\bM, \dGHP)$ is a Polish metric space: see Abraham, Delmas \&  Hoscheit \cite{ADHoscheit}, Theorem 2.5. 
We focus on metric spaces extending graph-trees that are called $\bbR$-trees and defined as follows. 
\begin{definition}
\label{Rtreedef} 
A metric space $(T, d)$ is a \emph{$\bbR$-tree} if it satisfies the following. 

\smallskip

\begin{compactenum}
\item[(a)] For any $\sigma_1, \sigma_2 \! \in\!  T$, there is a unique isometry 
$f:[0,d(\sigma_1,\sigma_2)] \! \rightarrow \! T$ such
that $f(0)\!=\! \sigma_1$ and $f(d(\sigma_1,\sigma_2))\! =\! \sigma_2$. We use the notation 
$\lgeo \sigma_1,\sigma_2\rgeo_T \! :=\! f([0,d (\sigma_1,\sigma_2)])$ for the geodesic arc. 

\smallskip

\item[(b)] $g([0,1]) \! = \! \lgeo g(0),g(1)\rgeo_T$, for any  
$g: [0, 1] \! \rightarrow \! T$, continuous and injective. 

\smallskip

\end{compactenum}
We denote by $\mathtt{Length}$ the $1$-dimensional Hausdorff measure on $T$, i.e., the (possibly infinite) Borel measure such that $\mathtt{Length}(\lgeo \sigma_1, \sigma_2\rgeo_T)\eqo d(\sigma_1, \sigma_2)$, for all $\sigma_1, \sigma_2\ino T$. \cq 
\end{definition}
\noi 
In other words, $(T,d)$ is a $\bbR$-tree if it is uniquely arcwise-connected and every arc is a geodesic. 
The following theorem shows that $\bbR$-trees are characterized metrically by the \emph{four point inequality}. \begin{theorem}
\label{4ptsth} A metric space $(T,d)$ is an $\bbR$-tree iff it is connected and satisfies the \emph{four points inequality}: 
for all $\sigma_1, \sigma_2, \sigma_3, \sigma_4 \ino T$,  
\begin{equation}
\label{4ptscondi}
d(\sigma_1 , \sigma_2 ) + d(\sigma_3 , \sigma_4 ) \leq \max \big( d(\sigma_1 , \sigma_3 ) + d(\sigma_2 , \sigma_4 ) 
 \, , \, d(\sigma_1 , \sigma_4 ) + d(\sigma_2 , \sigma_3 ) 
\big) .
\end{equation}
\end{theorem}
\noi
\textbf{Proof.} See Buneman \cite{Bun}, or Evans \cite{EvStF}, for a self-contained proof. \cqfd  

\smallskip

Instead of dealing directly with $\bbR$-tree metrics, 
we work with pseudo-metrics viewed as continuous two-parameters functions. 
More precisely, we set 
$\bC^*(\bbR_+^2, \bbR)\! :=\!  \{ (\zeta, f)\ino \bbR_+ \! \times  \bC(\bbR^2_+, \bbR)\! :\!  f(r,s) \eqo 0, (r,s) \ino [\zeta, \infty)^2 \}$, $\bC(\bbR_+^2, \bbR)$ being the space of continuous functions from $\bbR_+^2$ to $\bbR$.  
We easily see that 
$(\bC^*(\bbR_+^2, \bbR), \Delta)$ is Polish where for all $(\zeta, f)$, $(\zeta'\!, f')\ino \bC^*(\bbR_+^2, \bbR)$, we have set 
$$\Delta \big( (\zeta, f), (\zeta', f')\big)\! :=\!  |\zeta \! -\! \zeta'|+ \Delta_0 (f,f'),  \;  \textrm{with} \; 
\Delta_0 (f,f')\! :=\! \! \sum_{p\in \bbN^*} 2^{-p} \Big(1  \wedge \! \!\!\!\!\! \max_{\quad s,s'\in [0, p]} \!\!\!\!\!\!   |f(s,s')\! -\! f'(s,s')| \Big), $$ 
The space of continuous $\bbR$-tree pseudo-metrics is defined as follows. 
\begin{definition}
\label{psddef} $(a)$ We denote by $\MMT$ the space of $(\zeta, d)\ino \bC^*(\bbR^2_+ , \bbR)$ such that, for all 
$s_1$, $s_2$, $s_3$, $s_4\ino \bbR_+$, $d(s_1,s_2)\eqo  d(s_2,s_1)\geqo 0\eqo  d(s_1,s_1)$ and 
\begin{equation}
\label{fourpoints}
 d(s_1, s_2) + d(s_3, s_4) \leq \max \big( d(s_1, s_3) + d(s_2, s_4) \, , \,d(s_1, s_4) + d(s_2, s_3) \big) \; .
\end{equation}
We easily check that $\MMT$ is a closed subspace of $(\bC^*(\bbR_+^2, \bbR), \Delta)$ so that $(\MMT, \Delta)$ is Polish. 

\smallskip

\noi
$(b)$ Let $(\zeta, d)\ino \MMT$. Clearly, $d$ is a pseudo-metric on $\bbR_+$ and we write $s_1 \! \sim_d\!  s_2$ iff $d(s_1, s_2)\! = \! 0$, which defines an equivalence relation. We set  
$\mathtt{Tree} (\zeta, d)\! :=\! (T,d,r,\mu)$: here, $T$ is the quotient set $\bbR_+ / \! \!\! \sim_d$, $r\! :=\! \mathtt{proj}_d (0)$, where 
$\smash{\mathtt{proj}_d\! : \! \bbR_+ \! \rightarrow \! T}$ is the canonical projection, and $\mu$ is the pushforward measure 
of the Lebesgue measure on $[0, \zeta]$, i.e., $\smash{\int_{T} f(x) \, \mu (dx)\! = \! \int_0^\zeta   f(\mathtt{proj}_d (s))\,  ds}$, for any bounded Borel measurable function $f$. 
Observe that $\mathtt{proj}_d(s)\eqo r$, for all $s\ino [\zeta, \infty)$. Thus $\smash{T\eqo  \mathtt{proj}_d([0, \zeta])}$ and 
by Theorem \ref{4ptsth}, $(T,d)$ is a compact $\bbR$-tree. We view $r$ as a root and clearly $\mu \ino \cM_f(T)$.  \cq 
\end{definition}
We shall rely on the following continuity properties, simply adapted from \cite{DuKhLiTo22}. 
\begin{proposition}
\label{pseudoGHP} For all $(\zeta, d), (\zeta'\! , d') \ino \MMT$, 
$$1\wedge \dGHP\big(\mathtt{Tree} (\zeta,d),\mathtt{Tree} (\zeta'\! ,d') \big)\leq 3\! \cdot \!  2^{\zeta\vee \zeta'}\!  \Delta \big( (\zeta, d), (\zeta'\! , d')\big)\; .$$
\end{proposition}
\noi
\textbf{Proof.} W.l.o.g.~we can assume that $\zeta'\leqo \zeta$. We observe that $(\zeta, d')\ino \MMT$. Then, Proposition 4.14 in \cite{DuKhLiTo22} asserts that $\dGHP( \mathtt{Tree} (\zeta,d),\! \mathtt{Tree} (\zeta ,d'))\leqo \frac32 \! \max_{s,s'\in [0, \zeta]} \! |d(s,s')\! -\! d'(s,s')|$. Let $(T'\!\! , d'\! , r'\! , \mathrm{m'})$ stand for the $\bbR$-tree $\mathtt{Tree} (\zeta' \! ,d')$. We then check that $\mathtt{Tree} (\zeta ,d')\eqo \big(T'\! , d'\! , r'\! , \mathrm{m'} + (\zeta\! -\! \zeta')\delta_{r'} \big)$. Consequently,  $\dGHP ( \mathtt{Tree} (\zeta'\!  ,d'), 
\mathtt{Tree} (\zeta ,d'))\leqo \zeta\! -\! \zeta'$, which implies the desired inequality by a simple computation.  \cqfd 
\begin{proposition}
\label{submetric} Let $(E,d_E)$ be Polish and let $Z, Z_n$, $n\ino \bbN$, be $E$-valued r.v.s.  
Let $(\zeta, \bbd)$, $(\zeta_n, \bbd_n)$ and $(\zeta_n, \bbd_n^*)$, $n\ino \bbN$, be $\MMT$-valued r.v.s. 
We assume the following. 
\begin{compactenum}

\smallskip

\item[$(a)$] $(Z_n, (\zeta_n,\bbd_n)) \! \longrightarrow \! (Z, (\zeta, \bbd))$ in law on $E\! \times\! \MMT$.

\smallskip

\item[$(b)$] For all $n\ino \bbN$ and for all 
$ s, s^\prime\ino \bbR_+$, a.s.~$\bbd_n^*(s, s^\prime) \! \leq \! \bbd_n(s, s^\prime)$. 

\smallskip

\item[$(c)$] For all $s, s^\prime\ino \bbR_+$, $|\bbd_n^*(s, s^\prime) \! - \! \bbd_n(s, s^\prime)|\! \longrightarrow \! 0 $ in probability.

\smallskip

\end{compactenum}

\noi
Then, $(Z_n, (\zeta_n,\bbd^*_n),(\zeta_n,\bbd_n)) \! \longrightarrow \! (Z, (\zeta, \bbd), (\zeta,\bbd)) $ in law on 
$E\! \times\! \MMT^2$. 
\end{proposition} 
\noi
\textbf{Proof.} It is a straightforward adaptation of Proposition 4.4 in \cite{DuKhLiTo22}. \cqfd 

\smallskip

We next summarize from \cite{DuKhLiTo22} continuity properties of the $\bbR$-tree pseudo-distances that are obtained 
from continuous snakes as explained below. To that end, we recall from 
Definition \ref{snadef} that 
$\fSigma_1$ stands for the space of $1$-dimensional continuous snakes and from Remark \ref{easy} $(a)$ that 
$(\fSigma_1, D_1)$ is a Polish metric space, $D_1$ being defined in (\ref{snaspacedist}). 
\begin{definition}
\label{distsnadef} $(a)$ We set $\fSigma_1^*\!:=\! \big\{ (\zeta, (h,w))\ino \bbR_+ \!\! \times \! \fSigma_1: h(0)\eqo w_0(0)\eqo h(s)\eqo 0, s\ino [\zeta, \infty) \big\}$ and $D^*_1 ((\zeta, (h,w)), (\zeta'\! , (h'\! ,w')))\eqo |\zeta\! -\! \zeta'|+ D_1((h,w),(h'\! ,w'))$. We easily check that $(\fSigma_1^*, D_1^*)$ is Polish, its topology being the relative product topology of $\bbR_+$ and 
$\fSigma_1$.   

\smallskip

\noi
$(b)$ Let $(\zeta, (h,w))\ino \fSigma_1^*$ and $s,s'\ino \bbR_+$. We recall the notations $\widehat{w}_s\eqo w_s(h(s))$ and 
$m_h(s,s')\eqo \min_{r\in [s\wedge s', s\vee s']} h(r)$ and we define the following: 
\begin{equation}
\label{snadistdef}
d_h(s,s')\! := \! h(s)+h(s')-2m_h(s,s') \quad \textrm{and} \quad d_{h,w} (s, s')\! := \! \widehat{w}_{s}+  \widehat{w}_{s'} -2M_{h,w} (s, s')
\end{equation}
where $M_{h,w} (s,s')\eqo \big( \min_{r\in [m_h(s,s'), h(s)]} w_s(r) \big)\wedge \big( \min_{r\in [m_h(s,s'), h(s')]} w_{s'}(r)\big)$. \cq
\end{definition}
\begin{proposition}
\label{contsnadist}
 Let $(\zeta, (h,w)), (\zeta_n, (h_n,w_n))\ino \fSigma_1^*$, $n\ino \bbN$. Then the following holds true. 
\begin{compactenum}

\smallskip

\item[$(i)$] $(\zeta, d_h) $ and $(\zeta, d_{h,w})$ belong to $\MMT$.

\smallskip

\item[$(ii)$] If $\lim_{n\to 0} D^*_1 \big((\zeta_n, (h_n, w_n)), (\zeta, (h,w)) \big)\! = \! 0$, then 
$\lim_{n\to 0} \Delta \big( (\zeta_n, d_{h_n}), (\zeta, d_h) \big)\eqo 0$ and  $\lim_{n\to 0} \Delta \big( (\zeta_n, d_{h_n, w_n}), (\zeta, d_{h, w})\big) \! = \! 0$. 
\end{compactenum}
\end{proposition}
\noi\textbf{Proof.} $(i)$ is Lemma 4.22 \cite{DuKhLiTo22} and $(ii)$ is immediately derived from Lemma 4.20  \cite{DuKhLiTo22}.   \cqfd 
\begin{remark}
\label{scalingbehav} Let us discuss how scaling affects the previous spaces. Let $\lambda, b\ino \bbR_+^*$.

\noi
$(a)$  For $i\ino \{ 1, 2\}$, let $(E_i, d_i, \rho_i, \mu_i )$ be compact metric spaces where $\rho_i \! \in \! E_i$ and $\mu_i\ino \cM_f(E_i)$. We set $E^{_{\lambda,b}}_i $ $\! :=\!$ $  (E_i, \frac{1}{\lambda} d_i, \rho_i, \frac{1}{b} \mu_i)$, $i\ino \{ 1, 2\}$. Then $(\frac{1}{\lambda}\! \wedge\!  \frac{1}{b} ) \dGHP(E_1,E_2)\leqo  \dGHP(E^{_{\lambda,b}}_1,E^{_{\lambda,b}}_2)\leqo (\frac{1}{\lambda} \! \vee \! \frac{1}{b} ) \dGHP(E_1,E_2)$ (see Remark 4.12 in \cite{DuKhLiTo22} for details).

\smallskip

\noi
$(b)$ Let $(\zeta, d)\ino \MMT$. We set $(T,d,r,\mu) \! :=\! \mathtt{Tree} (\zeta, d)$ and $d'\! :=\! \frac{1}{\lambda} d(bs,bs')$, $s, s'\ino \bbR_+$. Then we easily check that $(\frac{1}{b} \zeta, d')\ino \MMT$ and $ \mathtt{Tree}(\frac{1}{b}\zeta, d')\eqo (T,\frac{1}{\lambda}d, r, \frac{1}{b} \mu)$. 

\smallskip

\noi
$(c)$ Let $(\zeta, (h,w))\ino \fSigma_1^*$ and $l\ino \bbR_+^*$. We set $h'(s)\! :=\! \frac{1}{l} h(bs)$ and $w'_s(r)\! :=\! \frac{1}{\lambda} w_{bs} (lr)$, $s,r \ino \bbR_+$. Then we easily check that $(\frac{1}{b} \zeta, (h',w'))\ino \fSigma_1^*$ and, for all $s,s'\ino\bbR_+$, that $d_{h'} (s,s')\! =\!  \frac{1}{l} d_h (bs,bs')$ and $d_{h'\! ,w'}(s,s')\eqo \frac{1}{\lambda} d_{h,w} (bs,bs')$. \cq 
\end{remark}
We next introduce the \emph{reflected Brownian cactus} as follows. To that end, we fix $\psi \ino \mathscr L$ and we assume 
$\texttt{Sheu} (\psi)$ (which implies $\texttt{Grey} (\psi)$), and as specified in $\mathbf{Case}(2)$, Section \ref{4casessec}, 
we denote by $\bN$ the excursion measure of the $\psi$-height process above $0$, we fix $c\ino \bbR_+^*$ and 
we consider the height process $(H_s)_{s\in \bbR_+}$ under $\bN (\, \cdot \, | \, \sup_{s\in \bbR_+}H_s \geko c)$, 
whose lifetime is denoted by $\zeta$. We also denote by $\mathbf W$ the continuous version of the $1$-dimensional 
Brownian snake with lifetime process $H$. We then set $C_s\eqo H_{s/2}$ 
and $W_{\! s}\eqo \mathbf W_{\! s/2}$, $s\ino \bbR_+$. Therefore $(2\zeta, (C,W))\ino \fSigma_1^*$. We easily observe that 
$|W|$ is also a snake with lifetime process $C$ too and that  $(2\zeta, (C,|W|))\ino \fSigma_1^*$. We next recall $d_{C}$ and 
$d_{C, |W|}$ from Definition \ref{distsnadef} and we use the following notation.  
\begin{equation}
\label{truecactdef}
\!\! \!\! \!   \mathtt{Tree} (2\zeta, d_C)\! =:\! (T_C, d_C, r_C, \mu_C) \; \, \textrm{and} \; \,  \mathtt{Tree} (2\zeta, d_{C, |W|})\! =:\! (T_{C, |W|}, d_{C, |W|}, r_{C, |W|}, \mu_{C, |W|}).
\end{equation}
Then $T_C$ is the \emph{$\psi$-Lévy tree conditioned to be higher than $c$} and we take $T_{C, |W|}$ as the definition of the \emph{reflected Brownian cactus with branching mechanism $\psi$}. To make the link with the definition given in Section \ref{traceapplsec}, we observe, thanks to the snake property, that it makes sense to define a $T_C$-index process 
$\smash{(\texttt{W}_{\sigma})_{\sigma \in T_C}}$ by setting $\smash{\texttt{W}_{\mathtt{proj}_{d_C } (s)}\!\!  :=\! \widehat{W}_{\! s}}$, for all 
$s\ino [0, 2\zeta]$. We easily check that $\smash{\sigma\ino T_C \! \mapsto \! \texttt{W}_\sigma}$ is continuous and conditionally given $C$, $\smash{(\texttt{W}_{\sigma})_{\sigma \in T_C}}$ is a centered Gaussian process such that 
$\smash{\bE \big [ | \, \texttt{W}_{\sigma}\!\!  -\! \texttt{W}_{\sigma' }|^2 | C\big]\eqo d_C(\sigma, \sigma')}$. We refer to D.~\& Le Gall \cite{DuLG05} for a detailled account on geometric properties of Lévy trees and to Lemmas 4.26 and 4.27 \cite{DuKhLiTo22} for basic geometric properties of $T_{C, |W|}$. 

\medskip

\noi
\textbf{Interpolation of $\bbT_{\! \mathtt b}$-valued BRWs, connections with snake metrics.}
We recall the definition of $\bbU$ from (\ref{UUlamdef}): we view it as a rooted ordered tree 
whose graph-distance is denoted by $d_{\mathtt{gr}}$. We embed it into a $\bbR$-tree by joining each adjacent vertices by a unit-length interval. More precisely, we set 
$\mathtt o\! :=\! \{ \varnothing\} \times \{ 0\}$ and we define the $\bbR$-tree $(\overline{\bbU}, \, \overline{\! d}_{\mathtt{gr}}, \mathtt o)$ by setting $\overline{\bbU}\! := \! \{ \mathtt o\} \cup  \bigcup_{u\in \bbU\backslash\{ \varnothing\}} \{ u\} \!  \times \!  (0, 1]$ and 
\begin{displaymath}
 \overline{\! d}_{\mathtt{gr}} \big( (u,s), (u'\! ,s')\big) \!\!  :=\!  
\left\{ \begin{array}{ll}
\! d_{\mathtt{gr}}(u,u') + s+ s'-2 \!\! &  \textrm{if $u \! \wedge\!  u'\!\!  \notin \! \{ u'\! ,u\}$,}\\
 \!\!\!   \big| (s+|u|\! -\! 1)_{_{\! +}} \!\! -\!  (s'\! + |u'|\! -\! 1)_{_{\! +}}\! \big| \!\! &  \textrm{if $u\! \wedge \! u'\!\!  \in \! \{ u'\! ,u\}$.}
\end{array} \right.
\end{displaymath}
We also define $\jmath \! :\! \bbU\! \to \! \overline{\bbU}$ by setting $\jmath (\varnothing)\! :=\! \mathtt  o$ 
and $\jmath(u)\! :=\! \{u\} \! \times \! \{ 1\}$.  

Let $t$ be a rooted ordered tree as in Definition \ref{Ulamtree}. 
We denote by $\overline{t}$ the $\bbR$-tree spanned by $\jmath(t)$ in $\overline{\bbU}$. Namely 
$ \overline{t}\! :=\! \bigcup_{u\in t} \lgeo \mathtt  o, \jmath(u)\rgeo_{\overline{\bbU}} $ (here, we recall from Definition \ref{Rtreedef} the notation $\lgeo\sigma_1, \sigma_2\rgeo_{\overline{\bbU}}$ for the geodesic arc joining $\sigma_1$ and $\sigma_2$ in the $\bbR$-tree $\overline{\bbU}$). Then $(\overline{t}, \, \overline{\! d}_{\mathtt{gr}}, \mathtt o)$ is a rooted complete $\bbR$-tree. 
If $t$ is finite, then $\overline{t}$ is compact and it is $\dGHP$-close to $t$ when the trees are equipped resp.~with their length and counting measures. More precisely, we easily check  that 
\begin{equation}
\label{ttbarclose}
\dGHP \big( (t,d_{\mathtt{gr}}, \varnothing, \#_t ) , (\overline{t}, \, \overline{\! d}_{\mathtt{gr}},\mathtt o, \mathtt{Length} ) \big) \leq 2, 
\end{equation} 
where $\#_t \eqo \sum_{u\in t} \delta_u$ stands for the counting measure on $t$ (see Lemma 5.2 in \cite{DuKhLiTo22} for details). 

We now assume that $t$ is finite. Its contour exploration $(v_k)_{0\leq k \leq 2(\# t-1)}$ (see Definition \ref{contsnadef}) extends continuously in  
$\overline{t}$ to a function $s\ino [0, 2(\# t \! -\! 1)] \! \mapsto \! v(\overline{t}, s) $ such that $ v(\overline{t},0)\eqo \mathtt  o$ and,  
for all $0\leqo k\leko  2(\# t-1)$ and all $s\ino (0, 1]$, 
\begin{equation}
\label{contcontexpl}
\textrm{$v(\overline{t},k+s)$ is the unique point of $\lgeo \jmath(v_k), \jmath(v_{k+1} )\rgeo_{\overline{t}}$ such that $\, \overline{\! d}_{\mathtt{gr}} (\jmath(v_k) , v(\overline{t},k+s) )\eqo s$.}
\end{equation} 
We easily check that $\, \overline{\! d}_{\mathtt{gr}} ( v(\overline{t},s), v(\overline{t},s'))\eqo d_{C(t)} (s,s')$ and in particular 
$\, \overline{\! d}_{\mathtt{gr}} ( \mathtt o, v(\overline{t},s))\eqo C_s(t)$, where $C_\cdot (t)$ stands for the contour process of $t$ (see Definition \ref{contsnadef}). This implies that 
\begin{equation}
\label{isometree}
\big( \overline{t}, \, \overline{\! d}_{\mathtt{gr}}, \mathtt o, 2\mathtt{Length} \big) \;\,  \textrm{and} \;\,  \mathtt{Tree} \big( 2(\# t\! -\! 1) , d_{C(t)}\big) \;\,  \textrm{are isometric.}
\end{equation}

We recall from the beginning of Section \ref{traceapplsec} that 
$ \bbT_{\! \mathtt b} \! :=\!  \bigcup_{n\in \bbN} \{ 1, \ldots, \mathtt b\}^n$ is the 
$\mathtt b$-ary rooted tree and that $\mathtt p(u,v)$, $u,v\ino \bbT_{\! \mathtt b}$, are the transition 
probabilities of the nearest neighbour null recurrent RW that is reflected at $\mathtt o$. 
Let $\btt$ be a random finite rooted ordered tree and $(\Upsilon_{\! u})_{u\in \btt}$ be a $\bbT_{\! \mathtt b}$-valued BRW whose law is given by 
 \begin{equation}
\label{lawBRWbis}
\bP \big( (\Upsilon_{\! u})_{u\in \mathtt t} \eqo (x_u)_{u\in \mathtt t} \, ;\,  \btt \eqo \mathtt t\big)= \bP ( \btt \eqo \mathtt t) 
\, \un_{\{ x_{\varnothing} = \mathtt o\}}\!\!\!  \prod_{u\in \mathtt t \backslash \{ \varnothing\} }\!\!  \mathtt p (x_{\overleftarrow{u}}, x_u)  
\end{equation}
for all finite rooted ordered tree $\mathtt t$ and all $x_u\ino \bbT_{\! \mathtt b}$, $u\ino \mathtt t$. Then, $(\Upsilon_{\! u})_{u\in \btt}$ 
extends continuously in $\overline{\bbT}_{\! \mathtt b}$ to a function $\sigma\ino \overline{\btt}\!  \mapsto 
 \! \overline{\Upsilon}_{\! \sigma}$ such that for all $u\ino \btt$, $\overline{\Upsilon}_{\! \jmath(u)}\!\! :=\! 
 \Upsilon_{\! u}$ and if $u\! \neq \! \varnothing$, for all $\sigma \ino \lgeo \jmath(\overleftarrow{u}), 
 \jmath(u)\rgeo_{\overline{\btt}}$, 
\begin{equation}
\label{interpolUpsi} \textrm{ $\overline{\Upsilon}_{\! \sigma}$ is the unique point of 
$\lgeo \overline{\Upsilon}_{\! \jmath (\overleftarrow{u})}, \overline{\Upsilon}_{\! \jmath(u)}\rgeo_{\overline{\bbT}_{\! \mathtt b}} $ 
such that  
$\, \overline{\! d}_{\mathtt{gr}} (\overline{\Upsilon}_{\! \jmath(\overleftarrow{u})}, \overline{\Upsilon}_{\! \sigma})\eqo 
\, \overline{\! d}_{\mathtt{gr}} (\jmath(\overleftarrow{u}) , \sigma)$. 
}
\end{equation}
As in (\ref{ttbarclose}), we easily get 
\begin{equation}
\label{Rbarclose}
\dGHP \big( (\mathcal R (\Upsilon) ,d_{\mathtt{gr}}, \varnothing, 2\mathbf m_{\mathtt{occ}} ) , 
(\mathcal R (\overline{\Upsilon}), \, \overline{\! d}_{\mathtt{gr}}, \mathtt o, \overline{\mathbf m}_{\mathtt{occ}}) \big) \leq 4, 
\end{equation}
where $\smash{\mathcal R (\Upsilon) \! :=\! \{ \Upsilon_{\! u}; u\ino \btt \}}$, 
$\smash{\mathcal R (\overline{\Upsilon}) \! :=\! \{ \overline{\Upsilon}_{\! \sigma}; \sigma\ino \overline{\btt}\}}$, 
$\smash{\mathbf m_{\mathtt{occ}}\! :=\!  \sum_{u\in \btt} \delta_{\Upsilon_{\! u}}}$ and 
$\smash{\int_{\mathcal R (\overline{\Upsilon}) } f(x) \,  
\overline{\mathbf m}_{\mathtt{occ}}(dx) \! := }$ $\smash{\! \int_{0}^{2(\# \btt-1)}\!  f(\overline{\Upsilon}_{\! v(\overline{\btt},s) })\, ds}$, 
for all bounded measurable $\smash{f\! : \! \overline{\bbT}_{\! \mathtt b}\! \to \! \bbR}$ (namely, 
$\smash{\overline{\mathbf m}_{\mathtt{occ}}}$ is the occupation measure yielded by the contour exploration of $\smash{\overline{\btt}}$). 
To simplify notation, we also set 
$$\zeta_\btt \! :=\! 2(\# \btt \! -\! 1) \quad \textrm{and} \quad d^*_{\overline{\Upsilon}} (s,s')\! :=\!   \, \overline{\! d}_{\mathtt{gr}} ( \overline{\Upsilon}_{\! v(\overline{\btt},s)} ,  \overline{\Upsilon}_{\! v(\overline{\btt},s')}), \; s,s'\ino [0, \zeta_{\btt} ].$$ 
Then, we easily check that a.s.~$(\zeta_{\btt} , d^*_{\overline{\Upsilon}})\ino \MMT$ and that 
\begin{equation}
\label{isometreesna}
\big(\mathcal R (\overline{\Upsilon}), \, \overline{\! d}_{\mathtt{gr}}, \mathtt o, \overline{\mathbf m}_{\mathtt{occ}} \big) \;\,  \textrm{and} \;\,  \mathtt{Tree} \big( \zeta_{\btt} , d^*_{\overline{\Upsilon}} \big) \;\,  \textrm{are isometric.}
\end{equation}

We recall that $|v|$ stands for the height of $v\ino \bbU$ and we observe that $(|\Upsilon_{\! u}|)_{u\in \btt}$ is a 
$1$-dimensional BRW whose jumps, conditionally given $\btt$, are those of a simple symmetric RW reflected at 
$0$. 
We next recall a result from \cite{DuKhLiTo22} stating that the graph-metric on 
$\mathcal R (\Upsilon)$ is close to the snake metric associated with  $(|\Upsilon_{\! u}|)_{u\in \btt}$: 
it is a geometric key point in the proof of Theorem \ref{applicactus}.  
\begin{lemma}
\label{snadistdomin}
We keep the previous notation. Then, there is a $\bbZ$-valued $\btt$-indexed BRW 
$ \smash{\mathbf S\! :=\! (S_u)_{u\in \btt}}$, whose jumps,  
conditionally given $\btt$, are independent with law $\smash{\frac{_{_1}}{^{^2}} \delta_{-1}\! + \frac{_{_1}}{^{^2}} \delta_{1}}$, and such that 
$\smash{(|\Upsilon_{\! u}|)_{u\in \btt}}$ $\smash{ \eqo (|S_{u}|)_{u\in \btt}}$ a.s. We set $\smash{\zeta_\btt \! :=\! 2(\# \btt \! -\! 1)}$, we recall from Definition \ref{contsnadef} that 
$\smash{C_\cdot (\btt)}$ and $\smash{W_\cdot (\bS , \cdot)}$ stand for resp.~the contour process of $\btt$ and the contour snake of $\bS$.
Then, a.s.~for all $s$, $s'\ino$  $[0, \zeta_{\btt}]$, 
\begin{equation}
\label{snadistdomina}
G_{s,s'}\! :=\! \tfrac12 d_{C(\btt), W(\bS, \cdot)} (s,s')-  \tfrac12 \, \overline{\! d}_{\mathtt{gr}} ( \overline{\Upsilon}_{\! v(\overline{\btt},s)} ,  \overline{\Upsilon}_{\! v(\overline{\btt},s')}) \geq 0 , 
\end{equation}
and $\bP (G_{s,s'} \geqo x) \leq \mathtt b^{2-x}$ for all $x\ino \bbR_+$.  
\end{lemma}
\noi\textbf{Proof.} See (5.19) p.~49 in \cite{DuKhLiTo22}. \cqfd 

\subsection{Proof of Theorem \ref{applicactus}}
\label{Pfapplicactussubsec}
 We argue under the assumptions of Theorem \ref{maincvsnake}, $\mathbf{Case}(2)$. Namely, we fix $c\ino \bbR_+^*$, 
$\psi\ino \mathscr L$, $(a_n)_{n\in \bbN}$, $(b_n)_{n\in \bbN}$, $(\mu_n)_{n\in \bbN}$, and we assume 
\texttt{\L{}uka}$_{^{\,}}$($\baa, \bbb, \bmu, \psi$) and \texttt{Sheu}$_{^{\,}}$($\baa, \bbb, \bmu $). This implies 
\texttt{Norm} ($\baa, \bbb$), i.e., $a_n \! \to \! \infty$ and $\lambda_n \! :=\! b_n/a_n \! \to \! \infty$ and \texttt{Sheu}($\psi$) (and thus \texttt{Grey}($\psi$), too). We denote by $\btt_n$ a random finite rooted ordered tree 
distributed as a single GW($\mu_n$)-tree conditionned to have total height $\, \geqo \lambda_n c$. Let 
 $\smash{(\Upsilon_{\! n,u})_{u\in \btt_n}}$ be a $\bbT_{\! \mathtt b}$-valued BRW distributed as in (\ref{lawBRW}) and let 
 $\smash{\bS_n \eqo (S_{n, u})_{u\in \btt_n}}$ be associated with $\smash{(\Upsilon_{\! n,u})_{u\in \btt_n}}$ as in Lemma \ref{snadistdomin}. Namely, $\bS_n$ is a $\bbZ$-valued BRW 
whose jumps, conditionally given $\btt_n$, are independent with law $\smash{\frac{_{_1}}{^{^2}} \delta_{-1}+ \frac{_{_1}}{^{^2}} \delta_{1}}$, and such that 
a.s.~$\smash{(|\Upsilon_{\! n,u}|)_{u\in \btt_n} \eqo (|S_{n,u}|)_{u\in \btt_n}}$. For all $s,r\ino \bbR_+$, we set 
\begin{equation}
\label{renormCWbis} 
\zeta_n:= b_n \zeta_{\btt_n} =2b_n(\# \btt_n \! -\! 1) , \quad   C^{_{(n)}}_{s}  \eqo \tfrac{1}{\lambda_n} C_{b_ns} (\btt_n) \quad \textrm{and} \quad   W^{_{(n)}}_{s} \! (r)  \eqo  \tfrac{1}{\sqrt{\lambda_n}} 
W_{\! b_ns} (\bS_n, \lambda_n r). 
\end{equation}  
Theorem \ref{maincvsnake} implies that 
\begin{equation}
\label{snakecvappli}
\Big( \zeta_n , \big(C^{_{(n)}}_{\cdot} \!\! ,  \, |W^{_{(n)}}_{\cdot}| \big)\Big) \xrightarrow[n\to \infty]{\; } \big( 2\zeta, (C_\cdot\, ,  |W_\cdot|)  \big)
\end{equation}
in law in $\bbR_+\! \times \!  \bC^{_0}_{^1} \! \times \! \mathbf C (\bbR_+, \bC^{_0}_{^1})$, i.e., in law in 
$(\fSigma^*_1, D^*_1)$. Here, $(C_{2s})_{s\in \bbR_+}$ is distribted as $H$ under $\bN(\, \cdot \, | \, \sup_{s\in \bbR_+} H_s\geko c)$, 
$\zeta$ is the lifetime of the excursion $H$ and $(W_{2s})_{s\in \bbR_+}$ is a continuous version of a Brownian snake with lifetime process $C$, as specified earlier. 

  We next recall from (\ref{contcontexpl}) the continous extension $\smash{(v(\overline{\btt}_n, s))_{s\in [0, \zeta_{\btt_n}] }}$  
of the contour exploration of $\btt_n$ and from (\ref{interpolUpsi}) the continuous extension 
$\smash{(\overline{\Upsilon}_{\! n, \sigma})_{\! \sigma \in \overline{\btt}_n}\! }$ of $\smash{(\Upsilon_{\! n,u})_{u\in \btt_n}}$. To simplify notation, for all $s,s'\ino [0, \zeta_n]$, we set 
$$\bbd^*_n (s,s')\! :=\!   \tfrac{1}{\sqrt{\lambda_n}} \, \overline{\! d}_{\mathtt{gr}} \big( \overline{\Upsilon}_{\! n, v(\overline{\btt}_n,b_ns)} ,  \overline{\Upsilon}_{\! n, v(\overline{\btt}_n,b_ns')} \big), \;   \fdelta_n \! :=\! d_{C^{_{(n)}}}   \; \textrm{and} \; \bdd_n \! :=\! d_{C^{_{(n)}}\! , |W^{_{(n)}}| } ,$$
where $\smash{d_{C^{_{(n)}}}} $ and $\smash{d_{C^{_{(n)}}\! , |W^{_{(n)}}| }}$ are as in Definition \ref{distsnadef}. By (\ref{snakecvappli}) and Proposition \ref{contsnadist}, the following convergence holds in law on $(\MMT, \Delta)^2$ 
\begin{equation}
\label{cvdistsimple}
\big( (\zeta_n, \fdelta_n), (\zeta_n, \bdd_n) \big)\xrightarrow[n\to \infty]{\;} \big( (2\zeta, d_C), (2\zeta, d_{C, |W|}) \big) \; .
\end{equation} 
By Lemma \ref{snadistdomin}, we then a.s.~get for all $s,s'\ino \bbR_+$ that $\bbd^*_n (s,s') \leqo \bbd_n(s,s')$ and also 
$\bP ( \bbd_n(s,s') \! -\! \bbd^*_n (s,s') \geqo 2\epp ) \leqo \mathtt b^{2-\epp \sqrt{\lambda_n}}$. Thus $ |\bbd_n(s,s') \! -\! \bbd^*_n (s,s')|\! \to \! 0$, in probability. Proposition \ref{submetric} applies and we get 
\begin{equation}
\label{cvdistmoinsimple}
\big( (\zeta_n, \fdelta_n), (\zeta_n, \bdd^*_n) \big)\xrightarrow[n\to \infty]{\;} \big( (2\zeta, d_C), (2\zeta, d_{C, |W|}) \big) \; .
\end{equation} 
 in law on $(\MMT, \Delta)^2$, which entails, by Proposition \ref{pseudoGHP} and notation (\ref{truecactdef}), 
 \begin{equation}
\label{cvTreedist}
\!\!\!\!\!\! \big( \mathtt{Tree}(\zeta_n, \fdelta_n), \mathtt{Tree} (\zeta_n, \bdd^*_n) \big)\xrightarrow[n\to \infty]{\;} \big( 
(T_C, d_C, r_C, \mu_C) , (T_{C, |W|}, d_{C, |W|}, r_{C, |W|}, \mu_{C, |W|}) \big)
\end{equation} 
in law on $(\bM, \dGHP)^2$. By (\ref{isometree}), (\ref{isometreesna}) and Remark \ref{scalingbehav} $(b)$, it implies that 
\begin{equation}
\label{cvTreenearly}
\!\!\!\!\!\! \Big( \mathcal R(\overline{\Upsilon_{\! n} }) , \tfrac{1}{\sqrt{\lambda_n}} \, \overline{\! d}_{\mathtt{gr}} , \mathtt o, \tfrac{1}{b_n} \overline{\mathbf m}^{(n)}_{\mathtt{occ}} \Big)\xrightarrow[n\to \infty]{\textrm{in law on $(\bM, \dGHP)$}} \big( 
 T_{C, |W|}, d_{C, |W|}, r_{C, |W|}, \mu_{C, |W|} \big)
\end{equation} 
jointly with $\big( \overline{\btt}_n, \tfrac{1}{\lambda_n} \, \overline{\! d}_{\mathtt{gr}} , \mathtt o, \tfrac{2}{b_n} \mathtt{Length}  \big)\! \to \! (T_C, d_C, r_C, \mu_C)$ in law on $(\bM, \dGHP)$. This implies the desired result by the bounds in 
(\ref{ttbarclose}), in (\ref{Rbarclose}) and in Remark \ref{scalingbehav} $(a)$. \cqfd

{\footnotesize

%\bibliographystyle{acm}
%\bibliography{Refsna}

}

\newpage

\appendix

\renewcommand{\theequation}{A.\arabic{equation}}

\section{The case of uniformly bounded i.i.d.~jumps.}
\label{unifboundedsec} 
\subsection{Statement of the main result}
Here we state a limit theorem for snakes of BRWs whose spatial displacements are $\bbR$-valued i.i.d.~uniformly bounded r.v.s. This result is slightly better than Theorem \ref{maincvsnake} which relies on a coupling argument and its assumptions  are in some sense optimal. 

\begin{theorem}
\label{unifbounded}
Let $\psi \in \mathscr L$ satisfy \emph{\texttt{Sheu}$_{^{\,}}$($\psi$)}. Let $(a_n)_{n\in \bbN}$, $(b_n)_{n\in \bbN}$, $(\mu_n)_{n\in \bbN}$, $\bS_n \eqo (S_{n, u})_{u\in \btt_n}$, $n\ino \bbN$, and $(C,Y,W)$ be as in \emph{$\textbf{Case (i)}$}, $0\leqo \mathbf i \leqo 2$. 
We recall from (\ref{psidisdef}) the definition of $\psi_n$. 
Recall that $C^{_{(n)}}_{\cdot}$, $V^{_{(n)}}_{\cdot}$ and $W^{_{(n)}}_{\cdot}$ are rescaled as in (\ref{renormCW}). 

\begin{compactenum}

\smallskip

\item[$(a)$] We assume \emph{\texttt{\L{}uka}$_{^{\,}}$($\baa, \bbb, \bmu, \psi$)} and 
\emph{\texttt{Hght}$_{^{\,}}$($\baa, \bbb, \bmu, \psi$)}.

\smallskip

\item[$(b)$] We assume for all $n\ino \bbN$, that conditionally given $\btt_n$, the jumps of $\bS_n$ are $\bbR$-valued and i.i.d.~random variables whose deterministic law is denoted by $\bgam_n$. We also assume that there are $c, \beta\ino \bbR_+^*$ such that for all $n\ino \bbN$, 
\begin{equation}
\label{cenuni}
\bgam_n ([-c,c]) \eqo 1, \quad 
\int_{\bbR}\!  y \bgam_n(dy)\eqo 0 \quad \textrm{and} \quad \beta_n \! :=\! \int_{\bbR}\! y^2 \bgam_n(dy) \xrightarrow[n\to \infty]{\; } \beta.
\end{equation}

\smallskip

\end{compactenum}
Then the following holds true. 

\begin{compactenum}

\smallskip

\item[$(i)$] \emph{\texttt{Sheu}$_{^{\,}}$($\baa, \bbb, \bmu$)} implies that the following convergence 
\begin{equation}
\label{cvsnake1}
\big( C^{_{(n)}}_\cdot, V^{_{(n)}}_\cdot ,  W^{_{(n)}}_\cdot \big)  \xrightarrow[n\to \infty]{\;} (C, Y , \sqrt{\beta}W) .
\end{equation}
holds weakly in $\bC^{_0}_{^1} \! \times \! \bD(\bbR_+, \bbR) \! \times \! \mathbf C (\bbR_+, \bC^{_0}_{^1})$. 

\smallskip

\item[$(ii)$] Suppose that we are in \emph{$\textbf{Case (0)}$}. We also assume for all $n\ino \bbN$ 
that $\bgam_n \eqo \frac{1}{2} (\delta_{-1} + \delta_{1})$ and that $\liminf_{n\to \infty} \mu_n (1) \geko 0$. Then \emph{\texttt{Sheu}$_{^{\,}}$($\baa, \bbb, \bmu$)} is equivalent to the tightness in $\bC^{_0}_{^1}$ of the laws of $(\widehat{W}^{_{(n)}}_s)_{s\in \bbR_+}$, $n\ino \bbN$. 
\end{compactenum}
\end{theorem}
\noi
\textbf{Proof:} see Section \ref{Thm1unpfsec}.  \cqfd 

\begin{remark}
\label{rererema} Theorem \ref{unifbounded} $(i)$ holds true in 
$\textbf{Case (3)}$: it is a specific case of 
a theorem due to Janson \& Marckert \cite{JanMar05} when $\alpha\eqo 2$ and to Marzouk \cite{Mar20} in the 
other stable cases. \cq  
\end{remark}
As already mentioned in Remark \ref{afterthm3} $\textbf{(d)}$, Theorem \ref{unifbounded} is better that Theorem \ref{maincvsnake} in the sense that (\ref{cvsnake1}) under the sole assumption \texttt{Sheu}$_{^{\,}}$($\baa, \bbb, \bmu$): namely, we do not require the assuumption 
$\lim_{n\to \infty} b_n / (a_n (\log a_n)^4)\eqo 0$. This unnecessary assumption is due to the coupling used in the proof of Theorem \ref{maincvsnake}. The proof of Theorem \ref{unifbounded} does not use coupling arguments and it proceeds along the same line of proof as Theorem \ref{Sheuexplain}: we study oscillation times via a discrete version of the functional equation of the exit processes of $1$-dimensional Brownian snakes stated in Lemma \ref{maxBrosnadiscr}. These estimates in the dicrete setting are the technical part of the proof of Theorem \ref{unifbounded}. 

\subsection{Oscillations of BRWs}
\label{osciBRWsec}
In this section, we study the oscilation times of real-valued BRWs: we prove an analogue of Lemma \ref{osc1lemma} and of Lemma \ref{renewosci}. 
More precisely, we introduce the following. 
\begin{definition}
\label{exittime} 
Let $(t(p))_{1\leq p\leq N}$ and let $ \big(S(p)\eqo (S_u(p))_{u\in t(p)}\big)_{1\leq p\leq N}$ be a (possibly infinite) sequence 
of $\bbR$-valued BRWs as in Definition \ref{forestdef}. We assume that $\# t(p)\leko \infty$, $1\leqo p\leq N$. Let $S\eqo (S_u)_{u\in t}$ be the BRW associated to 
the sequence $(S(p))_{1\leq p\leq N}$as in Remark \ref{contord}. Namely, $S_\varnothing \eqo 0$, for all $1\leqo p\leq N$, $t(p)\eqo \theta_{[p]} t$ and for all $u\ino t(p)$, $S_{[p]\ast u}\eqo S_u(p)$.    
We denote by $(v_k)$ (resp.~$(u_l)$ the vertices of $t$ listed in increasing contour order (resp.~depth first order) 
Let $y\ino \bbR_+^*$. We define recursively the sequence of times $(\sigma_q (y))_{q\in \bbN}$ and $(\bs_q (y))_{q\in \bbN}$ by setting $\sigma_0 (y) \eqo  \bs_0 (y) \eqo 0$ and 
\begin{eqnarray*}
\sigma_{q+1} (y) \eqo \inf \big\{ k \geko \sigma_{q} (y)   
\!\!\!\! \!\!\!  & : & \!\!\!\!  \!\!\!  | S_{v_{k+1}} \!\!  -\! S_{v_{k+1} \wedge v_{\sigma_q (y )+1}}| \geko  y \big\} \\
& \textrm{and} &  \bs_{q+1} (y) \eqo \inf \big\{ l \geko \bs_{q} (y)   
:  | S_{u_{l+1}} \!\!  -\! S_{u_{l+1} \wedge u_{\bs_q (y )+1}}| \geko  y \big\} , 
\end{eqnarray*}
with the convention that $\inf \emptyset \eqo \infty$. We call $(\sigma_q (y))_{q\in \bbN}$ (resp.~$(\bs_q (y))_{q\in \bbN}$) the \emph{contour} (resp.~\emph{depth-first}) \emph{$y$ oscillation times} of $S$. \cq 
\end{definition}

We see below in Lemma \ref{insert} that contour oscillation times have nice geometric properties and in Lemma \ref{timeincrease2} that depth-first oscillation times have convenient probabilistic properties. They are related as follows. 
\begin{lemma}
\label{timeincrease1} We keep the same notations as in Definition \ref{exittime} and we recall notation $(K_l)_{l\in \bbN}$ from (\ref{Klcontour}). Then, for all $q\ino \bbN$, $\bs_q (y)$ is finite iff $\sigma_q(y)$ is finite and in that case 
$K (\bs_q (y))\eqo \sigma_q(y)$. 
\end{lemma}
\noi
\textbf{Proof.} We argue recursively. First observe that the case where $q\eqo 0$ holds trivially. Then, we suppose that 
$\bs_q(y)$ and $\sigma_q(y)$ are finite and that $K (\bs_q (y))\eqo \sigma_q (y)$. We then set $w\eqo u_{\bs_{q} (y)+1}\eqo v_{\sigma_q(y)+1}$ and $J\eqo \{ u\ino t : w\!<_{\mathtt{lex}} \! u \; \textrm{and} \; |S_u \! -\! S_{u\wedge w}| \geko y \}$. We first observe that $\bs_{q+1} (y)$ $\leko$ $\infty$ iff $J\! \neq \! \emptyset$ and in that case $\bs_{q+1} (y) \eqo \min \{ l\ino \bbN: u_{l+1}\ino J \}$. 
We next prove similarly that $\sigma_{q+1} (y)$ $\leko$ $\infty$ iff $J\! \neq \! \emptyset$ and in that case $\sigma_{q+1} (y) \eqo \min \{ k\ino \bbN: v_{k+1}\ino J \}$. 

\noi
\emph{Indeed}, suppose that $u\ino J$. Since we deal with forests of finite trees, the contour exploration visits all vertices and there is $k\ino \bbN$ such that $v_{k+1}\eqo u$. Since $ w\! <_{\mathtt{lex}} u$, we necessarily get $\sigma_q(y) \leko k$ and therefore $\sigma_{q+1}(y) \leko \infty$. Conversely, let us suppose that $\sigma_{q+1}(y) \leko \infty$. By definition $v_{\sigma_{q+1} (y) +1} \! \notin \! \{ v_{k+1}\,  ; \, 0\leqo k \leko \sigma_{q+1} (y) \}$. Thus $w\! <_{\mathtt{lex}} v_{\sigma_{q+1} (y) +1}$ and we get $v_{\sigma_{q+1} (y) +1}\ino J$, which is not empty; it also entails immediately that $\sigma_{q+1} (y) \eqo \min \{ k\ino \bbN: v_{k+1}\ino J \}$.

We thus have proved that $\sigma_{q+1} (y) $ is finite iff $\bs_{q+1} (y)$ is finite. In that case, we also get 
$K (\bs_{q+1} (y)) \eqo \min \{ K_l : u_l \ino J \} \eqo \min \{ k \ino \bbN: v_{k+1} \ino J \}$ by definition of $(K_l)_{l\in \bbN}$. Thus $ K (\bs_{q+1} (y))\eqo \sigma_{q+1} (y)$, which completes the proof of the lemma by recursion. \cqfd 

\medskip

The following lemma is the analog of Lemma \ref{osc1lemma} and it shows how contour oscillation times can be used to control the oscillations of the snake. 
\begin{lemma}
\label{insert} We keep the same notations as in Definition \ref{exittime}. We fix $y_0\ino \bbR_+^*$ and we assume that
the jumps of $S$ are bounded by $y_0$: $|\xi_u| \leqo y_0$, $u\ino t$. 
We fix $k, k^\prime\ino \bbN$ such that $ k\leko k^\prime$, 
and $y, y_1, y_2\ino \bbR_+^*$ such that $y \geqo y_0$ and $y_2 \geko 3y + 2y_0 + y_1$. 
We also assume the following. 
\begin{compactenum}
\smallskip

\item[$(a)$] For all $v\ino \lgeo v_{k+1}\!  \wedge\!  v_{k^\prime+1}, v_{k+1} \rgeo $, $|S_{v_{k+1}} \! \! -\! S_v| \leqo y_1$. 

\smallskip

\smallskip

\item[$(b)$] There exists an integer $j$ such that $k\leko  j \leqo k^\prime$ and  $|S_{v_{k+1}} \! -\! S_{v_{j+1}}| \geko y_2$.

\end{compactenum}

\smallskip

\noi
Then, there exists $q \ino \bbN^*$ such that $k\leq \sigma_q (y) \leko \sigma_{q+1} (y) \leko k^\prime  $. 
\end{lemma}
\noi
\textbf{Proof.} To prove Lemma \ref{insert}, we first prove the following. 
 
\begin{lemma}
\label{elcontr} Let $\fftree\ino \bbT$ be the tree associated with a forest of finite trees. 
Let $(v_k)_{k\in \bbN}$ be its contour exploration. For all $k\ino \bbN$, set $I_{\leq k}\eqo \{ v_j\, ; \, j\leqo k\}$ and 
$I_{\geq k}\eqo \{ v_j\, ; \, j\geqo k\}$ and let $i_1, \ldots, i_4 \ino \bbN$ be such that $i_1\leqo i_2\leqo i_3\leqo i_4$. 
Then the following holds true.
\begin{compactenum}

\smallskip

\item[$(i)$] If $u\ino \lgeo v_{i_1} , v_{i_2} \rgeo$, then there exists $i\ino \bbN$ such that $i_1 \leqo i\leqo i_2$ and $u=v_i$.

\smallskip

\item[$(ii)$]  $\lgeo \varnothing , v_{i_1} \wedge v_{i_2} \rgeo \eqo I_{\leq i_1}\!  \cap I_{\geq i_2}$.

\smallskip

\item[$(iii)$] $v_{i_1} \wedge v_{i_3}\eqo (v_{i_1} \wedge v_{i_2})\wedge(v_{i_2} \wedge v_{i_3})$

\smallskip

\item[$(iv)$] $v_{i_1} \wedge v_{i_2}\ino \lgeo v_{i_1} \wedge v_{i_3}, v_{i_1}\rgeo$ and $v_{i_2} \wedge v_{i_3}\ino \lgeo v_{i_1} \wedge v_{i_3}, v_{i_3}\rgeo$.

\smallskip

\item[$(v)$] If $v_{i_3} \wedge v_{i_4} \ino \lgeo v_{i_1} \wedge v_{i_3}, v_{i_3}\rgeo$, then 
$v_{i_2} \wedge v_{i_4} \ino \lgeo v_{i_1} \wedge v_{i_3}, v_{i_3} \wedge v_{i_4}\rgeo$.

\smallskip

\item[$(vi)$] If $v_{i_3} \wedge v_{i_4} \ino \lgeo \varnothing , v_{i_1} \wedge v_{i_3}\rgeo$, then 
$v_{i_2} \wedge v_{i_4} =  v_{i_3} \wedge v_{i_4}$.

\end{compactenum}
\end{lemma}  
\noi
\textbf{Proof of Lemma \ref{elcontr}.} First note that $(i)$ is a simple consequence of Remark \ref{contord} \textbf{(a)}. Let us prouve $(ii)$: by Remark \ref{contord} \textbf{(b)}, $\lgeo \varnothing, v_k\rgeo\eqo I_{\leq k} \cap I_{\geq k}$. Thus, 
$\lgeo \varnothing , v_{i_1} \wedge v_{i_2} \rgeo= \lgeo \varnothing , v_{i_1}  \rgeo\cap \lgeo \varnothing , v_{i_2} \rgeo\eqo  (I_{\leq i_1} \cap I_{\geq i_1}) \cap (I_{\leq i_2} \cap I_{\geq i_2}) \eqo I_{\leq i_1}\!  \cap I_{\geq i_2}$.
Let us prove $(iii)$. First note that $v_{i_1} \wedge v_{i_2}$ and $ v_{i_2} \wedge v_{i_3}$ both 
belong to $\lgeo \varnothing, v_{i_2} \rgeo$. We set $u(v_{i_1} \wedge v_{i_2})\wedge(v_{i_2} \wedge v_{i_3})$. Then, 
$$ \lgeo \varnothing , u\rgeo \eqo \lgeo \varnothing , v_{i_1}\! \wedge\!  v_{i_2}\rgeo \cap  \lgeo \varnothing , v_{i_2}\! \wedge \! v_{i_3}\rgeo\eqo  (I_{\leq i_1} \cap I_{\geq i_2} )\cap (I_{\leq i_2} \cap I_{\geq i_3} )\eqo  I_{\leq i_1} \cap I_{\geq i_3} =  \lgeo \varnothing , v_{i_1} \! \wedge \! v_{i_3}\rgeo $$
by $(ii)$, which proves $(iii)$. Note that $(iv)$ is straightforward consequence of $(iii)$. 

Let us prove $(v)$ and $(vi)$. Since 
$v_{i_2} \wedge v_{i_4} \eqo (v_{i_2} \wedge v_{i_3}) \wedge ( v_{i_3} \wedge v_{i_4})$, we get 
$v_{i_2} \wedge v_{i_4}   \preceq v_{i_3} \wedge v_{i_4}$, 
where we recall that $\preceq$ stands for the genealogical order on $\bbU$. Similarly, we also get 
$v_{i_1} \wedge v_{i_4}   \preceq v_{i_2} \wedge v_{i_4}$. We thus have proved 
\begin{equation}
\label{Sandwirth}
v_{i_1} \wedge v_{i_4}\preceq v_{i_2} \wedge v_{i_4}  \preceq v_{i_3} \wedge v_{i_4}\; .
\end{equation}
Since $v_{i_1} \wedge v_{i_4} \eqo (v_{i_1} \wedge v_{i_3}) \wedge (v_{i_3} \wedge v_{i_4})$, 
the assumption in $(v)$ implies that $v_{i_1} \wedge v_{i_4} \eqo v_{i_1} \wedge v_{i_3}$ and we get $(v)$ by (\ref{Sandwirth}). Similarly, the assumption in $(vi)$ entails $v_{i_1} \wedge v_{i_4} \eqo v_{i_3} \wedge v_{i_4}$ which implies $(vi)$ by (\ref{Sandwirth}). \cqfd 

\smallskip

\noi
\textbf{Proof of Lemma \ref{insert}.} The assumptions of the jumps of $S$ first imply that 
 for all $q\ino \bbN^*$, 
\begin{equation}
\label{overshoot} 
\Big( \,  \sigma_{q} (y) \leko \infty \, \Big)  \; \Longrightarrow \; \Big( \,  \big| S_{v_{\sigma_{q} (y) +1}} \! \! - S_{v_{\sigma_{\! q-1} (y) +1} \wedge v_{\sigma_{q} (y) +1}} \big| \leq  y+ y_0 \, \Big)\; . 
\end{equation}

Next, let $j$ be as in Assumption $(b)$. Since $k\leqo j\leqo k^\prime$, Lemma \ref{elcontr} $(iv)$ implies that $v_{k+1}\wedge v_{j+1} \ino \lgeo v_{k+1}\wedge v_{k^\prime+1} , v_{k+1} \rgeo$ and Assumption $(a)$ entails that 
\begin{equation}  
\label{exhypc} 
|S_{v_{j+1}} \! -\! S_{v_{k+1} \wedge v_{j+1}}| \geko y_2-y_1 \; .
\end{equation}  
Let $q\ino \bbN^*$ be such that $\sigma_{\! q-1} (y) \leqo k \leko \sigma_{q} (y)$ (which is well-defined). Then, we prove that 
\begin{equation}
\label{insert1}
 \sigma_q (y) \leko j\; .
\end{equation}
\emph{Indeed}, let us suppose that $\sigma_{\! q-1} (y) \leqo k \leko j \leqo \sigma_q (y)$. By definition of  $\sigma_q(y)$ (and by (\ref{overshoot}) in the case where $j\eqo \sigma_q(y)$) we get 
\begin{equation}
\label{control1}
 \big| S_{v_{k+1}} \! - S_{v_{k+1} \wedge v_{\sigma_{\! q-1} (y)+1 }} \big| \leqo y \quad \textrm{and} \quad  \big| S_{v_{j+1}} \! - S_{v_{j+1} \wedge v_{\sigma_{\! q-1} (y) +1}} \big| \leqo  y+ y_0\; .
 \end{equation}
If $v_{k+1}\wedge v_{j+1}\eqo v_{j+1} \wedge v_{\sigma_{\! q-1} (y)+1}$, then 
$\big| S_{v_{j+1}} \!\!  - S_{v_{j+1}  \wedge v_{\sigma_{\! q-1} (y) +1}} \big|\eqo \big| S_{v_{j+1}} \!\!  - S_{v_{j+1} \wedge v_{k+1}} \big| \geko y_2\! -\! y_1$ which contradicts the second inequality in (\ref{control1}). Therefore 
$v_{k+1}\wedge v_{j+1}\! \neq \!  v_{j+1} \wedge v_{\sigma_{\! q-1} (y)+1}$. By Lemma \ref{elcontr} $(iv)$, it implies that $v_{k+1}\wedge v_{j+1} \! \in \,  \rgeo v_{j+1} \wedge v_{\sigma_{\! q-1} (y)+1}, v_{j+1}  \lgeo\, $ (note here that $v_{k+1}\wedge v_{j+1}$ must be distinct from $v_{j+1}$ by (\ref{exhypc})). Consequently by Lemma \ref{elcontr} $(i)$, there exists an integer $j_1$ such that $\sigma_{\! q-1} (y) \leko j_1 \leko j$  and 
$v_{k+1}\wedge v_{j+1} \eqo v_{j_1+1}$. Note that $v_{j_1+1} \wedge v_{\sigma_{\! q-1} (y)+1} \eqo v_{j+1} \wedge v_{\sigma_{\! q-1} (y)+1}$. By definition of $\sigma_q (y)$, it entails 
$$ \big| S_{v_{k+1}\wedge v_{j+1}} - S_{v_{j+1} \wedge v_{\sigma_{\! q-1} (y)+1}} \big|=  \big| S_{v_{j_1+1 }} - S_{v_{j_1+1} \wedge v_{\sigma_{\! q-1} (y) +1}} \big| \leqo y . $$
Thus, $| S_{v_{j+1}} \! -\! S_{v_{k+1} \wedge v_{j+1}}| \leqo |S_{v_{j+1}} \! - \! S_{v_{j+1} \wedge v_{\sigma_{\! q-1} (y)+1 } } | +  |S_{v_{k+1} \wedge v_{j+1}}\!  -\!  S_{v_{j+1} \wedge  v_{\sigma_{\! q-1} (y)+1}} | \leqo  2y + y_0$ by (\ref{control1}) which contradicts (\ref{exhypc}). This completes the proof of (\ref{insert1}).  \cq

\smallskip

We thus have proved that $\sigma_{\! q-1} (y) \leqo k\leko \sigma_q (y) \leko j$. We next prove
\begin{equation}
\label{depass1}
\big| S_{v_{j+1}} \! -\! S_{v_{j+1} \wedge v_{\sigma_q (y)+1} }  \big| \geko y_2 -y_1 -2y -y_0 . 
\end{equation}

Before proving (\ref{depass1}), let us show that it implies the desired result: since $y_2 \geko 3y+2y_0+ y_1$, (\ref{depass1}) entails that $\big| S_{v_{j+1}} \! -\! S_{v_{j+1} \wedge v_{\sigma_q (y)+1} }  \big| \geko y+ y_0$. 
It implies that $\sigma_{q+1} (y)\leko j$ by definition of $\sigma_{q+1} (y)$ and by 
(\ref{overshoot}). In particular, we get $\sigma_{q+1} (y)\leko k^\prime$, which is the desired result.

\smallskip

\smallskip

\noi
\textit{Proof of (\ref{depass1}).} Since $v_{j+1} \! \wedge \! v_{\sigma_{q} (y)+1}\ino 
\lgeo  \varnothing,  v_{\sigma_{q} (y)+1}\rgeo$, there are only two possible cases to consider: 
either $v_{j+1}\!  \wedge  v_{\sigma_{q} (y)+1}\!  \ino \lgeo  \varnothing,  v_{\sigma_{\! q-1} (y)+1} \! \wedge \! v_{\sigma_{q} (y)+1}\rgeo$ or 
$v_{j+1}\!  \wedge  v_{\sigma_{q} (y)+1} \!\!  \in\,  \rgeo  v_{\sigma_{\! q-1} (y)+1}\! \wedge \! v_{\sigma_{q} (y)+1}, v_{\sigma_{q} (y)+1} \rgeo$.

We first suppose that $v_{j+1}\!  \wedge  v_{\sigma_{q} (y)+1}\!  \ino \lgeo  \varnothing,  v_{\sigma_{\! q-1} (y)+1} \! \wedge \! v_{\sigma_{q} (y)+1}\rgeo$. Since $\sigma_{\! q-1} (y) \leqo k\leko \sigma_q (y) \leko j$, Lemma 
\ref{elcontr} $(vi)$ implies that 
$v_{k+1}\wedge v_{j+1} \eqo v_{j+1} \wedge v_{\sigma_{q} (y)+1}$, which entails (\ref{depass1}) by (\ref{exhypc}). 

We next assume that $v_{j+1}\!  \wedge  v_{\sigma_{q} (y)+1} \!\!  \in\,  \rgeo  v_{\sigma_{\! q-1} (y)+1}\! \wedge \! v_{\sigma_{q} (y)+1}, v_{\sigma_{q} (y)+1} \rgeo$. By Lemma \ref{elcontr} $(i)$, there exists $j_2$ such that 
$\sigma_{\! q-1} (y) \leko j_2 \leqo  \sigma_{q} (y)$ and $v_{j+1} \wedge v_{\sigma_{q} (y)+1} \eqo v_{j_2+1}$. 
Thus, $v_{j_2+1} \wedge v_{\sigma_{\! q-1} (y)+1}\eqo v_{\sigma_{\! q-1} (y)+1}\wedge v_{\sigma_{q} (y)+1}$. 
By definition of $\sigma_q (y)$ (and by (\ref{overshoot}) in the case where $j_2\eqo \sigma_q (y)$), we get 
\begin{equation}
\label{control2}
 \big| S_{v_{j+1} \wedge v_{\sigma_{q} (y)+1}} \! -\! S_{v_{\sigma_{\! q-1} (y)+1}\wedge v_{\sigma_{q} (y)+1}} \big|= \big| S_{v_{j_2+1}} \! -\! S_{v_{\sigma_{\! q-1} (y)+1}\wedge v_{j_2+1}} \big| \leq y+ y_0 \; .
\end{equation}
Since $v_{j+1} \wedge v_{\sigma_{q} (y)+1} \! \in\,  \rgeo  v_{\sigma_{\! q-1} (y)+1}\wedge v_{\sigma_{q} (y)+1}, v_{\sigma_{q} (y)+1} \rgeo$ and $\sigma_{\! q-1} (y) \leqo k\leko \sigma_q (y) \leko j$, Lemma  \ref{elcontr} $(v)$ implies $v_{k+1}\wedge v_{j+1} \ino \lgeo  v_{\sigma_{\! q-1} (y)+1}\wedge v_{\sigma_{q} (y)+1}, v_{\sigma_{q} (y)+1} \lgeo\, $. We next check that 
\begin{equation}
\label{control3}
\big| S_{v_{k+1} \wedge v_{j+1} } \! -\! S_{v_{\sigma_{\! q-1} (y)+1}\wedge v_{\sigma_{q} (y)+1}} \big| \leq y \; .
\end{equation}
\emph{Indeed}, it is immediate if 
$v_{k+1}\wedge v_{j+1}\eqo  v_{\sigma_{\! q-1} (y)+1}\wedge v_{\sigma_{q} (y)+1}$. Otherwise 
Lemma 
\ref{elcontr} $(i)$ implies that there is $j_3$ such that $\sigma_{\! q-1} (y) \leqo j_3 \leko  \sigma_{q} (y)$ and $v_{k+1}\wedge v_{j+1} \eqo v_{j_3+1}$. 
Therefore, $v_{j_3+1} \wedge v_{\sigma_{\! q-1} (y)+1}\eqo v_{\sigma_{\! q-1} (y)+1}\wedge v_{\sigma_{q} (y)+1}$ and thus 
$| S_{v_{k+1} \wedge v_{j+1} } \! -\! S_{v_{\sigma_{\! q-1} (y)+1}\wedge v_{\sigma_{q} (y)+1}} |\eqo | S_{v_{j_3+1}} \! -\! S_{v_{\sigma_{\! q-1} (y)+1}\wedge v_{j_3+1}} \big| \leqo y$,  
by definition of $\sigma_q (y)$, which proves (\ref{control3}).  

We now complete the proof of (\ref{depass1}) (and thus, the proof of the lemma as previously noticed) by observing first that 
$$ y_2 <  \big|  S_{v_{j+1}} \! -\! S_{v_{k+1}}\big| \leq  
\big|  S_{v_{j+1}} \! -\! S_{v_{k+1}\wedge v_{j+1} }\big| +
 \big|  S_{v_{k+1}} \! -\! S_{v_{k+1}\wedge v_{j+1} }\big| \leq  \big|  S_{v_{j+1}} \! -\! S_{v_{k+1}\wedge v_{j+1}}\big| + y_1$$
 by Assumptions $(b)$ and $(a)$. Then 
\begin{eqnarray*}
 \big|  S_{v_{j+1}} \! -\! S_{v_{k+1}\wedge v_{j+1}}\big| \leq  \big| S_{v_{j+1}}  \!\!  \!\! \!\!  & -&  \!\! \!\! \!  S_{v_{j+1} \wedge v_{\sigma_{q} (y)+1}}  \big|  + \big| S_{v_{j+1} \wedge v_{\sigma_{q} (y)+1}} \! -\! S_{v_{\sigma_{\! q-1} (y)+1}\wedge v_{\sigma_{q} (y)+1}} \big| \\ 
 &+&  \big| S_{v_{\sigma_{\! q-1} (y)+1}\wedge v_{\sigma_{q} (y)+1}}  \! -\! S_{v_{k+1} \wedge v_{j+1} } \big| \\
& \leq &   \big| S_{v_{j+1}}  \! -\!  S_{v_{j+1} \wedge v_{\sigma_{q} (y)+1}}  \big| \; + \; \overbrace{y+y_0}^{\textrm{by (\ref{control2})}} \; + \overbrace{ y }^{ \textrm{by (\ref{control3})} }.
\end{eqnarray*}
It completes the proof of (\ref{depass1}) and, as already mentionned, also the proof of the lemma. \cqfd 

\smallskip

We next consider BRWs which are indexed by infinite GW-forests and whose jumps are i.i.d.~and we show that their depth-first times of increase exhibit useful renewal properties. 
\begin{lemma}
\label{timeincrease2} Let $\mu$ be a probability on $\bbN$ such that $\mu(0)+ \mu(1)\leko 1$ and $\sum_{k\in \bbN} k\mu(k)\eqo 1$. Let $\bS (p)\eqo (S_u(p))_{u\in \tau(p)}$, $p\ino \bbN^*$, be a sequence of i.i.d.~real valued BRWs such that $\tau(p)$ is a GW($\mu$)-tree and conditionally given $\tau (p)$, the jumps of $\bS(p)$ are i.i.d.~centered and not a.s.~null. We denote by $\bS\eqo (S_u)_{u\in \tau}$ the BRW associated with $(\bS(p))_{p\in \bbN^*}$ as in Remark \ref{contord} $\textbf{(b)}$.
Let $y\ino \bbR_+^*$ and let $(\bs_q(y))_{q\in \bbN}$ the depth-first $y$ oscillation times of $\bS$ (see Definition \ref{exittime}). Then, a.s.~$\bs_q(y)\leko \infty$ and the r.v.s $(\bs_{q+1} (y) \! -\! \bs_q (y))_{q\in \bbN}$ are i.i.d. 
\end{lemma}
\noi
\textbf{Proof.} First observe that since $\mu (0) \! \neq \! 1$ and since 
the jumps of $\bS (p)$ are i.i.d.~centered random variables which are not a.s.~equal to $0$, we must have $\bP (\exists u\ino \tau (p): |S_u(p)| \geko y) \geko 0$. This implies a.s.~that $\bs_1(y)\leko \infty$. 

Let $(u_l)_{l\in \bbN}$ be the depth-first exploration of the vertices of $\tau$. 
Since $\tau$ is a forest of finite trees, $V_\cdot (\tau)$ completely encodes $\tau$ and we easily see that  
there is a measurable function $F\! : \! (\bbZ\! \times \! \bbR)^{\bbN^*}\!\!  \rightarrow \! \bbN\cup \{ \infty\} $ such that a.s.~$\bs_1 (y)\eqo F( (V_l(\tau), \xi_{u_{l+1}})_{l\in \bbN^*})$. 
Moreover, note that $\bs_1 (y)$ is a stopping time with respect to the filtration $(\mathscr F_l)_{l\in \bbN}$, where $\mathscr F_l$ stands for the sigma field generated by the r.v.~$V_{l^\prime}(\tau)$ and  $\xi_{u_{l^\prime +1}}$, $1\leqo  l^\prime \leqo l$. A simple recursive argument then entails that the $\bs_q(y)$ are $(\mathscr F_l)_{l\in \bbN}$-stopping times.

We next describe the forest of finite trees whose \L{}ukasiewicz path is 
 $(V_{l_0+l} (\tau)\! -\! V_{l_0} (\tau))_{l\in \bbN}$. 
To that end, we denote by $(w_p)_{p\in \bbN^*}$ the $<_{\mathtt{lex}}$-increasing sequence of vertices of the set 
$\{ w\ino \tau: u_{l_0+1}\!  \leq_{\mathtt{lex}}\!  w \; \textrm{and} \; \overleftarrow{w} \ino \lgeo \varnothing , u_{l_0+1} \rgeo \}$. 
We then denote by $\tau_{l_0}$ the forest formed by the trees $\theta_{w_p} \tau$, $p\ino \bbN^*$ as in Definition \ref{Ulamtree} \textbf{(c)}: namely, $k_{\varnothing} (\tau_{l_0}) \eqo \infty$, 
and $\theta_{[p]} \tau_{l_0} \eqo \theta_{w_p} \tau$ for all $p\ino \bbN^*$. 
Let $(u^\prime_l)_{l\in \bbN}$ be depth-first exploration of the vertices of $\tau_{l_0}$.  
Then, for all $l\ino \bbN^*$, there are unique $p\ino \bbN^*$ and $v\ino\bbU$ such that 
$u_{l_0+l}\eqo w_p\! \ast \! v$ and $u^\prime_{l}=[p]\! \ast \! v$. Consequently, there is a one-to-one 
$<_{\mathtt{lex}}$-increasing correspondence between $(u_{l_0+ l})_{l\in \bbN}$ and $(u^\prime_{l})_{l\in \bbN} $. 
Note that $k_{u_{l_0+l}} (\tau)\eqo k_{u^\prime_l} (\tau_{l_0})$ for all $l\ino \bbN^*$. 
For all $l\ino \bbN^*$, we then set $\xi^\prime_{u^\prime_l} \eqo \xi_{u_{l_0+l}}$, so we get 
$$\big( V_l (\tau_{l_0}) \, ,\,  \xi^\prime_{u^\prime_{l+1}}\big)_{l\in \bbN}= \big( V_{l_0+l} (\tau)\! -\! V_{l_0} (\tau)\, , \,  \xi_{u_{l_0+l+1}} \big)_{l\in \bbN}\; . $$ 
Then observe that the process on the right hand side has the same law as $(V_l (\tau), \xi_{u_{l+1}})_{l\in \bbN}$ and 
is independent from $\mathscr F_{l_0}$ (since $V_\cdot (\tau)$ is a random walk and since the 
jumps of $\bS$ are i.i.d.). 
 
 Then, denote by $(S^\prime_{u^\prime})_{u^\prime \in \tau_{l_0}}$ the BRW starting at $0$ whose jumps are $(\xi^\prime_{u^\prime})_{u^\prime\in \tau_{l_0}\backslash\{\varnothing\}}$. 
Observe that for all $l\ino \bbN$, $S^\prime_{u^\prime_{l+1}}\!\!\!\! =\!  S_{u_{l_0+l+1} } \! -\! S_{u_{l_0+1} \wedge u_{l_0+l+1}}$. Therefore  
$$ \inf \big\{ l\ino \bbN: 
|S^\prime_{u^\prime_{l+1}} | \! >\!  y \big\} \eqo  
 F \big( (V_{l+l_0}(\tau) \! -\! V_{l_0} (\tau), \xi_{u_{l_0+l+1}})_{l\in \bbN^*}\big) .$$
In particular (since the event $\{ \bs_q(y)\eqo l_0\}$ belongs to $\mathscr F_{l_0}$) if $\bs_q (y)\eqo l_0$, then
$\bs_{q+1} (y)\! -\! \bs_q (y)\eqo  
F \big( (V_{l+l_0}(\tau) \! -\! V_{l_0} (\tau), \xi_{u_{l_0+l+1}})_{l\in \bbN^*}\big)$ and the previous arguments imply that $\bs_{q+1} (y)\! -\! \bs_q (y)$ is independent from $\mathscr F_{\! \bs_q (y)}$ and that it has the same law as $\bs_1 (y)$, which completes the proof.  \cqfd

\subsection{Analytical estimates on BRWs indexed by GW-trees}
\label{reBRWGWsec}

In this section, we provide analytical estimates to study the extremal position of a $\bbR$-valued BRW $\bS\eqo (S_u)_{u\in \tau}$. More specifically, we prove in Lemma \ref{keylower} an analog for BRWs of Lemma \ref{maxBrosnadiscr}.
Here we overall assume the following:  
\begin{compactenum}

\smallskip

\item[\textbf{(a)}] 
The indexing tree $\tau$ is a GW($\mu$)-tree with $\sum_{k\in \bbN} k\mu (k) \eqo 1$ and $\, \mu (0)\! +\!  \mu (1)\leko 1$; 

\smallskip

\item[\textbf{(b)}] 
Conditionally given $\tau$, jumps $(\xi_u)_{u\in \tau \backslash \{ \varnothing\}}$ are i.i.d.~$\bbR$-valued r.v.~whose law is denoted by $\bgam $. 

\smallskip

\end{compactenum} 

\noi
Most of the time $S_\varnothing\eqo 0$ and $\bgam$ is assumed to be 
centered and uniformly bounded. We provide estimates on the law of $\max_{u\in \tau} S_u $ via analytical arguments involving the following functions that are derived from $\mu$ as follows. 
\begin{equation}
\label{regfpsiR}
\forall r\ino [0, 1], \qquad g_\mu (r)\eqo \sum_{k\in \bbN} r^k\mu (k), \quad f_\mu (r)\eqo 1\! -\! g_\mu (1\! -\! r) \quad \textrm{and} \quad  \Psimu (r)\eqo  r\! -\!  f_{\mu} (r) \; ,
\end{equation}
Recall from Lemma \ref{fungeneprop} that $f_\mu \! : \! [0, 1] \! \rightarrow \!  [0, 1\! -\! \mu (0)]$ is a strictly concave increasing bijection and the two functions $g_\mu\!  : \! [0, 1]\!  \rightarrow \! [\mu (0), 1]$ and $ \Psimu \! : \! [0, 1] \! \rightarrow \! [0, \mu (0)]$ are 
strictly convex increasing bijections. 
We next prove the following. 
\begin{lemma}
\label{fungenepropp} We keep the above notations and assumptions. 
Then, for all $s\ino [0, 1\! -\! \mu (0)]$, we set $R_\mu (s)\eqo f^{-1}_\mu (s)\! -\! s$, where $f^{-1}_\mu $ stands for the inverse of $f_\mu$. Then $R_\mu \! : \!  [0, 1\! -\! \mu (0)] \! \rightarrow [0, \mu (0)]$ is a strictly convex increasing bijection such that $R^\prime_\mu (0)\eqo 0$ and 
$R^\prime_\mu (1\! -\! \mu (0)) \eqo \tfrac{1}{\mu (1)} \! -\! 1$, this quantity being infinite if $\mu (1)\eqo 0$.  
\end{lemma}
\noi
\textbf{Proof.} Observe that for all $s\ino [0, 1 \! -\! \mu (0)]$,  
$$ R^{\prime }_\mu (s) \eqo \frac{1}{f^\prime_\mu (f^{-1}_\mu (s))} -1 \quad \textrm{and} \quad  R^{\prime \prime}_\mu (s) \eqo - \frac{f^{\prime \prime}_\mu (f^{-1}_\mu (s))}{ (f^\prime_\mu (f^{-1}_\mu (s))^3}  \; , $$
which entails the desired result. \cqfd

\smallskip

As mentioned at the begining of the section, we provide estimates of the law of $\max_{u\in \tau} S_u $ by studying the discrete \emph{exit process of $\bS$}, as we have done for the Brownian snake $\mathscr W$ in Section Lemma \ref{maxBrosnadiscr}. Of course here, these estimates are more technical to obtain. 
More precisely, we define the discrete exit process
$\mathcal Z_a$ of $\bS$ from the interval $(-\infty, a)$ by setting 
\begin{equation}
\label{ZZaadef}
\forall a\ino \bbR, \quad
 \mathcal Z_a \eqo \# \big\{ u\ino \tau : \; S_u \geqo a  \; \textrm{and} \;  \forall v \ino \lgeo \varnothing , u \lgeo \, , S_v \leko a \big\}  \; .
 \end{equation}
Observe that $\mathcal Z_a \eqo 1 $ if $a\ino (-\infty , S_\varnothing]$ and that $\mathcal Z_a \eqo 0$ iff $ \max_{u\in \tau} S_u \leko a$. 
For all $r\ino [0, 1)$, we also set 
\begin{equation}
\label{defvarphia}
\varphi_a (r) \eqo 1\! -\! \bE \big[(1\! -\!  r)^{\mathcal Z_a}\big] \quad \textrm{and} \quad  
\varphi_a (1)\eqo \lim_{r\rightarrow 1^-} \varphi_a (r)\eqo \bP \big(  \max_{u\in \tau} S_u \geqo a \big) \; .
\end{equation}

\vspace{-3mm}

\begin{lemma}
\label{mainequa} We keep the previous notations and assumptions and we suppose that $S_\varnothing \eqo 0$. We also assume that a.s.~$\limsup_{n\to \infty} \mathtt{U}_n \eqo \infty$, where $(\mathtt{U}_n)_{n\in \bbN}$ is a random walk with independent jumps whose law is $\bgam$. 
\begin{compactenum}

\smallskip

\item[$(i)$] 
 For all $a\ino \bbR_+^*$ and all $r\ino [0, 1]$, 
\begin{equation}
\label{equaZa}
\bE [\mathcal Z_a]\eqo 1 \quad \textrm{and} \quad R_\mu \big( \varphi_a (r) \big) \eqo \int_{\bbR}  \!\!  \bgam ( \mathrm dy) \, \big(\varphi_{a-y} (r) \! -\! \varphi_a (r) \big) \
\end{equation}
where $ R_\mu$ is defined in Lemma \ref{fungenepropp}. 

\smallskip

\item[$(ii)$] The function $a\ino \bbR \! \mapsto \! \varphi_a (1)$ is nonincreasing and 
\begin{equation}
\label{monovarphi}
\varphi_a (1) \geqo \varphi_b (\varphi_{a-b} (1)) , \quad a\ino \bbR_+^*, \; b \ino (0, a] \; .
\end{equation}

\item[$(iii)$] If  
$\bgam\eqo \tfrac{_1}{^2} \rho \delta_{-1}\! +\!  (1\! -\! \rho ) \delta_0 \! +\!  \tfrac{_1}{^2}  \rho \delta_1$, 
with $\rho \ino (0, 1]$, then $(\mathcal Z_{n})_{n\in \bbN}$ is a critical GW process and 
$\varphi_{n+m} (r)\eqo \varphi_n (\varphi_m (r)) $, $m, n\ino \bbN$ and $r\ino [0, 1]$. 
\end{compactenum}

\end{lemma}
\noi
\textbf{Proof.} 
By the 
many-to-one formula for critical GW($\mu$)-trees and the independence of the jumps of the branching random walk $(S_u)_{u\in \tau}$, we get 
$$ \bE \Big[\sum_{u\in \tau} F \big( S_v\, ; \, v\ino \lgeo \varnothing , u \rgeo   \big) \Big]= \sum_{n\in \bbN} \bE \big[ F \big( \mathtt{U}_k\, ; \, k\ino \{ 0, \ldots, n\}   \big) \big]  \; .$$
Let $T_a \eqo \inf\{ n\ino \bbN: \mathtt U_n \ino [a, \infty) \}$, with the convention that $\inf \emptyset \eqo \infty$. Then, the previous equality implies that 
$\bE [\mathcal Z_a]\eqo \sum_{n\in \bbN} \bP (T_a \eqo n ) \eqo \bP (T_a \leko \infty) \eqo 1$ by our assumption, which proves the first point of (\ref{equaZa}).

 Let us prove the second point of (\ref{equaZa}). For the need of the proof, we set 
 $\bP_{\! x} \eqo \bP (\, \cdot \, | S_\varnothing \eqo x)$, $x\ino \bbR$. We keep denoting $\bP_{\! 0}$ simply by $\bP$ and we also set $\phi_a (r)\eqo \bE [r^{\mathcal Z_a}]$. Then, for all $a\ino \bbR_+^*$ and all $r\ino [0, 1]$, we first get  
 $$ \phi_a(r) =  \bE \Big[ \bE \big[ r^{\mathcal Z_a} \big| k_\varnothing (\tau) , S_{[1]}, \ldots, S_{[k_\varnothing (\tau)]} \big] \Big]= \bE \Big[  \prod_{j=1}^{k_\varnothing (\tau)} \bE_{S_{[j]}} [r^{\mathcal Z_a}] \Big] =
\bE \Big[  \prod_{j=1}^{k_\varnothing (\tau)} \phi_{a-S_{[j]}} (r)  \Big]  \; , $$
because for all $x\ino \bbR$, $(x+S_u)_{u\in \tau}$ under $\bP$ has the same law as $(S_u)_{u\in \tau}$ under $\bP_{\! x}$ and therefore $\mathcal Z_{a-x}$ under $\bP $ has the same law as $\mathcal Z_{a}$ under $\bP_{\! x}$. Here, the products are taken equal to $1$ when $k_\varnothing (\tau)\eqo 0$. 
Consequently, we get
\begin{equation}
\label{gournigol}  \forall a\ino \bbR_+^*, \quad \phi_a (r) = \bE \Big[  \Big( \int_{\bbR}   \!\!  \bgam ( \mathrm dy) \, \phi_{a-y} (r) \Big)^{k_\varnothing (\tau)} \Big] = g_\mu \Big(  \int_{\bbR}  \!\! \bgam ( \mathrm dy)\,  \phi_{a-y} (r) \Big), 
\end{equation}
which easily entails the second point of (\ref{equaZa}). 

   We next prove $(ii)$. Let $a\ino \bbR_+^*$ and $b\ino (0, a]$. We set $G_b\eqo \big\{ u\ino \tau : \; S_u \geqo b  \; \textrm{and} \;  \forall v \ino \lgeo \varnothing , u \lgeo \, , S_v \leko b \big\}$. We denote by $\mathscr H$ the sigma field generated by $G_b$ and the $S_u$, $u\ino G_b$. By a direct computation, we see that conditional on $\mathscr H$, the branching random walks 
$(S_{u\ast v} \! -\! S_u)_{v\in \theta_u \tau}$, $u\ino G_b$ are independent copies of $(S_v)_{v\in \tau}$ 
under $\bP$. Thus we have 
$$ \phi_a (0)= \bP \big( \! \max_{v\in \tau} S_v \leko a\big) \eqo \bE \Big[ \prod_{u\in G_b} \phi_{a-S_u} (0)\Big] \leqo \bE \big[  (\phi_{a-b} (0))^{\# G_b}\big] = \phi_b (\phi_{a-b} (0))$$ 
since $S_u \geqo b$ for all $u\ino G_b$, since $c\! \mapsto \! \phi_c (0) $ is nondecreasing and since $\# G_b \eqo \mathcal Z_b$ by definition (here, the product is taken equal to $1$ if $G_b $ is empty). This immediately implies (\ref{monovarphi}). 

We now proceed to the proof of $(iii)$. 
With the same notations as above, we first get 
$$ \textrm{$\bP$-a.s.} \quad \bE \big[ r^{\mathcal Z_a} \big| \mathscr H\big]= \prod_{u\in G_b} \phi_{a-S_u} (r) \; .$$
Then, we fix $m, n\ino \bbN$ and we take $a\eqo m+n$ and $b\eqo n$. Since $S_u\eqo n$ for all $u\ino G_n$,
$$ \textrm{$\bP$-a.s.} \quad \bE \big[ r^{\mathcal Z_{n+m}} \big| \mathscr H\big]= \phi_{m} (r)^{\# G_n}= \phi_m (r)^{\mathcal Z_n} \; ,$$
which implies the desired result.  \cqfd

\medskip

The next lemma is the key analytical estimate used in the proof of Theorem \ref{unifbounded}. 
\begin{lemma}
\label{keylower} We assume that $\bS\eqo (S_u)_{u\in \tau}$ satisfies \emph{\textbf{(a)}} and \emph{\textbf{(b)}} as in the beginning of the section. We also suppose that $S_\varnothing\eqo 0$, that $\bgam$ is symmetric (i.e.~$\bgam (\!-\mathrm d y) \eqo \bgam (\mathrm d y)$), that $\bgam ([-1,1])\eqo 1$ and that $\gamun$ $ := $ $\bgam (\{ 1\})$ $\geko$ $0$. 
Recall the definition of $\Psimu$ from (\ref{gfpsiR}) and recall from (\ref{defvarphia}) that $\varphi_a (1)\eqo \bP \big(  \max_{u\in \tau} S_u \geqo a \big)$. Then, the following holds true. 

\smallskip

\begin{compactenum}
\item[$(i)$] For all $ a\ino [2, \infty)$
\end{compactenum}
\begin{equation}
\label{lowerbound} 
 \int_{\varphi_{a}(1)}^{\varphi_{a-2} (1)} \!\!\!\!  \frac{\mathrm d s}{\sqrt{\! \int_0^s  \Psimu (\gamun r) \, \mathrm d r}} \, \geq 2\sqrt{\gamun  } \quad \textrm{and thus} \quad  \int_{\varphi_{a}(1)}^{1}  \!\!  \frac{\mathrm d s}{\sqrt{\! \int_0^s \Psimu (\gamun r) \, \mathrm d r }}\,  \geq \tfrac{1}{2}\sqrt{\gamun   } a\; , 
\end{equation}
\begin{compactenum}
\item[$(ii)$] We assume that $\bgam \eqo \tfrac{_1}{^2} (\delta_{-1}+ \delta_1)$. We note that $\varphi^\prime_1 (1) \eqo \bP (\mathcal Z_1\eqo 1)\geqo \frac{_1}{^2} g'_\mu ( \frac{_1}{^2}) \geko 0$, and in particular  $\varphi^\prime_1 (1) \geqo \mu(1) /2$. 
Then, for all $n\ino \bbN$, 
\end{compactenum}
\begin{equation}
\label{upperbound} 
 \int_{\! \varphi_{n+1} (1)}^{\varphi_{n} (1)} \!\! \frac{ ds }{\sqrt{\int_0^{s} \Psimu \big( r \big) \mathrm dr\,  }\, } \leq \frac{2 }{\varphi^\prime_1 (1)^{5/2} } \quad \textrm{and thus} \quad  
\int^1_{\! \varphi_{n} (1)} \!\! \frac{ ds }{\sqrt{\int_0^{s} \Psimu \big( r \big) \mathrm dr \, }\, } \leq  \frac{2 n}{\varphi^\prime_1 (1)^{5/2} }\; , 
\end{equation}
\end{lemma}
\noi
\textbf{Proof.} Let $\xi$ be a r.v.~distributed according to $\bgam$. We set $\zeta\eqo |\xi|$. Since $\bgam$ is symmetric, there exists a r.v.~${\boldsymbol \epp}$, independent of $\zeta$ and such that 
$\bP({\boldsymbol \epp}\eqo \pm 1 ) \eqo \tfrac{1}{2}$ and $\xi \eqo {\boldsymbol \epp} \zeta$. We fix $r\ino [0, 1]$ and to simplify notation we simply denote $\varphi_a (r)$ by $\varphi_a$ for all $a\ino \bbR$. Then, (\ref{equaZa}) can be rewritten as 
$$R_\mu (\varphi_a) \eqo \tfrac{1}{2} \bE [\varphi_{a-\zeta} \! -\! \varphi_a ] \! -\! \tfrac{1}{2}  \bE [\varphi_{a} \! -\! \varphi_{a+\zeta} ] $$
that is a discrete version of the ordinary differential equation $2R_\mu (y(a))\eqo y^{\prime \prime} (a) $, $a\ino [0, \infty)$, where the function $y(\cdot)$ decreases to $0$. The solution of this ODE can be used as a guideline of the proof of Proposition \ref{keylower} $(i)$. It goes as follows:

\begin{compactenum}
\item[$-$]\texttt{Step 1}: $4y^\prime (a)R_\mu (y(a)) \eqo ((y^\prime(a) )^2)^\prime$.

\item[$-$]\texttt{Step 2}: thus, $4\! \int_{y(a)}^{y(a_1)} \! R_\mu (r) \, \mathrm d r\eqo (y^\prime(a_1))^2\! -\! (y^\prime(a))^2$ for all $a$, $a_1\ino [0, \infty)$ such that $a\leqo a_1$; we look at the solution such that 
$y$ and $y^\prime \! \rightarrow \! 0$ at $\infty$.

\item[$-$] \texttt{Step 3}: therefore we get 
$$-2\, \sqrt{\! \int_{0}^{y(a)} \! \!\! \!\! \! \! \!\! R_\mu (r) \, \mathrm d r} \, =   y^\prime(a) \quad \textrm{\emph{and thus}} \quad \int_{y(a)}^{y(0)} \!\! \frac{\mathrm ds}{\sqrt{\! \int_0^s \! R_\mu (r) \, \mathrm dr }}= 2a \; .$$ 
\end{compactenum}

To adapt \texttt{Step 1} to our discrete setting, we use $\zeta^\prime$, an independent copy of $\zeta$ and we get the following: 
\begin{eqnarray*}
\Delta_a (r) & := & 4 R_\mu (\varphi_a) \bE [\varphi_{a-\zeta} \! -\! \varphi_a ] + 4 R_\mu (\varphi_a) \bE [\varphi_{a} \! -\! \varphi_{a+\zeta} ] \\
& =&  2\big(\bE [\varphi_{a-\zeta} \! -\! \varphi_a ] \big)^2 - 2\big(\bE [\varphi_{a} \! -\! \varphi_{a+\zeta} ]  \big)^2 \\
&= & \bE \big[ 2(\varphi_{a-\zeta} \! -\! \varphi_a) (\varphi_{a-\zeta^\prime} \! -\! \varphi_a) \big] -\bE \big[ 2(\varphi_{a} \! -\! \varphi_{a+\zeta} )(\varphi_{a} \! -\! \varphi_{a+\zeta^\prime} ) \big] \\
& =&  \bE \big[  (\varphi_{a-\zeta} \! -\! \varphi_a)^2 \! +\!  (\varphi_{a-\zeta^\prime} \! -\! \varphi_a)^2 \! -\! (\varphi_{a-\zeta^\prime} \! -\! \varphi_{a-\zeta})^2 \big] \\
&  & \qquad  \qquad \!  -  \bE \big[  (\varphi_{a+\zeta} \! -\! \varphi_a)^2 \! +\!  (\varphi_{a+\zeta^\prime} \! -\! \varphi_a)^2 \! - (\varphi_{a+\zeta^\prime} \! -\! \varphi_{a+\zeta})^2 \big]\\
& =& 2 \, \underbrace{\bE \big[  (\varphi_{a-\zeta} \! -\! \varphi_a)^2 \! -\!  (\varphi_{a+\zeta} \! -\! \varphi_a)^2  \big]}_{\Delta^{*}_a(r)} +  \underbrace{\bE \big[ (\varphi_{a+\zeta^\prime} \! -\! \varphi_{a+\zeta})^2 \! -\! (\varphi_{a-\zeta^\prime} \! -\! \varphi_{a-\zeta})^2 \big]}_{\Delta^{**}_a(r)}. 
\end{eqnarray*}

To adapt \texttt{Step 2}, we proceed as follows. 
Let $c\ino (a+1, \infty)$. Since $\bP$-a.s.~$\zeta \ino [0, 1]$, we get 
\begin{eqnarray}
 \int_a^c\!\!   \Delta^{*}_b (r)\, \mathrm d b &\eqo &\bE \Big[\int_a^c \!\! (\varphi_{b-\zeta} \! -\! \varphi_b)^2  \mathrm d b \,  - \int_{a+\zeta}^{c+\zeta} \!\! (\varphi_{b-\zeta} \! -\! \varphi_b)^2  \mathrm d b  \, \Big] \nonumber \\
 \label{Delta1}&\eqo & \bE \Big[\int_a^{a+\zeta} \!\!\! (\varphi_{b-\zeta} \! -\! \varphi_b)^2  \mathrm d b\,  \Big]- \bE \Big[\int_c^{c+\zeta} \!\! \! (\varphi_{b-\zeta} \! -\! \varphi_b)^2  \mathrm d b \, \Big]. 
 \end{eqnarray}
We set $\zeta^\circ \eqo \tfrac{1}{2} (\zeta + \zeta^\prime)$ and $\zeta^\bullet\eqo  \tfrac{1}{2} (\zeta \! -\!  \zeta^\prime) $, then similarly: 
\begin{equation} \label{Delta2}
\int_a^c\!\!   \Delta^{**}_b (r)\, \mathrm d b = \bE \Big[\int_{c-\zeta^\circ}^{c+\zeta^\circ} \!\!\! (\varphi_{b-\zeta^\bullet} \! -\! \varphi_{b+\zeta^\bullet})^2  \mathrm d b\,  \Big]- \bE \Big[\int_{a-\zeta^\circ}^{a+\zeta^\circ} \!\!\! (\varphi_{b-\zeta^\bullet} \! -\! \varphi_{b+\zeta^\bullet})^2  \mathrm d b\,  \Big]. 
 \end{equation}
We now take $r\eqo 1$. Recall from Lemma \ref{mainequa} that 
$a\ino \bbR \! \mapsto \! \varphi_a (1)$ is nonincreasing. Thus $\Delta_a (1) \geqo 0$. Since clearly 
$\lim_{a\rightarrow \infty} \varphi_a (1)\eqo 0$, (\ref{Delta1}) and (\ref{Delta2}) imply that 
\begin{equation}
\label{voila0}
 \int_a^\infty \!\!\!\!    \Delta_b(1)\, \mathrm d b = 2 \, \bE \Big[\int_a^{a+\zeta} \!\!\!\!\!\!  \big(\varphi_{b-\zeta} (1) \! -\! \varphi_b(1) \big)^2  \mathrm d b\,  \Big]- 2\bE \Big[\int_{a-\zeta^\circ}^{a+\zeta^\circ} \!\!\!\!\!\! \!\!\! \big(\varphi_{b-\zeta^\bullet}(1) \! -\! \varphi_{b+\zeta^\bullet} (1) \big)^2  \mathrm d b\,  \Big]. 
\end{equation}
Since $\Delta_b(1) \geqo 4 R_\mu \big( \varphi_b(1) \big) \bE [\varphi_{b-\zeta} (1)  \! -\! \varphi_b(1)  ]$, (\ref{voila0}) implies 
\begin{equation}
\label{Deltaineq} 
2 \,  \int_a^\infty \!\!\!\!  R_\mu \big( \varphi_b(1) \big)\,  \bE [\varphi_{b-\zeta} (1)  \! -\! \varphi_b(1)  ]\, \mathrm d b\;  \leq 
 \, \bE \Big[\int_a^{a+\zeta} \!\!\!\!\!\!  \big(\varphi_{b-\zeta} (1) \! -\! \varphi_b(1) \big)^2  \mathrm d b\,  \Big]
 \end{equation}

We next adapt \texttt{Step 3}, which is more complicated than the two previous steps. We first need to modify the left member of (\ref{Deltaineq}).  
To that end, we introduce for all 
$z\ino (0, 1)$ and all $c\ino [0, 1] $
\begin{equation}
\label{defbeta}
\beta(z)\eqo \inf \big\{b\ino \bbR_+\, :  \; \varphi_b (1) \leko z   \big\} \quad \textrm{and} \quad I_c(z)\eqo \big\{ b\ino \bbR_+\, : \; \varphi_b (1) \leko z \leqo \varphi_{b-c} (1) \big\} \; , 
\end{equation}
that are well-defined since $\varphi_0 (1)\eqo 1$ and $\lim_{b\rightarrow \infty} \varphi_b (1)\eqo 0$ (note here that $I_c(z)$ may be possibly empty). We next prove the following. 
\begin{equation}
\label{girnogoul} 
2\gamun \!\! \int_0^{\varphi_{a-1} (1) }\!\!\!\!\!\!\!\! \!\!\!\!  \!\!\!\!   \mathrm d z \, \, 
 \Psimu \big( \gamun z \big)  \ell \big( (a, \infty) \! \cap \! I_1 (z) \big) \, \leq \,   \int_a^\infty \!\!\!\!  R_\mu \big( \varphi_b(1) \big)\,  \bE [\varphi_{b-\zeta} (1)  \! -\! \varphi_b(1)  ]\, \mathrm d b\; 
\end{equation}
where $\ell$ stands for the Lebesgue measure on the real line. 

\smallskip

\noi
\emph{Proof of (\ref{girnogoul}).} Let $z\ino (0, 1)$ and $c\ino [0, 1] $. First note that  
\begin{equation}
\label{Isimplif} (\beta (z), \beta (z) +c) \subset I_c(z) \subset [\beta (z), \beta (z) +c] \; .
\end{equation}
By definition of $I_\zeta (z)$, we next observe that for all $b\geq a$, 
\begin{equation}
\label{gornogol} 
\varphi_{b-\zeta}(1)   - \varphi_b (1) = \int_0^\infty \!\!\!\!  \un_{\{ b \in I_\zeta (z) \} }  
\, \mathrm d z \, = \, \int_0^{\varphi_{a-\zeta} (1) }\!\!\!\!\!\!\!\! \!\!\!\!\un_{\{ b \in I_\zeta (z) \} }  \, \mathrm d z\; .
\end{equation}
We thus get 
$\int_a^\infty R_\mu \big( \varphi_b (1) \big)  \big(\varphi_{b-\zeta} (1)  \! -\! \varphi_b(1)\big) \mathrm d b 
= \int_a^\infty \mathrm d b \int_0^{\varphi_{a-\zeta} (1) } \!   \mathrm d z \, R_\mu \big( \varphi_b (1) \big) \un_{\{ b\in I_\zeta (z)  \}}$. If $b\ino I_\zeta (z)$, then $\varphi_b (1) \geqo \varphi_\zeta (\varphi_{b-\zeta} (1))\geqo 
\varphi_\zeta (z)$ by (\ref{monovarphi}) and since $r\! \mapsto \! \varphi_\zeta (r)$ is nondecreasing. Since $R_\mu $ increases too (Lemma \ref{mainequa}), we get  
\begin{eqnarray}
\bE \Big[ \int_a^\infty R_\mu \big( \varphi_b (1) \big)  \big(\varphi_{b-\zeta} (1)  \! -\! \varphi_b(1)\big) \mathrm d b \Big]  
\!\!\!\! & \geq & \!\!\!\! \bE \Big[  \int_a^\infty \!\!\!\!  \mathrm db  \int_0^{\varphi_{a-\zeta} (1) }\!\!\!\!\!\!\!\! \!\!\!\! \!\!\!\!    \mathrm d z \, \, R_\mu \big( \varphi_\zeta (z) \big)\un_{\{ b \in I_\zeta (z) \} }  \Big]  \nonumber \\
\!\!\!\! & \geq & \!\!\!\! 
\bE \Big[ \int_0^{\varphi_{a-\zeta} (1) }\!\!\!\!\!\!\!\! \!\!\!\!  \!\!\!\!   \mathrm d z \, \, R_\mu \big( \varphi_\zeta (z) \big) 
\ell \big( (a, \infty) \! \cap \! I_\zeta (z) \big)  \Big] \nonumber \\
\label{gornigoul} \!\!\!\! & \geq & \!\!\!\! 
2\gamun \int_0^{\varphi_{a-1} (1) }\!\!\!\!\!\!\!\! \!\!\!\!  \!\!\!\!   \mathrm d z \, \, R_\mu \big( \varphi_1 (z) \big) 
\ell \big( (a, \infty) \! \cap \! I_1 (z) \big)
\end{eqnarray}
by Fubini and since  
$\bP (\zeta\eqo 1)\eqo 2\gamun$. Recall from (\ref{gfpsiR}) and from the definition of $R_\mu$ in Lemma \ref{fungenepropp} that $R_\mu \eqo \Psimu \circ f^{-1}_\mu $. Then, (\ref{gournigol}) easily implies that $\varphi_b (z) \eqo f_\mu \big(  \int_{[-1, 1]} \! \bgam (\mathrm d y ) \varphi_{b-y} (z) \big) $  
and thus 
\begin{equation}
\label{Rpsifourbi}
R_\mu \big( \varphi_1 (z) \big) \eqo \Psimu \big( f^{-1}_\mu  (\varphi_1 (z)\big) \big) = \Psimu \Big(\int_{[-1, 1]}  \!\!\!\! \!\!\!\! \!\!  \bgam ( \mathrm dy) \, \varphi_{1-y} (z)  \Big) \geq 
\Psimu \big( \gamun\varphi_0 (z)\big)= \Psimu \big(  \gamun z \big)  \; , 
\end{equation}
which implies (\ref{girnogoul}) by (\ref{gornigoul}). \cq 

\smallskip

To adapt \texttt{Step 3}, we next need to modify the right member of (\ref{Deltaineq}). 
To that end, we set $\overline{S}\eqo \max_{u\in \tau} S_u$ and we recall that $\varphi_b(1)\eqo \bP (\overline{S} \geqo b)$. 
Thus $b\mapsto \varphi_b(1)$ is left continuous and its right-limit $\varphi^+_b(1)$ is equal to $\bP (\overline{S}\geko b)$. 
We first prove that 
\begin{equation}
\label{phipm}
\forall a\ino [1, \infty) , \quad \varphi_{a-1}^+ (1)\!  -\! \varphi_{a+1}(1)= \bP \big( \, \overline{S}  \ino (a\! -\! 1, a+1) \big)> 0\; .
\end{equation}

\smallskip

\noi
\emph{Proof of (\ref{phipm})}. Since $a\ino [1, \infty)$, there is $p\ino \bbN^*\cap (a\! -\! 1, a+1)$. Clearly, we can find a finite tree 
$t\ino \bbT$ such that $\bP (\tau \eqo t)\geko 0$ and such that there exists $u_*\ino t$ with $|u_*|\eqo p$. 
Then, for all $u\ino t\backslash \{ \varnothing\}$, we set $\epp_u\eqo \un_{\{ u\in \lgeo \varnothing, u_* \rgeo\}} \! -\!  \un_{\{ u\notin \lgeo \varnothing ,u_* \rgeo\}} $. Then, on the event 
$A\eqo \{ \tau \eqo t\} \cap \{\forall u\ino t\backslash \{ \varnothing\},  \xi_u\eqo \epp_u \}$, we get $\overline{S}\eqo S_{u*}\eqo p$. Next observe that $\bP (A)\eqo ( \gamma_1)^{\# t-1} \bP (\tau \eqo t) \geko 0$, and we get the desired result.  \cq

\smallskip

We next prove that 
\begin{equation}
\label{Aupbound}
\forall a\ino [1, \infty) , \quad \Big( \bE \Big[\int_a^{a+\zeta} \!\!\!\!\!\!  \big(\varphi_{b-\zeta} (1) \! -\! \varphi_b(1) \big)^{\! 2}  \mathrm d b\,  \Big] \Big)^{\! \frac{1}{2}} \, \leq \! \int_{\varphi_{a+1}(1) }^{\varphi_{a-1}^+ (1)} \!\!\!\!\! \!\!\!\!\! \!\! \mathrm d s\,  \sqrt{\ell (I_1 (s)\!  \cap \! [a, a+1] ) } \; , 
\end{equation}
\noi
\emph{Proof of (\ref{Aupbound})}. To simplify notations, we denote by $A$ the left hand side of (\ref{Aupbound}). By (\ref{gornogol}), Fubini and a repeated use of Cauchy-Schwarz inequality, we get the following
\begin{eqnarray*}
A^2 & \overset{\textrm{by (\ref{gornogol})}}{=} & \bE \Big[\int_a^{a+ \zeta} \!\!\!\!  \!\!\!\! \mathrm d b \int_0^\infty \!\!  \!\!\!\!  \mathrm d z \!  \int_0^\infty \!\!  \!\!\!\!  \mathrm d z^\prime \, \un_{\{ b\in I_\zeta (z) \}} \un_{\{ b\in I_\zeta (z^\prime) \}}   \Big] \\
& =& \bE \Big[ \int_0^\infty  \!\!  \!\!\!\!   \mathrm d z \!  \int_0^\infty  \!\!  \!\!\!\!   \mathrm d z^\prime   \! \int_a^{a+ \zeta} \!\!\!\!  \!\!\!\! \mathrm d b  \, \un_{\{ b\in I_\zeta (z) \}} \un_{\{ b\in I_\zeta (z^\prime) \}}   \Big] \\
& \leq & \bE \left[ \int_0^\infty \!\! \!\!\!\!  \mathrm d z \!  \int_0^\infty \!\! \!\!\!\!   \mathrm d z^\prime  
\sqrt{\int_a^{a+ \zeta} \!\!\!\!  \!\!\!\!  \mathrm d b \,  \un_{\{ b\in I_\zeta (z) \}} } \, \sqrt{\int_a^{a+ \zeta} \!\!\!\! \!\!\!\!  \mathrm d b \,   \un_{\{ b\in I_\zeta (z^\prime) \}}  }\,  \right] \\
& =& \int_0^\infty \!\! \!\!\!\!  \mathrm d z \!  \int_0^\infty \!\! \!\!\!\!   \mathrm d z^\prime  \bE \Big[ \sqrt{\ell \big( [a, a\! +\!  \zeta] \! \cap \! I_\zeta (z) \big) } \sqrt{\ell \big( [a, a\! +\!  \zeta] \! \cap \! I_\zeta (z^\prime) \big)} \, \Big] \\
& \leq &  \int_0^\infty \!\! \!\!\!\!  \mathrm d z \!  \int_0^\infty \!\! \!\!\!\!   \mathrm d z^\prime   \sqrt{\bE \big[ \ell \big( [a, a\! +\!  \zeta] \! \cap \! I_\zeta (z) \big) \big] } \, \sqrt{\bE \big[ \ell \big( [a, a\! +\!  \zeta] \! \cap \! I_\zeta (z^\prime) \big) \big] } \\
& = & \left(  \,  \int_0^\infty \!\! \!\!\!\!  \mathrm d z \, \sqrt{\bE \big[ \ell \big( [a, a\! +\!  \zeta] \! \cap \! I_\zeta (z) \big) \big] }\,   \right)^2 \leq  \left(  \,  \int_0^\infty \!\! \!\!\!\!  \mathrm d z \, \sqrt{\ell \big( [a, a\! +\!  1] \! \cap \! I_1 (z) \big) }\,   \right)^2 
\end{eqnarray*}
since $I_\zeta(z) \! \subset \! I_1(z)$. 
Then (\ref{Aupbound}) is a consequence of the following:  
\begin{equation}
\label{supessen}
\big( \varphi_{a+1} (1) , \varphi^+_{a-1} (1) \big) \subset \big\{ z \ino \bbR_+ : \;  \ell \big( [a, a\! +\!  1] \! \cap \! I_1 (z) \big) \geko 0 \big\} \subset \big[ \varphi_{a+1} (1) , \varphi^+_{a-1} (1) \big]\; .
\end{equation}
\emph{Proof of (\ref{supessen})}. If $z\ino \big( \varphi_{a+1} (1) , \varphi^+_{a-1} (1) \big) $, then there exists $\epp \ino (0,\tfrac{_1}{^2})$ such that $\varphi_{a+ 1 -\epp} (1) \leko z \leko \varphi_{a-1 +\epp} (1)$. Thus, 
$a\! -\! 1 +\epp \leqo \beta (z) \leqo a+ 1 \! -\! \epp$. Since $(\beta (z), \beta (z) +1) \! \subset \! I_1 (z)$,  
we get 
$$ \ell \big( I_1 (z)\! \cap \! [a, a\! + \! 1]  \big) \geq \ell \big((\beta (z), \beta (z) \! +\! 1) \! \cap \! [a, a\! + \! 1]  \big) \geq \epp >0 \; ,$$
which proves the first inclusion of (\ref{supessen}). 

We next assume that $z\ino [0, \infty)$ is such that $\ell \big( [a, a\! +\!  1]  \cap  I_1 (z) \big) \geko 0$. Since $I_1 (z) $ is contained in $[\beta (z), \beta(z) +1]$, we get  $a\! -\! 1 \leko \beta  (z) \leko a\! + \! 1$. Since  $ \varphi^+_{\beta (z) } (1) \leqo z  \leqo  \varphi_{\beta (z) } (1)$ (this is a simple consequence of the definition of $\beta$), we get $\varphi_{a+1} (1) \leqo z \leqo \varphi^+_{a-1 } (1)$. 
It completes the proof of (\ref{supessen}) and, as already mentioned, it also completes the proof of (\ref{Aupbound}) \cq

\smallskip

By (\ref{Deltaineq}), (\ref{girnogoul}) and (\ref{Aupbound}), we thus have proved that for all $a\ino [1, \infty)$ 
\begin{equation}
\label{voila3}
\int_{\varphi_{a+1}(1) }^{\varphi_{a-1}^+ (1)} \!\!\!\!\! \!\!\!\!\! \!\! \mathrm d s\,  \sqrt{\ell (I_1 (s)\!  \cap \! [a, a+1] ) } \,  \geq \, 2 \sqrt{\gamun} \, \sqrt{\int_0^{\varphi_{a-1} (1) }\!\!\!\!\!\!\!\! \!\!\!\!  \!\!\!\!   \mathrm d z \, \, 
 \Psimu \big( \gamun z \big)  \ell \big( I_1 (z) \! \cap \!  (a, \infty)  \big)} \; .
\end{equation}
\noi
We are now able to adapt \texttt{Step 3} as follows. First note that $\varphi_{a-1}(1)\geko 0$ (by (\ref{phipm}) for instance)
and that the right member of (\ref{voila3}) is strictly positive by (\ref{supessen}) and (\ref{phipm}), which also imply that 
$\ell (I_1 (s)  \cap  [a, a+1] )\geko 0$ for all $s\ino \big( \varphi_{a+1} (1) , \varphi^+_{a-1} (1) \big) $. Thus, we get 
\begin{eqnarray}
\int_{\varphi_{(a+1)}^-(1) }^{\varphi_{a-1}^+ (1)} \!\!\!\!\! \!\!\!\!\!  \mathrm d s\, \left( \int_0^{s}\!\!\!   \mathrm d z \, \, 
 \Psimu \big( \gamun z \big) \frac{ \ell \big( I_1 (z) \! \cap \!  (a, \infty)  \big)}{\ell \big( I_1 (s) \! \cap \!  [a,a+1]  \big)} \right)^{\!\! \!  -\frac{1}{2}}   \!\!\!\!   \!\!\!\!   & \geq &  \!\!\!\!    \int_{\varphi_{(a+1)}^-(1) }^{\varphi_{a-1}^+ (1)} \!\!\!\!\! \!\!\!\!\!  \mathrm d s\, \left(\int_0^{\varphi_{a-1} (1) }\!\!\!\!\!\!\!\! \!\!\!\!  \!\!\!\!   \mathrm d z \, \, 
 \Psimu \big( \gamun z \big) \frac{ \ell \big( I_1 (z) \! \cap \!  (a, \infty)  \big)}{\ell \big( I_1 (s) \! \cap \!  [a, a+1]  \big)} \right)^{\!\!\!   -\frac{1}{2}} \nonumber \\
 \label{voila1}&\geq & 2 \sqrt{\gamun} .
 \end{eqnarray}

We next prove that for all $a\ino [1, \infty)$, for all $s\ino \big( \varphi_{a+1} (1) , \varphi^+_{a-1} (1) \big) $ and for all $z\ino [0, s]$, 
\begin{equation}
\label{voolume}
 \frac{ \ell \big( I_1 (z) \! \cap \!  (a, \infty)  \big)}{\ell \big( I_1 (s) \! \cap \!  [a,a+1]  \big)}\geq 1 \; .
\end{equation}
\emph{Proof of (\ref{voolume}).} Suppose that $z\leqo \varphi_a (1)$. By the definition (\ref{defbeta}) of $\beta$, $a\leqo \beta (z)$. Then, by (\ref{Isimplif}) we get $(\beta (z) , \beta (z) +1) \! \subset \! I_1 (z) \cap [a, \infty) \! \subset \! [\beta (z) , \beta (z) +1]$. Thus, $\ell  (I_1 (z) \cap [a, \infty) )\eqo 1$, which obviously implies (\ref{voolume}). 

  We next assume that $\varphi_a (1) \leko z \leqo s \leko \varphi^+_{a-1}(1)$. We first get 
$ a\! -\! 1 \leko  \beta(z) \leqo a$ by the definition (\ref{defbeta}) of $\beta$. By (\ref{Isimplif}), 
we get $(a, \beta (z) +1) \! \subset \! I_1 (z) \cap [a, \infty) \! \subset \! [a, \beta (z) +1]$. Similarly, we get  
$(a, \beta (s) +1) \! \subset \! I_1 (s) \cap [a, a+1] \! \subset \! [a, \beta (s) +1]$. Therefore 
$$  \frac{ \ell \big( I_1 (z) \! \cap \!  (a, \infty)  \big)}{\ell \big( I_1 (s) \! \cap \!  [a,a+1]  \big)}= \frac{\beta (z) +1\! -\! a}{\beta (s) +1 \! -\! a } \geq 1 $$
since $\beta$ is nonincreasing. It completes the proof of (\ref{voolume}). \cq

\smallskip

By (\ref{voila1}), we get for all $a\ino [1, \infty)$
$$  \int_{\varphi_{a+1}(1)}^{\varphi^+_{a-1} (1)} \!\!\!\!  \frac{\mathrm d s}{\sqrt{\! \int_0^s  \Psimu (\gamun r) \, \mathrm d r}} \, \geq 2\sqrt{\gamun\! } $$
which implies the first inequality in (\ref{lowerbound}) since  $\varphi^+_{a-1} (1) \leqo \varphi_{a-1} (1)$. The second inequality (\ref{lowerbound}) is an easy consequence of the first one: we leave the details to the reader.

\medskip

We next prove $(ii)$ and to that end we assume that 
$\bgam\eqo \tfrac{1}{2} (\delta_{-1} + \delta_{1}) $. We recall from Lemma \ref{mainequa} $(iii)$ that $(\mathcal Z_n)_{n\in \bbN}$ is a critical Galton-Watson process and that 
$\varphi_n (\varphi_m (z))\eqo \varphi_{n+m} (z)$, $m, n\ino \bbN$ and $z\ino [0, 1]$. 
We fix $n\ino \bbN$ and 
to simplify notation, we first set 
$\diffe_n \eqo \varphi_n (1) \! -\! \varphi_{n+1} (1)$. For all $a\ino (n, n+1]$, 
observe that $\mathcal Z_a  \eqo \mathcal Z_{n+1}$, $\varphi_a (z) \eqo \varphi_{n+1}(z) $ and thus 
$$ \Delta_a (1)=  4 R_\mu (\varphi_{n+1}(1) )\diffe_{n} 
+ 4 R_\mu (\varphi_{n+1} (1)) \diffe_{n+1} \; .
 $$
Furthermore, since $\bP$-a.s.~$\zeta\eqo \zeta^\prime\eqo \zeta^\circ  \eqo 1$ and  $\zeta^\bullet \eqo 0$, (\ref{voila0}) with $a\eqo n$ implies 
$$ 2 \sum_{k> n} \big( R_\mu (\varphi_{k}(1) )\diffe_{k-1} 
\! + \! R_\mu (\varphi_{k} (1)) \diffe_{k}\big)  \, = \frac{1}{2}\int_{n}^\infty \!\!\!\!\! \Delta_b (1) \, \mathrm db \eqo \diffe_n^2 \; .$$
 
 For all $n\ino \bbN$ and all $z\ino (\varphi_{n+1} (1), \varphi_n (1)]$, we next set $\sigma_-(z)\eqo \varphi_{n+1} (1)$ and $\sigma_+ (z)\eqo \varphi_n (1)$. Therefore, 
\begin{equation}
\label{voila2}
\diffe_n^2 = 2 \int_0^{\varphi_{n} (1)}  \!\!\!\!\!    \!\!\!\!\!   R_\mu (\sigma_- (z)) \, \mathrm d z  \, +\, 2 \int_0^{\varphi_{n+1} (1)}  \!\!\!\!\!  \!\!\!\!\!   \!\!\!\!\!   R_\mu (\sigma_+(z)) \, \mathrm d z   \; .
\end{equation}

We next modify the two integrals in the right hand side of (\ref{voila2}). To that end, 
we recall that $z\ino [0, 1] \! \mapsto \! \varphi_1 (z)$ is increasing and that $\varphi_1 ([0, 1]) \eqo [0, \varphi_1 (1))]$ where $\varphi_1 (1) \eqo \bP (\mathcal Z_1 \geko 0)$. We denote by 
$\varphi^{-1}_1\!  : \!  [0, \varphi_1 (1))] \! \mapsto \! [0, 1]$ the inverse of $ \varphi_1$. Then note that for all $z\ino [0, 1]$, $\sigma_- (z) \leqo \varphi_1 (1)$. Therefore $\varphi^{-1}_1 (\sigma_-(z))$ is well defined. We then prove that 
\begin{equation}
\label{fnurpss}  \sigma_+ (z) \eqo \varphi^{-1}_1 \big( \sigma_-(z) \big) \quad \textrm{and} \quad 
\sigma_- (z)  \eqo  \varphi^{-1}_1 \big( \sigma_-( \varphi_1 (z)) \big), \quad z\ino [0, 1] \; .
\end{equation}
\noi
\emph{Proof of (\ref{fnurpss}).} Assume that $z\ino (\varphi_{n+1} (1), \varphi_n (1)]$. Then $\sigma_+ (z)\eqo \varphi_{n} (1)$ and $\sigma_- (z)\eqo \varphi_{n+1} (1)\eqo \varphi_1 (\varphi_n (1))$, by 
Lemma \ref{mainequa} $(iii)$, which implies the first equality in (\ref{fnurpss}). We also get 
$$ \varphi_{n+2} (1) \eqo \varphi_{1} (\varphi_{n+1} (1)  ) \leko \varphi_1 (z)\leqo \varphi_{1} (\varphi_{n} (1)  ) \eqo \varphi_{n+1} (1). $$
Therefore, $\sigma_- (\varphi_1 (z))\eqo \varphi_{n+2} (1)\eqo \varphi_1 (\varphi_{n+1} (1))\eqo \varphi_1 (\sigma_-(z))$ and we get the second equality in (\ref{fnurpss}). \cq 

\smallskip

Thanks to (\ref{fnurpss}) and the change of variable $y\eqo \varphi_1 (z)$, we get 
$$\int_0^{\varphi_{n} (1)}  \!\!\!\!\!    \!\!\!\!\!   R_\mu (\sigma_- (z)) \, \mathrm d z=
\int_0^{\varphi_{n+1} (1)}  \!\!    \frac{R_\mu ( \varphi_1^{-1} (\sigma_- (y)))}{\varphi^\prime_1 (\varphi^{-1}_1 (y) )} \, \mathrm d y \leq \frac{1}{\varphi_1^\prime (1)} 
\int_0^{\varphi_{n+1} (1)}  \!\!\!\!\!    \!\!\!\!\! R_\mu \big(  \varphi_1^{-1}( \sigma_- (y)) \big) \, \mathrm d y  $$
since $z\! \mapsto\!  \varphi^\prime_1 (z)$ is decreasing (and note that $\varphi^\prime_1 (1)\eqo \bP (\mathcal Z_1 \eqo 1) \geko 0$). By (\ref{fnurpss}) again and (\ref{voila2}) we get 
$$ \diffe_n^2 \leq 2\left( 1+ \frac{1}{\varphi^\prime_1 (1)}\right) \! 
\int_0^{\varphi_{n+1} (1)}  \!\!\!\!\!\!\!    \!\!\!\!\! R_\mu \big(  \varphi_1^{-1}( \sigma_- (y)) \big) \, \mathrm d y  
\leq \frac{4}{\varphi^\prime_1 (1)} \! 
\int_0^{\varphi_{n+1} (1)}  \!\!\!\!\!\!\!    \!\!\!\!\! R_\mu \big(  \varphi_1^{-1}( y) \big) \, \mathrm d y  $$
since $\sigma_- (y)\leqo y \leqo \varphi_1 (1)$ for all $y\ino (0, \varphi_{n+1} (1)]$.

Since $(\mathcal Z_n)_{n\in \bbN}$ is a critical Galton-Watson process (by Lemma \ref{mainequa} $(iii)$), $\bE [r^{\mathcal Z_2}] \geqo r$ and thus, $\varphi_2 (r) \leqo r$, $r\ino [0, 1]$. 
We denote $\varphi_1^{-1} \! \circ \! \varphi_1^{-1}\! : \! [0, \varphi_2 (1) ] \! \rightarrow \! [0, 1] $ by 
$\varphi_1^{-2}$ that is the inverse of $\varphi_2$. Therefore, for all $y\ino [0, \varphi_2 (1)]$, we get $y \leqo \varphi_1^{-2} (y)$ and, by (\ref{Rpsifourbi}), 
$$ R_\mu \big(  \varphi_1^{-1}( y) \big) = R_\mu \big(  \varphi_1\big(  \varphi_1^{-2}( y) \big) \big) \eqo \Psimu \big( \tfrac{_1}{^2} (\varphi_1^{-2} (y) + y )\big) \leq \Psimu \big(\varphi_1^{-2} (y) \big) \; .$$
Thus, for all integers $n\geqo 2$, 
$$ \int_{\varphi_{n+1} (1)}^{\varphi_n (1)} \!\! \frac{\mathrm dz }{\sqrt{\int_0^{z} \Psimu \big( \varphi_1^{-2} (y) \big) \mathrm dy }} \; \leq \; \frac{\partial_n}{\sqrt{\int_0^{\varphi_{n+1} (1)} \!\! \Psimu \big( \varphi_1^{-2} (y) \big) \mathrm dy }} \; \leq \; \frac{2}{\sqrt{\varphi^\prime_1 (1)}} . $$
By an elementary change of variable we next get 
$$\int_{\varphi_{n+1} (1)}^{\varphi_n (1)} \!\! \frac{\mathrm dz }{\sqrt{\int_0^{z} \Psimu \big( \varphi_1^{-2} (y) \big) \mathrm dy }} =  \int_{\varphi_{n-1} (1)}^{\varphi_{n-2} (1)} \!\! \frac{\varphi_2^\prime (s) \mathrm ds }{\sqrt{\int_0^{s} \Psimu \big( r \big)\varphi_2^\prime (r)   \mathrm dr }} \geq (\varphi^\prime_1 (1))^2 
\int_{\varphi_{n-1} (1)}^{\varphi_{n-2} (1)} \!\! \frac{ ds }{\sqrt{\int_0^{s} \Psimu \big( r \big) \mathrm dr }}$$
since $\varphi^\prime_2 (z)\eqo \varphi^\prime_1 (z) \varphi^\prime_1 (\varphi_1 (z)) $ and 
$0 \leko \bP (\mathcal Z_1 \eqo 1) \eqo \varphi^\prime_1 (1) \leqo \varphi^\prime_1 (z) \leqo \varphi^\prime_1 (0)\eqo 1$. This implies (\ref{upperbound}), which completes the proof of the proposition.  \cqfd

\subsection{Estimates on the maximal displacement of BRWs with i.i.d jumps}
\label{genestisec}
In this section we provide general estimates on the maximal displacement of $\bbR$-valued
BRWs, which allow to remove the symmetry assumption in Lemma \ref{keylower} and also to remove the assumption $\gamma_1\geko 0$. We consider a $\bbR$-valued BRW $\bS\eqo (S_u)_{u\in t}$ where: 
\begin{compactenum}

\smallskip

\item[\textbf{(a)}] 
$t\ino \bbT$ is deterministic and finite and we recall that $\Gamma (\fftree)\eqo \max_{u\in \fftree} |u|$ is its \emph{total height};

\smallskip

\item[\textbf{(b)}] the jumps $(\xi_u)_{u\in t \backslash \{ \varnothing\}}$ are i.i.d.~$\bbR$-valued r.v.~whose law is denoted by $\bgam $; we assume that $\beta\! : = \! \int_{\bbR} y^2 \bgam (\mathrm dy)< \infty$ and that $ \int_{\bbR} y \bgam (\mathrm dy)\eqo 0$. 

\smallskip

\end{compactenum}

\noi
We also introduce the following. 
Let $\xi$ and $\xi^*$ be two independent r.v.~with law $\bgam$. Then, :  
\begin{equation}
\label{bgamdef} 
\textrm{$\bgam^{\mathtt{s}}$ is the law of $\tfrac{_1}{^2} (\xi\! -\! \xi^*)$,} \quad \bgam^\bullet \!\! := \tfrac{_1}{^4} (\delta_{-1} + \delta_1) + \tfrac{_1}{^2} \bgam^{\mathtt{s}} \quad \textrm{and} \quad \bgam^\circ\!\!  := \tfrac{_1}{^4} (\delta_{-1} + \delta_1) + \tfrac{_1}{^2} \delta_0\; .
\end{equation}

\begin{lemma}
\label{simplbound} We keep the notations and the assumption as above. We suppose that $S_\varnothing\eqo 0$. 

\begin{compactenum}

\smallskip

\item[$(i)$] Let $\bgam^{\mathtt{s}}$ be as in (\ref{bgamdef}) and let $(S^{\mathtt{s}}_u)_{u\in \fftree}$ be a $\bbR$-valued BRW whose jumps are i.i.d.~with law $\bgam^{\mathtt{s}}$. Suppose that $S^{\mathtt{s}}_\varnothing \eqo 0$. Then, 
\end{compactenum}
\begin{eqnarray}
\label{symcontrol}
\forall y_1, y_2 \ino \bbR_+^*, \qquad \tfrac{_1}{^2}\bP \big( \max_{u\in \fftree} |S^{\mathtt{s}}_u|  \geko y_1 \big) \!\!\!\!  &  \leqo &\!\!\!\! \bP \big( \max_{u\in \fftree} |S_u|  \geko y_1 \big) \nonumber \\
\textrm{and} \qquad \bP \big( \max_{u\in \fftree} |S_u | \geko 2y_1+ y_2 \big) \!\!\!\! & \leq& \!\!\!\! \bP \big( \max_{u\in \fftree} |S^{\mathtt{s}}_u|  \geko y_1 \big) +  \min \Big(1\, ,  \tfrac{\beta \Gamma (\fftree)}{y_2^2} \Big) .
\end{eqnarray}
\begin{compactenum}
\item[$(ii)$] Assume that $\bgam ([-1, 1])\eqo 1$ and let $\bgam^\bullet $ and $\bgam^\circ$ be as in (\ref{bgamdef}). Let $(S^\bullet_u)_{u\in \fftree}$ and $(S^\circ_u)_{u\in \fftree}$ 
be $\bbR$-valued BRWs  whose jumps are i.i.d.~with respective laws $\bgam^\bullet$ and $\bgam^\circ$. Suppose that 
$S^\bullet_\varnothing \eqo S^\circ_\varnothing \eqo 0$. 
Then, for all $y_1, y_2,y_3 \ino \bbR_+^*$
\end{compactenum}
\begin{equation}
\label{finpret} 
\bP \big( \max_{u\in \fftree} |S_u|  \geko 4(y_1+y_2)\! +\! y_3 \big) \leqo 2\bP \big( \max_{u\in \fftree} |S^\circ_u | \geko y_1 \big)\! +\!  2\bP \big( \max_{u\in \fftree} |S^\bullet_u|  \geko y_2 \big)\! +\!  \min \Big(1\, ,  \tfrac{\beta \Gamma (\fftree)}{y_3^2} \Big) .
\end{equation}
\end{lemma}
\noi
\textbf{Proof.} We first prove $(i)$. Let $(S^*_u)_{u\in \fftree}$ be an independent copy of $S$. For all $u\ino \fftree$ we set $S^{\mathtt{s}}_u\eqo \tfrac{_1}{^2}(S_{u}\! -\! S^*_u)$. Clearly, $S^{\mathtt{s}}$ is a BRW such that $S^{\mathtt{s}}_\varnothing \eqo 0$ and whose jumps are independent with law $\bgam^{\mathtt{s}}$ as defined in (\ref{bgamdef}). 
We easily observe that for all $y_1\ino \bbR_+^*$, 
$$ \bP \big( \max_{u\in \fftree} |S^{\mathtt{s}}_u|  \geko y_1 \big) \leqo  \bP \big( \max_{u\in \fftree} |S_u|  \! +\!  \max_{u\in \fftree} |S^*_u|  >  2y_1 \big) \leq 2\bP \big( \max_{u\in \fftree} |S_u|  \geko y_1 \big), $$
which prove the first inequality in (\ref{symcontrol}). 

Let us prove the second one. To that end, we introduce $u_0$, the $<_{\mathtt{lex}}$-minimal vertex $u\ino \fftree$ such that $|S_u|\eqo \max_{v\in \fftree} |S_v|$. Then for all $y_1, y_2 \ino \bbR_+^*$ 
$$\bP \big( \max_{u\in \fftree} |S_u| \geko  2y_1 \! +\! y_2 \big) \leq \bP \big(2|S^{\mathtt{s}}_{u_0}| \! +\!  |S^*_{u_0}| \geko 2y_1+y_2 \big)\leq \bP \big(|S^{\mathtt{s}}_{u_0}|  \geko y_1 \big)+\bP \big( |S^*_{u_0}| \geko y_2 \big). $$
First note that $ \bP \big(|S^{\mathtt{s}}_{u_0}|  \geko y_1 \big)\leqo  \bP \big( \max_{u\in \fftree} |S^{\mathtt{s}}_{u}|  \geko y_1 \! \big)$. Next, observe that $u_0$ only depends on $S$. It is therefore independent of $S^*$. Thus, 
$$ \bE \big[  \un_{\{|S^*_{u_0}| > y_2 \} }\, \big| \, S\big] \leq \min \Big(1 \, , \bE \Big[ \tfrac{(S^*_{u_0})^2}{y_2^2} \, \Big| \, S \Big] \Big) =  
\min \Big(1 \, , \bE \Big[ \tfrac{\beta |u_0|}{y_2^2} \, \Big| \, S \Big] \Big) \leq   \min \Big(1\, ,  \tfrac{\beta \Gamma (\fftree)}{y_2^2} \Big) , $$
which immediately entails the second inequality in (\ref{symcontrol}). 

\smallskip

We next prove $(ii)$. To that end, we assume that $\bgam ([-1, 1])\eqo 1$ and we introduce the independent r.v.~$(\eta_u, \epp_u, \xi^{\mathtt{s}}_u)_{u\in \fftree \backslash \{ \varnothing \}}$ that are distributed as follows: $\eta_u$ has law $\tfrac{_1}{^2} (\delta_0+ \delta_1)$, $\epp_u$ has law $\tfrac{_1}{^2} (\delta_{-1}+ \delta_1)$ and $\xi^{\mathtt{s}}_u$ has law $\bgam^{\mathtt{s}}$ (as defined in (\ref{bgamdef})). We denote by $(S^{\mathtt{s}}_u)_{u\in \fftree}$, $(S^\prime_u)_{u \in \fftree} $ and $(S^{\prime\prime}_u)_{u \in \fftree} $ the branching random walks such that $S^{\mathtt{s}}_\varnothing \eqo S^{\prime}_\varnothing \eqo S^{\prime\prime}_\varnothing \eqo 0$, with respective jumps 
$(\xi^{\mathtt{s}}_u)_{u \in \fftree \backslash \{ \varnothing \}}$ $(\eta_u \xi^{\mathtt{s}}_u)_{u \in \fftree \backslash \{ \varnothing \}}$ and $((1\! -\! \eta_u) \xi^{\mathtt{s}}_u)_{u \in \fftree \backslash \{ \varnothing \}}$. First note that for all $u\ino \fftree$, $S^{\mathtt{s}}_u\eqo S^{\prime}_u+ S^{\prime \prime}_u$. Then observe that 
$S^{\prime}$ and $S^{\prime \prime}$ have the same law. Therefore, for all $z\ino [0, \infty)$, 
\begin{equation}
\label{troutriv}
 \bP \big( \max_{u\in \fftree} |S^{\mathtt{s}}_u|  \geko 2z \big) \leq 
\bP \big( \max_{u\in \fftree} |S^{\prime}_u| \! + \!  \max_{u\in \fftree} |S^{\prime \prime}_u|  \geko 2z \big) \leq 2\bP \big( \max_{u\in \fftree} |S^{\prime}_u|  \geko z \big)\; .
\end{equation}

We next introduce the branching random walks $(S^\bullet_u)_{u\in \fftree}$ and $(S^\circ_u)_{u\in \fftree}$ 
such that $S^\bullet_\varnothing\eqo S^\circ_\varnothing \eqo 0$ 
and whose respective jumps are $((1\! -\! \eta_u) \epp_u + \eta_u\xi^{\mathtt{s}}_u)_{u\in \fftree \backslash \{ \varnothing \}}$ and 
$((1\! -\! \eta_u) \epp_u)_{u\in \fftree \backslash \{ \varnothing \}}$. First note that the jumps of $S^\bullet$ are independent with law $\bgam^\bullet$ and those of $S^\circ$ are also independent with law $\bgam^\circ$, 
where $\bgam^\bullet$ and $\bgam^\circ$ are defined in (\ref{bgamdef})). Next observe that for all $u\ino \fftree$, $S^\prime_u  \eqo  S^\bullet_u \! -\! S^\circ_u $. Thus, for all $y_1, y_2 \ino [0, \infty)$, we get 
\begin{eqnarray*}
 \bP \big( \max_{u\in \fftree} |S^{\prime}_u|  \geko y_1+ y_2 \big) &  \leq &  
\bP \big( \max_{u\in \fftree} |S^{\bullet}_u| \! + \!  \max_{u\in \fftree} |S^{\circ}_u|  \geko y_1+y_2 \big) \\
& \leq &\bP \big( \max_{u\in \fftree} |S^{\bullet}_u|  \geko y_1  \big)+
\bP \big( \max_{u\in \fftree} |S^{\circ}_u|  \geko y_2 \big), 
\end{eqnarray*}
which implies (\ref{finpret}) by (\ref{troutriv}) and (\ref{symcontrol}). \cqfd

\medskip

The following lemma provides a comparison of the tail distribution of the 
maximal displacement of any BRW whose i.i.d.~jumps have a third moment, 
with that of a BRW whose i.i.d.~jumps have law $\tfrac{_1}{^2} (\delta_{-1}+ \delta_1)$. This result is not used to proved Theorem \ref{unifbounded}. 
\begin{lemma}
\label{egoutploufa} Let $\bgam$ be a centered probability law on $\bbR$ such that $m_3 \! :=\!  \int_{\bbR} |y|^3 \bgam (\mathrm dy) \leko \infty$.
We also set $\beta \eqo  \int_{\bbR} y^2 \bgam (\mathrm dy)$ and $\varrho \eqo \beta^{3/2}/m_3 \ino (0, 1]$. 
Let $\fftree\ino \bbT$ be finite and let $(S_u)_{u\in \fftree}$ and $(S^o_u)_{u\in \fftree}$ be BRWs such that $S_\varnothing\eqo S^o_\varnothing \eqo 0$ and 
whose independent jumps have respective laws $\bgam$ and $\tfrac{_1}{^2} (\delta_{-1}+ \delta_1)$. 
 Then, for all $z \ino [\tfrac{_1}{^8}, \infty) $,
\begin{equation}
\label{bouclerla}
\bP \big(\!  \max_{u\in \fftree} |S_u| \geko \varrho \sqrt{\beta}  z \big) \, \geq \, \frac{1}{8} \Big( 1\! -\!  \frac{4}{ z\varrho^2 }\Big)\,  \bP \big(\!  \max_{u\in \fftree} |S^o_u| \geko 16 z \big) 
\end{equation}
\end{lemma}
\noi
\textbf{Proof.} We first deal with the case of a BRW $(S^{\mathtt{s}}_u)_{u\in \fftree}$ 
such that $S^{\mathtt{s}}_\varnothing \eqo 0$ and whose jumps $(\xi^{\mathtt{s}}_u)_{u\in t\backslash \{ \varnothing \}}$ are i.i.d~with law $\bgam^{\mathtt{s}}$ as in (\ref{bgamdef}). Since $\bgam^{\mathtt{s}}$ is symmetric, w.l.o.g.~we assume that there are independent r.v.~$\epp_u$, $\zeta_u$, 
$u\ino \fftree\backslash \{ \varnothing \}$ such that $\bP (\epp_u \eqo \pm 1)\eqo  \tfrac{1}{2}$, $\zeta_u \eqo |\xi^{\mathtt{s}}_u|$ and $\xi^{\mathtt{s}}_u\eqo \epp_u \zeta_u$. We denote by $(S^o_u)_{u\in \fftree}$ the BRW such that $S^o_\varnothing \eqo 0 $ and whose jumps are 
$(\epp_u)_{u\in t\backslash \{ \varnothing \}}$. Let $u^* \ino \fftree$ be the $<_{\mathtt{lex}}$-minimal vertex $u\ino \fftree$ such that $S^o_u \eqo \max_{v\in \fftree} S^o_v$. We set $n\eqo |u^*|$ and for all $j\ino \{ 0, \ldots, n\}$, we also set $u(j)\eqo u^*_{| j}$, which is the ancestor of $u^*$ at height $j$ in $\fftree$.
We then fix $k\ino \bbN^*$ and we work on the event $A_k\eqo \{ S^o_{\! u^*}\eqo k\}$. We next introduce 
the following integers 
$$ \big\{ j^{_+}_{^1}  \leko \ldots \leko j^{_+}_{^{\!\frac{n+k}{2}}} \big\} \eqo \big\{ j\ino \{ 1, \ldots, n\} \! : \epp_{u(j)}\eqo 1 \big\} $$ and 
$\big\{ j^{_-}_{^1} \leko \ldots \leko j^{_-}_{^{\! \frac{n-k}{2}}} \big\}\eqo \big\{1, \ldots, n \big\} \backslash \big\{ j^{_+}_{^1} \leko \ldots \leko j^{_+}_{^{\! \frac{n+k}{2}}} \big\}$. Note that 
$$ S^{\mathtt{s}}_{u^*}\eqo \zeta_{u(j^+_1)}+ \ldots  + \zeta_{u(j^+_k)} + Y \quad \textrm{where} \quad Y\eqo \sum_{1\leq m\leq \frac{n-k}{2}} \zeta_{u(j^+_{k+m})} \! -\! \zeta_{u(j^-_m )}. $$
Denote by $\bgam_{+}$ the law of the $\zeta_u$ and let $(\zeta^\prime_k)_{k\in \bbN^*}$ be an i.i.d.~sequence of r.v.~with law $\bgam_+$. 
Observe that conditionally given 
$S^o$ on $A_k$, the $(\zeta_{u(j)})_{1\leq j\leq n}$ are i.i.d.~with law $\bgam_+$ and $S^{\mathtt{s}}_{u^*}$ has thus the same law as $\zeta^\prime_1+ \ldots + \zeta^\prime_k + Y^\prime$ where $Y^\prime \eqo \sum_{1\leq m\leq \frac{n-k}{2}}  \zeta^\prime_{k+2j} \! -\! \zeta^\prime_{k+2j-1}$. 
Observe that $Y^\prime$ is independent from $\zeta^\prime_1+ \ldots + \zeta^\prime_k$ and that $Y^\prime$ has the same law as $-Y^\prime$, which implies that $\bP (Y^\prime \geqo 0) \geqo 1/2$. 
For all $z\ino [0, \infty)$, we therefore get the following. 
\begin{eqnarray}
\bE \big[ \un_{A_k \cap \{\max_{u\in \fftree}  S^{\mathtt{s}}_u \geq z\} } \big|  S^o\big] \!\!\!  &\geq & \!\! \!  \bE \big[ \un_{A_k \cap \{ S^{\mathtt{s}}_{u^*} \geq z\}} \big|  S^o\big]
\eqo   \bP \big(\zeta^\prime_1+ \ldots + \zeta^\prime_k + Y^\prime \geqo z \big) \un_{A_k}  \nonumber \\
\!\!\!  & \geq &\!\!\!   \bP \big(\zeta^\prime_1+ \ldots + \zeta^\prime_k  + Y^\prime \geqo z\, ; \, Y^\prime \geqo 0  \big) \un_{A_k} \nonumber  \\
\!\!\!  & \geq &\!\!\!   \bP \big(\zeta^\prime_1+ \ldots + \zeta^\prime_k  \geqo z\, ; \, Y^\prime \geqo 0  \big) \un_{A_k} \nonumber  \\
\label{espcondii} \!\!\!  & \geq &\!\!\! \tfrac{1}{2}\bP \big(\zeta^\prime_1+ \ldots + \zeta^\prime_k  \geqo z \big) \un_{A_k} .
\end{eqnarray}
For all $k\ino \bbN^*$ and $z\ino [0, \infty)$, we set $\Sigma (z,k)\eqo \bP \big(\zeta^\prime_1+ \ldots + \zeta^\prime_k  \geqo z \big)$. Then (\ref{espcondii}) easily entails for all $q\ino \bbN^*$ and $z\ino [0, \infty)$ that 
\begin{eqnarray*}
\bP \big(  \! \max_{u\in \fftree} S^{\mathtt{s}}_u\geq z\big) \!\!\! & \geq &\!\!\!  \bP \big( \!  \max_{u\in \fftree} S^{\mathtt{s}}_u\geqo z\; ; \; \max_{u\in \fftree} S^o_u \geqo q  \big) = \sum_{k\geq q} \bE \big[\un_{ A_k \cap \{ \max_{u\in \fftree} S^{\mathtt{s}}_u\geq z\}} \big]  \nonumber  \\
\!\!\! & \geq &\!\!\!  \sum_{k\geq q}  \tfrac{1}{2}\bE \big[ \Sigma (z, k)  \un_{A_k} \big]   \geq  
 \tfrac{1}{2}\Sigma (z,q) \bP \big(  \! \max_{u\in \fftree} S^o_u \geqo q \big) ,
\end{eqnarray*}
since $k\! \mapsto \! \Sigma (p,k)$ is nondecreasing. Next observe that 
$$ \bP \big(  \! \max_{u\in \fftree} |S^o_u| \geqo q \big) \leq \bP \big(  \! \max_{u\in \fftree} S^o_u \geqo q \big) +\bP \big(  \! -\! \min_{u\in \fftree} S^o_u \geqo q \big) \eqo 2\bP \big(  \! \max_{u\in \fftree} S^o_u \geqo q \big)$$
since the jumps are symmetric. Therefore, $\bP \big(  \! \max_{u\in \fftree} |S^{\mathtt{s}}_u|\geqo z\big) \geqo \tfrac{1}{4} \Sigma (z,q) \bP \big(  \! \max_{u\in \fftree} |S^o_u| \geqo q \big) $ and by the first inequality in (\ref{symcontrol}) we get 
\begin{equation}
\label{minokern}
\forall q\ino \bbN^*, \; \forall z\ino [0, \infty), \quad \bP \big(  \! \max_{u\in \fftree} |S_u|\geq z\big) \geq \tfrac{1}{8} \Sigma (z,q) \bP \big(  \! \max_{u\in \fftree} |S^o_u| \geqo q \big) \, .
\end{equation}  
We now provide a lower bound for $\Sigma (z,q)$ when $z$ is proportional to $q$. 
To that end, we set $m^{\prime}_k\eqo \int_{[0, \infty)} y^k\bgam_{\! +} (\mathrm d y)$ for all $k\ino \{ 1, 2, 3\}$. First observe that $(m^\prime_2)^2 \leqo m_3^\prime m_1^\prime$, by Cauchy-Schwarz. To simplify notation we set $c\eqo \tfrac{_1}{^2} (m^\prime_2)^{2}/m^\prime_3$. 
Then by Chebychev we get for all $q\ino \bbN^*$ 
$$ 1\! -\! \Sigma \big(cq, q \big) \leq 1\! -\! \Sigma \big(\tfrac{_1}{^2} m^\prime_1 q, q \big)\leq \bP \Big(\!\! \sum_{\; \,  1\leq k\leq q} \!\!\!\!  m^\prime_1\! -\! \zeta^\prime_k \, > \tfrac{_1}{^2} m^\prime_1 q\Big)\leq \frac{4\mathtt{var} (\zeta^\prime_1)}{q (m^\prime_1)^2} \leq \frac{4m_2^\prime}{q (m^\prime_1)^2} \leq 
\frac{4(m_3^\prime)^2}{q(m_2^\prime)^3}  . $$ 
We now connect these quantities to $\beta$ and $m_3$. 
To that end recall that $\bgam_{\! +}$ is the law of $\tfrac{1}{2}|\xi\! -\! \xi^*|$ where $\xi$ and $\xi^*$ are independent with law $\bgam$. We then observe that $m^\prime_2\eqo \tfrac{_1}{^2}\beta$ and $m_3^\prime \leqo m_3$, by convexity of $y\mapsto |y|^3$. Thus, if we set $\varrho \eqo \beta^{3/2}m_3^{-1}$, then for all $q\ino \bbN^*$, 
$$\Sigma \big(\tfrac{_1}{^8} \varrho \sqrt{\beta }q , q \big) =  \Sigma \big(\tfrac{_1}{^8} \beta^2 m_3^{-1} q ,  q \big) \geq 1- \frac{32\,  m_3^2}{q \, \beta^3} = 1- \frac{32}{q\varrho^2 },  $$
which easily implies (\ref{bouclerla}) by (\ref{minokern}). \cqfd

\smallskip

\subsection{Proof of Theorem \ref{unifbounded}}
\label{Thm1unpfsec}
We first prove  Theorem \ref{unifbounded} in $\textbf{Case (0)}$. We fix $\psi \ino \mathscr L$ that satisfies 
\texttt{Sheu}($\psi$). 
We assume that $(a_n)_{n\in \bbN}$, $(b_n)_{n\in \bbN}$ satisfy $\texttt{Norm} (\mathbf a, \mathbf b)$. For all $n\ino \bbN$, we fix $\mu_n $, a non-trivial critical 
offspring distribution (i.e.~it satisfies (\ref{nontricri})). We overall assume \texttt{\L{}uka}($\mathbf a , \mathbf b, \bmu, \psi$) and  \texttt{Hght}($\mathbf a , \mathbf b, \bmu$) . 
For all $n\ino \bbN$, $(\tau_n(p))_{p\in \bbN^*}$ 
is an infinite sequence of GW($\mu_n$)-tree. We denote by $\tau^{_\infty}_n$ the tree associted with 
$(\tau_n(p))_{p\in \bbN^*}$: namely $\theta_{[p]}\tau^{_\infty}_n \eqo \tau_n (p)$. We denote by $C^{_{(n)}}_\cdot$, $H^{_{(n)}}_\cdot$ and $V^{_{(n)}}_\cdot$ resp.~the contour process, the height process and the \L{}ukasiewicz path of the forest $(\tau_n(p))_{p\in \bbN^*}$, which are normalized as in (\ref{renormCW}). 
For all $x\ino \bbR_+^*$, we denote by $\tau^{_x}_n$ the tree associated with the finite forest $(\tau_n(p))_{1\leq p\leq \lfloor a_n x \rfloor}$. Note that $\tau^{_x}_n\subset \tau^{_\infty}_n$. 
 
 We then consider the infinite sequence of BRWs $\bS_n (p)\eqo (S_{n,u} (p))_{u\in \tau_n (p)}$, $p\ino \bbN^*$, such that $S_{n, \varnothing} (p)\eqo 0$ and 
conditionally given $(\tau_n(p))_{p\in \bbN^*}$, the jumps of $\bS_n(p)$ are $\bbR$-valued and i.i.d.~r.v.s whose deterministic law $\bgam_n$ satisfies (\ref{cenuni}). It is convenient to introduce $\bS_n \eqo (S_{n,u})_{u\in \tau^{_\infty}_n}$ such that $S_{n, \varnothing }\eqo 0$ and $S_{n, [p]\ast u}\eqo S_{n,u}(p)$, $p\ino \bbN^*$ and $u\ino \tau_n (p)$. 
We first prove the following. 
\begin{lemma}
\label{deviation} We keep the above notations. We assume \emph{\texttt{\L{}uka}($\mathbf a , \mathbf b, \bmu, \psi$)},  
\emph{\texttt{Hght}($\mathbf a , \mathbf b, \bmu$)} and \emph{\texttt{Sheu}$_{\,}$($\baa, \bbb, \bmu$)}, which imply \emph{\texttt{Sheu} ($\psi$)} and thus \emph{\texttt{Grey} ($\psi$)}. 
We only assume here that the $\bgam_n$ is such that $\bgam_n([-c,c])\eqo 1$, $\int_{\bbR} x\bgam_n(dx)\eqo 0$,
for all $n\ino \bbN$. Then, for all $x, y_1, y_2 \ino \bbR_+^*$, 
\begin{eqnarray}
\limsup_{n\to \infty} \bP \Big( \max_{u\in \tau^{_x}_{^n} } \, \lvert S_{n,u} \rvert \geq c\sqrt{\lambda_n} \big(8y_1+y_2 \big) \Big) 
\leq 8 \Big( 1 \!\!\!\! &- &\!\!\!\! e^{-6xw (y_1/12)} \Big) \nonumber \\
\label{deviaSenn} \!\!\!\! &+ &\!\!\!\! 
\bE \Big[ \min \Big(1, y_2^{-2}\Gamma_{\! x}\Big) \Big].
\end{eqnarray}
where $w$ stands for the inverse of $F$ as defined in (\ref{defFnF}) and where $\Gamma_{\! x}$ is defined in (\ref{szhghtGW}), \end{lemma}
\noi
\textbf{Proof.} To simplify notations, we set $\tau_n \! :=\! \tau_n (1)$, which is a single GW($\mu_n$)-tree and we set 
$\gamma_{n,1}\eqo \bgam_n(\{ c\})$. We first assume that 
\begin{equation}
\label{hypounn}
\textrm{$\bgam_n$ is symmetric} \quad \textrm{and}  \quad \liminf_{n\to \infty}\gamma_{n,1}\eqo \gamma \geko 0.
\end{equation}
For all $y\ino \bbR_+^*$ we also set 
$\varphi_{n,y} \eqo \bP \big(  \max_{u\in \tau_n} S_{n,[1]\ast u} \geqo c y\sqrt{\lambda_n} \,  \big)$. 
Then for all sufficiently large $n$ the BRW $\smash{(\frac{1}{c}S_{n, [1] \ast u})_{u\in \tau_n}}$ satisfies the assumptions of Proposition \ref{keylower}: by (\ref{lowerbound}) in Proposition \ref{keylower} $(i)$ entails 
$F_n \big( \gamma_{n,1}a_n  \varphi_{n,y} \big) \! -\! F_n \big( \gamma_{n,1} a_n \big) \geq \frac{_1}{^2} \gamma_{n, 1} y
$, which implies 
\begin{equation}
\label{Devstep1}
\gamma_{n,1}a_n  \varphi_{n,y} \leq w_n  \big( \, \frac{_{_1}}{^{^2}} \gamma_{n, 1} y+F_n ( \gamma_{n,1} a_n ) \big)\; .
\end{equation}
where $w_n$ is the inverse function of $F_n$ as defined in (\ref{defFnF}). By (\ref{hypounn}) we get $\frac{2}{3}\gamma \leko \gamma_{n,1}$ for all sufficiently large 
 $n\ino \bbN$. Therefore we first get $\lim_{n\to \infty}\gamma_{n,1} a_n\eqo \infty$ and \texttt{Sheu}$_{\,}$($\baa, \bbb, \bmu$) implies $\lim_{n\to \infty}F_n \big( \gamma_{n,1} a_n \big)\eqo 0$. Combined with Lemma \ref{controlpsi}, it easily entails 
\begin{equation}
\label{Devstep2}
\frac{_{_2}}{^{^3}}\gamma \, \limsup_{n\to \infty} a_n  \varphi_{n,y} \leq w \big(\frac{_{_1}}{^{^3}} \gamma  y\big) \; .
\end{equation}
We then get 
$$ \bP \big( \max_{u\in \tau^{_x}_{^n} } \, \lvert S_{n,u} \rvert \geq cy\sqrt{\lambda_n}  \big) \leq 2  \bP \big(  \max_{u\in \tau^{_x}_{^n} } \,  S_{n,u}  \geq cy\sqrt{\lambda_n} \big)= 2 \big( 1\! -\!  \big(1\! -\! \varphi_{n, y} \big)^{\lfloor a_n x\rfloor}\big), $$
(here we use the symmetry of $\bgam_n$ to get the first inequality). By (\ref{Devstep2}) we then get 
\begin{equation}
\label{Devstep3}
 \limsup_{n\to \infty} \bP \big( \max_{u\in \tau^{_x}_{^n} } \, \lvert S_{n,u} \rvert \geq cy\sqrt{\lambda_n}  \big) \leq  2 \big( 1\! -\!  e^{-\frac{3x}{2\gamma} w (\gamma y/3)} \big)
\end{equation}
for all $x,y\ino \bbR_+^*$ and for all $\bgam_n$ that satisfy $\bgam_n([-c,c])\eqo 1$, $\int_{\bbR} x\bgam_n(dx)\eqo 0$, and (\ref{hypounn}). 

   We now only assume that the $\bgam_n$ satisfy  (\ref{cenuni}). 
We denote by $(S^\circ_{u})_{n\in \tau^\infty_n}$ and $(S^{n,\bullet}_{n,u})_{n\in \tau^\infty_n}$ two BRWs such that conditionally given $\tau_{n}^{_\infty}$ their jumps are i.i.d.~with respective laws $\bgam^\circ\eqo \frac{_1}{^4} \delta_{-c} + \frac{_1}{^4} \delta_{c} +
\frac{_1}{^2} \delta_{0}$ and $\bgam_n^\bullet\eqo   \frac{_1}{^4} \delta_{-c} + \frac{_1}{^4} \delta_{c}+ \bgam^{\mathtt{s}}_n$, where $\bgam_n^{\mathtt{s}}$ is the symmetrized version of $\bgam_n$: namely, if $\xi$ and $\xi'$ are two independent r.v. with law $\bgam_n$, then $\frac{_1}{^2} (\xi\! -\! \xi')$ has law $\bgam^{\mathtt{s}}_n$. 

By (\ref{finpret}) in Lemma \ref{simplbound}, for all $n \ino \bbN$ and all $x, y_1, y'_2, y_2 \ino \bbR_+^*$ we get 
\begin{eqnarray}
\bP \big( \max_{u\in \tau_n^x} |S_{n,u}| \!\!\!  & \geko & \!\!\!   c\sqrt{\lambda_n}  \big( 4(y_1  +  y'_2)\! +\! y_2\big) \big) 
\leq \,  2\bP \big( \max_{u\in \tau_n^x} |S^\circ_u | \geko cy_1 \sqrt{\lambda_n} \,  \big) \nonumber  \\ 
\label{finpretbis}  
 \!\!\!\!  \!\!\!\! &  & \!\!\!\! \!\!\!\! \!\!\!\! + \,  
2\bP \big( \max_{u\in \tau_n^x} |S^\bullet_{n,u}|  \geko cy'_2 \sqrt{\lambda_n}\, \big)\! +\!  \bE \Big[ \min \big(1\, , c^{-2}\beta_n y_2^{-2} \lambda^{-1}_n \Gamma ( \tau_n^x)\big) \Big] .
\end{eqnarray}
where $\beta_n\eqo \int_{[-c,c]} z^2 \bgam_n (\mathrm d z)\leqo c^2$. Note that (\ref{Devstep3}) applies to $(S^\circ_{u})_{n\in \tau^\infty_n}$ and $(S^{n,\bullet}_{n,u})_{n\in \tau^\infty_n}$ with $\gamma\eqo 1/4$. 
This implies the desired result (\ref{deviaSenn}) since 
$$ \bE \Big[ \min \big(1\, , c^{-2}\beta_n y_2^{-2} \lambda^{-1}_n \Gamma ( \tau_n^x)\big) \Big] \leq  \bE \Big[ \min \big(1\, ,  y_2^{-2} \lambda^{-1}_n \Gamma ( \tau_n^x)\big) \Big] \xrightarrow[n\to \infty]{\; } \bE \Big[ \min \big(1\, ,  y_2^{-2}  \Gamma_{\! x} \big) \Big]$$
by the convergence in law (\ref{szhghtGW}). This completes the proof of the lemma. \cqfd

\smallskip

We next complete the proof Theorem \ref{unifbounded} $(i)$ in $\textbf{Case (0)}$. The proof follows the same line as  $\textbf{Step (8)}$ of the proof of Theorem \ref{Sheuexplain}. 
We first claim that we only need to prove the following. 
\begin{equation}
\label{Wendtight}
\forall N\ino \bbN, \; \forall \epp\ino  (0, 1), \qquad \lim_{p\to \infty} \limsup_{n\to \infty} \bP \big( \omega_{N,p} \big( \widehat{W}^{_{(n)}}_{\cdot} \big) \geqo \epp\big) = 0,
\end{equation}
where for all $N,p\ino \bbN$ and for all continuous function $f\ino \bC^{_0}_{^1}$ we have set 
$$ \omega_{N\! ,p} (f) \eqo \max_{0\leq j< N2^p} \max \big\{\,  | f(s)\! -\! f(j2^{-p}) |\, ; \, s\ino [j2^{-p} , (j+1)2^{-p}] \big\}  \; .$$ 
\noi
\emph{Indeed,} (\ref{Wendtight}) implies that the laws of $\widehat{W}^{_{(n)}}_{\cdot} $ are tight in $\bC^{_0}_{^1}$ by standard results (see e.g.~Billingsley \cite{Bil68} Thm 7.3 p.~82). 
By Lemma \ref{fdcvsna}, 
$$ \Big( C^{_{(n)}}_\cdot; V^{_{(n)}}_\cdot ;  \widehat{W}^{_{(n)}}_{\! s_0},  \ldots , \widehat{W}^{_{(n)}}_{\! s_p} \big) 
\xrightarrow[n\to \infty]{\;}
\big( C; Y ; \sqrt{\beta} \widehat{W}_{\! s_0}, \ldots, \sqrt{\beta} \widehat{W}_{\! s_p} \big) $$ 
weakly in $\bC^{_0}_{^1} \! \times \! \bD(\bbR_+, \bbR) \! \times \!\bbR^{p+1}$. 
Thus $\lim_{n\to \infty} ( C^{_{(n)}}_\cdot; V^{_{(n)}}_\cdot ;  \widehat{W}^{_{(n)}}_{\cdot} )
\eqo (C; Y ; \sqrt{\beta }\widehat{W}_{\cdot})$ weakly $\bC^{_0}_{^1} \! \times \! \bD(\bbR_+, \bbR) \! \times \! \bC^{_0}_{^1}$. This entails (\ref{cvsnake1}) by Proposition \ref{endsnake1} $(ii)$. \cq

\smallskip

It remains to prove (\ref{Wendtight}). 
To that end, we next show that 
\begin{equation}
\label{Wmgtight}
\forall N\ino \bbN, \; \forall \epp\ino  (0, 1), \qquad \lim_{p\to \infty} \limsup_{n\to \infty} \bP \big( \vartheta_{N,p} \big( W^{_{(n)}}_{\cdot} \big) \geqo \epp\big) = 0,
\end{equation}
where for all $n,p\ino \bbN$ and for all snake $(h, w)$ we recall that 
$$ \vartheta_{N\! ,\, p} (w) \eqo \!\! \max_{0\leq j< N2^p} \! \! \max \big\{\,  \big| \widehat{w}_{j2^{-p}} - w_{j2^{-p}} \big( h(j2^{-p})  - r\big) \big| \, ; \, r\ino \big[ 0\, , h(j2^{-p}) - m_h (j2^{-p} ,(j+1)2^{-p}) \big] \big\}  \; $$ 
(here $m_h(r_1,r_2)\eqo \min_{r\in [r_1, r_2]} h(r)$, for all real numbers $r_2\geqo r_1\geqo 0$). 

\smallskip

\noi
\emph{Proof of (\ref{Wmgtight})}. Recall that under the assumptions of  Theorem \ref{unifbounded} $(i)$, Lemma \ref{fdcvsna} applies and basic results on weak convergence in $\bC_{^0}^{_1}$ imply that $\vartheta_{N\! ,\, p} ( W^{_{(n)}}_{\cdot} )\! \to \! 
\sqrt{\beta } \vartheta_{N\! ,\, p} ( W_{\cdot} )$ in law in $\bbR_+$. Since $(C,W)$ is a continuous snake, we also get 
$\lim_{p\rightarrow \infty}  \vartheta_{N\! ,\, p} ( W_{\cdot} )\eqo 0$ in probability for all $N\ino \bbN$. This easily entails  (\ref{Wmgtight}).   \cq

\medskip

To complete the proof of  (\ref{Wendtight}), we recall that $\tau_{\! n}^{_\infty}\eqo (\tau_n(p))_{p\in \bbN}$ is a sequence of i.i.d.~GW($\mu_n$)-trees. We recall that  
$ (v_k)_{k\in \bbN}$ (resp.~$(u_l)_{l\in \bbN}$) are the vertices of $\tau_n^{_\infty}$ listed in increasing contour order (resp.~increasing depth-first order).  
For all $y\ino \bbR_+^*$ we recall from Definition \ref{exittime} the notation $(\sigma_{n,q} (y))_{q\in \bbN}$ and $(\bs_{n,q} (y))_{q\in \bbN}$ for the contour and depth-first $y$ oscillation times of increase of 
$(S_{n,u})_{u\in \tau_n^\infty}$ which are recursively defined as follows: 
we set $\sigma_{n,0} (y) \eqo  \bs_{n,0} (y) \eqo 0$ and 
\begin{eqnarray*}
\sigma_{n, q+1} (y) \eqo \inf \big\{ s \geko \sigma_{n,q} (y)   
\!\!\!\! \!\!\!  & : & \!\!\!\!  \!\!\!  | S_{n, v_{k+1}} \!\!  -\! S_{n, v_{k+1} \wedge v_{\sigma_{n,q} (y)+1}}| \geko  y  \big\} \\
& \textrm{and} &  \bs_{n,q+1} (y) \eqo \inf \big\{ l \geko \bs_{n,q} (y)   
:  | S_{n, u_{l+1}} \!\!  -\! S_{n, u_{l+1} \wedge u_{\bs_{n,q} (y )+1}}| \geko  y \big\} , 
\end{eqnarray*}
with the convention that $\inf \emptyset \eqo \infty$. To simplify notation we introduce the normalized versions of these times 
$$ \forall y\ino \bbR_+^*, \; \forall q,n\ino \bbN, \quad 
\sigma^{_{(n)}}_{q} (y)\eqo \tfrac{1}{b_n}\sigma_{n,q} \big( y\sqrt{\lambda_n} \, \big) \quad \textrm{and} \quad  \bs^{_{(n)}}_{q} (y)\eqo \tfrac{1}{b_n} \bs_{n,q} \big( y\sqrt{\lambda_n} \, \big) .$$
Let $N,p, n\ino \bbN$ and $\epp \ino (0, 1)$. We set $A_{n}\eqo \big\{ \omega_{N,p} \big(\widehat{W}^{_{(n)}}_\cdot \big) \geko \epp \big\}$ and $B_{n}  \eqo \big\{  \vartheta_{N,p} \big( W^{_{(n)}}_\cdot \big) \leqo \frac{_1}{^2} \epp\big\}$ and 
for all integers $0\leqo j\leko N2^{p}$, we also set $k_j\eqo \lfloor b_n j2^{-p}\rfloor$. We first observe that 
\begin{eqnarray*}
A_n\cap B_n \!\!\!\!\! & \subset &\!\!\!\!\! \Big\{ \, \forall j\ino \{ 0, \ldots,  N2^{p} \! -\! 1\}  , \, \forall v\ino \lgeo v_{k_j}, v_{k_j}\wedge v_{k_{j+1}}\rgeo \, : \, 
\lvert S_{n, v_{k_j}}  \!\! -\! S_{n,v}\rvert  \leq 2+ \frac{_{_1}}{^{^2}} \epp \sqrt{\lambda_n}  \\
 \textrm{and} \!\! \!\!& & \! \!\!\!\!\!\!\exists i\ino \{ 0, \ldots,  N2^{p} \! -\! 1\} , \, \exists l \ino \{ 1+k_j, \ldots k_{j+1} \} \, : \, 
  \lvert S_{n, v_{k_j}}  \!\! -\! S_{n,v_l}\rvert  \geko \epp \sqrt{\lambda_n} \! -\! 2 \Big\} 
\end{eqnarray*}
We choose $n_0 (\epp,p)\ino \bbN$ such that $\frac{_1}{^6} \epp \sqrt{\lambda_n} \geko 4$ and $2^{-p}b_n \geqo 1$. We apply Lemma \ref{insert} with $y_0\eqo 1$, $y_1\eqo 2+ \frac{_1}{^2} \epp \sqrt{\lambda_n}$, $y_2\eqo \epp \sqrt{\lambda_n} \! -\! 2$ and $y\eqo \frac{_1}{^{12}} \epp \sqrt{\lambda_n}$ and we get 
\begin{eqnarray*}
A_n \cap B_n  \!\!\! & \subset & \!\!\! \big\{ \exists q\ino \bbN, \, \exists i\ino \{ 0, \ldots,  N2^{p} \! -\! 1\} \, : \, k_i \leqo \sigma_{n,q} (y) \leko \sigma_{n, q+1} (y) \leqo k_{i+1} \big\} \\
\!\!\! & \subset & \Big\{ \min \big\{ \sigma^{_{(n)}}_{^{q+1}} (\frac{_{_1}}{^{^{12}}} \epp) \! -\!  \sigma^{_{(n)}}_{^q} (\frac{_{_1}}{^{^{12}}} \epp)\;   ; \;   q\ino \bbN\! :  \sigma^{_{(n)}}_{q}\!  (\frac{_{_1}}{^{^{12}}} \epp)\leqo N  \big\}  \leq 2.2^{-p} \Big\} 
\end{eqnarray*}
We now express the last event in terms of the $ \bs^{_{(n)}}_{^q} (\frac{_{_1}}{^{^{12}}} \epp) $. To that end, we 
recall from Lemma \ref{timeincrease1} that $K_{\bs_{n,q} (y)}\eqo \sigma_{n,q} (y)$ where as in (\ref{Klcontour}) for all $l\ino \bbN$, $K_l$ stands for $2l\! -\! H_l (\tau^{_\infty}_n)$. This implies that 
$$ \Big| \sigma^{_{(n)}}_{^{q+1}} (\frac{_{_1}}{^{^{12}}} \epp) \! -\!  \sigma^{_{(n)}}_{^q} (\frac{_{_1}}{^{^{12}}} \epp)  -2 \big( \bs^{_{(n)}}_{^{q+1}} (\frac{_{_1}}{^{^{12}}} \epp) \! -\!  \bs^{_{(n)}}_{^q} (\frac{_{_1}}{^{^{12}}} \epp)  \big) \Big| \leq \tfrac{2}{{a_n}}\, \max_{s\in [0, N]} H^{_{(n)}}_s \; .$$
We then fix $z\ino \bbR_+^*$ and we introduce the event $C_n \eqo \{ \max_{s\in [0, N]} H^{_{(n)}}_s \leqo z \}$ and the integer 
 $n_0 \eqo n_0(\epp, p,z)$ such that $2z/a_n \leqo 2^{-p}$ and for all integers $n\geqo n_0$.  Thus for all $n\geqo n_0$, we observe that 
$$ A_n \cap B_n \cap C_n \,  \subset  \Big\{ \min \big\{ \mathbf s^{_{(n)}}_{^{q+1}} (\frac{_{_1}}{^{^{12}}} \epp) \! -\!  \mathbf s^{_{(n)}}_{^q} (\frac{_{_1}}{^{^{12}}} \epp)\;   ; \;   q\ino \bbN\! :  \mathbf s^{_{(n)}}_{q}\!  (\frac{_{_1}}{^{^{12}}} \epp)\leqo \frac{_{_1}}{^{^2}}(N+1)  \big\}  \leq 2.2^{-p} \Big\} .$$
Then for all integers $N\geqo 3$ and $n\geqo n_0$, we get 
\begin{eqnarray}
\bP \big( \omega_{N,p} \!\!  \!\! \!\!  \!\! \!\!    & &\!\!  \!\! \!\!  \!\! \!\!    \big( \widehat{W}^{_{(n)}}_{\cdot} \big) \geqo \epp\big)\leq   \bP \big( \vartheta_{N,p} \big( W^{_{(n)}}_{\cdot} \big)  \geqo  \frac{_{_1}}{^{^2}} \epp\big) + \bP \big(  \max_{s\in [0, N]} H^{_{(n)}}_s \geko z  \big)  \nonumber  \\
\label{decompeven1}  +\,  \bP  \Big( \!\!\!\! \!\! &\min &   \!\!\!\!\!  \big\{ \mathbf s^{_{(n)}}_{^{q+1}} (\frac{_{_1}}{^{^{12}}} \epp) \! -\!  \mathbf s^{_{(n)}}_{^q} (\frac{_{_1}}{^{^{12}}} \epp)\;   ; \;   q\ino \bbN\! :  \mathbf s^{_{(n)}}_{q}\!  (\frac{_{_1}}{^{^{12}}} \epp)\leqo \frac{_{_1}}{^{^2}}(N+1)  \big\}  \leq 2.2^{-p}  \Big).
\end{eqnarray}
We next observe that for all $N\ino \bbN$ and all $\epp \ino \bbR_+^*$, 
\begin{equation}
\label{contrHetvarth}
\lim_{p\to \infty} \limsup_{n\to \infty}  \bP \big( \vartheta_{N,p} \big( W^{_{(n)}}_{\cdot} \big)\geqo \frac{_{_1}}{^{^2}} \epp\big) \eqo \lim_{z\to \infty} \limsup_{n\to \infty}  \bP \big(  \max_{s\in [0, N]} H^{_{(n)}}_s \geko z  \big)  =0
\end{equation}
By a general results stated in Ethier \& Kurtz \cite{EtKu86} Lemma 3.8.2 p.~134, we get that 
\begin{eqnarray}
\lim_{p\to \infty} \!\!\! \!\!\! &  &  \!\!\! \!\!\!  \limsup_{n\to \infty}   \bP  \Big(  \!\!\! \!\!\!\!\!\! \!\!\!   \!\!\!\!\!\!  \!\!\!\!\!\!  \!\!\!\!  \min_{\qquad \qquad \quad q\in \bbN\, : \, \bs^{_{(n)}}_{q}\!  (\frac{1}{12}\epp)\leq \frac{1}{2}(N+1) } \!\! \!\!\! \!\!\! \!\!\! \!\!\! \!\!\! \!\!\!\!\!\! \!\!\!\!  \! \bs^{_{(n)}}_{^{q+1}} (\frac{_{_1}}{^{^{12}}} \epp) \! -\!  \bs^{_{(n)}}_{^q} (\frac{_{_1}}{^{^{12}}} \epp) \leq  2.2^{-p} \Big)=0 \; \Longleftrightarrow \;  \nonumber \\
\label{EthKurlem}  \!\!\! \!\!\!  \!\!\!  & &  \!\!\! \!\!\!  \!\!\! \lim_{p\to \infty} \limsup_{n\to \infty} \, \sup_{q\in \bbN} \,  \bP \big(\,  \bs^{_{(n)}}_{q}\!  (\frac{_{_1}}{^{^{12}}}\epp)\leqo \frac{_{_1}}{^{^2}}(N+1) \, ; \, \bs^{_{(n)}}_{^{q+1}} (\frac{_{_1}}{^{^{12}}} \epp) \! -\!  \bs^{_{(n)}}_{^q} (\frac{_{_1}}{^{^{12}}} \epp) \, \leq 2.2^{-p} \big) \eqo 0 .
\end{eqnarray}
By Lemma \ref{timeincrease2} the r.v.s $(\bs^{_{(n)}}_{q+1} (\frac{_{_1}}{^{^{12}}} \epp)\! -\! \bs^{_{(n)}}_{q} (\frac{_{_1}}{^{^{12}}} \epp) )_{q\in \bbN^*}$ are i.i.d.~Thus 
(\ref{EthKurlem}) is equivalent to  
$$\lim_{p\to \infty} \limsup_{n\to \infty}   \bP \big(  \bs^{_{(n)}}_{1}\!  (\frac{_{_1}}{^{^{12}}}\epp) \leq 2.2^{-p} \big) \eqo 0 \; .$$
Namely, by (\ref{decompeven1}), (\ref{contrHetvarth}) and (\ref{EthKurlem}) we see that for all $N\ino \bbN$ and all $\epp \ino (0, 1)$, 
\begin{eqnarray}
\label{eeqquuiivv} \lim_{p\to \infty} \limsup_{n\to \infty} \bP \big( \omega_{N,p} \big( \widehat{W}^{_{(n)}}_{\cdot} \big)\!\! \!\!\!\!\!\! &  \geqo &\!\!\!\!\!\!\!\! \epp\big) \eqo 0  \\
 & \Longleftrightarrow & \lim_{p\to \infty} \limsup_{n\to \infty} \bP \big(\!\!\!\! \!\!\!\! \!\!\!\!  \max_{\quad \quad 0\leq l\leq 2b_n2^{-p}} \!\!\!\! \!\!\!\! \!\!\! \lvert S_{n,u_l} \rvert \geko \frac{_{_1}}{^{^{12}}}\epp \sqrt{\lambda_n} \, \big) \eqo 0. \nonumber
\end{eqnarray}
We fix $x\ino \bbR_+^*$ and we recall the notation $\varsigma_{n,x}\eqo \frac{1}{b_n}\sum_{1\leq p\leq \lfloor a_n x \rfloor} \# \tau_n (p)$ 
\begin{equation}
\label{decoup3}
\bP \big(\!\!\!\! \!\!\!\! \!\!\!\!  \max_{\quad \quad 0\leq l\leq 2b_n2^{-p}} \!\!\!\! \!\!\!\! \!\!\! \lvert S_{n,u_l} \rvert \geko \frac{_{_1}}{^{^{12}}}\epp \sqrt{\lambda_n} \, \big) \leq  \bP \big(\max_{u\in \tau^{_x}_n } \lvert S_{n,u} \rvert \geko \frac{_{_1}}{^{^{12}}}\epp \sqrt{\lambda_n} \, \big) + \bP \big(\varsigma_{n,x} \leq 22^{-p}    \big) 
\end{equation}
We recall from (\ref{szhghtGW}) that $\varsigma_{n,x}$ converges in law to $\varsigma_{-x}$.  
Then, as a consequence of (\ref{decoup3}) and of (\ref{deviaSenn}) in Lemma \ref{deviation} with $4y_1\eqo y_2\eqo \epp/(36c)$, we get 
\begin{eqnarray*} \limsup_{p\to \infty} \limsup_{n\to \infty} \bP \big(\!\!\!\! \!\!\!\! \!\!\!\!  \max_{\quad \quad 0\leq l\leq 2b_n2^{-p}} \!\!\!\! \!\!\!\! \!\!\! \lvert S_{n,u_l} \rvert \geko \frac{_{_1}}{^{^{12}}}\epp \sqrt{\lambda_n} \, \big) \leq 8 \big( 1 \!\!\!\! \! & -&  \!\!\!\! \!  e^{-6xw (\epp (12)^{-3} c^{-1}))} \big) \\
\!\!\!\!  &+ & \!\!\!\!  \bE \big[ \min \big(1 ,  (36c)^2\epp^{-2}\Gamma_{\! x}\big) \big] \; \underset{x\to 0^+}{-\!\!\! -\!\!\! \longrightarrow } \; 0, 
\end{eqnarray*}
by  (\ref{decoup3}). This in turn implies (\ref{Wendtight}) by (\ref{eeqquuiivv}). As already explained, it completes the proof Theorem \ref{unifbounded} $(i)$ in $\textbf{Case (0)}$. \cq

\smallskip

Theorem \ref{unifbounded} $(i)$ in $\textbf{Case (1)}$ and $\textbf{Case (2)}$ is derived from Theorem \ref{unifbounded} $(i)$ in $\textbf{Case (0)}$ exactly as Theorem \ref{Sheuexplain} $(i)$ in $\textbf{Case (1)}$ and $\textbf{Case (2)}$ is derived from Theorem \ref{Sheuexplain} $(i)$ \textbf{Case (0)}. \cq

\smallskip

We now proceed to the proof of Theorem \ref{unifbounded} $(ii)$, which is stated in $\textbf{Case (0)}$ only. 
We assume for all $n\ino \bbN$ 
that $\smash{\bgam_n \eqo \frac{1}{2} (\delta_{-1} + \delta_{1})}$ and that $\smash{\liminf_{n\to \infty} \mu_n (1) \geko 0}$.
By Theorem \ref{unifbounded} $(i)$ \texttt{Sheu} ($\mathbf a, \mathbf b, \bmu$) implies the tightness in $\bC^{_0}_{^1}$ of the laws of $(\widehat{W}^{_{(n)}}_s)_{s\in \bbR_+}$, $n\ino \bbN$. 

Conversely assume that the laws of the processes $\smash{(\widehat{W}^{_{(n)}}_\cdot )_{n\in \bbN}}$ are tight in $\bC^{_0}_{^1}$. 
By Proposition \ref{endsnake1} $(i)$, the laws of the whole snakes $\smash{(W^{_{(n)}}_{\cdot})_{n\in \bbN}}$, are tight in $\bC(\bbR_+, \bC^{_0}_{^1})$. This implies that the laws $\smash{(C^{_{(n)}}_{\cdot}\! ,V^{_{(n)}}_{\cdot} ,  W^{_{(n)}}_{\cdot}  )_{n\in \bbN}}$ are tight in $\bC^{_0}_{^1} \! \times \! \bD (\bbR_+, \bbR)\times \mathbf C (\bbR_+, \bC^{_0}_{^1})$. Let us consider a weak limiting process along a subsequence $(n_k)_{k\in \bbN}$. W.l.o.g.~by (\ref{basiccvtreee}) it can be written as $(C_\cdot, X_\cdot, W')$. Namely, 
$\lim_{k\to \infty} ( C^{_{(n_k)}}_{\cdot}\! ,V^{_{(n_k)}}_{\cdot} , W^{_{(n_k)}}_{\cdot})$ 
$\eqo $  $( C_\cdot, X_\cdot, W')$. Then Lemma \ref{fdcvsna} implies that $W'$ and $W$ have the same finite dimensional marginal laws. Since $W'$ is $\bC (\bbR_+, \bC^{_0}_{^1})$-valued, this shows that there is a continuous version of the $1$-dimensional Brownian snake with lifetime process $C_\cdot$. Thus \texttt{Sheu} $(\psi)$ holds true by Proposition \ref{BrosnapsiH} $(i)$. By taking $W$ continuous, Lemma \ref{fdcvsna} also implies that 
$(C_\cdot, X_\cdot, W')$ and $(C_\cdot, X_\cdot, W)$ have the same law. This proves existence and uniqueness of weak limits of the laws of the processes 
$( C^{_{(n)}}_{\cdot}\! ,V^{_{(n)}}_{\cdot} , W^{_{(n)}}_{\cdot}  )_{n\in \bbN}$, which therefore implies 
\begin{equation}
\label{cvsnake1BIS}
\big( C^{_{(n)}}_\cdot, V^{_{(n)}}_\cdot ,  W^{_{(n)}}_\cdot \big)  \xrightarrow[n\to \infty]{\; } (C, Y , W) .
\end{equation}
Let $\smash{x\ino \bbR_+^*}$. By general arguments (see e.g.~Jacod \& Shiryaev \cite{JaSh02} Proposition 2.11, Chapter VI, Section 2a p.~341 or more specifically Broutin, D.~\& Wang \cite{BrDuWa21} Lemma B.3 $(iv)$) 
$\varsigma_{n,x} \! \to \! \varsigma_{-x}$ holds jointly in law with (\ref{cvsnake1BIS}). This implies for all $\smash{x\ino \bbR_+^*}$ that 
\begin{equation}
\label{limmaxdisplBIS} 
M_{n, x}\! :=\! \! \! \! \! \! \! \!  \max_{\quad s\in [0, 2\varsigma_{n,x}] } \! \! \! \! \! \!  \widehat{W}^{_{(n)}}_{\! s}=\tfrac{1}{\sqrt{\lambda_n}} \max_{1\leq p\leq \lfloor a_nx \rfloor} \max_{u\in \tau_n (p)} \! S_{n,u} (p) \; \, 
 \xrightarrow[n\to \infty]{ \textrm{(law)}}  \! \! \! \! \! \!  \max_{\quad s\in [0, 2\varsigma_{-x}] }\! \! \! \! \! \!   \widehat{W}_{s} \; . 
\end{equation}
Since the r.v.s $\max_{u\in \tau_n (p)} S_{n,u} (p) $ are i.i.d., we get $\bP (M_{n,x} \leqo z)\eqo (1-\frac{1}{a_n} \varphi (n,z))^{\lfloor a_n x \rfloor} $ where we have set $\varphi (n, z)\eqo a_n \bP (\max_{u\in \tau_n (p)} S_{n,u} (p)  \geko \lambda_n z)$. We next recall from Proposition \ref{BrosnapsiH} $(ii)$  that $\bP (\max_{s\in [0, 2\varsigma_{-x}] } \widehat{W}_{s} \leqo z) \eqo e^{-xw(z)}$ where $w$ stands for the inverse of $F$ as defined in (\ref{defFnF}). Therefore (\ref{limmaxdisplBIS}) implies that $\lim_{n\to \infty} \varphi (n,z)\eqo w(z)$ for all $z\ino \bbR_+^*$. Then we use 
(\ref{upperbound}) in Lemma \ref{keylower} to get $F_n (\varphi (n,z)) \leqo 2^{7/2} \mu_n(1)^{-5/2} z$ 
(here $F_n$ is defined in (\ref{defFnF}) and in (\ref{upperbound}) we use the bound $\varphi'(1) \geqo \frac{_1}{^2}\mu(1)$). 
Let us set $q\eqo \frac{_1}{^2}\liminf_{n\to \infty} \mu_n (1) $ which is assumed to be positive. Then, there is an integer 
$n_1\eqo n_1(z)$ such that for all $n\geqo n_1$, $\varphi (n,z) \leko 2w(z)$ and $\mu_n(1) \geko q$ and thus 
$F_n (2w(z)) \leqo 2^{7/2} q^{-5/2} z$. This implies $\lim_{z\to 0^+}\limsup_{n\to \infty} F_n (2w(z)) \eqo 0$, which implies \texttt{Sheu}($ \mathbf a , \mathbf b, \bmu$) since $\lim_{z\to 0^+} w(z) \eqo \infty$. This completes the proof of Theorem \ref{unifbounded}. \cqfd

\end{document}